\documentclass{svmult}%
\usepackage{amsmath}%
\usepackage{amssymb}%
\usepackage[sort&compress,comma,numbers]{natbib}%

\usepackage{mathptmx}       
\usepackage{helvet}         
\usepackage{courier}        
\usepackage{makeidx}         
\usepackage{graphicx}        
\usepackage{multicol}        
\usepackage[bottom]{footmisc}

\begin{document}
\title*{Tarskian classical relevant logic}
\titlerunning{Tarskian classical relevant logic}
\author{Roger D.\ Maddux}
\institute{Roger D. Maddux \at Department of Mathematics, Iowa State
  University, Ames, Iowa 50011-2066, USA,
  \newline \email{maddux@iastate.edu}.
\newline Date: \today.}
%
%
\allowdisplaybreaks%
\def\alg#1#2{\text{$\mathsf{#1}_{#2}$}}%
\def\ra#1{\textup{\rm R$_{#1}$}}%
\def\conv#1{\setbox13\hbox{$#1$}\ifdim\wd13<1.501em{#1}^\smallsmile
  \else{(#1)^\smallsmile}\fi}%
\def\halfthinspace{\relax\ifmmode\mskip.5\thinmuskip\relax
  \else\kern.8888em\fi}%
\let\hts=\halfthinspace%
\def\rp{{\hts;\hts}}%
\def\rs{\mathop\dagger}%
\def\id{{1\kern-.08em{\text{\rm'}}}}%
\def\di{{0\kern-.04em{\text{\rm'}}}}%
\let\SS=\S%
\def\makecs#1#2{\makecsX {#1}#2,.}%
\def\makecsX#1#2#3.{\onecs{#1}{#2}%
  \ifx#3,\let\next\eatit\else\let\next\makecsX\fi\next{#1}#3.}%
\def\onecs#1#2{\expandafter\gdef\csname #2\endcsname%
  {{\csname #1\endcsname {#2}}}}%
\def\eatit#1#2.{\relax}%
\makecs{}{abcdefghijklmnopqrstuvwxyzABCDEFGHIJKLMNOPQRSTUVWXYZ}%
\def\ie{{\it i.e.}}%
\def\({\left(}%
\def\){\right)}%
\def\<{\left<}%
\def\>{\right>}%
\def\ex#1{\exists_{#1}}%
\def\all#1{\forall_{#1}}%
\def\gc#1{\mathfrak{#1}}
\def\NA{\mathsf{NA}}%
\def\SA{\mathsf{SA}}%
\def\RA{\mathsf{RA}}%
\def\RRA{\mathsf{RRA}}%
\def\MA{\mathsf{MA}}%
\def\VV{\mathsf{V}}%
\def\GG{\mathsf{G}}%
\def\HH{\mathsf{H}}%
\def\JJ{\mathsf{J}}%
\def\RR{\mathsf{R}}%
\def\RM{\mathsf{RM}}%
\def\KR{\mathsf{KR}}%
\def\CT{\mathsf{CT}}%
\def\TR{\mathsf{TR}}%
\def\TT{\mathsf{T}}%
\def\CR{\mathsf{CR^*}}%
\def\II{\mathcal{I}}%
\def\LL{\mathcal{L}}%
\def\Lp{{\mathcal L}^+}%
\def\Lx{{\mathcal L}^{\times}}%
\def\Lwx{{\mathcal L}\kern-.15em{w}^{\times}}%
\def\Ls{{\mathcal L}\kern-.15em{s}}%
\def\Mn{\mathcal{M}^{(n)}}%
\def\Cm#1{\setbox13\hbox{$#1$}\ifdim\wd13=0pt{\gc{Cm}}%
  \else{\gc{Cm}\({#1}\)}\fi}%
\def\mand{\mathrel{\text{ and }}}
\def\mor{\mathrel{\text{ or }}}
\def\MP{\text{\sf{MP}}}%
\def\iff{\mathrel\Leftrightarrow}%
\def\ifff{\text{iff }}%
\def\min#1{\overline{#1}}%
\def\rmin{{\sim}}%
\def\blank{\phantom{o}}%
\def\proves{\,\vdash}%
\def\0{\var0}%
\def\1{\var1}%
\def\2{\var2}%
\def\3{\var3}%
\def\fm#1#2#3{#2#1#3}%
\def\implies{\mathrel\Rightarrow}%
\def\stand{\hbox{\vrule height10pt depth5pt width0pt}}%
\def\rationals{\mathbb{Q}}%
\def\ra#1{\textup{\rm R$_{#1}$}}%
\def\display#1{$${\text{{#1}}}$$}%
\def\repl{{\sf Repl}}%
\def\tran{{\sf Trans}}%
\def\closure#1{\pmb[\,#1\,\pmb]}%
\def\eqsym{\id}%
\def\switch#1#2{\mathsf{S}_{#1#2}}%
\def\sep{\kern1pt\pmb|\kern1pt}%
\def\equals{\overset{\scriptscriptstyle\circ}=}%
\def\star#1{{#1}^*}%
\def\Re#1{\gc{Re}(#1)}%
\def\Sb#1{\gc{Sb}(#1)}%
\def\var#1{\mathsf v_{#1}}%
\def\ttt{\mathbf\t}%
\def\hom{{\rm End}(\gc\P)}%
\def\De#1{\mathsf{De}_{\gc\U}(#1)}%
\def\Dee{\mathsf{De}_{\gc\U}}%
\def\Dei{\mathsf{De}_{\gc\U_\i}}%
\def\holds{\scriptstyle\circ}%
\def\vi{\var\i}%
\def\vj{\var\j}%
\def\vk{\var\k}%
\def\prop#1{{\rm Pr}[\,#1\,]}%
\def\dee#1{\text{\bf\sf D}(#1)}%
\def\aut{\text{\sf Aut}}%
\def\autb{\text{\sf Aut}_1}%
\def\autf{\text{\sf Aut}_2}%
\def\autr{\text{\sf Aut}_{1,2}(R)}%
\def\perm{\text{\sf Perm}(K)}%
\def\atoms{\text{At}(\gc{A})}%
\setlength{\unitlength}{2mm}%
\maketitle%
\abstract{The Tarskian classical relevant logic $\TR$ arises from
  Tarski's work on the foundations of the calculus of relations and on
  first-order logic restricted to finitely many variables, presented
  by Tarski and Givant their book, \emph{A Formalization of Set Theory
    without Variables}, and summarized in first nine sections. $\TR$
  is closely related to the well-known logic $\KR$.  Every formula of
  relevance logic has a corresponding sentence in Tarski's extended
  first-order logic of binary relations with operators on the relation
  symbols.  A formula is in $\TR$ (by definition), or in $\KR$ (by a
  theorem), if and only if its corresponding sentence can be proved in
  first-order logic, using at most four variables, from the
  assumptions that all binary relations are dense and, for $\TR$,
  commute under composition, or, for $\KR$, are symmetric.  The
  vocabulary of $\TR$ is the same as the classical relevant logic
  $\CR$ proposed by Meyer and Routley but $\TR$ properly contains
  $\CR$. The frames characteristic for $\TR$ are the ones that are
  characteristic for $\CR$ and satisfy an extra frame condition.
  There are formulas in $\TR$ (but not in $\CR$) that correspond to
  this frame condition and provide a counterexample to a theorem of
  T.\ Kowalski.  The frames characteristic for $\TR$, or $\KR$, are
  the ones whose complex algebras are integral dense relation algebras
  that are commutative, or symmetric, respectively.  For both classes,
  the number of isomorphism types grows like the number of isomorphism
  types of ternary relations.  Asymptotic formulas are obtained for
  both classes. Similar results apply to a hierarchy of logics defined
  by the number of variables used in the first-order proofs of their
  corresponding sentences.}  \keywords{relevance logic, classical
  relevant logic, relation algebras, semi-associative relation
  algebras, provability in first-order logic with finitely many
  variables, sequent calculus}
\section{Introduction}\label{sect1}
In 1975, Alfred Tarski delivered a pair of lectures on relation
algebras at the University of Campinas. The videotaped lectures were
eventually transcribed and published in 2016 \cite{Tarski1975}. At the
end of his second lecture, Tarski said \cite[p.\,154]{Tarski1975},
\begin{quote}
  ``And finally, the last question, if it is so, you could ask me a
  question whether this definition of relation algebra which I have
  suggested and which I have founded---I suggested it many years
  ago---is justified in any intrinsic sense.  If we know that these
  are not all equations which are needed to obtain representation
  theorems, this means, to obtain the algebraic expression of
  first-order logic with two-place predicate, if we know that this is
  not an adequate expression of this logic, then why restrict oneself
  to these equations? Why not to add strictly some other equations
  which hold in representable relation algebras or maybe all?''
\end{quote}
Tarski defined relation algebras as those that satisfy the axioms
\eqref{BI}--\eqref{BX} in Table \ref{Lx-axioms}.  Each axiom is an
equation $\A\equals\B$ between predicates $\A,\B$ in Tarski's extended
system $\Lp$ of first-order logic (described in detail in
\SS\ref{sect3}).  In this system, $\id$ denotes the identity relation,
$+$ is an operation on predicates denoting union, $\min\blank$ denotes
complementation, $\rp$ denotes relative multiplication, and $\equals$
is a symbol denoting the equality of predicates, according to Tarski's
definitional axioms for $\Lp$ listed in Table \ref{L+defs}.  Since
$\Lp$ is a definitional extension of first-order logic $\LL$
(described in \SS\ref{sect2}) every equation $\A\equals\B$ in $\Lp$
can be translated into a sentence $\GG(\A\equals\B)$ of first-order
logic $\LL$ by eliminating predicate operators according to the
elimination mapping $\GG$, defined in Table \ref{G-defs}.  The answer
to Tarski's question ``whether this definition of relation algebra
\dots\ is justified in any intrinsic sense'' is Theorem
\ref{thm6}\eqref{3} in \SS\ref{sect9}: an equation is derivable from
the axioms for relation algebras iff its translation can be proved
with no more than four variables.

Half of this answer was known to Tarski already in the early 1940s.
The other half was proved thirty years later \cite{MR2628352}.  In a
manuscript started in 1942, Tarski created a system $\LL_3$ of logic
with only three variables (described in \SS\ref{sect7}) that is
equipollent with the equational theory $\Lx$ of relation algebras
(described in \SS\ref{sect5}). The equipollence of $\LL_3$ with $\Lx$
is stated as Theorem \ref{thm3} in \SS\ref{sect8}.

By 1953 Tarski had shown that set theory can be formalized in $\Lx$ as
equations between predicates of first-order logic, with proofs based
on just the axioms for relation algebras with substitution and modus
ponens as the only rules of inference. This result, announced in
\cite{Tarski1953a}, was eventually published in the book by Tarski and
Givant, \emph{A Formalization of Set Theory without Variables}, where
Theorem \ref{thm6}\eqref{3} is mentioned \cite[p.\,89,
  p.\,209]{MR920815}.

The characterization of the equations true in relation algebras as the
ones whose translations into first-order logic are provable with four
variables can be applied to the relevance logic $\RR$ of
\cite{MR0115902, MR0406756, MR1223997, MR141590, MR239942} and to the
classical relevant logic $\CR$ of \cite{MR0363789, MR0363789a}. The
connectives of $\RR$ and $\CR$ can be interpreted as operations on
binary relations according to Table \ref{defs-ops}. Define a predicate
$\A$ of $\Lp$ to be valid under density and commutativity if $\A$
denotes a relation containing the identity relation whenever the
connectives in $\A$ are interpreted as operations on a set $\S$ of
dense binary relations, where $\S$ is closed under the operations in
Table \ref{defs-ops} that correspond to connectives occurring in $\A$
and $\S$ is closed and commutative under relative multiplication. When
interpreted this way, every predicate in $\RR$ or $\CR$ is valid under
density and commutativity.
\begin{table}
  \begin{align*}
    &\text{Name of connective}&&\text{Interpretation as an operation}\\
    &\text{disjunction}&\A\lor\B&=\{\<\x,\y\>:\<\x,\y\>\in\A
    \mor\<\x,\y\>\in\B\}\\
    &\text{conjunction}&\A\land\B&=\{\<\x,\y\>:\<\x,\y\>\in\A
    \mand\<\x,\y\>\in\B\}\\
    &\text{Boolean negation}&\neg\A
    &=\{\<\x,\y\>:\x,\y\in\U\mand\<\x,\y\>\notin\A\}\\
    &\text{De Morgan negation}&\rmin\A&=\{\<\x,\y\>:\x,\y\in\U\mand
    \<\y,\x\>\notin\A\}\\
    &\text{implication}&\A\to\B&=\{\<\x,\y\>:\text{$\x,\y\in\U$, and
      for all $\z\in\U$,}\\
    &&&\qquad\qquad\text{if $\<\z,\x\>\in\A$ then $\<\z,\y\>\in\B$}\}\\
    &\text{fusion}&\A\circ\B&=\{\<\x,\y\>: \text{for some $\z$, }
    \<\x,\z\>\in\B\text{ and }\<\z,\y\>\in\A\}\\
    &\text{Routley star}&\star\A&=\{\<\x,\y\>:\<\y,\x\>\in\A\}\\
    &\text{truth}&\ttt&=\{\<\x,\x\>:\x\in\U\} \\
  \end{align*}
  \caption{Operations for interpreting formulas as relations on a 
    base set $U$.}
  \label{defs-ops}
\end{table}

For predicates $\A$ and $\B$ let $\A\leq\B$ be the equation
$\A+\B\equals\B$. The equation $\id\leq\A$ asserts that $\A$ contains
the identity relation.  The density of $\A$ is expressed by the
equation $\A\leq\A\rp\A$ and commutativity by
$\A\rp\B\equals\B\rp\A$. By the completeness theorem for first-order
logic, if $\A$ is in $\RR$ or $\CR$ then the translation
$\GG(\id\leq\A)$ is provable from the sets of equations expressing
density and commutativity, as defined in \eqref{dens-eqs} and
\eqref{comm-eqs} in \SS\ref{sect10}. Any such proof will involve some
finite number of variables but it turns out that if $\A$ is a theorem
of $\RR$ or $\CR$ then the translation of $\id\leq\A$ into $\LL$ can
be actually be proved with no more than four variables (Theorem
\ref{thm11}).

In \SS\ref{sect19} there are two examples, \eqref{L''} and
\eqref{M''}, of a predicate $\A$ with the property that $\id\leq\A$
translates to a logically valid sentence that cannot be proved with
four variables. In both cases the translation $\GG(\id\leq\A)$ can be
proved with five variables and requires no appeal to density or
commutativity. Theorem \ref{5var} shows that they are 5-provable but
not 4-provable. Predicates in the vocabulary of $\CR$ that are
5-provable but not 4-provable have been known for a long time but
\cite{MR2496334} was the first to find examples in the vocabulary of
$\RR$. By creating infinitely many such predicates, \cite{MR2496334}
proved that $\TT_\omega$ is not finitely axiomatizable.  The two
examples in \SS\ref{sect19} were created later \cite[\SS8]{MR2641636}.

Tarski's classical relevant logic $\TR$ is defined in \eqref{TR-def}
of \SS\ref{sect10} as the set of predicates $\A$ in the vocabulary of
$\CR$ such that $\GG(\id\leq\A)$ is 4-provable from density and
commutativity.  Consequently $\CR\subseteq\TR$. Although $\TR$ does
not contain any 5-provable formulas, equality still fails.  The frame
conditions \eqref{center reflection}--\eqref{right reflection} in
\SS\ref{sect12} hold in the frames characteristic for $\TR$ (Theorem
\ref{thm10-}\eqref{charTR}) but they do not hold in all the frames
that are characteristic for $\CR$ (Theorem \ref{thm15}\eqref{thm15i}).
The frame conditions \eqref{center reflection}--\eqref{right
  reflection} correspond to predicates \eqref{reflection1} and
\eqref{reflection1a}.  These predicates were created using the same
device by which \eqref{L''} and \eqref{M''} were obtained from
predicates in the vocabulary of $\CR$ that are 5-provable but not
4-provable.  They are confined to the vocabulary of $\RR$ and belong
to $\TR$, but they are not theorems of $\RR$ (Theorem
\ref{thm15}\eqref{thm15i}). In Theorem \ref{thm15}\eqref{thm15ii} they
are shown to be valid in a frame satisfying \eqref{identity} iff it
satisfies \eqref{left reflection}.  Furthermore, \eqref{reflection1}
and \eqref{reflection1a} are 3-provable without assuming density or
commutativity (Theorem \ref{thm11}).  These observations show in
\SS\ref{sect18} that \cite[Thm 8.1]{MR3289545} is incorrect.

Although $\CR$ and $\RR$ cannot be characterized as the formulas that
are 4-provable from density and commutativity, $\TR$ does have that
characterization simply because it is defined that way. The logic
$\KR$, which figures prominently in the research of Alasdair Urquhart
\cite[\SS65]{MR1223997}, \cite{MR771777, MR1234403, MR1720799,
  Urquhart2017, MR4038537} also has such a characterization despite
being defined in a completely different way.  Both logics $\TR$ and
$\KR$ can be correlated with classes of relation algebras.  A
predicate $\A$ is in $\TR$ iff the equation $\id\leq\A$ is true in
every dense commutative relation algebra iff $\GG(\id\leq\A)$ is
4-provable from density and commutativity (Theorem
\ref{thm10-}\eqref{charTR}).  Similarly, $\A$ is in $\KR$ iff the
equation $\id\leq\A$ is true in every dense symmetric relation algebra
iff $\GG(\id\leq\A)$ is 4-provable from density and symmetry (Theorem
\ref{thm10-}\eqref{charKR}).

Theorem \ref{thm11} shows that dozens of formulas and rules are
provable with one to four variables, with or without additional
non-logical assumptions selected from density, commutativity, or
symmetry.  For example, the permutation axiom $(\A\to(\B\to\C)) \to
(\B\to(\A\to\C))$ is 4-provable from commutativity (Theorem
\ref{thm11}\eqref{perm}) and the contraction axiom $(\A\to(\A\to\B))
\to (\A\to\B)$ is 4-provable from density (Theorem
\ref{thm11}\eqref{contract2}).  It was recognized long ago that
density and commutativity are optional hypotheses.  For example,
permutation \eqref{perm} and contraction \eqref{contract2} are not
taken as axioms of Basic Logic \cite{Routleyetal1982, MR3728341}.

Many formulas of $\RR$ and $\CR$ depend on density and commutativity.
The number of variables required for a proof of validity is another
classificatory principle.  For example, Theorem \ref{thm12} shows that
permutation \eqref{perm} is not 3-provable from density and symmetry
(which implies commutativity by Lemma \ref{SA:symm->comm}). In fact,
Theorem \ref{thm13} shows that permutation is not even
$\omega$-provable from density alone without commutativity.
Similarly, Theorem \ref{thm12} shows that contraction
\eqref{contract2} is not 3-provable from density and symmetry while
Theorem \ref{thm-new} shows that it is not $\omega$-provable from
commutativity alone without density.
\begin{table}
  \centerline{\bf Formalisms of \cite{MR920815}}\medskip
  \centerline{
    \begin{tabular}{|l|c|l|l|l|c|}\hline  
      \stand Formalism&Section&Sentences&Provability
      &\qquad Axioms&Rules\\ \hline 
      \stand\quad $\LL$&\SS\ref{sect2}&\quad$\Sigma$
      &\quad$\proves$&\eqref{AI}--\eqref{AIX}&$\MP$\\
      \stand\quad $\Lp$&\SS\ref{sect3}&\quad$\Sigma^+$
      &\quad$\proves^+$
      &\quad"\quad plus~\eqref{DI}--\eqref{DV}&$\MP$\\\hline
      \stand\quad $\Lx$&\SS\ref{sect5}&\quad$\Sigma^\times$
      &\quad$\proves^\times$
      &\eqref{BI}--\eqref{BIII}\eqref{BIV}\eqref{BV}--\eqref{BX}
      &\repl,~\tran\\
      \stand\quad $\LL_3$&\SS\ref{sect7}&\quad$\Sigma_3$
      &\quad$\proves_3$
      &\eqref{AI}--\eqref{AVIII}\eqref{AIX'}\eqref{AX}&$\MP$\\
      \stand\quad $\Lp_3$&\SS\ref{sect7}&\quad$\Sigma^+_3$
      &\quad$\proves^+_3$
      &\quad"\quad plus~\eqref{DI}--\eqref{DV}&$\MP$\\\hline
      \stand\quad $\Lwx$&\SS\ref{sect9}&\quad$\Sigma^\times$
      &\quad$\proves^\times_s$
      &\eqref{BI}--\eqref{BIII}\eqref{BIV'}\eqref{BV}--\eqref{BX} 
      &\repl,~\tran\\
      \stand\quad $\Ls_3$&\SS\ref{sect9}&\quad$\Sigma_3$
      &\quad$\proves_s$&\eqref{AI}--\eqref{AVIII}\eqref{AIX'}&$\MP$\\
      \stand\quad $\Ls_3^+$&\SS\ref{sect9}&\quad$\Sigma^+_3$
      &\quad$\proves^+_s$
      &\quad"\quad plus~\eqref{DI}--\eqref{DV}&$\MP$\\\hline
      \stand\quad $\LL_4$&\SS\ref{sect9}&\quad$\Sigma_4$
      &\quad$\proves_4$&\eqref{AI}--\eqref{AVIII}\eqref{AIX'}&$\MP$\\
      \stand\quad $\Lp_4$&\SS\ref{sect9}&\quad$\Sigma^+_4$
      &\quad$\proves^+_4$
      &\quad"\quad plus~\eqref{DI}--\eqref{DV}&$\MP$\\\hline
    \end{tabular}}
  \bigskip
  \caption{
    For a given system $\mathcal F$, the axioms for $\mathcal F^+$ 
    are those of $\mathcal F$ plus \eqref{DI}--\eqref{DV}. 
    \newline
    Axioms for $\mathcal L^\times$ and $\RA$ coincide.
    Axioms for $\mathcal Lw^\times$ and $\SA$ coincide.
    \newline
    Formalisms between horizontal lines are equipollent. }
  \label{table1}
\end{table}

The remainder of this introduction is a detailed review of the
contents of each section.  \SS\SS\ref{sect2}--\ref{sect9} present
Tarski's work on logic with finitely many variables.  These sections
summarize the first 100 pages of \cite{MR920815}.  Table \ref{table1}
lists the formalisms treated in \SS\SS\ref{sect2}--\ref{sect9}.
First-order logic $\LL$ is presented in \SS\ref{sect2}.  Axioms for
first-logic are listed in Table \ref{L-axioms}.  \SS\ref{sect3}
presents the definitional extension $\Lp$ of $\LL$.  Axioms that
define the predicate operators are listed in Table \ref{L+defs}.  In
\SS\ref{sect4}, the translation mapping $\GG$, defined in Table
\ref{G-defs}, shows how to translate equations in the first-order
sentences.  Theorem \ref{thm1} states that $\LL$ and $\Lp$ are
equipollent in means of expression and proof.
\SS\SS\ref{sect2}--\ref{sect4} summarize \cite[Ch.\ 1--2]{MR920815}.

The equational formalism $\Lx$ of \cite[Ch.\ 3]{MR920815} is defined
in \SS\ref{sect5}.  Axioms for $\Lx$ and $\Lwx$ are listed in Table
\ref{Lx-axioms}. The axioms for $\Lx$ and $\Lwx$ also axiomatize the
variety $\RA$ of relation algebras and the variety $\SA$ of
semi-associative relation algebras, respectively. If both associative
laws \eqref{BIV} and \eqref{BIV'} are excluded, the remaining
equations in Table \ref{Lx-axioms} axiomatize the variety $\NA$ of
non-associative relation algebras.  Theorem \ref{thm1a} says that
$\Lx$ and $\Lwx$ are subformalisms of $\Lp$. \SS\ref{sect6} contains
basic definitions and facts about relation algebras, semi-associative
relation algebras, non-associative relation algebras, representable
relation algebras, and the rules of equational logic.  Dense,
commutative, symmetric, simple, and integral algebras are defined in
\SS\ref{sect6}.  Key facts, presented in Lemmas \ref{SA:symm->comm}
and \ref{commSA}, are that symmetric semi-associative relation
algebras are commutative and simple commutative semi-associative
relation algebras are integral.  The predicate algebra $\gc\P$, proper
relation algebras, the variety $\RRA$ of representable relation
algebras, algebras of binary relations, the satisfaction relation, and
the denotation function are all defined in \SS\ref{sect6}.  Theorems
\ref{eq} and \ref{sem} relate provability in $\Lx$, $\Lwx$, and $\Lp$
to truth in $\RA$, $\SA$, and $\RRA$, respectively. The free $\RA$,
$\SA$, and $\RRA$ are constructed as quotients of the predicate
algebra $\gc\P$.

$\Lx$ is compared to $\Lp$ in \SS\ref{sect7}.  $\Lx$ is weaker than
$\Lp$ because there is a 5-provable but not 4-provable equation and a
sentence not equivalent to any equation.  A study of the translation
mapping $\GG$ leads to Tarski's idea for a 3-variable formalism.
Tarski's proposal is realized in the construction of the 3-variable
formalisms $\LL_3$ and $\Lp_3$.  The axioms for these formalisms
include the associative law \eqref{BIV}, which requires four variables
to prove.  Theorems \ref{thm2} and \ref{thm3} in \SS\ref{sect8} state
that $\LL_3$, $\Lp_3$, and $\Lx$ are equipollent in means of
expression and proof.  An alternative way to express 3-variable
sentences as equations is also presented.

\citet[\SS3.10]{MR920815} introduced the \emph{standardized}
3-variable formalisms $\LL_s$ and $\Lp_s$.  These formalisms (which
should have been called $\LL_3$ and $\Lp_3$, but the names were
already taken) are shown to be equipollent to $\Lwx$ in
\SS\ref{sect9}.  Weakening \eqref{BIV} to the semi-associative law
\eqref{BIV'}, which only requires three variables to prove, produces
the equational formalism $\Lwx$ and the standardized 3-variable
formalisms $\Ls_3$ and $\Ls_3^+$. Because it is 4-provable, the
associative law \eqref{BIV} is part of the standard 4-variable
formalism $\Lp_4$.  Theorem \ref{thm5} states that the formalisms
$\Lwx$, $\Ls_3$, and $\Ls_3^+$ are equipollent in means of expression
and proof, Theorem \ref{thm5+} gives the connections between
provability in $\Lp_4$ $\Lp_3$, and $\Lx$, and Theorem \ref{thm6}
links theories in $\Lp_4$, $\Lp_3$, and $\Lx$.

In \SS\ref{sect10} the logics $\TT^\Psi_n$, $\CT^\Psi_n$, $\TT_n$,
$\CT_n$ for $3\leq n\leq\omega$, and $\TR$ are officially defined,
where $\Psi$ is a set of equations that serve as non-logical
assumptions.  For this purpose, equations of density $\Xi^d$,
commutativity $\Xi^c$, and symmetry $\Xi^s$ are defined in
\eqref{dens-eqs}, \eqref{comm-eqs}, and \eqref{symm-eqs},
respectively, and some predicate operators are defined in
\eqref{or}--\eqref{t=id} for use as connectives in relevance
logic. The definitions match the interpretations in Table
\ref{defs-ops}. For example, the Routley star is converse, truth
$\ttt$ is the identity predicate $\id$, $\neg$ is used as Boolean
negation (as well as negation in $\LL$), and $\circ$ is defined as
reversed relative multiplication.  Two groups of operators are
distinguished, the ``relevance logic operators'' and ``classical
relevant logic operators''.  The difference between logics $\TT_n$ and
$\CT_n$ lies only in their vocabularies: $\TT_\n$ uses the relevance
logic operators while $\CT_\n$, its ``classical'' counterpart, uses
the classical relevant logic operators.

\SS\ref{sect12} contains a review of material on frames, including the
definitions of $\KR$-frames, $\KR$ itself, the $\CR$-frames
characteristic for $\CR$, complex algebras of frames, the pair-frame
on a set, validity in a frame, and 12 frame conditions.  Lemmas
\ref{five} and \ref{lem8} present some basic connections between frame
conditions.  Theorem \ref{thm7-} relates conditions on complex
algebras to frame conditions.  Theorem \ref{thm7} characterizes frames
whose complex algebras are in $\NA$, $\SA$, and $\RA$.  Theorem
\ref{thm8} is the Representation Theorem for $\NA$, $\SA$, and $\RA$:
every algebra in $\NA$, $\SA$, or $\RA$ is embeddable in the complex
algebra of a frame satisfying the characteristic conditions in Theorem
\ref{thm7}.  Theorem \ref{lem8+} says that every group can be viewed
as a frame.  By Lemma \ref{lem8a}, all $\CR$-frames are commutative.
Theorem \ref{thm8a} is the key connection between frames and relation
algebras: the complex algebra of a frame $\gc\K$ is a dense
commutative relation algebra iff $\gc\K$ is a $\CR$-frame satisfying
\eqref{left reflection} and the complex algebra of $\gc\K$ is a dense
symmetric relation algebra iff $\gc\K$ is a $\KR$-frame.

The sequent calculus from \cite{MR722170} is presented in
\SS\ref{sect11} with definitions of $\n$-provability in the sequent
calculus, $\n$-dimensional relational basis (Definition \ref{basis}),
and the variety $\RA_n$ of $\n$-dimensional relation algebras.
Theorem \ref{key} gives the key connections between algebras and
provability in the sequent calculus: 3-provability matches up with
$\SA$, 4-provability with $\RA$, $\omega$-provability with $\RRA$, and
the $\n$-provability of a sequent is characterized as a satisfaction
relation on an algebra.  Lemma \ref{equiv} has several derived rules
of inference for the sequent calculus.

In \SS\ref{sect13} the results of of Tarski and Givant are combined
with the sequent calculus and frame characterizations to characterize
$\KR$, $\TR$, and Tarski's relevance logics of 3, 4, and $\omega$
variables.  Lemma \ref{lem5} says that an equation is true in every
algebra satisfying the equations in $\Psi$ iff a certain condition on
homomorphisms holds. This lemma is used for the major characterization
theorems.  Theorem \ref{thm9-} characterizes $\CT^\Psi_3$ and $\CT_3$
in six ways, Theorem \ref{thm10-} characterizes $\KR$, $\TR$,
$\CT^\Psi_4$, and $\CT_4$ in eight ways, Theorem \ref{KT}
characterizes $\KR$ and $\TR$ in two more ways, and Theorem
\ref{thm10} characterizes $\CT^\Psi_\omega$ and $\CT_\omega$ in four
ways. Theorem \ref{two-obs-2} extends the characterizations to cover
$\TT^\Psi_3$, $\TT_3$, $\TT^\Psi_4$, $\TT_4$, $\TT^\Psi_\omega$, and
$\TT_\omega$.

For Theorem \ref{thm11} in \SS\ref{sect14}, dozens of predicates and
rules are provided with $\n$-proofs in the sequent calculus, where
$\n$ ranges from 1 to 4, sometimes under various non-logical
assumptions.  Theorem \ref{thm12} in \SS\ref{sect15} shows that the
eleven predicates of Theorem \ref{thm11} that are 4-provable
(sometimes from density, commutativity, or both) are not 3-provable
from density and symmetry.  The proof uses a frame in Table \ref{K1}
whose complex algebra is a dense symmetric semi-associative relation
algebra that is not a relation algebra.  Theorem \ref{thm13} in
\SS\ref{sect16} shows the five predicates of Theorem \ref{thm11} that
use commutativity are not $\omega$-provable from density alone. The
proof uses a frame in Table \ref{K2} whose complex algebra is a
non-commutative dense representable relation algebra.  Theorem
\ref{thm-new} in \SS\ref{sect16a} shows six predicates of Theorem
\ref{thm11} that rely on density are not $\omega$-provable from
symmetry. The proof uses the frame in Table \ref{K9} of the 2-element
group.

Theorem \ref{thm15} of \SS\ref{sect17} shows that $\TR$ exceeds $\CR$.
The frames characteristic for $\TR$ satisfy the frame condition
\eqref{left reflection}, which is needed to insure axiom \eqref{BIX}
holds.  However, according to \cite[p.\ 104]{MR2067967},
\begin{quote}
  ``But (56) [$(\A\circ\B)^{-1}=\B^{-1}\circ\A^{-1}$] does not
  correspond to any formula in the primitive vocabulary of $\RR$, nor
  do I know of any such formula that it implies which is not also a
  theorem of $\RR$.  So we are left with a nagging question.''
\end{quote}
It turns out that axiom \eqref{BIX}, $(\A\circ\B)^{-1} =
\B^{-1}\circ\A^{-1}$, corresponds to \eqref{reflection1}. The
predicates \eqref{reflection1} and \eqref{reflection1a} are in the
vocabulary of $\RR$ but they are not theorems of $\RR$ because they
fail in the frame $\gc\K_4$ in Table \ref{K3} whose complex algebra is
Dunn monoid that cannot be embedded in a relation algebra because it
fails to satisfy \eqref{left reflection}.  Consequently predicates
\eqref{reflection1} and \eqref{reflection1a} are in $\TT_3$ but not
$\RR$. Predicate \eqref{dedekind}, which uses the Routley star, is in
$\CT_3$ but not $\CR$.  \SS\ref{sect18} points out that
\eqref{reflection1} is a counterexample to \cite[Thm 8.1]{MR3289545}
and the complex algebra of $\gc\K_4$ is a counterexample to \cite[Thm
  7.1]{MR3289545}.  \SS\ref{sect19} presents the two examples of
predicates in the vocabulry of $\RR$ that are 5-provable but not
4-provable.  Asymptotic formulas for the numbers of $\TR$-frames and
$\KR$-frames are obtained in \SS\ref{sect20}. Their numbers grow like
$\c^{\n^3}$ for some $\c>1$. For any fixed $3\leq\n\in\omega$, the
probability that a randomly selected $\TR$-frame or $\KR$-frame
validates every $\n$-provable predicate approaches $1$ as the number
of elements in the frame grows.  Some questions are raised in
\SS\ref{sect21}.
\section{First-order logic $\LL$ of binary relations}
\label{sect2}
Tarski and Givant let $\LL$ be a first-order language with equality
symbol $\overset{\text{${}_\circ$}}{\mathbf1}$ and exactly one binary
relation symbol $\mathbf{E}$ \cite[p.\,5]{MR920815}, while $\Mn$ is a
first-order language with equality symbol $\overset{\text{${}_\circ$}}
{\mathbf1}$ and exactly $\n+1$ binary relation symbols, where
$\n\geq0$ \cite[p.\,191]{MR920815}.  They also consider formalisms
$\mathcal{M}$ and $\mathcal{M}^\times$ with any cardinality of binary
relation symbols \cite[p.\,237]{MR920815}.  Tarski and Givant used a
single relation symbol because they were presenting Tarski's
formalization of set theory without variables.  For set theory it is
usually sufficient to have just one relation symbol intended to denote
the relation of membership.  There is no need here for such a
restriction.  Changing notation and the number of relation symbols, we
assume instead that $\id$ is the {\bf equality symbol of $\LL$} and
that $\LL$ has a countable infinite set $\Pi$ of binary relation
symbols (including $\id$), but no function symbols and no constants.
(Tarski and Givant let $\pmb\Pi$ be what we call $\Pi^+$ in
\SS\ref{sect3}.)

The relation symbols in $\Pi$ are called {\bf atomic predicates}, and
those that are distinct from $\id$ are also called {\bf propositional
  variables} because they will play the r\^ole of variables in
formulas of relevance logic.  The {\bf connectives of $\LL$} are {\bf
  implication $\implies$} and {\bf negation $\neg$}, and $\forall$ is
the {\bf universal quantifier}.  $\LL$ has a countable set of {\bf
  variables} $\Upsilon=\{\var\i:\i\in\omega\}$, ordered in the natural
way so that $\var\i$ precedes $\var\j$ if $\i<\j$.  Thus, $\0$ and
$\var1$ are the first and second variables.  For every $\n\in\omega$,
let $\Upsilon_\n=\{\var\i:\i<\n\}$.  The {\bf atomic formulas of
  $\LL$} are the ones of the form $\x\A\y$, where $\x,\y\in\Upsilon$
are variables and $\A\in\Pi$ is an atomic predicate. For example,
$\x\id\y$ is an atomic formula since $\id\in\Pi$.  The set $\Phi$ of
{\bf formulas of $\LL$} is the intersection of every set that contains
the atomic formulas and includes $\varphi\implies\psi$, $\neg\varphi$,
and $\all\x\varphi$ for every variable $\x\in\Upsilon$ whenever it
contains $\varphi$ and $\psi$.  The set of {\bf sentences of $\LL$}
(formulas with no free variables) is $\Sigma$. The connectives $\lor$,
$\land$, and $\iff$, and the existential quantifier $\exists$ are
defined for all $\varphi,\psi\in\Phi$ by $\varphi\lor\psi =
\neg\varphi\implies\psi$, $\varphi\land\psi =
\neg(\varphi\implies\neg\psi)$, $\varphi\iff\psi = \neg((\varphi
\implies\psi) \implies \neg(\psi\implies\varphi))$, and $\ex\x\varphi
= \neg\all\x\neg\varphi$ for every $\x\in\Upsilon$. When a connective
is used more than once without parentheses, we restore them by
association to the left. For example, $\varphi\lor\psi\lor\xi=
(\varphi\lor\psi)\lor\xi$.  When parentheses are omitted from a
formula, the unary connective $\neg$ should be applied first, followed
by $\land$, $\lor$, $\implies$, and $\iff$, in that order.

In formulating axioms and deductive rules for $\LL$,
\citet[p.\,8]{MR920815} adopt the system $\mathcal\S_1$ of
\cite{Tarski1965}, which provides axioms for the logically valid
sentences and requires only the rule $\MP$ of {\bf modus ponens} (to
infer $\B$ from $\A\to\B$ and $\A$).  Tarski's system $\mathcal\S_2$
provides axioms for the logically valid formulas (not just the
sentences), and uses the rule of {\bf generalization} (to infer
$\all\x\varphi$ from $\varphi$) as well as $\MP$.  The systems
$\mathcal\S_1$ and $\mathcal\S_2$ in \cite{Tarski1965} were obtained
by modifying a system of \cite{MR0002508, Quine1951, MR0136536,
  MR695499} which also uses only $\MP$.  Tarski's systems avoid the
notion of substitution. \citet{Henkin1949, MR97c:03005} proved
G\"odel's completeness theorem for the case in which there are
relation symbols of arbitrary finite rank but no constants and no
function symbols.  He used $\MP$ and a restricted form of
generalization as rules of inference. \citet[Thms 1 and 5]{Tarski1965}
proved that his systems $\mathcal\S_1$ and $\mathcal\S_2$ are complete
by deriving Henkin's axioms and noting that both systems are
semantically sound.

For every formula $\varphi\in\Phi$, the {\bf closure} $\closure
\varphi$ of $\varphi$ is a sentence obtained by universally
quantifying $\varphi$ with respect to every free variable in
$\varphi$. The closure operator is determined by the following
conditions: $\closure\varphi=\varphi$ for every sentence
$\varphi\in\Sigma$, and if $\x$ is the last variable (in the ordering
of the variables) that occurs free in $\varphi$, then $\closure\varphi
= \closure{ \all\x \varphi}$.

The set $\Lambda[\LL]$ of {\bf logical axioms for $\LL$}, or simply
$\Lambda$, is the set of sentences that coincide with one of the
sentences \eqref{AI}--\eqref{AIX} shown in Table \ref{L-axioms}, where
$\varphi,\psi,\xi\in\Phi$.
\begin{table}
\begin{align}
  &\label{AI}\tag{AI} \index{AI@(AI)}
  \closure{(\varphi\implies\psi)\implies((\psi\implies\xi)\implies
    (\varphi\implies\xi))}\\
  &\label{AII}\tag{AII} \index{AII@(AII)}
  \closure{(\neg\varphi\implies\varphi)\implies\varphi}\\
  &\label{AIII}\tag{AIII} \index{AIII@(AIII)}
  \closure{\varphi\implies(\neg\varphi\implies \psi)}\\
  &\label{AIV}\tag{AIV} \index{AIV@(AIV)}
  \closure{\all\x\all\y\varphi\implies\all\y\all\x\varphi}\\
  &\label{AV}\tag{AV} \index{AV@(AV)}
  \closure{\all\x(\varphi\implies\psi)\implies(\all\x\varphi\implies
    \all\x\psi)}\\
  &\label{AVI}\tag{AVI} \index{AVI@(AVI)}
  \closure{\all\x\varphi\implies\varphi}\\
  &\label{AVII}\tag{AVII} \index{AVII@(AVII)}
  \closure{\varphi\implies\all\x\varphi}\text{ where
    $\x$ is not free in  $\varphi$}\\
  &\label{AVIII}\tag{AVIII} \index{AVIII@(AVIII)}
  \closure{\ex\x(\x\eqsym\y)}\text{ where $\x\neq\y$}\\
  &\label{AIX}\tag{AIX} \index{AIX@(AIX)}
  \closure{\x\eqsym\y\implies(\varphi\implies\psi)} \text{ where
    $\varphi$ is atomic, $\x$ occurs in $\varphi$, and $\psi$}\\
  &\notag\text{is obtained from $\varphi$ by replacing a single
    occurrence of $\x$ by $\y$}
\end{align}
\caption{Axioms for first-order logic $\mathcal{L}$, where
  $x,y\in\Upsilon$, $\varphi,\psi\in\Phi$.}
\label{L-axioms}
\end{table}
If $\Psi\subseteq\Sigma$ then a sentence $\varphi\in\Sigma$ is {\bf
  provable in $\LL$ from $\Psi$}, written $\Psi\proves\varphi$ or
$\proves\varphi$ if $\Psi=\emptyset$, if $\varphi$ is in every set
that is closed under $\MP$ and contains $\Psi\cup\Lambda$.  The {\bf
  theory generated by $\Psi$} in $\LL$ is
\begin{equation*}
  \Theta\eta\,\Psi=\{\varphi:\varphi\in\Sigma,\,\Psi\proves\varphi\}.
\end{equation*}
Two formulas $\varphi,\psi\in\Phi$ are {\bf logically equivalent} in
$\LL$, written $\varphi\equiv\psi$, if
$\proves\closure{\varphi\iff\psi}$.
\section{Extending $\LL$ to $\Lp$}\label{sect3} 
Tarski and Givant extend $\LL$ to $\Lp$ by adding a second equality
symbol $\equals$ and four operators $+$, $\min\blank$, $\rp$, and
$\conv{}$ that act on relation symbols and produce new relation
symbols.  The set $\Pi^+$ of {\bf predicates} of $\Lp$ is the
intersection of every set containing $\Pi$ that also contains $\A+\B$,
$\min\A$, $\A\rp\B$, and $\conv\A$ whenever it contains $\A$ and $\B$.
Predicates obtained in distinct ways are distinct, so, for example, if
$\A+\B=\C+\D$ then $\A=\C$ and $\B=\D$.  Three predicates in $\Pi^+$
are defined by
\begin{align}\label{10di}
  1&=\id+\min\id,& 0&=\min{\id+\min\id}, &\di&=\min\id,
\end{align}
and two additional predicate operators are defined for all
$\A,\B\in\Pi^+$ by
\begin{align}\label{dotmin}
  \A\cdot\B&=\min{\min\A+\min\B},&\A\rs\B&=\min{\min\A\rp\min\B}.
\end{align}
When parentheses are omitted, the unary operators should be evaluated
first, followed by $\rp$, $\cdot$, $\rs$, and then $+$, in that order.
For example, $\A\rs\B+\C\rp\D\cdot\E =(\A\rs\B)+((\C\rp\D)\cdot\E)$.
Tarski and Givant add a formula $\A\equals\B$, called an {\bf
  equation}, for any predicates $\A,\B\in\Pi^+$.  The set of {\bf
  equations of $\Lp$} is $\Sigma^\times$.  For all $\A,\B\in\Pi^+$,
let $\A\leq\B$ be the equation $\A+\B\equals\B \in \Sigma^\times$
\cite[p.\,236]{MR920815}, which we call an {\bf inclusion}.

Our notation for the equality symbol is derived from
\citet{MR0192999}, who denoted the identity relation on a set by
$\id$.  Schr\"oder obtained his notation from the Boolean unit $1$ by
adding an apostrophe.  Similarly, Schr\"oder added a comma to the
symbol $\cdot$ for intersection to obtain his symbol $\rp$ for
relative product. In \cite{MR920815}, Tarski altered Schr\"oder's
system by using a circle instead of an apostrophe or comma, as well as
making many symbols boldface. For example, instead of $\rp$ he used
$\pmb\odot$.  Tarski and Givant originally used a boldface equality
symbol instead of $\equals$, but $\pmb=$ is not as easily
distinguished from the usual equality symbol $=$ as is $\equals$. The
notation for $\equals$ used here was inspired by Tarski's device of
adding circles to Boolean notation.

The {\bf atomic formulas of $\Lp$} are $\x\A\y$ and $\A\equals\B$,
where $\x,\y\in\Upsilon$ and $\A,\B\in\Pi^+$.  The set $\Phi^+$ of
{\bf formulas of $\Lp$} is the intersection of every set containing
the atomic formulas of $\Lp$ that also contains $\varphi\implies\psi$,
$\neg\varphi$, and $\all\x\varphi$ for every $\x\in\Upsilon$ whenever
it contains $\varphi$ and $\psi$.  The set $\Sigma^+$ of {\bf
  sentences of $\Lp$} is the set of formulas that have no free
variables.  Equations have no free variables so $\Sigma^\times
\subseteq \Sigma^+$.

The set $\Lambda[\Lp]$ of {\bf logical axioms of $\Lp$}, or simply
$\Lambda^+$ \cite[p.\,25]{MR920815}, is the union of $\Lambda[\LL]$
with the set of sentences that coincide with one of the sentences in
Table \ref{L+defs} for some $\A,\B\in\Pi^+$.
\begin{table}
  \begin{align}
    \label{DI}\tag{DI}
    \closure{\0\A+\B\var1&\iff\0\A\var1\lor\0\B\var1}\\
    \label{DII}\tag{DII}
    \closure{\0\min\A\var1&\iff \neg\0\A\var1}\\
    \label{DIII}\tag{DIII}
    \closure{\0\A\rp\B\var1&\iff\ex\z(\0\A\z\land\z\B\var1)}\\
    \label{DIV}\tag{DIV}
    \closure{\0\conv\A\var1&\iff\var1\A\0}\\
    \label{DV}\tag{DV}
    \A\equals\B&\iff\closure{\0\A\var1\iff\0\B\var1}
  \end{align}
  \caption{Definitional axioms for the extension $\Lp$ of $\LL$,
    where $A,B\in\Pi^+$.}
  \label{L+defs}
\end{table}
If $\Psi\subseteq\Sigma^+$, then a sentence $\varphi\in\Sigma^+$ is
{\bf provable in $\Lp$ from $\Psi$}, written $\Psi\proves^+\varphi$ or
$\proves^+\varphi$ if $\Psi=\emptyset$, if $\varphi$ is in every set
that contains $\Psi\cup\Lambda^+$ and is closed under $\MP$. The {\bf
  theory generated by $\Psi$} in $\Lp$ is
\begin{equation*}
  \Theta\eta^+\Psi=\{\varphi:\varphi\in\Sigma^+,\,\Psi\proves\varphi\}.
\end{equation*}
Two formulas $\varphi, \psi\in \Pi^+$ of $\Lp$ are {\bf logically
  equivalent} in $\Lp$, written $\varphi\equiv^+\psi$, if $\proves^+
\closure{\varphi\iff\psi}$.  The {\bf calculus of relations} may be
defined as $\Theta\eta^+\emptyset$. One may also consider it to be the
closure of $\Theta\eta^+\emptyset$ under the connectives $\neg$ and
$\implies$, since Schr\"oder and Tarski showed that every
propositional combination of equations is logically equivalent to an
equation \cite[2.2(vi)]{MR920815}.

\section{Equipollence of $\LL$ and $\Lp$}\label{sect4}
$\LL$ and $\Lp$ are expressively and deductively equipollent. To prove
this, Tarski defined a translation mapping $\GG$ from formulas of
$\Lp$ to formulas of $\LL$.  See \cite[2.3(iii)]{MR920815} for the
definition of $\GG$ and \cite[2.4(iii)]{MR920815} for the definition
of {\bf translation mapping} from one formalism to another.  $\GG$
eliminates operators in accordance with the definitional axioms
\eqref{DI}--\eqref{DV}. If $\varphi,\psi\in\Phi^+$,
$\x,\y\in\Upsilon$, and $\A,\B\in\Pi^+$, then the conditions
determining $\GG$ are shown in Table \ref{G-defs}.
\begin{table}
\begin{align*}
  \GG(\x\A\y)&=\x\A\y\quad\text{if $\A\in\Pi$}\\
  \GG(\varphi\implies\psi)&=\GG(\varphi)\implies\GG(\psi)\\
  \GG(\neg\varphi)&=\neg\GG(\varphi)\\
  \GG(\all\x\varphi)&=\all\x\GG(\varphi)\\
  \GG(\x\A+\B\y)&=\GG(\x\A\y)\lor\GG(\x\B\y)\\
  \GG(\x\min\A\y)&=\neg\GG(\x\A\y)\\
  \GG(\x\A\rp\B\y)&=\ex\z(\GG(\x\A\z)\land\GG(\z\B\y))\\
  &\text{ where $\z$ is the first variable distinct from $\x$ and $\y$}\\
  \GG(\x\conv\A\y)&=\GG(\y\A\x)\\
  \GG(\A\equals\B)&=\closure{\GG(\0\A\var1)\iff\GG(\0\B\var1)}
\end{align*}
\caption{Definition of translation mapping $\GG\colon\Phi^+\to\Phi$,
  where $x,y,z\in\Upsilon$, $A,B\in\Pi^+$, and $\varphi,\psi\in\Phi^+$.}
\label{G-defs}
\end{table}
From the first four conditions it follows that $\GG$ leaves formulas
of $\LL$ unchanged.  The next result states that $\LL$ is a
subformalism of $\Lp$ and $\LL$ is expressively and deductively
equipollent with $\Lp$. Part \eqref{thm1iv} is the {\bf main mapping
  theorem for $\LL$ and $\Lp$}.
\begin{theorem}\label{thm1}{\rm\cite[\SS2.3]{MR920815}}
  Formalisms $\LL$ and $\Lp$ are equipollent.
  \begin{enumerate}
  \item\label{thm1i} $\Phi\subseteq\Phi^+$ and
    $\Sigma\subseteq\Sigma^+$ {\rm[2.3(i)]}.
  \item\label{thm1ii} $\GG$ maps $\Phi^+$ onto $\Phi$ and $\Sigma^+$
    onto $\Sigma$ {\rm[2.3(iv)($\delta$)]}.
  \item\label{thm1iii} $\varphi\equiv^+\GG(\varphi)$ if
    $\varphi\in\Phi^+$ {\rm[2.3(iv)($\varepsilon$)]}.
  \item\label{thm1iv} $\Psi\proves^+\varphi$ iff
    $\{\GG(\psi):\psi\in\Psi\} \proves \GG(\varphi)$ if
    $\Psi\subseteq\Sigma^+$ and $\varphi\in\Sigma^+$ {\rm[2.3(v)]}.
  \item\label{thm1v} $\Psi\proves^+\varphi$ iff
    $\Psi\proves\varphi$, if $\Psi\subseteq\Sigma$ and
    $\varphi\in\Sigma$ {\rm[2.3(ii)(ix)]}.
  \item\label{thm1vi} $\Theta\eta\,\Psi=\Theta\eta^+\Psi\cap\Sigma$ if
    $\Psi\subseteq\Sigma$ {\rm[2.3(x)]}.
  \end{enumerate}
\end{theorem}
\section{Equational formalisms $\Lx$ and $\Lwx$}
\label{sect5}
The equational formalisms $\Lx$ and $\Lwx$ are defined by
\citet[\SS3.1 and p.\,89]{MR920815}. $\Lx$ is the primary subject of
their book but $\Lwx$ makes only an incidental appearance as a
weakening of $\Lx$.  $\Lx$ is equipollent with the 3-variable
formalisms $\LL_3$ and $\Lp_3$ described in \SS\ref{sect7} while
$\Lwx$ is equipollent with the ``(\emph{standardized})
\emph{formalisms}'' $\Ls_3$ and $\Ls_3^+$ described in \SS\ref{sect9}.
Tarski and Givant said, ``These standardized formalisms are
undoubtedly more natural and more interesting in their own right than
$\LL_3$ and $\Lp_3$'' \cite[p.\,89]{MR920815}.

The axioms of $\Lx$ and $\Lwx$ are certain equations in
$\Sigma^\times$ and their deductive rules apply to equations.
$\Lambda[\Lx]$, or simply $\Lambda^\times$, is the set of {\bf axioms
  of $\Lx$} and $\Lambda[\Lwx]$, or simply $\Lambda^\times_s$, is the
set of {\bf axioms of $\Lwx$}.  An equation $\varepsilon\in
\Sigma^\times$ belongs to $\Lambda^\times$ if there are predicates
$\A,\B,\C\in\Pi^+$ such that $\varepsilon$ coincides with one of the
equations \eqref{BI}--\eqref{BX} listed in Table \ref{Lx-axioms}, and
$\varepsilon$ belongs to $\Lambda^\times_s$ if $\varepsilon$ coincides
with one of the equations \eqref{BI}--\eqref{BIII}, \eqref{BIV'},
\eqref{BV}--\eqref{BX} in Table \ref{Lx-axioms}.
\begin{table}
  \begin{align*}
    \label{BI}   \tag{\ra1}\A+\B&\equals\B+\A\\
    \label{BII}  \tag{\ra2}\A+(\B+\C)&\equals(\A+\B)+\C\\
    \label{BIII} \tag{\ra3}\min{\min\A+\min\B}+\min{\min\A+\B}
    &\equals\A\\
    \label{BIV}  \tag{\ra4}\A\rp(\B\rp\C)&\equals(\A\rp\B)\rp\C\\
    \label{BIV'} \tag{\ra4$'$} \A\rp(\B\rp1)&\equals(\A\rp\B)\rp1\\
    \label{BV}   \tag{\ra5}(\A+\B)\rp\C&\equals\A\rp\C+\B\rp\C\\
    \label{BVI}  \tag{\ra6}\A\rp\id&\equals\A\\
    \label{BVII} \tag{\ra7}\conv{\conv\A}&\equals\A\\
    \label{BVIII}\tag{\ra8}\conv{\A+\B}&\equals\conv\A+\conv\B\\
    \label{BIX}  \tag{\ra9}\conv{\A\rp\B}&\equals\conv\B\rp\conv\A\\
    \label{BX} \tag{\ra{10}}\conv\A\rp\min{\A\rp\B}+\min\B
    &\equals\min\B
  \end{align*}
  \caption{Axioms for the equational formalisms $\Lx$ and $\Lwx$, 
    where $A,B,C\in\Pi^+$.}
  \label{Lx-axioms}
\end{table}
Note that $\Lambda^\times_s\subseteq \Lambda^\times$ because
$1\in\Pi^+$, hence every instance of \eqref{BIV'} is also an instance
of \eqref{BIV}. Deducibility in $\Lx$ and $\Lwx$ is defined as it is
in equational logic.  The {\bf transitivity rule} \tran\ is to infer
$\B\equals\C$ from $\A\equals\B$ and $\A\equals\C$, and the {\bf
  replacement rule} \repl\ is to infer $\min\A\equals\min\B$,
$\conv\A\equals\conv\B$, $\A+\C\equals\B+\C$, and
$\A\rp\C\equals\B\rp\C$ from $\A\equals\B$.  For every
$\Psi\subseteq\Sigma^\times$, an equation $\varepsilon\in
\Sigma^\times$ is {\bf provable in $\Lx$ from $\Psi$}, written
$\Psi\proves^\times\varepsilon$ or $\proves^\times\varepsilon$ when
$\Psi=\emptyset$, iff $\varepsilon$ is in every set that contains
$\Psi\cup\Lambda^\times$ and is closed under \tran\ and \repl.
Similarly, $\varepsilon$ is {\bf provable in $\Lwx$ from $\Psi$},
written $\Psi\proves^\times_s\varepsilon$ or $\proves^\times_s
\varepsilon$ when $\Psi=\emptyset$, iff $\varepsilon$ belongs to
every set that contains $\Psi\cup\Lambda^\times_s$ and is closed under
\tran\ and \repl.  For every $\Psi\subseteq\Sigma^\times$, the {\bf
  theory generated by $\Psi$} in $\Lx$ is
\begin{equation*}
  \Theta\eta^\times\Psi=\{\varepsilon:\varepsilon\in\Sigma^\times,\,
  \Psi\proves^\times\varepsilon\},
\end{equation*}
and the {\bf theory generated by $\Psi$} in $\Lwx$ is
\begin{equation*}
  \Theta\eta^\times_s\Psi=\{\varphi:\varepsilon\in\Sigma^\times,\,
  \Psi\proves^\times_s\varepsilon\}.
\end{equation*}
The rules stated here employ simplifications (mentioned but not proved
by \citet[p.\,47]{MR920815}) made possible by the presence of certain
equations in $\Lambda^\times_s$.  The equation $\A\equals\A$ is
deducible in $\Lx$ and $\Lwx$ for every predicate $\A\in\Pi^+$ because
it follows by \tran\ from two instances of \eqref{BIII}, \eqref{BVI},
or \eqref{BVII}.  To derive $\B\equals\A$ from $\A\equals\B$, first
derive $\A\equals\A$ using \tran\ and one of \eqref{BIII},
\eqref{BVI}, or \eqref{BVII} and then apply \tran\ to $\A\equals\B$
and $\A\equals\A$.  To derive $\A\equals\C$ from $\A\equals\B$ and
$\B\equals\C$, first derive $\B\equals\A$ from $\A\equals\B$ and apply
\tran\ to $\B\equals\A$ and $\B\equals\C$.

For equational logic in general the replacement rule would include the
equations $\C+\A\equals\C+\B$ and $\C\rp\A\equals\C\rp\B$ as equations
derivable from $\A\equals\B$, but they can be derived.  To get
$\C\rp\A\equals\C\rp\B$ from $\A\equals\B$, first derive
$\conv\A\rp\conv\C\equals \conv\B\rp\conv\C$ by applying \repl\ twice.
Two instances of \eqref{BIX} are $\conv{\C\rp\A} \equals
\conv\A\rp\conv\C$ and $\conv{\C\rp\B} \equals \conv\B\rp\conv\C$.
Use the laws of equality proved above to get $\conv{\C\rp\A} \equals
\conv{\C\rp\B}$ from these last three equations.  Next, obtain
$\conv{\conv{\C\rp\A}}\equals\conv{\conv{\C\rp\B}}$ by \repl\ and
complete the proof using two instances of \eqref{BVII} and the laws of
equality.  Thus, the presence of \eqref{BVII} and \eqref{BIX} in
$\Lambda^\times$ and $\Lambda^\times_s$ is enough to make
$\C\rp\A\equals\C\rp\B$ derivable from $\A\equals\B$.  It is easier to
derive $\C+\A\equals\C+\B$ from $\A\equals\B$ using \eqref{BI} and one
of \eqref{BIII}, \eqref{BVI}, or \eqref{BVII}.
\begin{theorem}{\rm\cite[\SS3.4]{MR920815}}\label{thm1a}
  $\Lwx$ is a subformalism of $\Lx$ and $\Lx$ is a subformalism of
  $\Lp$.
  \begin{enumerate}
  \item\label{thm1ai} $\Sigma^\times \subseteq\Sigma^+$ {\rm[3.4(i)]}.
  \item\label{thm1aii} if $\Psi\proves^\times_s\varepsilon$ then
    $\Psi\proves^\times\varepsilon$, for every
    $\Psi\subseteq\Sigma^\times$ and $\varepsilon\in\Sigma^\times$.
  \item\label{thm1aiii}
    $\Theta\eta^\times_s\Psi\subseteq\Theta\eta^\times\Psi$ for every
    $\Psi\subseteq\Sigma^\times$.
  \item\label{thm1aiv} if $\Psi\proves^\times\varepsilon$ then
    $\Psi\proves^+\varepsilon$, for every $\Psi\subseteq\Sigma^\times$
    and $\varepsilon\in\Sigma^\times$ {\rm[3.4(ii)]}.
  \item\label{thm1av}
    $\Theta\eta^\times\Psi\subseteq\Theta\eta^+\Psi\cap\Sigma^\times$
    for every $\Psi\subseteq\Sigma^\times${\rm[3.4(vii)]}.
  \end{enumerate}
\end{theorem}
\proof Part \eqref{thm1ai} is the observation made in the previous
section that equations have no free variables.  Part \eqref{thm1aii}
follows from $\Lambda^\times_s\subseteq\Lambda^\times$ and part
\eqref{thm1aiii} follows from part \eqref{thm1aii}.  Part
\eqref{thm1aiv} can be proved by induction on provability in $\Lp$.
One shows that the axioms of $\Lx$ are provable in $\Lp$ and that if
the hypotheses of \tran\ or \repl\ are provable in $\Lp$ then so are
their conclusions.  This would be tedious to carry out according to
the definitions of the notions involved.  For example, if $\varphi$
were $\A+\B\equals\B+\A$, an instance of axiom \eqref{AI} where
$\A,\B\in\Pi^+$, one would have to provide a sequence of sentences in
$\Sigma^+$, each of which is either an instance of
\eqref{AI}--\eqref{AIX} or \eqref{DI}--\eqref{DV} or follows from two
previous sentences by $\MP$, ending with $\A+\B\equals\B+\A$. It is
much easier to proceed semantically, taking advantage of G\"odel's
completeness theorem for $\LL$ (see \cite[\SS1.4]{MR920815}) and its
implications for $\Lp$ (see \cite[\SS2.2]{MR920815}). It then becomes
clear that \eqref{AI} expresses the fact that the operation of forming
the union of two binary relations is commutative and that this fact
can be proved in $\Lp$.  Similarly, all the other axioms of $\Lx$ can
be seen as logically valid (and therefore provable in $\Lp$) when they
are interpreted according to \eqref{DI}--\eqref{DV}.  \endproof
\section{Relation algebras, semi-associative and representable}
\label{sect6}
Since $\Pi^+$ is closed under the predicate operators $+$,
$\min\blank$, $\rp$, $\conv{}$, and contains $\id$, we may define the
{\bf predicate algebra} $\gc\P$ by
\begin{align*}
  \gc\P&=\<\Pi^+,+,\min\blank,\rp,\conv{},\id\>.
\end{align*}
Then $\gc\P$ is an absolutely free algebra that is freely generated by
the propositional variables $\{\A:\id\neq\A\in\Pi\}$
\cite[p.\ 238]{MR920815}. This means that any function mapping the
propositional variables into an algebra of the same similarity type as
$\gc\P$ has a unique extension to a homomorphism from $\gc\P$ into
that algebra.  Consider an algebra $\gc\A$ having the same similarity
type as $\gc\P$, say
\begin{equation*}
  \gc\A=\<\U,+,\min\blank,\rp,\conv{},\id\>,
\end{equation*}
where $\U$ is a set called the {\bf universe} of $\gc\A$, $\id\in\U$,
$+$ and $\rp$ are binary operations on $\U$, and $\min\blank$ and
$\conv{}$ are unary operations on $\U$.  Define three additional
elements of $\U$ by $\di=\min\id$, $1=\id+\di$, and $0=\min1$. Define
the binary operation $\cdot$ on $\U$ by $\x\cdot\y =
\min{\min\x+\min\y}$ and the binary relation $\leq$ on $\U$ by
$\x\leq\y$ if $\x+\y=\y$.  The algebra $\gc\A$ is {\bf dense} if
$\x\leq\x\rp\x$ for every $\x\in\U$, {\bf commutative} if
$\x\rp\y=\y\rp\x$ for all $\x,\y\in\U$, {\bf symmetric} if
$\conv\x=\x$ for every $\x\in\U$, and {\bf integral} if $0\neq1$ and
$\x\rp\y=0$ imply $\x=0$ or $\y=0$, for all $\x,\y\in\U$.  An element
$\x\in\U$ is an {\bf atom} of $\gc\A$ if $\x\neq0$ and for every
$\y\in\U$ either $\x\leq\y$ or $\y=0$, and $\atoms$ is the set of
atoms of $\gc\A$.  The algebra $\gc\A$ {\bf atomic} if $0\neq\y\in\U$
implies $\x\leq\y$ for some atom $\x\in\atoms$.  The algebra $\gc\A$
is {\bf simple} if $0\neq1$ and has no non-trivial homomorphic images,
meaning that every homomorphic image of $\gc\A$ is either a 1-element
algebra or is isomorphic to $\gc\A$.  The algebra $\gc\A$ is {\bf
  semi-simple} if it is isomorphic to a subdirect product of simple
homomorphic images of $\gc\A$.  It follows from the definition of
subdirect product that every semi-simple algebra is isomorphic to a
subalgebra of a direct product of simple homomorphic images of
$\gc\A$.

For every homomorphism $\h\colon \gc\P\to\gc\A$ that maps the
predicate algebra $\gc\P$ into $\gc\A$, let $\models_\h$ be the
relation that holds between the algebra $\gc\A$ and a set of equations
$\Psi\subseteq\Sigma^\times$ if $\h(\A)=\h(\B)$ whenever
$\A,\B\in\Pi^+$ and $\A\equals\B\in\Psi$.  For any
$\varepsilon\in\Sigma^\times$ let $\gc\A\models_\h\varepsilon$ mean
the same as $\gc\A\models_\h\{\varepsilon\}$.  An equation
$\A\equals\B\in\Sigma^\times$ is {\bf true in $\gc\A$} if
$\gc\A\models_\h\A\equals\B$ for every homomorphism
$\h\colon\gc\P\to\gc\A$ from the predicate algebra $\gc\P$ to $\gc\A$.
For example, if $\id\neq\A,\B\in\Pi$, then $\gc\A$ is commutative iff
the equation $\A\rp\B\equals\B\rp\A$ is true in $\gc\A$, dense iff the
inclusion $\A\leq\A\rp\A$ is true in $\gc\A$, and symmetric iff
$\conv\A\equals\A$ is true in $\gc\A$. The way these equivalences are
established is illustrated by the proof of Lemma \ref{SA:symm->comm}
below.

The algebra $\gc\A$ is a {\bf relation algebra} if the equations
\eqref{BI}--\eqref{BX} in Table \ref{Lx-axioms} are true in $\gc\A$.
$\gc\A$ is a {\bf semi-associative relation algebra} if the equations
\eqref{BI}--\eqref{BIII}, \eqref{BIV'}, and \eqref{BV}--\eqref{BX} are
true in $\gc\A$.  $\gc\A$ is a {\bf non-associative relation algebra}
if the equations \eqref{BI}--\eqref{BIII} and \eqref{BV}--\eqref{BX}
are true in $\gc\A$.  Let $\RA$ be the class of relation algebras, let
$\SA$ be the class of semi-associative relation algebras, and let
$\NA$ be the class of non-associative relation algebras.  It follows
immediately from their definitions that $\NA\subseteq\SA\subseteq\RA$.
\begin{lemma}\label{SA:symm->comm}
  Assume $\gc\A=\<\U,+,\min\blank,\rp,\conv{},\id\>$ is a symmetric
  algebra in which \eqref{BIX} is true. Then $\gc\A$ is commutative.
  In particular, every symmetric semi-associative relation algebra is
  commutative.
\end{lemma}
\proof Assume $\gc\A=\<\U,+,\min\blank,\rp,\conv{},\id\>$, $\gc\A$ is
symmetric, and \eqref{BIX} is true in $\gc\A$.  Let $\x,\y\in\U$ and
let $\id\neq\A,\B\in\Pi$ be distinct propositional variables.  Since
$\gc\P$ is absolutely freely generated by the propositional variables,
there is a homomorphism $\h\colon\gc\P\to\gc\A$ such that $\h(\A)=\x$
and $\h(\B)=\y$.  Since \eqref{BIX} is true in $\gc\A$,
$\h(\conv{\A\rp\B})=\h(\conv\B\rp\conv\A)$, hence
\begin{align*}
  \x\rp\y&=(\x\rp\y)\conv{}&&\text{$\gc\A$ is symmetric}
  \\&=\conv{\h(\A)\rp\h(\B)}&&\text{choice of $\h$}
  \\&=\h(\conv{\A\rp\B})&&\text{$\h$ is a homomorphism}
  \\&=\h(\conv\B\rp\conv\A)&&\text{\eqref{BIX} is true in $\gc\A$}
  \\&=\conv{\h(\B)}\rp\conv{\h(\A)}&&\text{$\h$ is a homomorphism}
  \\&=\conv\y\rp\conv\x&&\text{choice of $\h$}
  \\&=\y\rp\x&&\text{$\gc\A$ is symmetric}.
\end{align*}
\endproof \citet[Th.\ 4.15]{MR0045086} proved that every relation
algebra is semi-simple; see \cite[Thm 12.10]{MR3699801}. By \cite[Cor
  8(7)]{MR2628352} or \cite[Thm 388]{MR2269199} it is also true that
every semi-associative relation algebra is semi-simple.  By \cite[Thm
  379(iii)]{MR2269199}, \cite[Thm 7(20)]{MR2628352}, or \cite[Thm
  29]{MR1049616}, a semi-associative relation algebra
$\gc\A=\<\U,+,\min\blank,\rp,\conv{},\id\>$ is simple iff $0\neq1$ and
for all $\x,\y\in\U$, if $(\x\rp1)\rp\y=0$ then $\x=0$ or $\y=0$.
\citet[4.17]{MR0045086} proved that a relation algebra is integral iff
its identity element is an atom.  This result also extends to $\SA$.
By \cite[Thm 353]{MR2269199} or \cite[Thm 4]{MR1387901} a
semi-associative relation algebra $\gc\A$ is integral iff $\id$ is an
atom of $\gc\A$.  These facts and some basic observations from
universal algebra are used in the proof of the following lemma.
\begin{lemma}\label{commSA}\ 
 \begin{enumerate}
  \item\label{commSAi} Assume
    $\gc\A=\<\U,+,\min\blank,\rp,\conv{},\id\>$ is a simple
    commutative semi-associative relation algebra.  Then $\gc\A$ is
    integral and $\id$ is an atom of $\gc\A$.
  \item\label{commSAii} An equation is true in every commutative
    semi-associative relation algebra iff it is true in every
    integral commutative semi-associative relation algebra.
  \item\label{commSAiii} An equation is true in every commutative
    semi-associative relation algebra iff it is true in every
    commutative semi-associative relation algebra in which the
    identity element is an atom.
 \end{enumerate}
\end{lemma}
\proof For part \eqref{commSAi}, suppose $\gc\A\in\SA$ is commutative
and simple.  First note that $0\neq1$ since $\gc\A$ is simple.  To
show $\gc\A$ is integral, suppose that $\x\rp\y=0$.  We have
$(\x\rp1)\rp\y= (1\rp\x)\rp\y$ since $\gc\A$ is commutative,
$(1\rp\x)\rp\y= 1\rp(\x\rp\y)=1\rp0$ by \cite[Thm 13]{MR1049616} or
\cite[Thm 354]{MR2269199} and the assumption that $\x\rp\y=0$, and
$1\rp0=0$ by \cite[Thm 287]{MR2269199}, so $(\x\rp1)\rp\y=0$.  Since
$\gc\A$ is a simple semi-associative relation algebra, it follows, as
noted above, that either $\x=0$ or $\y=0$.  This shows $\gc\A$ is
integral, so we conclude that $\id$ is an atom of $\gc\A$, also noted
above.

Parts \eqref{commSAii} and \eqref{commSAiii} are equivalent because a
semi-associative relation algebra is integral iff its identity element
is an atom. One direction of each part is also trivially true.  It
suffices therefore to assume that an equation $\varepsilon$ is true in
every integral commutative semi-associative relation algebra and show
that it is true in every commutative semi-assciative relation algebra.

Suppose $\gc\A$ is a commutative semi-associative relation algebra.
As noted above, $\gc\A$ is isomorphic to a subalgebra of a direct
product of simple semi-associative relation algebras that are
homomorphic images of $\gc\A$. Homomorphic images of commutative
algebras are commutative, so all these simple semi-associative
relation algebras are commutative. By part \eqref{commSAi}, they are
also integral, so by hypothesis the equation $\varepsilon$ is true in
all of them. If an equation is true in a collection of algebras, then
it is also true in their direct product, and if it is true in an
algebra then it is also true in all the subalgebras of that
algebra. These two facts combine to show that $\varepsilon$ must
therefore be true in $\gc\A$, as desired.  \endproof The next theorem
provides a link between relation algebras and deducibility in $\Lx$
and between semi-associative relation algebras and deducibility in
$\Lwx$.  Following their proof of \cite[8.2(x)]{MR920815}, which is
part \eqref{eq3}, Tarski and Givant say, ``It may be noticed that the
theorem just proved could be given a stronger form by using the notion
of a free algebra with defining relations. \dots\ A precise
formulation of the improved Theorem (x) would be rather involved, and
we leave it to the reader.''  This stronger form is part \eqref{eq1}.

For every $\Psi\subseteq\Sigma^\times$, let $\simeq^\times_\Psi$ and
$\simeq^s_\Psi$ be the binary relations defined for any
$\A,\B\in\Pi^+$ by $\A\simeq^\times_\Psi\B$ iff
$\Psi\proves^\times\A\equals\B$ \cite[p.\,238]{MR920815}, and
$\A\simeq^s_\Psi\B$ iff $\Psi\proves^\times_s\A\equals\B$.  For every
$\Psi\subseteq\Sigma^+$, let $\simeq^+_\Psi$ be the binary relation
defined for any $\A,\B\in\Pi^+$ by $\A\simeq^+_\Psi\B$ iff
$\Psi\proves^+\A\equals\B$ \cite[p.\,240]{MR920815}.  In all three
relations, reference to $\Psi$ is omitted when $\Psi=\emptyset$.
Theorem \ref{eq} concerns $\simeq^\times_\Psi$ and $\simeq^s_\Psi$,
while Theorem \ref{sem} deals with $\simeq^+_\Psi$.
\begin{theorem}\label{eq} Assume $\Psi\subseteq\Sigma^\times$. 
  \begin{enumerate}
  \item\label{eq5} $\simeq^\times_\Psi$ is a congruence relation on
    $\gc\P$ and the quotient algebra $\gc\P/{\simeq^\times_\Psi}$ is a
    relation algebra {\rm\cite[8.2(ix)]{MR920815}}.
  \item\label{eq6} $\simeq^s_\Psi$ is a congruence relation on
    $\gc\P$ and the quotient algebra $\gc\P/{\simeq^s_\Psi}$ is a
    semi-associative relation algebra.
  \item\label{eq1} For every $\varepsilon \in \Sigma^\times$,
    $\Psi\proves^\times\varepsilon$ iff for every $\gc\A\in\RA$ and
    every homomorphism $\h\colon\gc\P\to\gc\A$, if
    $\gc\A\models_\h\Psi$ then $\gc\A\models_\h\varepsilon$.
  \item\label{eq2} For every $\varepsilon \in \Sigma^\times$,
    $\Psi\proves^\times_s\varepsilon$ iff for every $\gc\A\in\SA$ and
    every homomorphism $\h\colon\gc\P\to\gc\A$, if
    $\gc\A\models_\h\Psi$ then $\gc\A\models_\h\varepsilon$.
  \item\label{eq3} $\gc\P/{\simeq^\times}$ is a relation algebra that is
    $\RA$-freely generated by $\{\A/{\simeq^\times}:\id\neq\A\in\Pi\}$
    {\rm\cite[8.2(x)]{MR920815}}.
  \item\label{eq4} $\gc\P/{\simeq^\times}$ is a semi-associative relation
    algebra that is $\SA$-freely generated by
    $\{\A/{\simeq^s}:\id\neq\A\in\Pi\}$.
  \end{enumerate}
\end{theorem}
\proof For part \eqref{eq5}, note that $\simeq^\times_\Psi$ is a
congruence relation on $\gc\P$ because of the rules \tran\ and \repl\
and their consequences.  It follows that $\simeq^\times_\Psi$
determines a quotient homomorphism $\q\colon \Pi^+\to\{\A/
{\simeq^\times_\Psi}: \A\in\Pi^+\}$ that carries each predicate
$\A\in\Pi^+$ to its equivalence class $\A/{\simeq^\times_\Psi}$ under
$\simeq^\times_\Psi$.  To show that the quotient algebra
$\gc\P/{\simeq^\times_\Psi}$ is in $\RA$, we must prove that every
axiom in $\Lambda^\times$ is true in $\gc\P/{\simeq^\times_\Psi}$. For
that we assume $\h\colon\gc\P\to\gc\P/{\simeq^\times_\Psi}$ is a
homomorphism.  We must show that if $\lambda\in\Lambda^\times$ then
$\gc\P/{\simeq^\times_\Psi}\models_\h\lambda$.  We do just one
example, say an instance $\A+\B\equals\B+\A$ of \eqref{BI}.  Since
$\h(\A)$ and $\h(\B)$ are equivalence classes of predicates, we may
choose $\C,\D\in\Pi^+$ such that $\h(\A)=\C/{\simeq^\times_\Psi} =
\q(\C)$ and $\h(\B)=\D/{\simeq^\times_\Psi} = \q(\D)$.  Then
\begin{align*}
  \h(\A+\B) &=\h(\A)+\h(\B)&&\text{$\h$ is a homomorphism}
  \\&=\q(\C)+\q(\D)&&\text{choice of $\C,\D$}
  \\&=\q(\C+\D)&&\text{$\q$ is a homomorphism}
  \\&=\q(\D+\C)&&\Psi\proves^\times\C+\D\equals\D+\C
  \\&=\q(\D)+\q(\C)&&\text{$q$ is a homomorphism}
  \\&=\h(\B)+\h(\A)&&\text{choice of $\C,\D$}
  \\&=\h(\B+\A)&&\text{$\h$ is a homomorphism}
\end{align*}
Proofs for the other axioms are similar. Thus,
$\gc\P/{\simeq^\times_\Psi}\in\RA$.  Furthermore,
$\gc\P/{\simeq^\times_\Psi}\models_\q\Psi$ just by the definitions of
$\simeq^\times_\Psi$ and $\q$.  The proof of part \eqref{eq6} is the
same, but with $\simeq^\times_\Psi$, $\RA$, and $\Lambda^\times$
replaced by $\simeq^s_\Psi$, $\SA$, and $\Lambda^\times_s$,
respectively.

For part \eqref{eq1}, assume $\Psi\proves^\times\varepsilon$,
$\gc\A\in\RA$, $\h\colon\gc\P\to\gc\A$ is a homomorphism, and
$\gc\A\models_\h\Psi$.  We wish to show $\gc\A\models_\h\varepsilon$.
Let $\Theta=\{\theta: \theta\in\Sigma^\times,\,
\gc\A\models_\h\theta\}$.  We will show
$\Theta\eta^\times\Psi\subseteq\Theta$.  We have
$\Lambda^\times\subseteq\Theta$ by the definition of $\RA$ and
$\Psi\subseteq\Theta$ by the hypothesis $\gc\A\models_\h\Psi$ and the
definition of $\models_\h$.  Next, we show that $\Theta$ is closed
under the rules \tran\ and \repl.  To see this for \tran, assume
$\A\equals\B,\A\equals\C\in\Theta$.  Then $\h(\A)=\h(\B)$ and
$\h(\A)=\h(\C)$ by the definition of $\Theta$.  It follows that
$\h(\B)=\h(\C)$, so $\B\equals\C\in\Theta$.  For \repl, we assume
$\A\equals\B\in\Theta$ and wish to show the conclusions of \repl\ are
in $\Theta$.  We have $\h(\A)=\h(\B)$ since $\A\equals\B\in\Theta$.
This implies $\conv{\h(\A)}=\conv{\h(\B)}$,
$\min{\h(\A)}=\min{\h(\B)}$, $\h(\A)\rp\h(\C)=\h(\B)\rp\h(\C)$, and
$\h(\A)+\h(\C)=\h(\B)+\h(\C)$ for every $\C\in\Pi^+$.  Then
$\h(\conv\A)=\h(\conv\B)$, $\h(\min\A)=\h(\min\B)$,
$\h(\A\rp\C)=\h(\B\rp\C)$, and $\h(\A+\C)=\h(\B+\C)$ since $\h$ is a
homomorphism.  Thus, the conclusions of \repl\ are also in $\Theta$.
Since $\Theta$ is a set containing $\Lambda^\times\cup\Psi$ and is
closed under \tran\ and \repl, it contains $\Theta\eta^\times\Psi$.
By hypothesis, we have $\varepsilon\in\Theta\eta^\times\Psi$, hence
$\varepsilon\in\Theta$, as desired. This completes the proof of one
direction of part \eqref{eq1}.

Now suppose that $\varepsilon=(\A\equals\B)$ for some $\A,\B\in\Pi^+$
and $\Psi\not\proves^\times\varepsilon$.  We wish to show
$\gc\A\not\models_\h\varepsilon$ for some $\gc\A\in\RA$ and some
homomorphism $\h:\gc\P\to\gc\A$ such that $\gc\A\models_\h\Psi$.  It
suffices to let $\gc\A=\gc\P/{\simeq^\times_\Psi}$ and $\h=\q$.  Since
$\A\equals\B$ is not provable in $\Lx$ from $\Psi$, the equivalence
classes $\q(\A)=\A/{\simeq^\times_\Psi}$ and $\q(\B)=
\B/{\simeq^\times_\Psi}$ are distinct, hence $\q(A)\neq\q(\B)$.  Thus,
we have an algebra $\gc\A\in\RA$ and a homomorphism $\h\colon
\gc\P\to\gc\A$ such that $\gc\A\models_\h\Psi$ but not
$\gc\A\models_\h\varepsilon$, as desired.  This completes the proof of
part \eqref{eq1}.  For part \eqref{eq2}, repeat the proof of part
\eqref{eq1} using $\SA$, $\Lambda^\times_s$, and $\simeq^s_\Psi$ in
place of $\RA$, $\Lambda^\times$, and $\simeq^\times_\Psi$.

Parts \eqref{eq3} and \eqref{eq4} follow from parts \eqref{eq1} and
\eqref{eq2} when $\Psi=\emptyset$.  We show that the quotient algebra
$\gc\P/{\simeq^\times}$ is $\RA$-freely generated by the
$\simeq^\times$-equivalence classes of the propositional variables.
One starts with an arbitrary map $\f$ from $\{\A/{\simeq^\times}:
\id\neq\A\in\Pi\}$ into a relation algebra $\gc\A$.  Then $\f$
determines a map on the propositional variables $\{\A:\id\neq
\A\in\Pi\}$ that sends $\A$ to $\f(\A/{\simeq^\times})$.  Since
$\gc\P$ is absolutely freely generated by the propositional variables,
this map has a unique extension to a homomorphism $\h\colon
\gc\P\to\gc\A$ that vacuously satisfies the condition $\gc\A\models_\h
\emptyset$.  Now $\h$ sends any two equivalent predicates to the same
thing, for if $\A\simeq^\times\B$, then $\proves^\times\A\equals\B$,
hence $\gc\A\models_\h\A\equals\B$ by part \eqref{eq1}, \ie,
$\h(\A)=\h(\B)$. This means that $\h$ determines a map $\g$ from
$\{\A/{\simeq^\times} :\A\in\Pi^+\}$ into $\gc\A$ that satisfies the
condition $\g(\A/{\simeq^\times})=\h(\A)$ for every $\A\in\Pi^+$.
This condition implies $\g$ is a homomorphism because, for example, if
$\A,\B\in\Pi^+$ then
\begin{align*}
  \g(\A/{\simeq^\times}+\B/{\simeq^\times})
  &=\g((\A+\B)/{\simeq^\times})=\h(\A+\B)
  \\&=\h(\A)+\h(\B)=\g(\A/{\simeq^\times})+\g(\B/{\simeq^\times}).
\end{align*}
This proves part \eqref{eq3}.  The proof of part \eqref{eq4} is
essentially the same.
\endproof
An algebra $\gc\A= \<\U, +, \min\blank, \rp, \conv{}, \id\>$ is {\bf
  proper relation algebra} if there is an equivalence relation
$\E\in\U$ such that $\U$ is a set of binary relations included in $\E$
and the following conditions hold.
\begin{align}
  \label{union}
  \A+\B&=\{\<\x,\y\>:\<\x,\y\>\in\A\mor\<\x,\y\>\in\B\}=\A\cup\B,\\
  \label{comp}
  \min\A&=\{\<\x,\y\>:\<\x,\y\>\in\E\mand\<\x,\y\>\notin\A\}
  =\E\setminus\A,\\
  \label{relprod}
  \A\rp\B&=\{\<\x,\y\>: \text{for some $\z$, }
  \<\x,\z\>\in\A\text{ and }\<\z,\y\>\in\B\},\\
  \label{converse}
  \conv\A&=\{\<\x,\y\>:\<\y,\x\>\in\A\},\\
  \label{id}
  \id&=\{\<\x,\x\>:\<\x,\x\>\in\E\}.
\end{align}
For any equivalence relation $\E$, $\Sb\E$ is the proper relation
algebra whose universe $\wp(\E)$ consists of all relations included in
$\E$, called the {\bf algebra of subrelations of $\E$}.  For any set
$\U$, $\gc\A$ is {\bf proper relation algebra on} $\U$ if conditions
\eqref{union}--\eqref{id} hold when $\E=\U\times\U$.  For example,
$\Re\U$ is defined as the proper relation algebra on $\U$ whose
universe $\wp(\U\times\U)$ is the set of all binary relations on $\U$,
called the {\bf algebra of binary relations on $\U$}.  The algebra
$\gc\A$ is a {\bf representable relation algebra} if it is isomorphic
to a proper relation algebra.  Let $\RRA$ be the class of
representable relation algebras.  Straightforward computations show
that the axioms of $\Lx$ and $\Lwx$ are true in every representable
relation algebra.  Therefore, $\RRA\subseteq\RA\subseteq\SA$.

Suppose $\U$ is a set and $\R(\A)$ is a binary relation on $\U$
whenever $\id\neq\A\in\Pi$. This determines a relational structure
$\gc\U= \<\U,\R(\A)\>_{\id\neq\A\in\Pi}$.  The {\bf denotation
  function} $\mathsf{De}_{\gc\U}$ of $\gc\U$ is defined for all
$\A,\B\in\Pi^+$ by the following conditions.
\begin{align}
  \De\A&=\R(\A)\text{ if }\id\neq\A\in\Pi,\label{De1}\\
  \De{\A+\B}&=\De\A\cup\De\B,\label{De2}\\
  \De{\min\A}&=(\U\times\U)\setminus\De\A,\label{De3}\\
  \De{\A\rp\B}&=\De\A\rp\De\B,\label{De4}\\
  \De{\conv\A}&=\conv{\De\A},\label{De5}\\
  \De\id&=\{\<\u,\u\>:\u\in\U\},\label{De6}
\end{align}
where the symbols on the right in \eqref{De4} and \eqref{De5} denote
the operations on binary relations defined in \eqref{relprod} and
\eqref{converse}.  This definition is used by \cite[p.\,26, p.\,47,
  6.1(i)]{MR920815}.  Note that $\Dee$ is a homomorphism from $\gc\P$
into $\Re\U$. Since $\gc\P$ is absolutely freely generated by
$\{\A:\id\neq\A\in\Pi\}$, $\Dee$ could also have been defined as the
unique homomorphism from $\gc\P$ into $\Re\U$ determined by condition
\eqref{De1}.

For any sequence $\s\colon\omega\to\U$ of elements of $\U$ and any
formula $\varphi\in\Phi^+$ the {\bf satisfaction relation} is defined
by induction on the complexity of formulas in $\Phi^+$ as follows.
\begin{align*}
  \gc\U&\models\A\equals\B[\s]&&\ifff\De\A=\De\B, \text{ for
    $\A,\B\in\Pi^+$,}\\
  \gc\U&\models\var\i\A\var\j[\s]&&\ifff\<\s_\i,\s_\j\>\in\De\A,
  \text{ for $\A\in\Pi^+$ and
    $\i,\j\in\omega$,}\\ \gc\U&\models\varphi\implies\psi[\s] &&\ifff
  \gc\U\models\psi[\s] \text{ or not }
  \gc\U\models\varphi[\s],\\ \gc\U&\models\neg\varphi[\s] &&\ifff
  \text{not }
  \gc\U\models\varphi[\s],\\ \gc\U&\models\all{\var\i}\varphi[\s]&&\ifff
  \gc\U\models\varphi[\s^{\i/\u}] \text{ for every $\u\in\U$ and
    $\i\in\omega$, where}\\ &&&\s^{\i/\u}_\j=\begin{cases}
  \s_\j&\text{ if }\i\neq\j\in\omega,\\ \u&\text{ if }\j=\i.
  \end{cases}
\end{align*}
Finally, $\gc\U\models\varphi$ iff $\gc\U\models\varphi[\s]$ for all
$\s:\omega\to\I$, and for every $\Psi\subseteq\Phi^+$,
$\gc\U\models\Psi$ if $\gc\U\models\varphi$ for every
$\varphi\in\Psi$.  When $\gc\U\models\varphi$ or $\gc\U\models\Psi$,
we say $\gc\U$ is a {\bf model of} $\varphi$ or $\Psi$, respectively.
\begin{theorem}\label{sem} Assume $\Psi\subseteq\Sigma^+$. 
  \begin{enumerate}
  \item\label{sem1} $\simeq^+_\Psi$ is a congruence relation on
    $\gc\P$ and the quotient algebra $\gc\P/{\simeq^+_\Psi}$ is a
    representable relation algebra {\rm\cite[8.3(vii)]{MR920815}}.
  \item\label{sem3} For every $\varphi \in \Sigma^+$,
    $\Psi\proves^+\varphi$ iff for every proper relation algebra
    $\gc\A$ on a set $\U$ and every homomorphism
    $\h\colon\gc\P\to\gc\A$, if $\gc\U= \<\U,\h(\A)\>_{\A\in\Pi}$ is a
    model of $\Psi$ then $\gc\U$ is a model of $\varphi$.
  \item\label{sem2} $\gc\P/{\simeq^+}$ is a representable relation
    algebra that is $\RRA$-freely generated by
    $\{\A/{\simeq^+}:\id\neq\A\in\Pi\}$
    {\rm\cite[8.3(viii)]{MR920815}, \cite[Thm 11(4)]{MR2628352},
      \cite[Thm 553]{MR2269199}, \cite[Thm 4.3]{MR1720820}}.
  \end{enumerate}
\end{theorem}
\proof
The proof of Theorem \ref{sem} makes use of G\"odel's completeness
theorem for $\LL$ and its extension to $\Lp$ \cite[\SS1.4,
\SS2.2]{MR920815}.  For a detailed proof of the completeness theorem
for $\Lp$ in a more general setting that allows predicates of
arbitrary finite rank, see \cite[Thm 170]{MR2269199}.  By the
extension of the completeness theorem to $\Lp$, $\A\simeq^+_\Psi\B$
iff $\De\A=\De\B$ whenever $\gc\U$ is a model of $\Psi$. It is
apparent from this that $\simeq^+_\Psi$ is a congruence relation, but
for a detailed proof see \cite[Thm 130]{MR2269199}.

Tarski and Givant present the following proof of part \eqref{sem1}.
First construct an indexed system $\<\gc\U_\i:\i\in\I\>$ of structures
$\U_\i=\<\U_\i,\R_\i(\A)\>_{\id\neq\A\in\Pi}$ such that every model of
$\Psi$ is elementarily equivalent to (satisfies the same sentences as)
one of the indexed structures.  Let $\gc\A$ be the the direct product
of the system $\<\Re{\U_\i}:\i\in\I\>$.  For each $\i\in\I$ let
$\h_\i$ be the unique homomorphism from $\gc\P$ into $\Re{\U_\i}$ that
extends $\R_\i\colon\{\id\neq\A\in\Pi\}\to\wp(\U_\i\times\U_\i)$.
Since $\h_\i$ is a homomorphism that agrees with $\Dei$ on the
propositional variables it follows from the definition of $\Dei$ by
induction on predicates that $\h_\i=\Dei$.  If $\A\simeq^+_\Psi\B$,
\ie, $\Psi\proves^+\A\equals\B$, then $\gc\U_\i\models\A\equals\B$ by
the soundness part of the completeness theorem, hence
$\h_\i(\A)=\Dei\A=\Dei\B=\h_\i(\B)$. It follows that $\h_\i$
determines a map $\g_\i$ from $\{\A/{\simeq^+_\Psi} :\A\in\Pi^+\}$
into $\gc\A$ that satisfies the condition $\g(\A/{\simeq^+_\Psi})
=\h(\A)$ for every $\A\in\Pi^+$.  This condition implies $\g_\i$ is a
homomorphism as in the proof of Theorem \ref{eq}\eqref{eq3}.  Define
$\k$ by setting
\begin{equation*}
  \k(\A/{\simeq^+_\Psi})=\<\g_\i(\A/{\simeq^+_\Psi}):\i\in\I\>
\text{ for every }\A\in\Pi^+.
\end{equation*}
Then $\k$ is a homomorphism from $\gc\P/{\simeq^+_\Psi}$ into $\gc\A$.
Suppose $\A/{\simeq^+_\Psi}\neq\B/{\simeq^+_\Psi}$ for some
$\A,\B\in\Pi^+$. By the completeness theorem for $\Lp$, there is a
model of $\Psi$ in which $\A\equals\B$ fails, so by the definition of
the system $\<\gc\U_\i:\i\in\I\>$ there is some $\j\in\I$ such that
$\g_\j(\A/{\simeq^+_\Psi})\neq\g_\j(\B/{\simeq^+_\Psi})$ and hence
$\k(\A/{\simeq^+_\Psi})\neq\k(\B/{\simeq^+_\Psi})$. This shows that
$\k$ is one-to-one and therefore an isomorphism of
$\gc\P/{\simeq^+_\Psi}$ into the proper relation algebra $\gc\A$.
\endproof
\section{3-variable formalisms $\Lp_3$ and $\LL_3$}\label{sect7}
$\Lx$ is weaker than $\LL^+$ and $\LL$.  It follows from Theorem
\ref{eq}\eqref{eq1} and Theorem \ref{sem}\eqref{sem2} that if an
equation is provable in $\Lp$ but is not true in some $\gc\A\in\RA$,
then it is not provable in $\Lx$ and $\Lx$ is therefore weaker than
$\Lp$ in means of proof.  To get the following equation, which is
provable in $\Lp$ but not $\Lx$, Tarski and Givant used a
non-representable relation algebra found by \citet{MR0286735} that is
generated by a single element.  Givant constructed the equation and it
was later simplified by George McNulty and Tarski
\cite[3.4(vi)]{MR920815}.
\begin{align*}
  1=1\rp(\A\rs\A+(\A\rp\A+\id+
  (\min{\A+\conv\A})\rp(\min{\A+\conv\A})
  \cdot\min{\conv\A})\cdot\min\A+\min{\A\rp\conv\A})\rp1.
\end{align*}
$\Lx$ is also weaker than $\LL$ and $\Lp$ in means of expression due
to Korselt's result, reported by \citet{MR1511835}, that no equation
in $\Sigma^\times$ is logically equivalent to any sentence asserting
the existence of four distinct objects, such as
\begin{equation*}
  \closure{\ex\w\neg(\w\id\x\lor\w\id\y\lor\w\id\z)}.
\end{equation*}
Tarski greatly generalized Korselt's theorem
\cite[3.5(viii)]{MR920815}.

$\Lx$ seems to be correlated with the logic of three variables because
Korselt's sentence uses four variables while $\GG(\A\equals\B)$
contains at most three.  Indeed, it is apparent from the definition of
$\GG$ that if neither $\A$ nor $\B$ contains an occurrence of~$\rp$,
then only the first two variables occur in $\GG(\A\equals\B)$ but
if~$\rp$ occurs in $\A$ or $\B$ then $\GG(\A\equals\B)$ does have the
first three variables in it.  This suggests that perhaps \emph{every}
sentence containing only the first three variables is logically
equivalent to an equation in $\Sigma^\times$.  Tarski was able to show
that this is actually the case. He proposed the construction of
3-variable formalisms $\Lp_3$ and $\LL_3$ that would be equipollent
with $\Lx$ in means of expression and proof.

For every finite $\n\geq3$ let $\Phi^+_\n$ be the set of formulas in
$\Phi^+$ that contain only variables in $\Upsilon_\n$ and let
\begin{align*}
  \Phi_\n&=\Phi\cap\Phi^+_\n,&
  \Sigma_\n&=\Sigma\cap\Phi^+_\n,&
  \Sigma^+_\n&=\Sigma^+\cap\Phi^+_\n.
\end{align*}
Tarski's theorem that every sentence in $\Sigma^+_3$ is logically
equivalent to an equation in $\Sigma^\times$ suggests that
$\Sigma^+_3$ and $\Sigma_3$ should be the sets of sentences of $\Lp_3$
and $\LL_3$. The restriction of $\GG$ to $\Sigma^+_3$ could serve as
the translation mapping from $\Lp_3$ to $\LL_3$.

Tarski's initial proposal came in two parts.  First, Tarski proposed
restricting the axioms \eqref{AI}--\eqref{AIX} and the rule $\MP$ to
those instances that belong to $\Sigma^+_3$.  Givant found these
restricted axioms were too weak and suggested replacing \eqref{AIX}
with \eqref{AIX'}, called the {\bf general Leibniz law}, which is
formulated in terms of a variant type of substitution defined by
\citet[pp.\,66--67]{MR920815}.
\begin{align}
  &\label{AIX'}\tag{AIX$'$}
  \closure{\x\eqsym\y\implies(\varphi\implies\varphi[\x/\y])}.
\end{align}
The variant substitution is complicated so Tarski and Givant borrowed
an idea from \cite{MR2628352} to formulate an alternate axiom
\eqref{AIX''}. If $\x,\y\in\Upsilon$ and $\varphi\in\Phi^+$ then
$\switch\x\y\varphi$ is the result of interchanging $\x$ and $\y$
throughout the formula $\varphi$.  The function $\switch\x\y \colon
\Phi^+\to\Phi^+$ is determined by these rules, in which $\hat\x=\y$,
$\hat\y=\x$, and $\hat\v=\v$ if $\v\neq\x,\y$ for every
$\v\in\Upsilon$.
\begin{align*}
  \switch\x\y(\v\A\w) &=\hat\v\A\hat\w,\\
  \switch\x\y(\varphi\implies\psi)&=\switch\x\y(\varphi)
  \implies\switch\x\y(\psi),\\
  \switch\x\y(\neg\varphi)&=\neg\switch\x\y(\varphi),\\
  \switch\x\y(\all\v\varphi)&=\all{\hat\v}\switch\x\y(\varphi).
\end{align*}
Givant proved that \eqref{AIX''} can be used instead of \eqref{AIX'}
in the axiomatization of $\Lp_3$,
\begin{align}
  &\label{AIX''}\tag{AIX$''$}
  \closure{\x\eqsym\y\implies(\varphi\implies\switch\x\y\varphi)}.
\end{align}
Tarski knew by the early 1940s that \eqref{BIV} can not be proved with
only three variables and would have to be included in the
axiomatization of $\Lp_3$ by {\it fiat}. The second part of Tarski's
proposal was to include the general associativity axiom
\begin{align}
  \label{AX}\tag{AX}
  \pmb[\,&\ex\z(\ex{\y}(\varphi[\x,\y] {\land} \psi[\y,\z]) {\land}
  \xi[\z,\y])\iff \ex\z(\varphi[\x,\z] {\land} \ex\x(\psi[\z,\x]
  {\land} \xi[\x,\y]))\pmb].
\end{align} 
This axiom involves the complicated substitution but Givant proved it
could be replaced by \eqref{AX'}, in which the free variables of
formulas $\varphi,\psi,\xi\in\Phi^+$ are just $\x$ and $\y$,
\begin{align}
  &\label{AX'}\tag{AX$'$}
  \closure{\ex\z(\ex\y(\varphi\land\switch\x\z\psi)\land\switch\x\z\xi)
    \iff\ex\z(\switch\y\z\varphi\land\ex\x(\switch\y\z\psi\land\xi)}.
\end{align} 
The sets of sentences of $\Lp_3$ and $\LL_3$ are $\Sigma^+_3$ and
$\Sigma_3$, their axioms are the sentences in $\Sigma^+_3$ and
$\Sigma_3$ that are instances of \eqref{AI}--\eqref{AVIII},
\eqref{AIX'}, or \eqref{AX}, and their rule of inference is \MP.  For
simpler axiom sets use \eqref{AIX''} and \eqref{AX'} instead of
\eqref{AIX'} and \eqref{AX}.  For every $\Psi\subseteq\Sigma^+_3$ a
sentence $\varphi\in\Sigma^+_3$ is {\bf provable in $\Lp_3$ from
  $\Psi$}, written $\Psi\proves^+_3\varphi$ or $\proves^+_3\varphi$ if
$\Psi=\emptyset$, if $\varphi$ is in every set that contains $\Psi$
and the axioms of $\Lp_3$ and is closed under $\MP$.  The {\bf theory
  generated by $\Psi$} in $\Lp_3$ is
\begin{equation*}
  \Theta\eta^+_3\,\Psi
  =\{\varphi:\varphi\in\Sigma^+_3,\,\Psi\proves^+_3\varphi\}.
\end{equation*}
The notions of $\Psi\proves_3\varphi$, $\varphi$ is {\bf provable in
  $\LL_3$ from $\Psi$}, and $\Theta\eta_3\,\Psi$, the {\bf theory
  generated by $\Psi$} in $\LL_3$, are defined similarly for every
$\Psi\subseteq\Sigma_3$ and sentence $\varphi\in\Sigma_3$.
\section{Equipollence of $\Lx$, $\LL_3$, and $\Lp_3$}\label{sect8}
With Givant's changes Tarski's proposal worked. It provided a
3-variable restriction $\Lp_3$ of $\LL^+$ and a 3-variable restriction
$\LL_3$ of $\LL$. Both are equipollent with $\Lx$ in means of
expression and proof.  For the equipollence of $\LL_3$ and $\Lp_3$ the
appropriate translation mapping from $\Lp_3$ to $\LL_3$ is simply the
restriction of $\GG$ to $\Sigma^+_3$. Part \eqref{thm2iv} of the
following theorem is the {\bf main mapping theorem} for $\LL_3$ and
$\Lp_3$.
\begin{theorem}{\rm\cite[\SS3.8]{MR920815}}\label{thm2}
  Formalisms $\LL_3$ and $\Lp_3$ are equipollent.
  \begin{enumerate}
  \item\label{thm2i} $\Phi_3\subseteq\Phi^+_3$ and
    $\Sigma_3\subseteq\Sigma^+_3$  {\rm[3.8(viii)($\alpha$)]}.
  \item\label{thm2ii} $\GG$ maps $\Phi^+_3$ onto $\Phi_3$ and
    $\Sigma^+_3$ onto $\Sigma_3$ {\rm[3.8(ix)($\delta$)]}.
  \item\label{thm2iii} $\varphi \equiv^+_3\GG(\varphi)$ if $\varphi\in
    \Phi^+_3$ {\rm[3.8(ix)($\varepsilon$)]}.
  \item\label{thm2iv} $\Psi\proves^+_3\varphi$ iff
    $\{\GG(\psi):\psi\in\Psi\} \proves_3 \GG(\varphi)$, for
    $\Psi\subseteq\Sigma^+_3$ and $\varphi\in\Sigma^+_3$
    {\rm[3.8(xi)]}.
  \item\label{thm2v} $\Psi\proves^+_3\varphi$ iff
    $\Psi\proves_3\varphi$, for $\Psi\subseteq\Sigma_3$ and
    $\varphi\in\Sigma_3$, {\rm[3.8(viii)($\beta$),
      3.8(xii)($\beta$)]}.
  \end{enumerate}
\end{theorem}
For the equipollence of $\Lx$ and $\Lp_3$, Tarski and Givant construct
a function $\HH$ on $\Phi^+_3$ whose restriction to $\Sigma^+_3$ is an
appropriate translation mapping from $\Lp_3$ to $\Lx$ \cite[pp.\
77--79]{MR920815}.  They begin, ``Its definition is complicated and
must be formulated with care. We give here enough hints for
constructing such a definition, without formulating it precisely in
all details,'' and end with, ``We hope the above outline gives an
adequate idea of the definition of $\HH$''.  After obtaining the main
mapping theorem for $\Lx$ and $\Lp_3$ they say,
\begin{quote}
  ``The construction used here to establish these equipollence results
  has clearly some serious defects, if only from the point of view of
  mathematical elegance.  Actually, this applies to the proof of the
  equipollence of $\Lp_3$ and $\Lx$.  The splintered character of the
  definition of the translation mapping $\HH$, with its many cases, is
  a principal cause of the fragmented nature of certain portions of
  the argument; the involved notion of substitution (which we have to
  use because of the restricted number of variables in our formalisms)
  is another detrimental factor. As a final result, the construction
  is so cumbersome in some of its parts---culminating in the proofs of
  (iv) and (v)---that we did not even attempt to present them in full.
  A different construction that would remove most of the present
  defects would be very desirable indeed.'' \cite[p.\,87]{MR920815}
\end{quote}
Another description of $\HH$ was given by \citet{MR2310682} and
simpler alternative constructions appear in
\cite[p.\ 192--3]{MR2628352} and \cite[p.\,543--4,
  p.\,548--9]{MR2269199}.  One of these is presented here as a
response to Tarski and Givant. It can be precisely defined in one page
instead of outlined in three.  We start with an auxiliary map
$\JJ\colon\Phi^+_3\to\Phi^+_3$.  The map $\JJ$ has the property that
for every formula $\varphi\in \Phi^+_3$ there are $\k\in\omega$ and
finite sequences of atomic predicates $\R,\S,\T\in{}^\k\Pi$ such that
\begin{equation*}
  \JJ(\varphi)=\bigwedge_{\i<\k}(\0\R_\i\var2\lor\var2\S_\i\var1
  \lor\0\T_\i\var1)\in\Phi^+_3.
\end{equation*}
$\JJ$ is defined by induction on the complexity of formulas.  For the
meanings of $0$, $1$, $\cdot$, and $\rs$, recall definitions
\eqref{10di} and \eqref{dotmin}.  If $\A,\B\in\Pi$ then
\begin{align*}
  \JJ(\0\A\var1)&=\00\var2\lor\var20\var1\lor\0\A\var1,\\
  \JJ(\var1\A\0)&=\00\var2\lor\var20\var1\lor\0\conv\A\var1,\\
  \JJ(\var1\A\var2)&=\00\var2\lor\var2\conv\A\var1\lor\00\var1,\\
  \JJ(\var2\A\var1)&=\00\var2\lor\var2\A\var1\lor\00\var1,\\
  \JJ(\0\A\var2)&=\0\A\var2\lor\var20\var1\lor\00\var1,\\
  \JJ(\var2\A\0)&=\0\conv\A\var2\lor\var20\var1\lor\00\var1,\\
  \JJ(\0\A\0)&=\00\var2\lor\var20\var1\lor
  \0(\A\cdot\id)\rp1\var1,\\
  \JJ(\var1\A\var1)&=\00\var2\lor\var20\var1\lor
  \01\rp(\A\cdot\id)\var1,\\
  \JJ(\var2\A\var2)&=\0(1\rp(\A\cdot\id))\var2\lor
  \var20\var1\lor\00\var1,\\
  \JJ(\A\equals\B)&=\00\var2\lor\var20\var1\lor
  \00\rs(\A\cdot\B+\min\A\cdot\min\B)\rs0\var1.
\end{align*}
This completes the cases in which $\varphi$ is atomic.  If
$\varphi,\varphi'\in\Phi^+_3$ and there are $\k,\k'\in\omega$,
$\R,\S,\T\in{}^\k\Pi$, and $\R',\S',\T'\in{}^{\k'}\Pi$ such that
\begin{align*}
  \JJ(\varphi)&=\bigwedge_{\i<\k}(\0\R_\i\var2\lor
  \var2\S_\i\var1\lor\0\T_\i\var1),\\
  \JJ(\varphi')&=\bigwedge_{\j<\k'}(\0\R'_\j\var2\lor
  \var2\S'_\j\var1\lor\0\T'_\j\var1),
\end{align*}
then
\begin{align*}
  &\JJ(\neg\varphi)=\bigwedge_{f\in{}^\k3} \Bigg(
  \0\Big(\sum_{f(\i)=0}\min{\R_\i}\Big) \var2\lor
  \var2\Big(\sum_{f(\i)=1}\min{\S_\i}\Big) \var1\lor
  \0\Big(\sum_{f(\i)=2}\min{\T_\i}\Big)\var1\Bigg),\\
  &\JJ(\varphi\implies\varphi')= \bigwedge_{f\in{}^\k3,\,\j<\k'}
  \Bigg(\0\Big(\sum_{f(\i)=0}\min{\R_\i}+\R'_\j\Big)\var2\lor
  \var2\Big(\sum_{f(\i)=1}\min{\S_\i}+\S'_\j\Big)\var1\\
  &\hskip2.77in\lor\0\Big(\sum_{f(\i)=2}\min{\T_\i}+
  \T'_\j\Big)\var1\Bigg),\\
  &\JJ(\all\0\varphi)=\bigwedge_{\i<\k}\Big(\00\var2\lor
  \var2(\conv\R_\i\rs\T_\i)+\S_\i\var1\lor\00\var1\Big),\\
  &\JJ(\all{\var1}\varphi)=\bigwedge_{\i<\k}
  \Big(\0(\T_\i\rs\conv\S_\i)+\R_\i\var2\lor
  \var20\var1\lor\00\var1\Big),\\
  &\JJ(\all{\var2}\varphi)=
  \bigwedge_{\i<\k}\Big(\00\var2\lor\var20\var1\lor
  \0(\R_\i\rs\S_\i)+\T_\i\var1\Big).
\end{align*}
Now we can define $\HH$ on sentences $\varphi\in\Sigma_3$. Apply $\JJ$
to $\all{\var2}\varphi$, obtaining $\k<\omega$ and finite sequences
$\R,\S,\T\in{}^\k\Pi$ such that
\begin{equation*}
  \JJ(\all{\var2}\varphi)=
  \bigwedge_{\i<\k}\Big(\var00\var2\lor\var20\var1
  \lor\var0(\R_\i\rs\S_\i)+\T_\i\var1\Big),
\end{equation*}
and set $\HH(\varphi)$ equal to
\begin{align*}
  1\equals
  \big((\R_0\rs\S_0)+\T_0\big)\cdot\big((\R_1\rs\S_1)+\T_1\big)\cdot\ldots
  \cdot\big((\R_{\k-1}\rs\S_{\k-1})+\T_{\k-1}\big).
\end{align*}
The original construction of $\HH$ in \cite[\SS3.9]{MR920815} has
desirable properties not shared by the mapping $\HH$ defined above.
Tarski's $\HH$ is defined on all formulas but the definition given
here is restricted to sentences. By
\cite[3.9(iii)($\delta$)]{MR920815} the sets of free variables of
$\varphi$ and $\HH(\varphi)$ are the same for every formula
$\varphi\in\Phi^+_3$ and $\HH$ maps $\Sigma^+_3$ onto $\Sigma^\times$
(not just into, as is the case for the $\HH$ constructed here).
Furthermore, $\HH$ produces simpler output in certain cases.  For
example, $\HH(\neg\A\equals\B) = 1\rp\min{\A\cdot\B+\min\A\cdot\min\B}
\rp1\equals1$, $\HH(\neg\A\equals1)=1\rp\min\A\rp1\equals1$, and
$\HH(\A\equals\B)=\A\equals\B$ whenever $\A,\B\in\Pi$.  On the other
hand, to insure that the output of $\HH$ has the same set of free
variables as the input, it is necessary in the definition of
$\HH(\neg\varphi)$ and $\HH(\varphi\implies\psi)$ to consider many
cases that depend on the free variables of $\HH(\varphi)$ and
$\HH(\psi)$.  This causes the very long proof by cases of the main
mapping theorem for $\Lx$ and $\Lp_3$ encountered by Tarski and
Givant.

The following theorem summarizes the principal parts of the
equipollence of $\Lx$ and $\Lp_3$. It is stated for the versions of
$\JJ$ and $\HH$ constructed here, so part \eqref{thm3ii} only says
``into'' instead of ``onto''.  The other parts are the same as the
corresponding versions in \cite[\SS3.9]{MR920815}.  Part
\eqref{thm3iv} is the {\bf main mapping theorem} for $\Lx$ and
$\Lp_3$, while part \eqref{thm3v}, a corollary of part \eqref{thm3iv},
is the equipollence of $\Lx$ and $\Lp_3$ in means of proof.
\begin{theorem}{\rm\cite[\SS3.9]{MR920815}}\label{thm3}
  Formalisms $\Lx$ and $\Lp_3$ are equipollent in means of expression
  and proof.
  \begin{enumerate}
  \item\label{thm3i} $\Sigma^\times\subseteq\Sigma^+_3$ {\rm[3.9(i)]}.
  \item\label{thm3ii} $\HH$ maps $\Sigma^+_3$ into $\Sigma^\times$.
  \item\label{thm3iii} $\varphi\equiv^+_3\JJ(\varphi)$ if
    $\varphi\in\Phi^+_3$ {\rm[3.9(iii)($\varepsilon$)]} and
    $\JJ(\varphi)\equiv\HH(\varphi)$ if $\varphi\in\Sigma^+_3$.
  \item\label{thm3iv} $\Psi\proves^+_3\varphi$ iff
    $\{\HH(\psi):\psi\in\Psi\}\proves^\times\HH(\varphi)$, for
    $\Psi\subseteq\Sigma^+_3$ and $\varphi\in\Sigma^+_3$
    {\rm[3.9(vii)]}.
  \item\label{thm3v} $\Psi\proves^+_3\varepsilon$ iff
    $\Psi\proves^\times\varepsilon$, for $\Psi\subseteq\Sigma^\times$
    and $\varepsilon\in\Sigma^\times$ {\rm[3.9(ix)]}.
  \item\label{thm3vi}
    $\Theta\eta_3\,\Psi=\Theta\eta^+_3\,\Psi\cap\Sigma_3$, for
    $\Psi\subseteq\Sigma_3$ {\rm[3.9(x)]}.
  \item\label{thm3vii}
    $\Theta\eta^\times\,\Psi=\Theta\eta^+_3\,\Psi\cap\Sigma^\times$,
    for $\Psi\subseteq\Sigma^\times$ {\rm[3.9(xi)]}.
  \end{enumerate}
\end{theorem}
\section{Equipollent 3-variable formalisms $\Ls$, $\Ls^+$, $\Lwx$}
\label{sect9}
Since Tarski and Givant included the axiom \eqref{AX} only to achieve
the equipollence of $\LL_3$ and $\Lp_3$ with $\Lx$, they defined the
``(\emph{standardized}) \emph{formalisms}'' $\Ls_3$ and $\Ls_3^+$ by
deleting \eqref{AX} from the axiom sets of $\LL_3$ and $\Lp_3$
\cite[p.\,89]{MR920815}.  They did not introduce any special notation
for provability in these formalisms. We use $\proves_s$ and
$\proves^+_s$ for $\Ls_3$ and $\Ls_3^+$.  Since \eqref{BIV} is the
axiom of $\Lx$ that corresponds to \eqref{AX}, Tarski and Givant asked
whether simply deleting \eqref{BIV} from the axioms of $\Lx$ would
produce a formalism equipollent with the standardized formalisms
$\Ls_3$ and $\Ls_3^+$.  The answer is ``no''.  For example, the {\bf
  semi-associative law} \eqref{BIV'} in Table \ref{Lx-axioms} is
provable in $\Ls_3^+$ but cannot be derived in $\Lx$ from just the
axioms \eqref{BI}--\eqref{BIII} and \eqref{BV}--\eqref{BX} when
$\id\neq\A,\B\in\Pi$. Another equation with the same property is
\begin{align*}
  &\A\rp1\equals(\A\rp1)\rp1.
\end{align*}
Adding either one of these as an axiom produces a formalism
equipollent with $\Ls_3$ and $\Ls_3^+$. Therefore, Tarski and Givant
defined a weakened equational formalism $\Lwx$ by replacing
\eqref{BIV} with \eqref{BIV'} in the axiomatization of $\Lx$.  The
equipollence of $\Lwx$ with $\Ls_3$ and $\Ls_3^+$ is stated in the
next theorem and was noted by \citet[p.\,89, p.\,209]{MR920815}.  Part
\eqref{thm5iv} tells us that the axioms for $\SA$ characterize the
equations provable in standard first order logic of three variables
without associativity.
\begin{theorem}\label{thm5}{\rm\cite[Thm 11(30)]{MR2628352},
    \cite[Thm 6.3]{MR1720820}, \cite[Thm 569(i)(ii)]{MR2269199}}
  Formalisms $\Lwx$, $\Ls_3$, and $\Ls_3^+$ are equipollent in means
  of expression and proof.
  \begin{enumerate}
  \item\label{thm5i} $\varphi\equiv^+_s\GG(\varphi)
    \equiv^+_s\JJ(\varphi)$, for $\varphi\in\Phi^+_3$.
  \item\label{thm5ii} $\Psi\proves^+_s\varphi$ iff
    $\{\GG(\psi):\psi\in\Psi\} \proves_s \GG(\varphi)$ iff
    $\{\HH(\psi):\psi\in\Psi\}\proves^\times_s\HH(\varphi)$, for
    $\Psi\subseteq\Sigma^+_3$ and $\varphi\in\Sigma^+_3$.
  \item\label{thm5iii} $\Psi\proves^+_s\varphi$ iff
    $\Psi\proves_s\varphi$, for $\Psi\subseteq\Sigma_3$ and
    $\varphi\in\Sigma_3$.
  \item\label{thm5iv} $\Psi\proves^+_s\varepsilon$ iff
    $\Psi\proves^\times_s\varepsilon$, for $\Psi\subseteq
    \Sigma^\times$ and $\varepsilon\in\Sigma^\times$.
  \end{enumerate}
\end{theorem}
\citet[p.\,91]{MR920815} also define formalisms $\LL_\n$ and
$\LL_\n^+$ for every finite $\n\geq4$, imitating the definitions of
$\Ls_3$ and $\Ls_3^+$, but with $\n$ in place of $3$.  The sets of
formulas of $\LL_\n$ and $\LL_\n^+$ are $\Phi_\n$ and $\Phi^+_\n$, the
sets of sentences are $\Sigma_\n$ and $\Sigma^+_\n$, the sets of
axioms are those instances of axioms \eqref{AI}--\eqref{AVIII} and
\eqref{AIX'} that lie in $\Sigma_\n$ and $\Sigma^+_\n$, and the only
rule of inference is \MP.  The next two theorems involve the first of
these formalisms, when $\n=4$.  They are a precise expression of the
fact, noted by \citet[p.\ 92]{MR920815}, that a sentence in
$\Sigma^+_3$ (containing only three variables) can be proved using
four variables iff it can be proved using just three variables
together with the assumption that relative multiplication is
associative. (The notationally peculiar equivalence of $\proves_3$
with $\proves_4$ and $\proves^+_3$ with $\proves^+_4$ is due to the
inclusion of associativity in the definitions of $\proves_3$ and
$\proves^+_3$.)
\begin{theorem}\label{thm5+}{\rm\cite[Thm 11(31)]{MR2628352},
    \cite[Thm 6.4]{MR1720820}, \cite[Thm 569(iii)(iv)]{MR2269199}}
  \begin{enumerate}
  \item\label{thm5+i} $\Psi\proves^+_4\varphi$ iff
    $\{\GG(\psi):\psi\in\Psi\} \proves_4 \GG(\varphi)$ iff
    $\{\HH(\psi):\psi\in\Psi\}\proves^\times\HH(\varphi)$, for
    $\Psi\subseteq\Sigma^+_3$ and $\varphi\in\Sigma^+_3$.
  \item\label{thm5+iv} $\Psi\proves^+_4\varphi$ iff
    $\Psi\proves^+_3\varphi$, for $\Psi\subseteq\Sigma^+_3$ and
    $\varphi\in\Sigma^+_3$.
  \item\label{thm5+ii} $\Psi\proves^+_4\varphi$ iff
    $\Psi\proves_4\varphi$ iff $\Psi\proves_3\varphi$, for
    $\Psi\subseteq\Sigma_3$ and $\varphi\in\Sigma_3$.
  \item\label{thm5+iii} $\Psi\proves^+_4\varepsilon$ iff
    $\Psi\proves^\times\varepsilon$, for $\Psi\subseteq\Sigma^\times$
    and $\varepsilon\in\Sigma^\times$.
  \end{enumerate}
\end{theorem}
Theorem \ref{thm5+} has the following consequences when
$\Psi=\emptyset$.
\begin{theorem}\label{thm6}{\rm\cite[Thm 24]{MR1011183}}
  \begin{enumerate}
  \item\label{1} $\Theta\eta_3\,\emptyset=\Theta\eta_4
    \,\emptyset\cap\Sigma_3.$
  \item\label{2} $\Theta\eta^+_3\,\emptyset=\Theta\eta^+_4
    \,\emptyset\cap\Sigma^+_3.$
  \item\label{3} $\Theta\eta^\times\,\emptyset=\Theta\eta^+_4
    \,\emptyset\cap\Sigma^\times.$
  \end{enumerate}
\end{theorem}
Theorem \ref{thm6}\eqref{3} asserts that an equation is true in every
relation algebra iff its translation into a sentence containing only
three variables can be proved with four variables.  Theorem
\ref{thm6}\eqref{3} provides an answer to Tarski's question, ``whether
this definition of relation algebra \dots\ is justified in any
intrinsic sense''. The answer is that Tarski's definition
characterizes the equations provable with four variables.  Theorem
\ref{thm6} says, \display{\sf true in relation algebras $\iff$
  3-provable with associativity $\iff$ 4-provable.}
\section{Relevance logics $\TT_n$, $\CT_n$, $3\leq n\leq\omega$, and $\TR$}
\label{sect10}
The connectives of relevance logic that are interpreted in Table
\ref{defs-ops} as operations on binary relations are $\lor$, $\land$,
$\neg$, $\rmin$, $\to$, $\circ$, ${}^*$, and $\ttt$. To connect these
interpretations with the formalisms of Tarski and Givant we use the
symbols for these connectives to denote operators on $\Pi^+$.
Although the symbols $\lor$, $\land$, and $\neg$ have already appeared
as connectives in $\LL$, they will also denote operators on $\Pi^+$.
For their classical relevant logic $\CR$, Routley and Meyer introduce
{\bf Boolean negation $\neg$} and the {\bf Routley star $\star{}$}.
Let $\lor$, $\land$, $\neg$, $\rmin$, $\to$, $\circ$, ${}^*$, and
$\ttt$ denote the operators on $\Pi^+$ defined by
\begin{align}
  \label{or}
  \A\lor\B&=\A+\B,\\
  \label{and}
  \A\land\B&=\A\cdot\B,\\
  \label{neg}
  \neg\A&=\min\A,\\
  \label{rmin}
  \rmin\A&=\min{\conv\A},\\
  \label{to}
  \A\to\B&=\min{\conv\A\rp\min\B},\\
  \label{circ}
  \A\circ\B&=\B\rp\A,\\
  \label{star}
  \star\A&=\conv\A,\\
  \label{t=id}
  \ttt&=\id.
\end{align}
When parentheses are omitted, unary operators should be evaluated
first, followed by $\rp$, $\circ$, $\cdot$, $\land$, $\rs$, $+$,
$\lor$, and then $\to$, in that order.  Repeated binary operators of
equal precedence should be evaluated from left to right.  Definitions
\eqref{rmin} and \eqref{to} produce the standard connection between
$\to$ and $\circ$,
\begin{align*}
  \rmin(\A\to\rmin\B)
  &=\min{\bigg(\min{\conv\A\rp\min{\min{\conv\B}}}\bigg)\conv{}}
  \simeq^\times \B\rp\A = \A\circ\B.
\end{align*}
Besides the laws of double negation and the commutativity of converse
and complementation, this computation involves the axiom \eqref{BIX}.
The definitions also match the interpretations in Table
\ref{defs-ops}. To show this, we compute $\GG$ for \eqref{or},
\eqref{and}, \eqref{rmin}, and \eqref{to}. In three cases we convert
to simpler equivalent formulas.  The results agree with Table
\ref{defs-ops}. Let $\x=\0$, $\y=\var1$, and $\z=\var2$.  Assume
$\A,\B$ are propositional variables, $\id\neq\A,\B\in\Pi$, so that,
for example, $\GG(\x\A\y)=\x\A\y$ and $\GG(\x\B\y)=\x\B\y$. Then
\begin{align*}
  \GG(\x\A\lor\B\y)&=\GG(\x\A+\B\y)=\GG(\x\A\y)\lor\GG(\x\B\y)
  =\x\A\y\lor\x\B\y,\\
  \GG(\x\A\land\B\y)&=\GG(\x\A\cdot\B\y) =\GG(\x\min{\min\A+\min\B}\y)
  =\neg\GG(\x\min\A+\min\B\y)\\
  &=\neg(\GG(\x\min\A\y)\lor\GG(\x\min\B\y))
  =\neg(\neg\GG(\x\A\y)\lor\neg\GG(\x\B\y))\\
  &=\neg(\neg\x\A\y\lor\neg\x\B\y) \equiv^+\x\A\y\land\x\B\y,\\
  \GG(\x\rmin\A\y) &=\GG(\x\min{\conv\A}\y) =\neg\GG(\x\conv\A\y)
  =\neg\GG(\y\A\x) =\neg(\y\A\x),\\
  \GG(\x\A\to\B\y) &=\GG(\x\min{\conv\A\rp\min\B}\y)
  =\neg\GG(\x\conv\A\rp\min\B\y)
  =\neg\ex\z(\GG(\x\conv\A\z)\land\GG(\z\min\B\y))\\
  &=\neg\neg\all\z\neg(\GG(\z\A\x)\land\neg\GG(\z\B\y))
  =\neg\neg\all\z\neg(\z\A\x\land\neg\z\B\y)\\
  &=\neg\neg\all\z\neg\neg\big(\z\A\x\implies\neg\neg\z\B\y\big)
  \equiv^+\all\z\big(\z\A\x\implies\z\B\y\big).
\end{align*}
The {\bf relevance logic operators} are $\lor$, $\land$, $\rmin$,
$\to$, $\circ$, and $\ttt$. The {\bf classical relevant logic
  operators} are $\lor$, $\land$, $\neg$, $\rmin$, $\to$, $\circ$,
$\star{}$, and $\ttt$. The classical relevant logic operators include
$\lor$, $\neg$, $\star{}$, $\ttt$, and $\circ$. The first four of
these coincide with Tarski's $+$, $\min\blank$, $\conv{}$, $\id$, and
$\rp$ can be defined from $\circ$, so the closure of $\Pi$ under the
classical relevant logic operators is $\Pi^+$.  Let $\Pi^r$ be the
closure of $\Pi$ under just the relevance logic operators. Then
$\Pi^r$ is a proper subset of $\Pi^+$ since $\neg$ and $\star{}$
cannot be defined from the relevance logic operators.  

For the next definitions, recall that for any $\A,\B\in\Pi^+$,
$\A\leq\B$ is the equation $\A+\B\equals\B\in\Sigma^\times$.  Suppose
$3\leq\n\leq\omega$ and $\Psi\subseteq\Sigma^\times$. The equations in
$\Psi$ are the {\bf non-logical assumptions} of the logics
$\CT^\Psi_\n$ and $\TT^\Psi_\n$ defined next. We omit reference to
$\Psi$ when $\Psi=\emptyset$.
\begin{align}
  \label{CTnPsi}
  \CT^\Psi_\n&=\{\A\colon\A\in\Pi^+,\, \Psi\proves^+_\n\id\leq\A\},\\
  \label{TTnPsi}
  \TT^\Psi_\n&=\CT^\Psi_\n\cap\Pi^r,\\
  \label{CTn}
  \CT_\n&=\{\A\colon\A\in\Pi^+,\,\proves^+_\n\id\leq\A\},\\
  \label{TTn}
  \TT_\n&=\CT_\n\cap\Pi^r.
\end{align}
$\TT_\n$ is {\bf Tarski's basic $\n$-variable relevance logic} with no
non-logical assumptions. It uses only the relevance logic operators.
Adding ``$\sf\C$'' to get $\CT_\n$ indicates its classical
counterpart, in which $\neg$ and $\star{}$ are admitted, thus allowing
the full range of classical relevant logic operators.  Next, we define
some special sets of equations.
\begin{align}
  \Xi^d&=\{\A\leq\A\rp\A:\A\in\Pi^+\},
  \label{dens-eqs}\\
  \Xi^c&=\{\A\rp\B\equals\B\rp\A:\A,\B\in\Pi^+\},
  \label{comm-eqs}\\
  \Xi^s&=\{\A\equals\conv\A:\A\in\Pi^+\}.
  \label{symm-eqs}
\end{align} 
We refer to $\Xi^d$, $\Xi^c$, and $\Xi^s$ as the {\bf equations of
  density}, {\bf commutativity}, and {\bf symmetry}, respectively.
The equations of density and commutativity are used to define $\TR$,
{\bf Tarski's classical relevant logic}, by
\begin{align}\label{TR-def}
  \TR&=\CT^\Psi_4\text{ where }\Psi=\Xi^d\cup\Xi^c.
\end{align}
\section{Frames and the relevance logics $\CR$ and $\KR$}
\label{sect12}
A {\bf frame} is a quadruple $\gc\K=\<\K,\R,\star{},\II\>$ consisting
of a set $\K$, a ternary relation $\R\subseteq\K^3$, a unary operation
$\star{}\colon\K\to\K$, and a subset $\II\subseteq\K$. The associated
{\bf complex algebra of $\gc\K$} is $\Cm{\gc\K} = \<\wp(\K), \cup,
\min\blank, \rp, \conv{}, \II\>$, where $\wp(\K)$ is the set of
subsets of $\K$, and the operations $\cup$, $\min\blank$, $\rp$, and
$\conv{}$ are defined on subsets $\X,\Y\subseteq\K$ by
\begin{align}
  \label{cup}
  \X\cup\Y&=\{\x:\x\in\X\mor\x\in\Y\},\\
  \label{min}
  \min\X&=\K\setminus\X=\{\x:\x\in\K,\,\x\notin\X\},\\
  \label{rp}
  \X\rp\Y&=\{\z:\<\x,\y,\z\>\in\R
  \text{ for some $\x\in\X$, $\y\in\Y$}\},\\
  \label{conv}
  \conv\X&=\{\star\z:\z\in\X\}.
\end{align}
A predicate $\A\in\Pi^+$ is {\bf valid in} $\gc\K$, and $\gc\K$ {\bf
  validates} $\A$, if the equation $\id\leq\A$ is true in
$\Cm{\gc\K}$. $\gc\K$ {\bf invalidates} $\A$ if $\A$ is not valid in
$\gc\K$.  Any homomorphism from the predicate algebra into
$\Cm{\gc\K}$ must send $\id$ to $\II$, so $\A\in\Pi^+$ is valid in
$\gc\K$ iff $\II\subseteq\h(\A)$ for every homomorphism
$\h:\gc\P\to\Cm{\gc\K}$.  A homomorphism $\h:\gc\P\to\Cm{\gc\K}$ {\bf
  validates} $\A$ if $\II\subseteq\h(\A)$ and $\h$ {\bf invalidates}
$\A$ otherwise.  For every set $\U$, $\gc\U^2=\<\U^2,\R,\star{},\II\>$
is the {\bf pair-frame on $\U$} where
\begin{align}
  \label{pairs}
  \U^2&=\{\<\x,\y\>:\x,\y\in\U\},\\
  \label{triples}
  \R&=\{\<\<\x,\y\>,\<\y,\z\>,\<\x,\z\>\>:\x,\y,\z\in\U\},\\
  \label{pair-star}
  \star{\<\x,\y\>}&=\<\y,\x\>\text{ for all $\x,\y\in\U$},\\
  \label{id-rel} 
  \II&=\{\<\x,\x\>:\x\in\U\}.
\end{align}
Note that $\Re\U$, the algebra of binary relations on $\U$, is the
complex algebra of the pair-frame on $\U$, \ie, $\Re\U=\Cm{\gc\U^2}$.

The following conditions on a frame $\gc\K$ are written in a
first-order language with equality symbol $=$, ternary relation symbol
$\R$, unary function symbol $\star{}$, and unary relation symbol
$\II$, but we will also frequently use an atomic formula like
$\R\x\y\z$ as an abbreviation for $\<\x,\y,\z\>\in\R$ when $\R$ is a
ternary relation. Each condition should be read as holding for all
$\v,\w,\x,\y,\z\in\K$.
\begin{align}
  \label{left rotation}&\R\x\y\z\implies\R\y\star\z\star\x,\\
  \label{right rotation}&\R\x\y\z\implies\R\star\z\x\star\y,\\
  \label{center reflection}&\R\x\y\z\implies\R\star\y\star\x\star\z,\\
  \label{left reflection}&\R\x\y\z\implies\R\star\x\z\y,\\
  \label{right reflection}&\R\x\y\z\implies\R\z\star\y\x,\\
  \label{identity}&\x=\y\iff\ex\u(\II\u\land\R\x\u\y),\\
  \label{involution}&\star\x\star{}=\x,\\
  \label{semi-Pasch}
  &\ex\x(\R\v\w\x\land\R\x\y\z)\implies\ex\u\R\v\u\z,\\
  \label{Pasch}
  &\ex\x(\R\v\w\x\land\R\x\y\z)\implies\ex\u(\R\v\u\z\land\R\w\y\u),\\
  \label{dense} &\R\x\x\x,\\
  \label{comm} &\R\x\y\z\iff\R\y\x\z,\\
  \label{symm} &\star\x=\x.
\end{align} 
When \eqref{involution} holds, the implication $\implies$ in
\eqref{left rotation}--\eqref{right reflection} can be replaced by an
equivalence $\iff$, since applying each of these frame conditions to
$\R\x\y\z$ two or three times produces $\R\star\x\star{}
\star\y\star{} \star\z\star{}$.  See Lemma \ref{lem8}\eqref{lem8i}
below.

Condition \eqref{Pasch} holds in a frame $\gc\K$ iff the operation
$\rp$ in the complex algebra $\Cm{\gc\K}$ satisfies the associative
law \eqref{BIV}, and \eqref{semi-Pasch} holds iff the semi-associative
law \eqref{BIV'} is true in $\Cm{\gc\K}$.  Note that \eqref{Pasch}
implies \eqref{semi-Pasch}, reflecting the fact that \eqref{BIV'} is a
special case of \eqref{BIV}. See Theorem \ref{thm7} below.

Conditions \eqref{left rotation}--\eqref{Pasch} hold in the pair-frame
on any set. Table \ref{triangles} illustrates some triples in the
ternary relation $\R$ of the pair-frame on the set $\U=\{1,2,3\}$ when
$x=\<1,2\>$, $y=\<2,3\>$, and $z=\<1,3\>$. Table \ref{triangles} can
be used to correlate each of the first five frame conditions with
permutations of $\{1,2,3\}$.  The triangle containing $\R\x\y\z$
illustrates the hypothesis of each frame condition. The vertices of
the triangle containing the conclusion of each condition are also
labelled with 1, 2, and 3. By matching up the vertices of the
conclusion with the vertices of the hypothesis one obtains a
permutation of $\{1,2,3\}$.  In cycle notation, the permutations match
up with the conditions in this way: $(1,2,3)$ with \eqref{left
  rotation}, $(1,3,2)$ with \eqref{right rotation}, $(1,3)$ with
\eqref{center reflection}, $(1,2)$ with \eqref{left reflection}, and
$(2,3)$ with \eqref{right reflection}. Successively applying the frame
conditions is the same as composing the correlated permutations.
\begin{table}
\setlength{\unitlength}{1.3mm}
\begin{align*}\thicklines
&\begin{picture}(30,20)(-5,-5)
 \put(    -1,-1    ){1}
 \put(   9.2,14.5  ){2}
 \put(  19.5,-1    ){3}
 \put(   1  , 1.5  ){\line( 2, 3){ 8}}
 \put(   1.5, 0    ){\line( 1, 0){17}}
 \put(  11  ,13.5  ){\line( 2,-3){ 8}}
 \put(  -5.5, 8    ){$\<1,2\>=\x$}
 \put(  15.5, 8    ){$\y=\<2,3\>$}
 \put(   6  ,-3    ){$\z=\<1,3\>$}
 \put(   7  , 3    ){$\R\x\y\z$}
\end{picture}
&
&\begin{picture}(30,20)(-5,-5)
 \put(   -1 ,-1    ){2}
 \put(  9.2 ,14.5  ){1}
 \put( 19.5 ,-1    ){3}
 \put(   1  , 1.5  ){\line( 2, 3){ 8}}
 \put(   1.5, 0    ){\line( 1, 0){17}}
 \put(  11  ,13.5  ){\line( 2,-3){ 8}}
 \put(  -5.5, 8    ){$\<2,1\>=\star\x$}
 \put(  15.5, 8    ){$\z=\<1,3\>$}
 \put(   6  ,-3    ){$\y=\<2,3\>$}
 \put(   7  , 3    ){$\R\star\x\z\y$}
\end{picture}
\\\thicklines
&\begin{picture}(30,20)(-5,-5)
 \put(    -1,-1    ){3}
 \put(   9.2,14.5  ){1}
 \put(  19.5,-1    ){2}
 \put(     1, 1.5  ){\line( 2, 3){ 8}}
 \put(   1.5, 0    ){\line( 1, 0){17}}
 \put(    11,13.5  ){\line( 2,-3){ 8}}
 \put(  -5.5, 8    ){$\<3,1\>=\star\z$}
 \put(  15.5, 8    ){$\x=\<1,2\>$}
 \put(   5  ,-3    ){$\star\y=\<3,2\>$}
 \put(   6.5, 3    ){$\R\star\z\x\star\y$}
\end{picture}
&
&\begin{picture}(30,20)(-5,-5)
 \put(    -1,-1    ){1}
 \put(   9.2,14.5  ){3}
 \put(  19.5,-1    ){2}
 \put(   1  , 1.5){\line( 2, 3){ 8}}
 \put(   1.5, 0  ){\line( 1, 0){17}}
 \put(  11  ,13.5){\line( 2,-3){ 8}}
 \put(  -5.5, 8  ){$\<1,3\>=\z$}
 \put(  15.5, 8  ){$\star\y=\<3,2\>$}
 \put(   6  ,-3  ){$\x=\<1,2\>$}
 \put(   7  , 3  ){$\R\z\star\y\x$}
\end{picture}
\\\thicklines
&\begin{picture}(30,20)(-5,-5)
 \put(    -1,-1    ){2}
 \put(   9.2,14.5  ){3}
 \put(  19.5,-1    ){1}
 \put(   1  , 1.5  ){\line( 2, 3){ 8}}
 \put(   1.5, 0    ){\line( 1, 0){17}}
 \put(  11  ,13.5  ){\line( 2,-3){ 8}}
 \put(  -5.5, 8    ){$\<2,3\>=\y$}
 \put(  15.5, 8    ){$\star\z=\<3,1\>$}
 \put(   5  ,-3    ){$\star\x=\<2,1\>$}
 \put(   6.5, 3    ){$\R\y\star\z\star\x$}
\end{picture}
&
&\begin{picture}(30,20)(-5,-5)
 \put(    -1,-1    ){3}
 \put(   9.2,14.5  ){2}
 \put(  19.5,-1    ){1}
 \put(   1  , 1.5  ){\line( 2, 3){ 8}}
 \put(   1.5, 0    ){\line( 1, 0){17}}
 \put(  11  ,13.5  ){\line( 2,-3){ 8}}
 \put(  -5.5, 8    ){$\<3,2\>=\star\y$}
 \put(  15.5, 8    ){$\star\x=\<2,1\>$}
 \put(   5  ,-3    ){$\star\z=\<3,1\>$}
 \put(   6  , 3    ){$\R\star\y\star\x\star\z$}
\end{picture}
\end{align*}
\caption{Triples in the ternary relation of the pair-frame 
  on $\{1,2,3\}$.}
\label{triangles}
\end{table}
The next lemma is an expression of the fact that the symmetric group
on a 3-element set is generated by a permutation of order 3 (2
choices) together with a permutation of order 2 (3 choices), and it is
also generated by any two permutations of order 2.
\begin{lemma}\label{five}
  Suppose a frame $\gc\K=\<\K,\R,\star{},\II\>$ satisfies
  \eqref{involution}.
 \begin{enumerate}
  \item\label{five1} If $\gc\K$ satisfies either \eqref{left rotation}
    or \eqref{right rotation} and any one of \eqref{center
      reflection}--\eqref{right reflection}, then it satisfies all
    five conditions \eqref{left rotation}--\eqref{right reflection}.
  \item\label{five2} If $\gc\K$ satisfies any two of \eqref{center
    reflection}--\eqref{right reflection} then it satisfies
    \eqref{left rotation}--\eqref{right reflection}.
 \end{enumerate}
\end{lemma}
The following lemma is useful in proving Theorems \ref{thm7} and
\ref{thm8} below.
\begin{lemma}\label{lem8}
  Assume $\gc\K=\<\K,\R,\star{},\II\>$ is a frame.
  \begin{enumerate}
  \item\label{lem8i} If \eqref{left reflection} and \eqref{identity}
    hold then \eqref{involution} holds and $\R\x\y\z\iff\R\star\x\y\z$
    for all $\x,\y,\z\in\K$.
  \item\label{lem8ii} If \eqref{left reflection}, \eqref{right
    reflection}, and \eqref{identity} hold, then $\v=\star\v$ for all
    $\v\in\II$.
  \item\label{lem8iii} Assume \eqref{left reflection}, \eqref{right
      reflection}, \eqref{identity}, and \eqref{semi-Pasch}. If
    $\x\in\K$, $\u,\v\in\II$, $\R\x\u\x$, and $\R\x\v\x$, then
    $\u=\v$.
  \end{enumerate}
\end{lemma}
\proof[\ref{lem8i}] Since $\x=\x$, by \eqref{identity} there must be
some $\u\in\II$ such that $\R\x\u\x$. Applying \eqref{left reflection}
twice yields $\R\star\x\star{}\u\x$, so we obtain $\star\x\star{}=\x$
by \eqref{identity}.  Since $\star{}$ is an involution, \eqref{left
  reflection} implies that its converse also holds.
\proof[\ref{lem8ii}] Assume $\v\in\II$. By \eqref{identity} there is
some $\u\in\II$ such that $\R\v\u\v$.  Then $\R\star\v\v\u$ by
\eqref{left reflection} and $\R\v\star\u\v$ by \eqref{right
  reflection}.  From $\R\v\star\u\v$ we also have $\R\star\v\v\star\u$
by \eqref{left reflection}.  From $\R\star\v\v\u$ and
$\R\star\v\v\star\u$ we obtain $\star\v=\u$ and $\star\v=\star\u$ by
\eqref{identity} since $\v\in\II$.  From $\star\v=\star\u$ we get
$\v=\u$ by part \eqref{lem8i}, so from $\star\v=\u$ we get
$\star\v=\v$, as desired.  \proof[\ref{lem8iii}] Assume $\u,\v\in\II$,
$\x\in\K$, $\R\x\u\x$, and $\R\x\v\x$.  Then $\R\star\x\x\u$ and
$\R\star\x\x\v$ by \eqref{left reflection}.  From $\R\star\x\x\u$ we
get $\R\u\star\x\star\x$ by \eqref{right reflection}.  Apply
\eqref{semi-Pasch} to $\R\u\star\x\star\x$ and $\R\star\x\x\v$,
obtaining some $\y\in\K$ such that $\R\u\y\v$. Then $\R\star\u\v\y$ by
\eqref{left reflection} so $\star\u=\y$ by \eqref{identity} since
$\v\in\II$.  We therefore have $\R\u\star\u\v$.  But $\star\u\in\II$
by part \eqref{lem8ii} since $\u\in\II$ by assumption, hence $\u=\v$
by \eqref{identity}.  \endproof The next theorem asserts that
conditions \eqref{dense}, \eqref{comm}, or \eqref{symm} hold in a
frame iff its complex algebra is dense, commutative, or symmetric,
respectively
\begin{theorem}\label{thm7-}
  Let $\gc\K=\<\K,\R,\star{},\II\>$ be a frame.
  \begin{enumerate}
  \item\label{thm7-ii} $\Cm{\gc\K}$ is dense iff $\gc\K$ satisfies
    \eqref{dense}.
  \item\label{thm7-i} $\Cm{\gc\K}$ is commutative iff $\gc\K$
    satisfies \eqref{comm}.
  \item\label{thm7-iii}$\Cm{\gc\K}$ is symmetric iff $\gc\K$ satisfies
    \eqref{symm}.
  \end{enumerate}
\end{theorem}
\proof For part \eqref{thm7-ii}, assume $\Cm{\gc\K}$ is dense and
$\x\in\K$.  Then $\{\x\}\leq\{\x\}\rp\{\x\}$ by density, hence
$\R\x\x\x$ and \eqref{dense} holds. Assume \eqref{dense} holds and
$\X\subseteq\K$.  For every $\x\in\X$, $\R\x\x\x$ implies
$\{\x\}\subseteq\{\x\}\rp\{\x\}\subseteq\X\rp\X$, hence
$\X\subseteq\X\rp\X$, which shows $\Cm{\gc\K}$ is dense.  Part
\eqref{thm7-iii} is equally easy. For part \eqref{thm7-i}, assume
$\Cm{\gc\K}$ is commutative, $\x,\y,\z\in\K$, and $\R\x\y\z$. Then
$\z\in\{\x\}\rp\{\y\}$, but $\{\x\}\rp\{\y\}=\{\y\}\rp\{\x\}$ since
$\Cm{\gc\K}$ is commutative, so $\R\y\x\z$. This shows that
\eqref{comm} holds. For the converse, assume \eqref{comm} holds and
let $\X,\Y\subseteq\K$ be elements of the complex algebra. By
\eqref{comm} and the definition of $\rp$ we have
\begin{align*}
  \X\rp\Y& =\{\z:\z\in\K, \R\x\y\z\text{ for some }\x\in\X,\y\in\Y\}\\
  &=\{\z:\z\in\K, \R\y\x\z\text{ for some }\x\in\X,\y\in\Y\} =\Y\rp\X,
\end{align*} 
so $\Cm{\gc\K}$ is commutative. \endproof The frames whose complex
algebras are semi-associative relation algebras or relation algebras
are characterized next. Because of Lemmas \ref{five} and \ref{lem8},
Theorems \ref{thm7} and \ref{thm8} remain true if \eqref{left
  reflection}--\eqref{identity} are replaced by \eqref{left
  rotation}--\eqref{involution}.
\begin{theorem}\label{thm7}{\rm\cite[Thm 2.2]{MR662049}}
  Let $\gc\K=\<\K,\R,\star{},\II\>$ be a frame.
  \begin{enumerate}
  \item\label{thm7iii} $\Cm{\gc\K}\in\NA$ iff $\gc\K$ satisfies
    \eqref{left reflection}--\eqref{identity}.
  \item\label{thm7i} $\Cm{\gc\K}\in\SA$ iff $\gc\K$ satisfies
    \eqref{left reflection}--\eqref{identity} and \eqref{semi-Pasch}.
  \item\label{thm7ii} $\Cm{\gc\K}\in\RA$ iff $\gc\K$ satisfies
    \eqref{left reflection}--\eqref{identity} and \eqref{Pasch}.
  \end{enumerate}
\end{theorem}
The J\'onsson-Tarski Representation Theorem \cite[Thm 3.10]{MR0044502}
in combination with Theorems \ref{thm7-} and \ref{thm7} produces a
representation theorem for $\NA$ and the subvarieties that can be
obtained by imposing semi-associativity, associativity, density,
commutativity, or symmetry.  For example, $\gc\A$ is a dense
commutative relation algebra iff $\gc\A$ is isomorphic to a subalgebra
of the complex algebra of a frame satisfying \eqref{left
  reflection}--\eqref{identity} and \eqref{Pasch}--\eqref{dense}.
\begin{theorem}\label{thm8}{\rm\cite[Thm 4.3]{MR662049}}
  Assume $\gc\A=\<\U,+,\min\blank,\rp,\conv{},\id\>$ is an algebra
  satisfying axioms \eqref{BI}--\eqref{BIII}, \eqref{BV},
  \eqref{BVII}--\eqref{BX}.  Then there is a frame $\gc\K =
  \<\K,\R,\star{},\II\>$ such that the following statements hold.
  \begin{enumerate}
  \item\label{thm8i} $\gc\A\cong\gc\A'\subseteq\Cm{\gc\K}$ for some
    subalgebra $\gc\A'$ of the complex algebra $\Cm{\gc\K}$.
  \item\label{thm8i.a} $\Cm{\gc\K}$ is an algebra satisfying axioms
    \eqref{BI}--\eqref{BIII}, \eqref{BV}, \eqref{BVII}--\eqref{BX}.
  \item\label{thm8ii-} $\Cm{\gc\K}\in\NA$ iff $\gc\A\in\NA$.
  \item\label{thm8iii-} $\gc\A\in\NA$ iff $\gc\K$ satisfies
    \eqref{left reflection}--\eqref{identity}.
  \item\label{thm8ii} $\Cm{\gc\K}\in\SA$ iff $\gc\A\in\SA$.
  \item\label{thm8iii} $\gc\A\in\SA$ iff $\gc\K$ satisfies \eqref{left
    reflection}--\eqref{identity} and \eqref{semi-Pasch}.
  \item\label{thm8iv} $\Cm{\gc\K}\in\RA$ iff $\gc\A\in\RA$.
  \item\label{thm8v} $\gc\A\in\RA$ iff $\gc\K$ satisfies \eqref{left
    reflection}--\eqref{identity} and \eqref{Pasch}.
  \item\label{thm8vii} $\gc\A$ is dense iff $\gc\K$ satisfies
    \eqref{dense}.
  \item\label{thm8vi} $\gc\A$ is commutative iff $\gc\K$ satisfies
    \eqref{comm}.
  \item\label{thm8viii}$\gc\A$ is symmetric iff $\gc\K$ satisfies
    \eqref{symm}.
  \end{enumerate}
\end{theorem}
\proof This theorem was originally derived from \cite[Thms 2.15, 2.17,
  2.18]{MR0044502}, but the desired frame $\gc\K = \<\K,\R,\star{},
\II\>$ can be obtained directly from the algebra $\gc\A$.  A subset
$\X\subseteq\U$ is an {\bf ultrafilter} of $\gc\A$ if $\X\neq\U$ and
for all $\x,\y\in\U$, if $\x\in\X$ then $\x+\y\in\X$ and if
$\x,\y\in\X$ then $\x\cdot\y\in\X$.  Let $\K$ be the set of
ultrafilters of $\gc\A$, and define $\R\subseteq\K^3$,
$\star{}\colon\K\to\K$, and $\II\subseteq\K$ by
\begin{align*}
  \R&=\{\<\X,\Y,\Z\>:\X,\Y,\Z\in\K,
  \{\x\rp\y:\x\in\X,\y\in\X\}\subseteq\Z\},\\
  \star\X&=\{\conv\x:\x\in\X\}\text{ for every }\X\in\K,\\
  \II&=\{\X:\id\in\X\in\K\}.
\end{align*}
Define $\varepsilon\colon\A\to\wp(\K)$ by $\varepsilon(\x) =
\{\X:\x\in\X\in\K\}$ for every $\x\in\A$.  Then $\varepsilon$ is an
isomorphic embedding of $\gc\A$ onto a subalgebra $\gc\A'$ of
$\Cm{\gc\K}$, as required by part \eqref{thm8i}.  This construction
was first described in R.\ McKenzie's dissertation \cite[Thm
  2.11]{MR2616328}.  \endproof The next theorem shows how frames arise
from groups.  Note that $|\X|$ is the number of elements in $\X$.
\begin{theorem}\label{lem8+}
  Suppose $\gc\G=\<\K,\rp,\star{},\e\>$ is a group.  Define a ternary
  relation $\R$ on $\K$ by $\R=\{\<\x,\y,\z\>:\x\rp\y=z\}$.  Then
  $\gc\K=\<\K,\R,\star{},\{\e\}\>$ is a frame satisfying \eqref{left
    rotation}--\eqref{Pasch}.
\end{theorem}
\proof We use the properties of groups, that $\rp$ is associative, and
for all $\x,y\in\K$, $\x\rp\e = \x = \e\rp\x$, $\star\x\star{} = \x$,
$\x\rp\star\x = \e = \star\x\rp\x$, and $\star{(\x\rp\y)} =
\star\y\rp\star\x$.  Frame conditions \eqref{left
  rotation}--\eqref{right reflection} all have the same assumption,
$\R\x\y\z$, and their conclusions are $\R\y\star\z\star\x$,
$\R\star\z\x\star\y$, $\R\star\y\star\x\star\z$, $\R\star\x\z\y$, and
$\R\z\star\y\x$.  According to the definition of $\R$, we assume
$\z=\x\rp\y$ and prove $\y\rp\star\z = \star\x $ for \eqref{left
  rotation} by $\y\rp\star\z = \y\rp(\x\rp\y)\star{}
=\y\rp(\star\y\rp\star\x) = (\y\rp\star\y)\rp\star\x = \e\rp\star\x
=\star\x$, $\star\z\rp\x =\star\y$ for \eqref{right rotation} by
$\star\z\rp\x = \star{(\x\rp\y)}\rp\x =(\star\y\rp\star\x)\rp\x =
\star\y\rp(\star\x\rp\x) = \star\y\rp\e =\star\y$,
$\star\y\rp\star\x=\star\z$ for \eqref{center reflection} by
$\star\y\rp\star\x= \star{(\x\rp\y)} =\star\z$, $\star\x\rp\z=\y$ for
\eqref{left reflection} by $\star\x\rp\z= \star\x\rp(\x\rp\y) =
(\star\x\rp\x)\rp\y = \e\rp\y =\y$, and $\z\rp\star\y=\x $ for
\eqref{right reflection} by $\z\rp\star\y= (\x\rp\y)\rp\star\y =
\x\rp(\y\rp\star\y) = \x\rp\e = \x$.

For \eqref{identity} in one direction, let $\x\in\K$.  We want
$\R\x\u\x$ for some $\u\in\K$. Take $\u=\e$, and get $\x\rp\e=\x$,
\ie, $\R\x\e\x$, as desired.  For the other direction, assume
$\R\x\e\y$. This implies $\x\rp\e=\y$, but $\x\rp\e=\x$, so
$\x=\y$. We have \eqref{involution} since $\star\x\star{} = \x$.  For
\eqref{Pasch}, we assume $\R\v\w\x$ and $\R\x\y\z$ and wish to show
$\R\v\u\z$ and $\R\w\y\u$.  From our hypotheses we have $\v\rp\w=\x$
and $\x\rp\y=\z$. Let $\u=\w\rp\y$.  Then $\R\w\y\u$ and
$\z=\x\rp\y=(\v\rp\w)\rp\y = \v\rp(\w\rp\y) = \v\rp\u$, hence
$\R\v\u\z$, as desired. Note that \eqref{semi-Pasch} is a trivial
consequence of \eqref{Pasch}.  \endproof The frame of a group
satisfies \eqref{dense} iff it has only one element.  Groups (treated
as frames) that satisfy \eqref{comm} are usually called {\bf Abelian
  groups}. The frame of a group satisfies \eqref{symm} iff every
element is its own inverse.  Such groups are called {\bf Boolean
  groups} \cite{MR1503469}.

$\CR$-frames were introduced by \citet[p.\,184]{MR0363789a}. A frame
$\gc\K=\<\K,\R,\star{},\II\>$ is a {\bf $\CR$-frame} if $\II=\{0\}$
and for all $\a,\b,\c,\d\in\K$,
\begin{enumerate}
\item[]p1. $\R0\a\b\iff\a=\b$,
\item[]p2. $\ex\x(R\a\b\x\land\R\x\c\d)\iff\ex\y(\R\a\c\y\land\R\y\b\d)$,
\item[]p3. $R\a\a\a$,
\item[]p4. $\star\a\star{}=\a$,
\item[]p5. $\R\a\b\c\iff\R\a\c^*\b^*$.
\end{enumerate}
Meyer and Routley define the logic $\CR$ as the sets of predicates
that are valid in all $\CR$-frames \cite[p.\,187]{MR0363789a} and they
characterize Anderson and Belnap's relevance logic $\RR$
\cite{MR0115902, MR0406756, MR1223997, MR141590, MR239942} as those
predicates in $\Pi^r$ that are in $\CR$ when \eqref{rmin} is taken as
the definition of $\rmin$ \cite[Translation Theorem,
p.\,190]{MR0363789a}. Thus,
\begin{align}
  \label{CR-def}
  \CR&=\{\A:\A\in\Pi^+,\text{$\A$ is valid in every $\CR$-frame}\},\\
  \label{R-char}
  \RR&=\CR\cap\Pi^r.
\end{align}
The logic $\KR$ can be obtained from $\RR$ by adding \eqref{symm-ax}
in Theorem \ref{thm11}; see \cite[\SS65.1.2]{MR1223997} (by Urquhart).
A frame $\gc\K$ is a {\bf $\KR$-frame} if it is a $\CR$-frame
satisfying \eqref{symm} \cite[p.\,350]{MR1223997}.  Urquhart added,
``A slight modification of the usual completeness proof for $\RR$
shows that $\KR$ is complete with respect to the class of all $\KR$
model structures.'' Thus,
\begin{align}\label{KR-def}
  \KR&=\{\A:\A\in\Pi^+,\text{$\A$ is valid in every $\KR$-frame}\}.
\end{align}
\begin{lemma}\label{lem8a}
  If $\gc\K=\<\K,\R,\star{},\II\>$ is a frame such that {\rm p1} holds
  for some $0\in\K$ and {\rm p2} holds then \eqref{comm} holds. Every
  $\CR$-frame satisfies \eqref{comm}.
\end{lemma}
\proof Assume $\R\a\b\c$.  We have $\R0\a\a$ by p1. Applying p2 to
$\R0\a\a$ and $\R\a\b\c$, we obtain some $\y\in\K$ such that $\R0\b\y$
and $\R\y\a\c$.  Then $\b=\y$ by p1 so $\R\b\a\c$.  \endproof When
\eqref{comm} holds we can restate p5 in two ways. By switching the
order of the first two entries in $\R\a\star\c\star\b$ in p5 one gets
\eqref{right rotation}. Switching the first two entries in $\R\a\b\c$
in p5 and interchanging $\a$ and $\b$ in the entire statement gives
\eqref{left rotation}.  In the presence of p1 and p2, the three
postulates p5, \eqref{left rotation}, and \eqref{right rotation} are
equivalent because of Lemma \ref{lem8a}.
\begin{theorem}\label{thm8a}
  Let $\gc\K=\<\K,\R,\star{},\{0\}\>$ be a frame.
  \begin{enumerate}
  \item\label{thm8ai} $\gc\K$ is a $\CR$-frame satisfying \eqref{left
      reflection} iff $\Cm{\gc\K}$ is a dense commutative relation
    algebra.
  \item\label{thm8aii} $\gc\K$ is a $\KR$-frame iff $\Cm{\gc\K}$ is a
    dense symmetric relation algebra.
  \end{enumerate}
\end{theorem}
\proof Assume $\Cm{\gc\K}$ is a dense commutative relation algebra. We
will show $\gc\K$ is a $\CR$-frame.  By Theorem
\ref{thm7}\eqref{thm7ii}, $\gc\K$ satisfies \eqref{left reflection},
\eqref{right reflection}, \eqref{identity}, and \eqref{Pasch}.
Postulate p1 is what \eqref{identity} reduces to when $\II=\{0\}$.
Since $\Cm{\gc\K}$ is dense, $\gc\K$ satisfies \eqref{dense} by
Theorem \ref{thm7-}\eqref{thm7-ii} and \eqref{dense} coincides with
p3.  Condition p4 follows from \eqref{identity} and \eqref{left
  reflection} by Lemma \ref{lem8}\eqref{lem8i}.  For p5, assume
$\R\a\b\c$.  By Theorem \ref{thm7-}\eqref{thm7-i}, $\gc\K$ satisfies
\eqref{comm} since $\Cm{\gc\K}$ is commutative, so $\R\b\a\c$.  Then
$\R\star\b\c\a$ by \eqref{left reflection}, so $\R\a\star\c\star\b$ by
\eqref{right reflection}.  Thus, p5 holds in one direction and the
other direction follows from this by p4.  Note that the implication in
p2 from right to left is formally identical to the implication from
left to right, so we need only prove the latter. Assume $\R\a\b\x$ and
$\R\x\c\d$.  By \eqref{Pasch}, $\R\a\u\d$ and $\R\b\c\u$ for some
$\u\in\K$, so $\R\c\b\u$ since \eqref{comm} holds.  We get
$\R\star\a\d\u$ from $\R\a\u\d$ by \eqref{left reflection} and
$\R\u\star\b\c$ from $\R\c\b\u$ by \eqref{right reflection}.  By
\eqref{Pasch}, there is some $\y\in\K$ such that $\R\star\a\y\c$ and
$\R\d\star\b\y$.  By \eqref{left reflection}, \eqref{right
  reflection}, and p4, $\R\a\c\y$ and $\R\y\b\d$.

Assume $\gc\K$ is a $\CR$-frame satisfying \eqref{left reflection}.
Then \eqref{comm} holds by Lemma \ref{lem8a} so $\Cm{\gc\K}$ is
commutative by Lemma \ref{thm7-}\eqref{thm7-i}.  $\Cm{\gc\K}$ is dense
by p3 and Lemma \ref{thm7-}\eqref{thm7-ii}, \eqref{identity} holds by
p1 since $\II=\{0\}$, and \eqref{right reflection} follows immediately
from \eqref{left reflection} and \eqref{comm}.  To prove \eqref{Pasch}
assume $\R\v\w\x$ and $\R\x\y\z$.  Then $\R\w\v\x$ by \eqref{comm}, so
by p2 applied to $\R\w\v\x$ and $\R\x\y\z$, there is some $\u\in\K$
such that $\R\w\y\u$ and $\R\u\v\z$, hence $\R\v\u\z$ by \eqref{comm}.
Since \eqref{left reflection}, \eqref{right reflection},
\eqref{identity}, and \eqref{Pasch} hold, we conclude that
$\Cm{\gc\K}\in\RA$ by Theorem \ref{thm7}\eqref{thm7ii}.

Suppose $\gc\K$ is a $\KR$-frame.  Since $\gc\K$ is a $\CR$-frame and
satisfies \eqref{dense} and \eqref{symm} by definition, $\Cm{\gc\K}$
is dense and symmetric by Theorem
\ref{thm7-}\eqref{thm7-ii}\eqref{thm7-iii}.  If $\R\a\b\c$ then
$\R\a\star\c\star\b$ by p5, hence $\R\star\a\c\b$ by \eqref{symm}.
Thus, \eqref{left reflection} holds and we conclude that $\Cm{\gc\A}$
is a relation algebra by part \eqref{thm8ai}.

For the converse, assume $\Cm{\gc\K}$ is a dense symmetric relation
algebra. Then $\gc\K$ satisfies \eqref{symm} by Theorem
\ref{thm7-}\eqref{thm7-iii}.  $\Cm{\gc\K}$ is commutative by Lemma
\ref{SA:symm->comm}, hence $\gc\K$ satisfies \eqref{comm} by Theorem
\ref{thm7-}\eqref{thm7-i}.  By part \eqref{thm8ai}, $\gc\K$ is a
$\CR$-frame. Because $\gc\K$ satisfies \eqref{symm}, $\gc\K$ is also a
$\KR$-frame.
\endproof
\section{The $n$-variable sequent calculus}
\label{sect11}
Assume $1\leq\n\in\omega$.  An {\bf $\n$-sequent} is an ordered pair
$\<\Gamma,\Delta\>$, written $\Gamma\sep\Delta$, of sets
$\Gamma,\Delta \subseteq \{\x\A\y : \x,\y\in\Upsilon_\n,
\,\A\in\Pi^+\}$. An $\n$-sequent $\Gamma\sep\Delta$ is an {\bf axiom}
if $\Gamma\cap\Delta \neq\emptyset$ or $\x\id\x\in\Delta$ for some
$\x\in\Upsilon_\n$.

Let $\Psi$ be a set of $\n$-sequents. A sequent is {\bf $\n$-provable
  from $\Psi$} (just {\bf $\n$-provable} when $\Psi=\emptyset$) if it
is contained in every set of $\n$-sequents that includes $\Psi$ and
the axioms and is closed under the rules of inference in Table
\ref{rules}.
\begin{table}
  \begin{align*}
    &\boxed{\text{Axiom}}\quad
  \begin{aligned}
    &\\\hline\stand \x\A\y,\Gamma&\sep\x\A\y,\Delta
  \end{aligned} 
  &&\boxed{\text{Cut}}\quad
  \begin{aligned}
    \Gamma&\sep\Delta,\x\A\y\\
    \x\A\y,\Gamma'&\sep\Delta'\\
    \hline\stand \Gamma,\Gamma'&\sep\Delta,\Delta'
  \end{aligned}
  \\\\
    &\boxed{\text{$\id|$}}\quad
  \begin{aligned}
    \x\A\y,\Gamma&\sep\Delta
    \\\hline\stand \x\A\z,\z\id\y,\Gamma&\sep\Delta
  \end{aligned} 
  &&\boxed{|\id}\quad
  \begin{aligned}
    &\\\hline\stand \Gamma&\sep\Delta,\x\id\x
  \end{aligned}
  \\\\
  &\boxed{+|}\quad
  \begin{aligned}
    \x\A\y,\Gamma&\sep\Delta\\
    \x\B\y,\Gamma'&\sep\Delta'\\
    \hline\stand\x{\A+\B}\y,\Gamma,\Gamma' &\sep\Delta,
    \Delta'\end{aligned} && \boxed{|+}\quad
  \begin{aligned}
    \Gamma&\sep\Delta,\x\A\y,\x\B\y\\
    \hline\stand \Gamma&\sep\Delta,\x{\A+\B}\y
  \end{aligned}
  \\\\
  &\boxed{\cdot|}\quad
  \begin{aligned}
    \Gamma,\x\A\y,\x\B\y&\sep\Delta\\
    \hline\stand \Gamma,\x{\A\cdot\B}\y&\sep\Delta
  \end{aligned}
  &&\boxed{|\cdot}\quad
  \begin{aligned}
    \Gamma&\sep\Delta,\x\A\y\\
    \Gamma'&\sep\Delta',\x\B\y\\
    \hline\stand\Gamma,\Gamma'
    &\sep\Delta,\Delta',\x{\A\cdot\B}\y
  \end{aligned}
  \\\\
  &\boxed{\min\blank|}\quad
  \begin{aligned}
    \Gamma&\sep\Delta,\x\A\y\\\hline\stand
    \x{\min\A}\y,\Gamma&\sep\Delta
  \end{aligned}
  &&\boxed{|\min\blank}\quad
  \begin{aligned}
    \x\A\y,\Gamma&\sep\Delta\\\hline\stand
    \Gamma&\sep\Delta,\x{\min\A}\y
  \end{aligned}
  \\\\
  &\boxed{\rp|}\quad
  \begin{aligned}
    \x\A\y,\y\B\z,\Gamma&\sep\Delta\\\hline\stand
    \x\A\rp\B\z,\Gamma&\sep\Delta,\text{no $\y$}
  \end{aligned}
  &&\boxed{|\rp}\quad
  \begin{aligned}
    \Gamma&\sep\Delta,\x\A\y\\
    \Gamma'&\sep\Delta',\y\B\z\\\hline\stand
    \Gamma,\Gamma'&\sep\Delta,\Delta',\x\A\rp\B\z
  \end{aligned}
\\\\
  &\boxed{\conv{}|}\quad
  \begin{aligned}
    \x\A\y,\Gamma&\sep\Delta\\\hline\stand
    \y\conv\A\x,\Gamma&\sep\Delta
  \end{aligned}
  &&\boxed{|\conv{}}\quad
  \begin{aligned}
    \Gamma&\sep\Delta,\x\A\y\\\hline\stand
    \Gamma&\sep\Delta,\y\conv\A\x
  \end{aligned}
\end{align*}
\caption{Axioms and rules of inference for the $n$-variable sequent
  calculus, where \newline $1\leq n\leq\omega$, $x,y,z\in\Upsilon_n$,
  $A,B\in\Pi^+$, and $\Gamma,\Delta\subseteq\{xAy:x,y\in\Upsilon_n,\,
  A\in\Pi^+\}$.}
\label{rules}
\end{table}
In the rules in Table \ref{rules}, $\Gamma$, $\Gamma'$, $\Delta$, and
$\Delta'$ are sets of formulas in $\{\x\A\y : \x,\y\in\Upsilon_\n,
\,\A\in\Pi^+\}$, $\A,\B\in\Pi^+$ are predicates, and $\x,\y,\z
\in\Upsilon_\n$.  The notation ``no $\y$'' in rule $\rp|$ means that
$\y\neq\x,\z$ and $\y$ does not occur in any formula in $\Gamma$ or
$\Delta$.

The rules are taken from \cite{MR722170}.  The rules $|\cdot$ and
$\cdot|$ are derived from the rules for $\min\blank$ and $+$ through
the definition of $\cdot$ in \eqref{dotmin}.  Braces and union symbols
are frequently omitted from the notation for sequents in favor of
commas. For example, we write $\Gamma,\x\A\y,\x\B\y$ instead of
$\Gamma \cup \{\x\A\y\} \cup \{\x\B\y\}$.

An {\bf $\n$-proof from $\Psi$} is a sequence of sequents in which
every sequent is either in $\Psi$, or is an axiom, or follows from one
or two previous sequents in the sequence by one of the rules of
inference in Table \ref{rules}.  Whenever a rule is applied in an
$\n$-proof we include the numbers for the previous sequents used by
the rule and the name of the rule. For every application of the rule
$|\rp$ we also include notation of the form ``no $\y$'' as a reminder
that the eliminated variable $\y$ must not occur in the conclusion of
$\rp|$.

If $\Gamma\sep\Delta$ is an $\n$-sequent, then an {\bf $\n$-proof of
  $\Gamma\sep\Delta$ from $\Psi$} is an $\n$-proof from $\Psi$ in
which $\Gamma\sep\Delta$ occurs. It is straightforward to prove that
an $\n$-sequent is $\n$-provable from $\Psi$ iff it has an $\n$-proof
from $\Psi$.

The sequent calculus is connected with classes of algebras defined in
\cite{MR722170} and originally called $\MA_\n$, later renamed $\RA_\n$
\cite[Def 4]{MR1011183}.  
\begin{definition}\label{basis}
  Assume $\gc\A=\<\U, +, \min\blank, \rp, \conv{}, \id\>\in\NA$ is a
  non-associative relation algebra.  For $3\leq\n\leq\omega$, an {\bf
    $\n$-dimensional relational basis} for $\gc\A$ is a set
  $\B\subseteq\atoms^{\n\times\n}$ of $\n\times\n$ matrices of atoms
  of $\gc\A$ such that
  \begin{enumerate}
  \item\label{basis1} for all $\i,\j,\k<\n$ and all $\x,\y,\z\in\B$,
    $\x_{\i\i}\leq\id$, $\conv{\x_{\i\j}}=\x_{\j\i}$, and
    $\x_{\i\k}\leq\x_{\i\j}\rp\x_{\j\k}$,
  \item\label{basis2} for every $\a\in\atoms$ there is some $\x\in\B$
    such that $\x_{01}=\a$,
  \item\label{basis3} if $\i,\j<\n$, $\x\in\B$, $\a,\b\in\atoms$,
    $\x_{\i\j}\leq\a\rp\b$, and $\i,\j\neq\k<\n$, then there is some
    $\y\in\B$ such that $\y_{\i\k}=\a$, $\y_{\k\j}=\b$, and
    $\x_{\ell\m}=\y_{\ell\m}$ whenever $\k\neq\ell,\m<\n$.
  \end{enumerate}
  The algebra $\gc\A$ is an {\bf $\n$-dimensional relation algebra} if
  $\gc\A$ is a subalgebra of an atomic semi-associative relation
  algebra that has an $\n$-dimensional relational basis. $\RA_\n$ is
  the class of $\n$-dimensional relation algebras.  
\end{definition}
Part \eqref{key0.1} of the next theorem is proved by the equational
axiomatizations of $\RA_n$ in \cite[Thms 414, 419]{MR2269199} and
\cite[\SS13.8]{MR1935083}.
\begin{theorem}\label{key}{\rm\cite[Thms 2, 3, 6, 9, 10]{MR722170}}
  \begin{enumerate}
  \item\label{key0} If $3\leq\m\leq\n\leq\omega$ then
    $\RA_\m\supseteq\RA_\n$.
  \item\label{key0.1} If $3\leq\n\leq\omega$ then $\RA_\n$ is a
    variety.
  \item\label{key1} $\RA_3=\SA$.
  \item\label{key2} $\RA_4=\RA$.
  \item\label{key3} $\RA_\omega=\RRA=\bigcap_{3\leq\n<\omega}\RA_\n$.
  \item\label{key4} Assume $3\leq\n\leq\omega$,
    $\mathcal\R\subseteq(\Pi^+)^2$, and $\A,\B\in\Pi^+$.  Then the
    following conditions are equivalent.
    \begin{enumerate}
    \item $\fm\A\0\1\sep\fm\B\0\1$ is $\n$-provable from
      $\{\fm\X\0\1\sep\fm\Y\0\1:\<\X,\Y\>\in\mathcal\R\}$.
    \item For every $\gc\A\in\RA_\n$ and homomorphism
      $\h\colon\gc\P\to\gc\A$, if $\h(\X)\leq\h(\Y)$ whenever
      $\<\X,\Y\>\in\mathcal\R$, then $\h(\A)\leq\h(\B)$.
    \end{enumerate}
  \end{enumerate}
\end{theorem}
The following lemma is used extensively in the proof of Theorem
\ref{thm11} below.
\begin{lemma}\label{equiv}
  Assume $\A,\B\in\Pi^+$, $1\leq\n\in\omega$, $\i,\j,\k\leq\n$, and
  $\Psi$ is a set of $\n$-sequents. In the six derived rules below,
  any $\n$-proof that contains the sequent above the horizontal line
  can be extended to an $\n$-proof that contains the one below the
  line.
  \begin{align}
    \label{equiv-1}
    &\begin{aligned} \fm\A\vi\vj&\sep\fm\B\vi\vj\\ \hline\stand
       &\sep\fm{\A\to\B}\vj\vj
    \end{aligned}\text{ if }\i\neq\j\\
  \label{equiv-2}
    &\begin{aligned} &\sep\fm{\A\to\B}\vj\vj\\ \hline\stand
       \fm\A\vi\vj&\sep\fm\B\vi\vj
  \end{aligned}\\
  \label{equiv-3} 
  &\begin{aligned} &\sep\fm\A\vi\vi\\ \hline\stand
     \fm\id\vj\vk&\sep\fm\A\vj\vk
  \end{aligned}\text{ if }\i\neq\j,\k\\
  \label{equiv-4}
  &\phantom{\fm\A\vi\vi\,}
  \begin{aligned}
    &\sep\fm\A\vi\vi\\ \hline\stand &\sep\fm\A\vj\vj
  \end{aligned}\\
  \label{equiv-5}
  &\begin{aligned} \fm\A\vi\vj&\sep\fm\B\vi\vj\\ \hline\stand
     \fm\A\vj\vi&\sep\fm\B\vj\vi
  \end{aligned}\\
  \label{equiv-6}
  &\begin{aligned} \fm\id\vi\vj&\sep\fm\A\vi\vj\\ \hline\stand
     &\sep\fm\A\vi\vi
  \end{aligned}\text{ if }\i\neq\j\\
  \label{equiv-7} 
  &\begin{aligned} \fm\A\vi\vj&\sep\fm\B\vi\vj\\ \hline\stand
     \fm\A\vi\vi&\sep\fm\B\vi\vi
  \end{aligned}
\end{align}
\end{lemma}
\proof[\ref{equiv-1}]
\begin{align*}
  1.&&\fm\A\vi\vj&\sep\fm\B\vi\vj
  &&\text{in an $\n$-proof from $\Psi$}\\
  2.&&\fm{\conv\A}\vj\vi&\sep\fm\B\vi\vj
  &&\text{1, $\conv{}|$}\\
  3.&&\fm{\conv\A}\vj\vi,\fm{\min\B}\vi\vj&\sep
  &&\text{2, $\min\blank|$}\\
  4.&&\fm{\conv\A\rp\min\B}\vj\vj&\sep
  &&\text{3, $\rp|$, no $\vi$, $\i\neq\j$}\\
  5.&&&\sep\fm{\min{\conv\A\rp\min\B}}\vj\vj
  &&\text{4, $|\min\blank$}\\
  6.&&&\sep\fm{\A\to\B}\vj\vj &&\text{5, \eqref{to}}
\end{align*}
\proof[\ref{equiv-2}]
\begin{align*}
  1.&&&\sep\fm{\A\to\B}\vj\vj
  &&\text{in an $\n$-proof from $\Psi$}\\
  2.&&&\sep\fm{\min{\conv\A\rp\min\B}}\vj\vj
  &&\text{1, \eqref{to}}\\
  3.&&\fm{\conv\A\rp\min\B}\vj\vj&\sep\fm{\conv\A\rp\min\B}\vj\vj
  &&\text{axiom}\\
  4.&&\fm{\conv\A\rp\min\B}\vj\vj,
  \fm{\min{\conv\A\rp\min\B}}\vj\vj&\sep
  &&\text{3, $\min\blank|$}\\
  5.&&\fm{\conv\A\rp\min\B}\vj\vj&\sep
  &&\text{2, 4, Cut}\\
  6.&&\fm\A\vi\vj&\sep\fm\A\vi\vj
  &&\text{axiom}\\
  7.&&\fm\A\vi\vj&\sep\fm{\conv\A}\vj\vi
  &&\text{6, $|\conv{}$}\\
  8.&&\fm\B\vi\vj&\sep\fm\B\vi\vj
  &&\text{axiom}\\
  9.&&&\sep\fm\B\vi\vj,\fm{\min\B}\vi\vj
  &&\text{8, $|\min\blank$}\\
  10.&&\fm\A\vi\vj&\sep\fm\B\vi\vj,\fm{\conv\A\rp\min\B}\vj\vj
  &&\text{7, 9, $|\rp$}\\
  11.&&\fm\A\vi\vj&\sep\fm\B\vi\vj &&\text{5, 10, Cut}
\end{align*}
\proof[\ref{equiv-3}]
\begin{align*}
  1.&&&\sep\fm\A\vi\vi&&\text{in an $\n$-proof from $\Psi$}\\
  2.&&&\sep\fm{\conv\A}\vi\vi&&\text{1, $|\conv{}$}\\
  3.&&\fm\A\vj\vk&\sep\fm\A\vj\vk&&\text{axiom}\\
  4.&&\fm\A\vj\vi,\fm\id\vi\vk&\sep\fm\A\vj\vk&&\text{3, $\id|$}\\
  5.&&\fm{\conv\A}\vi\vj,\fm\id\vi\vk&\sep\fm\A\vj\vk&&\text{4, $\conv{}|$}\\
  6.&&\fm{\conv\A}\vi\vi,\fm\id\vi\vj,\fm\id\vi\vk&\sep\fm\A\vj\vk
  &&\text{5, $\id|$}\\
  7.&&\fm\id\vi\vj,\fm\id\vi\vk&\sep\fm\A\vj\vk&&\text{2, 6, Cut}\\
  8.&&\fm{\conv\id}\vj\vi,\fm\id\vi\vk&\sep\fm\A\vj\vk&&\text{7, $\conv{}|$}\\
  9.&&\fm{\conv\id\rp\id}\vj\vk&\sep\fm\A\vj\vk
  &&\text{8, $\rp|$, no $\vi$, $\i\neq\j,\k$}\\
  10.&&\fm\id\vj\vk&\sep\fm\id\vj\vk&&\text{axiom}\\
  11.&&&\sep\fm\id\vj\vj&&\text{$|\id$}\\
  12.&&&\sep\fm{\conv\id}\vj\vj&&\text{11, $|\conv{}$}\\
  13.&&\fm\id\vj\vk&\sep\fm{\conv\id\rp\id}\vj\vk&&\text{10, 12, $|\rp$}\\
  14.&&\fm\id\vj\vk&\sep\fm\A\vj\vk&&\text{9, 13, Cut}
\end{align*}
\proof[\ref{equiv-4}] This part is trivial if $\i=\j$ so assume
$\i\neq\j$.
\begin{align*}
  1.&&&\sep\fm\A\vi\vi&&\text{in an $\n$-proof from $\Psi$}\\
  2.&&\fm\id\vj\vj&\sep\fm\A\vj\vj
  &&\text{\eqref{equiv-3}, $\i\neq\j$}\\
  3.&&&\sep\fm\id\vj\vj&&\text{$|\id$}\\
  4.&&&\sep\fm\A\vj\vj&&\text{2, 3, Cut}
\end{align*}
\proof[\ref{equiv-5}] This part is trivial if $\i=\j$, so assume
$\i\neq\j$.
\begin{align*}
  1.&&\fm\A\vi\vj&\sep\fm\B\vi\vj&&\text{in an $\n$-proof from $\Psi$}\\
  2.&&&\sep\fm{\A\to\B}\vj\vj&&\text{\eqref{equiv-1}, $\i\neq\j$}\\
  2.&&&\sep\fm{\A\to\B}\vi\vi&&\text{\eqref{equiv-4}}\\
  3.&&\fm\A\vj\vi&\sep\fm\B\vj\vi&&\text{\eqref{equiv-2}}
\end{align*}
\proof[\ref{equiv-6}]
\begin{align*}
  1.&&\fm\id\vi\vj&\sep\fm\A\vi\vj&&\text{in an $\n$-proof}\\
  2.&&\fm\id\vj\vi&\sep\fm\id\vj\vi&&\text{axiom}\\
  3.&&\fm\id\vi\vj,\fm\id\vj\vi&\sep\fm{\A\rp\id}\vi\vi&&\text{1, 2, $|\rp$}\\
  4.&&\fm{\id\rp\id}\vi\vi&\sep\fm{\A\rp\id}\vi\vi
  &&\text{3, $\rp|$, no $\vj$}\\
  5.&&&\sep\fm\id\vi\vi&&\text{$|\id$}\\
  6.&&&\sep\fm{\id\rp\id}\vi\vi&&\text{5, $|\rp$}\\
  7.&&&\sep\fm{\A\rp\id}\vi\vi&&\text{4, 6, Cut}\\
  8.&&\fm\A\vi\vi&\sep\fm\A\vi\vi&&\text{axiom}\\
  9.&&\fm\A\vi\vj,\fm\id\vj\vi&\sep\fm\A\vi\vi&&\text{8, $\id|$}\\
  10.&&\fm{\A\rp\id}\vi\vi&\sep\fm\A\vi\vi&&\text{9, $|\rp$, no $\vj$}\\
  11.&&&\sep\fm\A\vi\vi&&\text{7, 10, Cut}
\end{align*}
\proof[\ref{equiv-7}]
\begin{align*} 
  1.&&\fm\A\vi\vj&\sep\fm\B\vi\vj
  &&\text{in an $\n$-proof}\\
  2.&&\fm\A\vi\vi,\fm\id\vi\vj&\sep\fm\B\vi\vj
  &&\text{1, $\id|$}\\
  3.&&\fm\id\vj\vi&\sep\fm\id\vj\vi
  &&\text{axiom}\\
  4.&&\fm\A\vi\vi,\fm\id\vi\vj,\fm\id\vj\vi&\sep\fm{\B\rp\id}\vi\vi
  &&\text{2, 3, $|\rp$}\\
  5.&&\fm\A\vi\vi,\fm{\id\rp\id}\vi\vi&\sep\fm{\B\rp\id}\vi\vi
  &&\text{4, $\rp|$, no $\vj$}\\
  6.&&&\sep\fm\id\vi\vi
  &&\text{$|\id$}\\
  7.&&&\sep\fm{\id\rp\id}\vi\vi
  &&\text{6, $\rp$}\\
  8.&&\fm\A\vi\vi&\sep\fm{\B\rp\id}\vi\vi
  &&\text{5, 7, Cut}\\
  9.&&\fm\B\vi\vi&\sep\fm\B\vi\vi
  &&\text{axiom}\\
  10.&&\fm\B\vi\vj,\fm\id\vj\vi&\sep\fm\B\vi\vi
  &&\text{9, $\id|$}\\
  11.&&\fm{\B\rp\id}\vi\vi&\sep\fm\B\vi\vi
  &&\text{10, $\rp|$, no $\vj$}\\
  12.&&\fm\A\vi\vi&\sep\fm\B\vi\vi
  &&\text{8, 11, Cut}\\
\end{align*}
\endproof Rules \eqref{equiv-1} and \eqref{equiv-2} do not involve
$\id$ in their statements and and do not use rules $|\id$ and $\id|$
in their proofs.  Rules \eqref{equiv-4}, \eqref{equiv-5}, and
\eqref{equiv-7} do not involve $\id$ in their statements but rules
$|\id$ and $\id|$ are used in their proofs.  There are alternate
proofs of \eqref{equiv-4} and \eqref{equiv-5} that avoid the use of
rules $|\id$ and $\id|$.  Observe that if the variables $\vi$ and
$\vj$ are interchanged in an axiom or a rule in Table \ref{rules}, the
result is still an axiom or rule. It follows that if variables $\vi$
and $\vj$ are interchanged throughout an $\n$-proof, the result is
still an $\n$-proof.  Consequently, if $\sep\fm\A\vi\vi$,
$\fm\A\vi\vj\sep \fm\B\vi\vj$, or $\fm\A\vi\vj\sep\fm\B\vj\vi$ is
$\n$-provable, so is $\sep\fm\A\vj\vj$, $\fm\A\vj\vi\sep\fm\B\vj\vi$,
or $\fm\A\vj\vi\sep\fm\B\vi\vj$, respectively.
\section{Characterizing $\CT^\Psi_3$, $\CT^\Psi_4$, $\TR$, $\KR$,
and $\CT^\Psi_\omega$}
\label{sect13}
Let $\hom$ be the set of endomorphisms of the predicate algebra
$\gc\P$, that is, homomorphisms from $\gc\P$ to itself.  For every
$\Psi\subseteq\Sigma^\times$, let
\begin{equation*}
  \Xi(\Psi) = \{\g(\P)\equals\g(\Q) : \P\equals\Q\in\Psi, 
  \P,\Q\in\Pi^+, \g\in\hom\}.
\end{equation*}
A {\bf $\Psi$-algebra} is an algebra similar to $\gc\P$ in which every
equation in $\Psi$ is true.
\begin{lemma}\label{lem5}
  Let $\A,\B\in\Pi^+$ and $\Psi\subseteq\Sigma^\times$. Let $\VV$ a
  class of algebras similar $\gc\P$ that is closed under subalgebras.
  Then the following statements are equivalent.
  \begin{enumerate}
  \item\label{l1} If $\h\colon\gc\P\to\gc\A$ is a homomorphism and
    $\gc\A\in\VV$ then $\gc\A\models_\h\Xi(\Psi)$ implies
    $\gc\A\models_\h\A\equals\B$.
  \item\label{l2} $\A\equals\B$ is true in every $\Psi$-algebra in
    $\VV$.
  \end{enumerate}
\end{lemma}
\proof Assume \eqref{l2}. To prove \eqref{l1}, assume $\gc\A\in\VV$,
$\h\colon\gc\P\to\gc\A$ is a homomorphism, and $\gc\A\models_\h
\Xi(\Psi)$.  We must show $\h(\A)=\h(\B)$.  Let $\gc\B$ be the
subalgebra of $\gc\A$ whose universe is $\h(\Pi^+) = \{\h(\X) :
\X\in\Pi^+\}$.  Note that $\gc\B\in\VV$ by our assumption on $\VV$.
We will show $\gc\B$ is a $\Psi$-algebra in $\VV$. Consider any
homomorphism $\f\colon\gc\P\to\gc\B$.  Construct an endomorphism
$\g\colon\gc\P\to\gc\P$ as follows.  For every propositional variable
$\id\neq\C\in\Pi$, since $\f(\C)\in\h(\Pi^+)$ we may choose
$\D\in\Pi^+$ such that $\h(\D)=\f(\C)$. Set $\g(\C)=\D$ so that
$\h(\g(\C))=\f(\C)$.  Make such a choice for every propositional
variable.  Since $\gc\P$ is absolutely freely generated by the
propositional variables, these choices extend to the desired
endomorphism $\g\in\hom$.  Since $\h(\g(\C))=\f(\C)$ for every
propositional variable, the properties of the homomorphisms $\f,\g,\h$
imply that $\h(\g(\E))=\f(\E)$ for every $\E\in\Pi^+$.  For every
equation $\C\equals\D\in\Psi$, $\g(\C)\equals\g(\D) \in \Xi(\Psi)$, so
our assumption $\gc\A\models_\h \Xi(\Psi)$ implies $\h(\g(\C)) =
\h(\g(\D))$, hence $\f(\C)=\f(\D)$. This completes the proof that
$\gc\B$ is a $\Psi$-algebra in $\VV$. Since $\A\equals\B$ is true in
every $\Psi$-algebra in $\VV$ and the homomorphism $\h$ maps $\gc\P$
into $\gc\B\in\VV$, it follows that $\h(\A)=\h(\B)$.

For the converse we assume \eqref{l1} and wish to show $\A\equals\B$
is true in every $\Psi$-algebra in $\VV$.  Assume $\gc\A\in\VV$ is a
$\Psi$-algebra. To show $\A\equals\B$ is true in $\gc\A$, we consider
an arbitrary homomorphism $\h:\gc\P\to\gc\A$ and want to show
$\h(\A)=\h(\B)$.  Consider an equation in $\Xi(\Psi)$. It has the form
$\g(\C)\equals\g(\D)$ for some $\g\in\hom$ and some $\C\equals\D
\in\Psi$.  The composition $\h\circ\g$ of $\h$ and $\g$ is also a
homomorphism from $\gc\P$ to $\gc\A$, but $\gc\A$ is a $\Psi$-algebra
in $\VV$, so $\h(\g(\C)) = (\h\circ\g)(\C) = (\h\circ\g)(\D) =
\h(\g(\D))$, \ie, $\gc\A \models_\h \g(\C) \equals \g(\D)$. This shows
$\gc\A\models_\h\Xi(\Psi)$, so by our hypothesis that \eqref{l1}
holds, we conclude that $\gc\A \models_\h \A\equals\B$, \ie,
$\h(\A)=\h(\B)$.  Since $\h$ was an arbitrary homomorphism, we have
shown that $\A\equals\B$ is true in $\gc\A$.
\endproof
If $\Psi\subseteq\Sigma^\times$ is a set of equations and
$1\leq\n\leq\omega$, then the set of {\bf $\n$-sequents corresponding
  to} $\Psi$ consists of all sequents of the form
$\fm\A{\var\i}{\var\j},\Gamma \sep \fm\B{\var\i}{\var\j},\Delta$ or
$\fm\B{\var\i}{\var\j},\Gamma \sep \fm\A{\var\i}{\var\j},\Delta$, where
$\i,\j\leq\n$, $\A,\B\in\Pi^+$, $\A\equals\B\in\Psi$, and $\Gamma$ and
$\Delta$ are sets of $\n$-sequents.
\begin{theorem}\label{thm9-} Characterizations of $\CT^\Psi_3$ and $\CT_3$.
  \begin{enumerate}
  \item\label{thm9-1} If $\Psi\subseteq\Sigma^\times$, then the
    following statements are equivalent for every $\A\in\Pi^+$.
    \begin{enumerate}
    \item\label{thm9-i} $\A\in\CT^\Psi_3$.
    \item\label{thm9-ii} $\Psi\proves^+_s\id\leq\A$.
    \item\label{thm9-iii} $\Psi\proves^\times_s\id\leq\A$.
    \item\label{thm9-iv} If $\h\colon\gc\P\to\gc\A$ is a homomorphism
      and $\gc\A\in\SA$ then $\gc\A\models_\h\Psi$ implies
      $\gc\A\models_\h\id\leq\A$.
    \item\label{thm9-v} If $\gc\K$ is a frame satisfying \eqref{left
        reflection}, \eqref{right reflection}, \eqref{identity}, and
      \eqref{semi-Pasch}, $\h\colon\gc\P\to\gc\A$ is a homomorphism,
      and $\Cm{\gc\K} \models_\h \Psi$, then $\h(\id)\subseteq\h(\A)$.
    \item\label{thm9-viii} $\sep\fm\A\0\0$ is 3-provable from the
      sequents corresponding to $\Psi$.
    \item\label{thm9-vi} $\{\GG(\varepsilon):\varepsilon\in\Psi\}
      \proves_s\all\0\GG(\0\A\0)$.
    \end{enumerate}
  \item\label{thm9-2} The following statements are equivalent for
    every $\A\in\Pi^+$.
    \begin{enumerate}
    \item\label{thm9-i-e} $\A\in\CT_3$.
    \item\label{thm9-ii-e} $\proves^+_s\id\leq\A$.
    \item\label{thm9-iii-e} $\proves^\times_s\id\leq\A$.
    \item\label{thm9-iv-e} $\id\leq\A$ is true in every
      semi-associative relation algebra.
    \item\label{thm9-v-e} $\id\leq\A$ is valid in every frame
      satisfying \eqref{left reflection}, \eqref{right reflection},
      \eqref{identity}, and \eqref{semi-Pasch}.
    \item\label{thm9-viii-e} $\sep\fm\A\0\0$ is 3-provable.
    \item\label{thm9-vi-e} $\proves_s\all\0\GG(\0\A\0)$.
    \end{enumerate}
  \end{enumerate}
\end{theorem}
\proof Parts \eqref{thm9-i} and \eqref{thm9-ii} are equivalent by
definition.  Parts \eqref{thm9-ii} and \eqref{thm9-iii} are equivalent
by Theorem \ref{thm5}\eqref{thm5iv}.  Parts \eqref{thm9-iii} and
\eqref{thm9-iv} are equivalent by Theorem \ref{eq}\eqref{eq2}.

To prove that parts \eqref{thm9-iv} and \eqref{thm9-v} are equivalent,
first assume \eqref{thm9-iv}, that $\gc\A\models_\h\Psi$ implies
$\gc\A\models_\h\id\leq\A$ whenever $\gc\A\in\SA$ and
$\h\colon\gc\P\to\gc\A$ is a homomorphism.  Suppose the frame $\gc\K$
satisfies the four conditions and that $\Cm{\gc\K} \models_\h \Psi$
for some homomorphism $\h\colon \gc\P \to \Cm{\gc\K}$.  We wish to
show $\h(\id)\subseteq\h(\A)$.  Since $\Cm{\gc\K}\in\SA$ by Theorem
\ref{thm7}\eqref{thm7i}, our hypothesis on $\A$ tells us that
$\gc\A\models_\h\id\leq\A$. By the definition of $\models_\h$, we get
$\h(\id)\subseteq\h(\A)$, as desired. Thus, \eqref{thm9-v} holds.

For the converse, assume \eqref{thm9-v}, that $\k(\id)\subseteq\k(\A)$
whenever $\gc\K$ is a frame satisfying \eqref{left reflection},
\eqref{right reflection}, \eqref{identity}, and \eqref{semi-Pasch},
and $\k\colon\gc\P\to\Cm{\gc\K}$ is a homomorphism such that
$\Cm{\gc\K} \models_\k \Psi$.  Assume the hypotheses of
\eqref{thm9-iv}, that $\gc\A\in\SA$ and $\gc\A\models_\h\Psi$ for some
homomorphism $\h\colon\gc\P\to\gc\A$.  We wish to show
$\gc\A\models_\h\id\leq\A$.  By Theorem
\ref{thm8}\eqref{thm8i}\eqref{thm8ii}\eqref{thm8iii} there is a frame
$\gc\K$ such that $\Cm{\gc\K}\in\SA$, $\gc\K$ satisfies the four
conditions, and $\gc\A$ is isomorphic to a subalgebra of $\Cm{\gc\K}$.
By composing $\h$ with the isomorphism from $\gc\A$ into $\Cm{\gc\K}$
we get a homomorphism $\k\colon\gc\P\to\Cm{\gc\K}$.  If
$\B\equals\C\in\Psi$ then $\h(\B)=\h(\C)$ since $\gc\A\models_\h\Psi$,
and this equality is preserved under the isomorphism from $\gc\A$ into
$\Cm{\gc\K}$, hence $\k(\B)=\k(\C)$.  This shows that
$\Cm{\gc\K}\models_\k\Psi$.  All the conditions are now met for
concluding from our hypotheses on $\A$ that $\k(\id)\subseteq\k(\A)$.
By applying the inverse of the isomorphism from $\gc\A$ into
$\Cm{\gc\K}$ to both sides of this equation, we get back to
$\h(\id)\subseteq\h(\A)$.  By definition of $\models$, this means that
$\gc\A\models_\h\id\leq\A$, as desired.

Parts \eqref{thm9-iv} and \eqref{thm9-viii} are equivalent by Theorem
\ref{key}\eqref{key1}\eqref{key4} and Lemma \ref{equiv}.  The
equivalence of parts \eqref{thm9-ii} and \eqref{thm9-vi} follows from
Theorem \ref{thm5}\eqref{thm5ii} when $\varphi$ is $\id\leq\A$,
together with the observation that $\GG(\id\leq\A) \equiv^+_s
\all\0\GG(\0\A\0)$ by Theorem \ref{thm5}\eqref{thm5i}. This completes
the proof of part \eqref{thm9-1}.

For part \eqref{thm9-2} it is enough to note that the statements
\eqref{thm9-i}--\eqref{thm9-vi} are equivalent to the corresponding
statements \eqref{thm9-i-e}--\eqref{thm9-vi-e} when $\Psi=\emptyset$.
This is true by notational convention in all but two cases. Once we
know that \eqref{thm9-iv} and \eqref{thm9-iv-e} are equivalent and
that \eqref{thm9-v} and \eqref{thm9-v-e} are equivalent when
$\Psi=\emptyset$, we get the equivalence of
\eqref{thm9-i-e}--\eqref{thm9-vi-e} from the equivalence of
\eqref{thm9-i}--\eqref{thm9-vi}.

Part \eqref{thm9-iv} coincides with part \eqref{l1} in Lemma
\ref{lem5} when $\VV=\SA$ and $\id\leq\A$ replaces $\A\equals\B$.
Applying Lemma \ref{lem5} with $\psi=\emptyset$, we conclude that
\eqref{thm9-iv} holds iff $\id\leq\A$ is true in every
$\emptyset$-algebra in $\SA$. Every equation in $\emptyset$ is true in
every algebra, so the latter statement simply says that $\id\leq\A$ is
true in every semi-associative relation algebra, that is,
\eqref{thm9-iv-e} holds. Thus, \eqref{thm9-iv} and \eqref{thm9-iv-e}
are equivalent when $\Psi=\emptyset$.

To see that \eqref{thm9-v} and \eqref{thm9-v-e} are equivalent when
$\Psi=\emptyset$, note that since $\Cm{\gc\K}\models_\h\emptyset$ is
vacuously true, \eqref{thm9-v} asserts that $\h(\id)\subseteq\h(\A)$
for every homomorphism $\h\colon\gc\P\to\Cm{\gc\K}$, \ie, $\id\leq\A$
is true in $\Cm{\gc\K}$, \ie, $\A$ is valid in $\gc\K$, whenever
$\gc\K$ is a frame satisfying \eqref{left reflection}, \eqref{right
  reflection}, \eqref{identity}, and \eqref{semi-Pasch}.  But that is
exactly what \eqref{thm9-v-e} says.  \endproof For the last two parts
of the following theorem, recall that $\Xi^d$, $\Xi^c$, and $\Xi^s$
are the equations of density \eqref{dens-eqs}, commutativity
\eqref{comm-eqs}, and symmetry \eqref{symm-eqs}, respectively.
\begin{theorem}\label{thm10-} Characterizations of $\CT^\Psi_4$, $\CT_4$,
  $\TR$, and $\KR$.
  \begin{enumerate}
  \item\label{thm10-1} If $\Psi\subseteq\Sigma^\times$, then the
    following statements are equivalent for every $\A\in\Pi^+$.
    \begin{enumerate}
    \item\label{thm10-i} $\A\in\CT^\Psi_4$.
    \item\label{thm10-ii} $\Psi\proves^+_4\id\leq\A$.
    \item\label{thm10-iii} $\Psi\proves^\times\id\leq\A$.
    \item\label{thm10-iv} If $\gc\A\in\RA$, $\h\colon\gc\P\to\gc\A$ is
      a homomorphism, and $\gc\A\models_\h\Psi$ then
      $\gc\A\models_\h\id\leq\A$.
    \item\label{thm10-v} If $\gc\K$ is a frame satisfying \eqref{left
        reflection}, \eqref{right reflection}, \eqref{identity}, and
      \eqref{Pasch}, $\h\colon\gc\P\to\gc\A$ is a homomorphism, and
      $\Cm{\gc\K} \models_\h \Psi$, then $\h(\id)\subseteq\h(\A)$.
    \item\label{thm10-vii} The sequent $\sep\fm\A\0\0$ is 4-provable
      from the sequents corresponding to $\Psi$.
    \item\label{thm10-vi} $\{\GG(\varepsilon):\varepsilon\in\Psi\}
      \proves_4\all\0\GG(\0\A\0)$.
    \item\label{thm10-viii} $\Psi\proves^+_3\id\leq\A$.
    \item\label{thm10-ix} $\{\GG(\varepsilon):\varepsilon\in\Psi\}
      \proves_3\all\0\GG(\0\A\0)$.
    \end{enumerate}
  \item\label{thm10-2} The following statements are equivalent for
    every $\A\in\Pi^+$.
    \begin{enumerate}
    \item\label{thm10-i-e} $\A\in\CT_4$.
    \item\label{thm10-ii-e} $\proves^+_4\id\leq\A$.
    \item\label{thm10-iii-e} $\proves^\times\id\leq\A$.
    \item\label{thm10-iv-e} $\id\leq\A$ is true in every relation
      algebra.
    \item\label{thm10-v-e} $\id\leq\A$ is valid in every frame
      satisfying \eqref{left reflection}, \eqref{right reflection},
      \eqref{identity}, and \eqref{Pasch}.
    \item\label{thm10-vii-e} The sequent $\sep\fm\A\0\0$ is
      4-provable.
    \item\label{thm10-vi-e} $\proves_4\all\0\GG(\0\A\0)$.
    \item\label{thm10-viii-e} $\proves^+_3\id\leq\A$.
    \item\label{thm10-ix-e} $\proves_3\all\0\GG(\0\A\0)$.
    \end{enumerate}
  \item\label{charTR} The following statements are equivalent for
    every $\A\in\Pi^+$.
    \begin{enumerate}
    \item\label{charTRi} $\A\in\TR$.
    \item\label{charTRii} $\Xi^d\cup\Xi^c\proves^+_4\id\leq\A$.
    \item\label{charTRiii} $\Xi^d\cup\Xi^c\proves^\times\id\leq\A$.
    \item\label{charTRiv} $\id\leq\A$ is true in every dense
      commutative relation algebra.
    \item\label{charTRv} $\A$ is valid in every frame satisfying
      \eqref{left reflection}, \eqref{right reflection},
      \eqref{identity}, \eqref{Pasch}, \eqref{dense}, and
      \eqref{comm}.
    \item\label{charTRvii} The sequent $\sep\0\A\0$ is 4-provable from
      the sequents corresponding to $\Xi^d\cup\Xi^c$.
    \item\label{charTRvi} $\Xi^d\cup\Xi^c\proves_4\all\0\GG(\0\A\0)$.
    \item\label{charTRviii} $\Xi^d\cup\Xi^c\proves^+_3\id\leq\A$.
    \item\label{charTRix} $\Xi^d\cup\Xi^c\proves_3\all\0\GG(\0\A\0)$.
    \end{enumerate}
  \item\label{charKR} The following statements are equivalent for
    every $\A\in\Pi^+$.
    \begin{enumerate}
    \item\label{charKRi} $\A\in\KR$.
    \item\label{charKRii} $\Xi^d\cup\Xi^s\proves^+_4\id\leq\A$.
    \item\label{charKRiii} $\Xi^d\cup\Xi^s\proves^\times\id\leq\A$.
    \item\label{charKRiv} $\id\leq\A$ is true in every dense symmetric
      relation algebra.
    \item\label{charKRv} $\A$ is valid in every frame satisfying
      \eqref{left reflection}, \eqref{right reflection},
      \eqref{identity}, \eqref{Pasch}, \eqref{dense}, and
      \eqref{symm}.
    \item\label{charKRvii} The sequent $\sep\0\A\0$ is 4-provable from
      the sequents corresponding to $\Xi^d\cup\Xi^s$.
    \item\label{charKRvi} $\Xi^d\cup\Xi^s\proves_4\all\0\GG(\0\A\0)$.
    \item\label{charKRviii} $\Xi^d\cup\Xi^s\proves^+_3\id\leq\A$.
    \item\label{charKRix} $\Xi^d\cup\Xi^s\proves_3\all\0\GG(\0\A\0)$.
    \end{enumerate}
  \end{enumerate}
\end{theorem}
\proof The proof of part \eqref{thm10-1} differs slightly from the
proof of Theorem \ref{thm9-}\eqref{thm9-1}.  Parts \eqref{thm10-i} and
\eqref{thm10-ii} are equivalent by definition, \eqref{thm10-ii} and
\eqref{thm10-iii} are equivalent by Theorem
\ref{thm5+}\eqref{thm5+iii}, and \eqref{thm10-iii} is equivalent to
\eqref{thm10-iv} by Theorem \ref{eq}\eqref{eq1}.  The proof that parts
\eqref{thm10-iv} and \eqref{thm10-v} are equivalent is the same as the
proof that parts \eqref{thm9-iv} and \eqref{thm9-v} of Theorem
\ref{thm9-} are equivalent, except that one uses Theorems
\ref{thm7}\eqref{thm7ii}, \ref{thm8}\eqref{thm8iv}, and
\ref{thm8}\eqref{thm8v} in place of Theorems \ref{thm7}\eqref{thm7i},
\ref{thm8}\eqref{thm8ii} and \ref{thm8}\eqref{thm8iii}, respectively.
Parts \eqref{thm10-iv} and \eqref{thm10-vii} are equivalent by Theorem
\ref{key}\eqref{key2}\eqref{key4} and Lemma \ref{equiv}.  Part
\eqref{thm10-viii} is equivalent to part \eqref{thm10-ii} by Theorem
\ref{thm5+}\eqref{thm5+iv}, equivalent to part \eqref{thm10-iii} by
Theorem \ref{thm3}\eqref{thm3v}, and equivalent to part
\eqref{thm10-ix} by Theorem \ref{thm2}\eqref{thm2iv} together with the
observation that $\id\leq\A \equiv^+_3 \all\0\GG(\0\A\0)$ by Theorem
\ref{thm2}\eqref{thm2iii}.  Finally, \eqref{thm10-vi} and
\eqref{thm10-ix} are equivalent by Theorem \ref{thm5+}\eqref{thm5+ii}.

The proof of part \eqref{thm10-2} is the same as the proof of Theorem
\ref{thm9-}\eqref{thm9-2} except that $\RA$ and \eqref{Pasch} replace
$\SA$ and \eqref{semi-Pasch}.  Parts \eqref{charTR} and \eqref{charKR}
follow mostly from part \eqref{thm10-1} by taking
$\Psi=\Xi^d\cup\Xi^c$ and $\Psi=\Xi^d\cup\Xi^s$, respectively.  In all
three cases one uses various instances of Lemma \ref{lem5}, taking
$\Psi$ to be $\emptyset$, $\Xi^d\cup\Xi^c$, or $\Xi^d\cup\Xi^s$, and
arguing as in the proof of Theorem \ref{thm9-}\eqref{thm9-2}.

There are some differences between the proofs of parts \eqref{thm10-1}
and \eqref{thm10-2} and the proofs of parts \eqref{charTR} and
\eqref{charKR}.  While the definition of $\TR$ still produces the
equivalence of \eqref{charTRi} and \eqref{charTRii}, the definition of
$\KR$ yields the equivalence of \eqref{charKRi} and \eqref{charKRv}.
In the proof that \eqref{charTRiv} and \eqref{charTRv} are equivalent,
one needs to observe that $\Cm{\gc\K}$ is a dense commutative relation
algebra iff $\gc\K$ satisfies \eqref{left reflection}, \eqref{right
  reflection}, \eqref{identity}, \eqref{Pasch}, \eqref{dense}, and
\eqref{comm}.  In the proof that \eqref{charKRiv} and \eqref{charKRv}
are equivalent, one must observe that $\Cm{\gc\K}$ is a dense symmetic
relation algebra iff $\gc\K$ satisfies \eqref{left reflection},
\eqref{right reflection}, \eqref{identity}, \eqref{Pasch},
\eqref{dense}, and \eqref{symm}.  These observations follow from Lemma
\ref{SA:symm->comm} and Theorems \ref{thm7-},
\ref{thm7}\eqref{thm7ii}, and \ref{thm8a}\eqref{thm8aii}.  \endproof
The inclusion of commutativity in the definitions of $\TR$ and $\KR$
allows some further characterizations as a consequence of Lemma
\ref{commSA}.
\begin{theorem}\label{KT} Integral characterizations of $\TR$ and $\KR$.
  \begin{enumerate}
  \item\label{KT1} The following statements are equivalent for every
    $\A\in\Pi^+$.
    \begin{enumerate}
    \item\label{KT1i} $\A\in\TR$.
    \item\label{KT1ii} $\id\leq\A$ is true in every integral
      dense commutative relation algebra.
    \item\label{KT1iii} $\A$ is valid in every frame
      $\gc\K=\<\K,\R,\star{},\II\>$ satisfying $|\II|=1$, \eqref{left
      reflection}, \eqref{right reflection}, \eqref{identity},
      \eqref{Pasch}, \eqref{dense}, and \eqref{comm}.
    \end{enumerate}
  \item\label{KT2} The following statements are equivalent for
    every $\A\in\Pi^+$.
    \begin{enumerate}
    \item\label{KT2i} $\A\in\KR$.
    \item\label{KT2ii} $\id\leq\A$ is true in every integral
      dense symmetric relation algebra.
    \item\label{KT2iii} $\A$ is valid in every frame
      $\gc\K=\<\K,\R,\star{},\II\>$ satisfying $|\II|=1$, \eqref{left
      reflection}, \eqref{right reflection}, \eqref{identity},
      \eqref{Pasch}, \eqref{dense}, and \eqref{symm}.
    \end{enumerate}
  \end{enumerate}
\end{theorem}
Theorems \ref{thm9-}, \ref{thm10-}, and \ref{KT} link the logics of
3-provability with semi-associative relation algebras and the logics
of 4-provability with relation algebras. The logics of
$\omega$-provability are linked with representable relation algebras.
\begin{theorem}\label{thm10} Characterizations of $\CT^\Psi_\omega$
  and $\CT_\omega$.
  \begin{enumerate}
  \item\label{thm10-a} If $\Psi\subseteq\Sigma^\times$, then the
    following statements are equivalent for every $\A\in\Pi^+$.
    \begin{enumerate}
    \item\label{thm10i} $\A\in\CT^\Psi_\omega$.
    \item\label{thm10ii} $\Psi\proves^+\id\leq\A$.
    \item\label{thm10iv} If $\h\colon\gc\P\to\gc\A$ is a homomorphism
      and $\gc\A\in\RRA$ then $\gc\A\models_\h\Psi$ implies
      $\gc\A\models_\h\id\leq\A$.
    \item\label{thm10viii} $\sep\fm\A\0\0$ is $\omega$-provable from
      the sequents corresponding to $\Psi$.
    \item\label{thm10vi} $\{\GG(\varepsilon):\varepsilon\in\Psi\}
      \proves\all\0\GG(\0\A\0)$.
    \end{enumerate}
  \item\label{thm10-b} The following statements are equivalent for
    every $\A\in\Pi^+$.
    \begin{enumerate}
    \item\label{thm10i-e} $\A\in\CT_\omega$.
    \item\label{thm10ii-e} $\proves^+\id\leq\A$.
    \item\label{thm10iv-e} $\id\leq\A$ is true in every
      representable relation algebra.
    \item\label{thm10viii-e} $\sep\fm\A\0\0$ is $\omega$-provable.
    \item\label{thm10vi-e} $\proves\all\0\GG(\0\A\0)$.
    \end{enumerate}
  \end{enumerate}
\end{theorem}
\proof Parts \eqref{thm10i} and \eqref{thm10ii} are equivalent by
definition, \eqref{thm10ii} and \eqref{thm10iv} by Theorem \ref{sem},
and \eqref{thm10ii} and \eqref{thm10vi} by Theorem
\ref{thm1}\eqref{thm1iii}\eqref{thm1iv} and $\GG(\id\leq\A) \equiv^+
\all\0\GG(\0\A\0)$, and \eqref{thm10iv} and \eqref{thm10viii} by
Theorem \ref{key}\eqref{key3}\eqref{key4} and Lemma \ref{equiv}.
\endproof Two observations close this section. For $5\leq\n<\omega$,
membership in $\CT_\n$ as defined in \eqref{CTn} cannot be
characterized as $\n$-provability in the sequent calculus.  The
situation is quite complicated and full of non-finite axiomatizibility
results; see \cite{MR1935083}. However, the implication in one
direction still holds.
\begin{theorem}\label{two-obs-1} If $\Psi\subseteq\Sigma^\times$,
    $\A\in\Pi^+$, $5\leq\n<\omega$, and $\sep\fm\A\0\0$ is
  $\n$-provable in the sequent calculus from $\Psi$, then
  $\A\in\CT^\Psi_\n$.
\end{theorem}
Several results in this section also apply to $\TT^\Psi_\n$.
\begin{theorem}\label{two-obs-2} Theorems \ref{thm9-}, \ref{thm10-},
    \ref{thm10}, and \ref{two-obs-1} continue to hold if $\A$ is
    assumed to be in $\Pi^r$ instead of $\Pi^+$, $\Psi$ is a set of
    equations between predicates in $\Pi^r$, and $\CT^\Psi_3$,
    $\CT^\Psi_4$, $\CT^\Psi_\n$, and $\CT^\Psi_\omega$ are replaced by
    $\TT^\Psi_3$, $\TT^\Psi_4$, $\TT^\Psi_\n$, and $\TT^\Psi_\omega$.
\end{theorem}
\section{Theorems and derived rules of $\CT_3$, $\CT_4$, and $\TR$}
\label{sect14}
Theorem \ref{thm9-} provides several ways to show $\A\in\CT_3$.  One
can prove $\id\leq\A$ in $\Lwx$ or $\Ls_3^+$, or prove $\id\leq\A$ is
true in every semi-associative relation algebra, or prove
$\all\0\GG(\0\A\0)$ in $\Ls_3$, or prove $\A$ is valid in every frame
satisfying \eqref{left reflection}, \eqref{right reflection},
\eqref{identity}, and \eqref{Pasch}, or show that $\sep\fm\A\0\0$ is
3-provable.

Theorem \ref{thm11} below shows that $\CT_3$ includes truth
\eqref{ttt}, laws of the excluded middle for Boolean negation
\eqref{lem} and De Morgan negation \eqref{lem1}, self-implication
\eqref{self}, basic laws of disjunction and conjunction
\eqref{or1}--\eqref{and-assoc} and \eqref{to1}--\eqref{to2}, laws of
distributivity \eqref{dist-or-and}--\eqref{dist-and-or}, laws of
double Boolean negation \eqref{doubleBA1}--\eqref{doubleBA2} and
double De Morgan negation \eqref{doubleDM1}--\eqref{doubleDM2}, De
Morgan laws for De Morgan negation \eqref{DeM1}--\eqref{DeM4} and for
Boolean negation \eqref{DeM5}--\eqref{DeM8}, explosion for Boolean
negation \eqref{explosion}, and some basic laws of fusion
\eqref{fus1}--\eqref{fus4}. Furthermore, $\CT_3$ has many derived
rules of inference, including adjunction \eqref{adj}, modus ponens
\eqref{modus-p}, disjunctive syllogism \eqref{disj-syll},
contraposition for De Morgan negation
\eqref{contra-rule1}--\eqref{contra-rule2}, a cut rule
\eqref{cut-rule}, the suffixing rule \eqref{suff-rule}, the prefixing
rule \eqref{pref-rule}, and monotonicity for fusion \eqref{fus-rule}.
The prefixing axiom \eqref{no.1} is not in $\CT_3$ but it is in
$\CT_4$, as shown by Theorem \ref{thm12} in \SS\ref{sect15}.  The
suffixing axiom \eqref{suff} is not even $\omega$-provable and the
same applies to the axioms of permutation \eqref{perm} and
contraposition \eqref{contra}, as shown by Theorem \ref{thm13} in
\SS\ref{sect16}.

Although explosion for Boolean negation \eqref{explosion} is in
$\CT_3$, explosion for De Morgan negation $\A\land\rmin\A\to\B$ and
positive paradox $\A\to(\B\to\A)$ are not even $\omega$-provable,
because as relations they need not contain the identity relation. For
example, if $\<\y,\x\>\in\A$, $\<\x,\y\>\notin\A$, and
$\<\y,\x\>\notin\B$, then $\<\x,\x\>\notin(\A\land\rmin\A)\to\B$, and
if $\<\y,\x\>\in\A$, $\<\z,\y\>\in\B$, and $\<\z,\x\>\notin\A$, then
$\<\x,\x\>\notin\A\to(\B\to\A)$.

The terminology of $\n$-proofs and $\n$-provability extends from
sequents to predicates.  A predicate $\A\in\Pi^+$ is {\bf
  $\n$-provable from $\Psi$} if the sequent $\sep\fm\A\0\0$ is
$\n$-provable from $\Psi$ and an {\bf $\n$-proof of $\A$ from $\Psi$}
is an $\n$-proof of the sequent $\sep\fm\A\0\0$ from $\Psi$.  Thus,
$\A$ is $\n$-provable from $\Psi$ iff there is an $\n$-proof of $\A$
from $\Psi$. By \eqref{equiv-1}, \eqref{equiv-2}, and \eqref{equiv-4},
$\A\to\B$ is $\n$-provable from $\Psi$ iff there is an $\n$-proof from
$\Psi$ of the sequent $\fm\A{\var\i}{\var\j}\sep\fm\B{\var\i}{\var\j}$
for distinct $\i,\j<\n$.  For this reason, we refer to an $\n$-proof
from $\Psi$ of the sequent $\fm\A{\var\i}{\var\j} \sep
\fm\B{\var\i}{\var\j}$ as an $\n$-proof of $\A\to\B$ from $\Psi$.
Accordingly, $\A$ (or $\A\to\B$) is {\bf $\n$-provable from density},
{\bf $\n$-provable from commutativity}, or {\bf $\n$-provable from
  symmetry} if the sequent $\sep\fm\A\vi\vi$ (or the sequent
$\fm\A\vi\vj \sep \fm\B\vi\vj$) is $\n$-provable from the equations of
the equations of density $\Xi^d$ \eqref{dens-eqs}, commutativity
$\Xi^c$ \eqref{comm-eqs}, or the equations of symmetry $\Xi^s$
\eqref{symm-eqs}, respectively, for any distinct $\i,\j<\n$.

Derived rules of inference are stated in the form $\A\proves\B$ or
$\A,\B\proves\C$.  We say $\A\proves\B$ is {\bf $\n$-provable} if
every proof of $\sep\fm\A\0\0$ can be extended by an $\n$-proof to a
proof of $\sep\fm\B\0\0$.  Similarly, we say $\A,\B\proves\C$ is {\bf
  $\n$-provable} if every proof of $\sep\fm\A\0\0$ and every proof of
$\sep\fm\B\0\0$ can be concatenated and extended by an $\n$-proof to
obtain a proof of $\sep\fm\C\0\0$.  For example, if $\sep\fm\A\0\0$ is
$\m$-provable and $\A\proves\B$ is $\n$-provable, then $\B$ is
$\max(\n,\m)$-provable.
\begin{theorem}\label{thm11}
  Assume $\A,\B,\C,\D,\E,\F\in\Pi^+$.
  \begin{enumerate}
  \item Predicates \eqref{ttt}--\eqref{lem1} are 1-provable.
  \item Predicates \eqref{self}--\eqref{explosion} are 2-provable.
  \item Predicates \eqref{to1}--\eqref{min-th} are 3-provable.
  \item Rules \eqref{adj}--\eqref{star-rule} are 1-provable.
  \item Rules \eqref{tran-rule}--\eqref{red-rule2} are 2-provable.
  \item Rules \eqref{suff-rule}--\eqref{cycling} are 3-provable.
  \item Predicates \eqref{no.1}--\eqref{no.6} are 4-provable.
  \item Predicates \eqref{reductio}--\eqref{contract4} are 3-provable
    from density.
  \item Predicates \eqref{contract2}--\eqref{contract3} are 4-provable
    from density.
  \item Predicates \eqref{mp}--\eqref{contra} are 3-provable from
    commutativity.
  \item Predicates \eqref{perm}--\eqref{suff} are 4-provable from
    commutativity.
  \item Predicate \eqref{self-dist} is 4-provable from density
    and commutativity.
  \item Predicate \eqref{symm-ax} is 2-provable from symmetry.
  \item Predicate \eqref{comm-ax} is 3-provable from symmetry.
  \item The predicates and rules of $\CT_3$, $\CT_4$, and $\TR$
    include\display{\begin{tabular}{|l|l|l|}\hline 
        \rm\,Logic&\rm\,Predicates&\rm\,Rules\\\hline 
        \,$\CT_3$ 
        &\,\eqref{ttt}--\eqref{min-th} 
        &\,\eqref{adj}--\eqref{cycling} 
        \\
        \,$\CT_4$ 
        &\,\eqref{ttt}--\eqref{min-th}\eqref{no.1}--\eqref{no.6} 
        &\,\eqref{adj}--\eqref{cycling} 
        \\
        \,$\TR$ 
        &\,\eqref{ttt}--\eqref{min-th}\eqref{no.1}--\eqref{no.6}
        \eqref{mp}--\eqref{self-dist} 
        &\,\eqref{adj}--\eqref{cycling} 
        \\
        \hline
      \end{tabular}} 
  \end{enumerate}
  The predicates and rules that have been provided with $\n$-proofs
  are marked with asterisks. 
  \begin{align}
    \intertext{1-provable predicates:}
    \ttt& \label{ttt}&&*\\
    \A&\lor\neg\A\label{lem}&&*\\
    \A&\lor\rmin\A\label{lem1}&&*
    \intertext{2-provable predicates:}
    \A&\to\A\label{self}&&*\\
    ((\A\to\A)\to\B)&\to\B\label{E-ax}&&*\\
    \A\lor\A&\to\A\label{or1}&&\\
    \A&\to\A\lor\B\label{or2}&&\\
    \B&\to\A\lor\B\label{or3}&&\\
    \A\lor\B&\to\B\lor\A\label{or-comm}&&\\
    \A&\to\A\land\A\label{and1}&&\\
    \A\land\B&\to\A\label{and2}&&\\
    \A\land\B&\to\B\label{and3}&&\\
    \A\land\B&\to\B\land\A\label{and-comm}&&\\
    (\A\lor\B)\lor\C&\to\A\lor(\B\lor\C)\label{or-assoc}&&\\
    (\A\land\B)\land\C&\to\A\land(\B\land\C)\label{and-assoc}&&\\
    (\A\lor\B)\land\C&\to(\A\land\C)\lor(\B\land\C)\label{dist-or-and}&&\\
    (\A\land\B)\lor\C&\to(\A\lor\C)\land(\B\lor\C)\label{dist-and-or}&&\\
    \rmin\rmin\A&\to\A\label{doubleDM1}&&\\
    \A&\to\rmin\rmin\A\label{doubleDM2}&&\\
    \rmin(\A\lor\B)&\to\rmin\A\land\rmin\B\label{DeM1}&&\\
    \rmin(\A\land\B)&\to\rmin\A\lor\rmin\B\label{DeM2}&&\\
    \rmin\A\land\rmin\B&\to\rmin(\A\lor\B)\label{DeM3}&&\\
    \rmin\A\lor\rmin\B&\to\rmin(\A\land\B)\label{DeM4}&&\\
    \neg\neg\A&\to\A\label{doubleBA1}&&\\
    \A&\to\neg\neg\A\label{doubleBA2}&&\\
    \neg(\A\lor\B)&\to\neg\A\land\neg\B\label{DeM5}&&\\
    \neg(\A\land\B)&\to\neg\A\lor\neg\B\label{DeM6}&&\\
    \neg\A\land\neg\B&\to\neg(\A\lor\B)\label{DeM7}&&\\
    \neg\A\lor\neg\B&\to\neg(\A\land\B)\label{DeM8}&&\\
    \A&\to\B\lor\neg\B\label{top}&&*\\
    \neg\A\land\A&\to\B\label{explosion}&&*
    \intertext{3-provable predicates:}
    (\A\to\B)\land(\A\to\C)&\to(\A\to\B\land\C)\label{to1}&&\\
    (\A\to\C)\land(\B\to\C)&\to(\A\lor\B\to\C)\label{to2}&&\\
    (\A\to\B)\land(\C\to\D)&\to(\A\land\C\to\B\land\D)\label{to3}&&\\
    (\A\to\B)\land(\C\to\D)&\to(\A\lor\C\to\B\lor\D)\label{to4}&&*\\
    (\A\to\B)\lor(\C\to\D)&\to(\A\land\C\to\B\lor\D)\label{to5}&&\\
    \A\circ\B&\to\rmin(\A\to\rmin\B)\label{fus1}&&\\
    \rmin(\A\to\rmin\B)&\to\A\circ\B\label{fus2}&&\\
    (\A\to\B)\circ\A&\to\B\label{fus3}&&\\
    \A&\to(\B\to\A\circ\B)\label{fus4}&&*\\
    \A&\to((\B\to\rmin\A)\to\rmin\B)\label{mp-what}&&*\\
    \A\circ\B\land\C
    &\to\A\circ(\B\land\rmin\D)\lor(\A\land\C\circ\D)\circ\B
    \label{reflection1}&&*\\
    \A\circ\B\land\C
    &\to(\A\land\rmin\D)\circ\B\lor\A\circ(\B\land\D\circ\C)
    \label{reflection1a}&&\\
    \A\circ\B\land\C
    &\to(\A\land\C\circ\star\B)\circ(\B\land\star\A\circ\C)
    \label{dedekind}&&\\
    \ttt\circ\A&\to\A
    \label{right-id}&&*\\
    \A\circ\ttt&\to\A
    \label{left-id}&&*\\
    \A&\to\A\circ\star\ttt
    \label{A->Aot*}&&*\\
    \ttt&\to\star\ttt
    \label{t->t*}&&*\\
    \star\ttt&\to\ttt
    \label{t*->t}&&*\\
    \ttt\land\rmin\ttt&\to\A
    \label{min-th}&&*
    \intertext{1-provable rules:}
    \A,\B&\proves\A\land\B\label{adj}&&*\\
    \A\to\B,\A&\proves\B\label{modus-p}&&*\\
    \A\lor\B,\rmin\A&\proves\B\label{disj-syll}&&*\\
    \A&\proves\A^*\label{star-rule}&&*
    \intertext{2-provable rules:}
    \A\to\B,\B\to\C&\proves\A\to\C\label{tran-rule}&&*\\
    \A\to\B&\proves\rmin\B\to\rmin\A\label{contra-rule1}&&*\\
    \A\to\rmin\B&\proves\B\to\rmin\A\label{contra-rule2}&&*\\
    \A\land\B\to\C,\B\to\C\lor\A&\proves\B\to\C\label{cut-rule}&&*\\
    \A&\proves(\A\to\B)\to\B\label{E-rule}&&*\\
    \A\land\B\to\neg\C&\proves\A\land\C\to\neg\B\label{antilogism}&&*\\
    \B\land\neg\C\to\A\land\neg\A&\proves\B\to\C\label{red-rule1}&&\\
    \B\land\C^*\to\A\land\neg\A&\proves\B\to\rmin\C\label{red-rule2}&&
    \intertext{3-provable rules:}
    \A\to\B&\proves(\B\to\C)\to(\A\to\C)\label{suff-rule}&&*\\
    \A\to\B&\proves(\C\to\A)\to(\C\to\B)\label{pref-rule}&&*\\
    \A\to\B,\C\to\D&\proves(\B\to\C)\to(\A\to\D)\label{affixing}&&*\\
    \A\to\B,\C\to\D&\proves\A\circ\C\to\B\circ\D\label{fus-rule}&&*\\
    \A\to(\B\to\C)&\proves\B\to(\rmin\C\to\rmin\A)\label{cycling}&&*
    \intertext{4-provable predicates:}
    (\A\to\B)&\to((\C\to\A)\to(\C\to\B))    \label{no.1}&&*\\
    (\A\to(\B\to\C))&\to(\A\circ\B\to\C)    \label{no.2}&&*\\
    (\A\circ\B\to\C)&\to(\A\to(\B\to\C))    \label{no.3}&&*\\
    (\A\to\B)&\to(\A\circ\C\to\B\circ\C)    \label{no.4}&&*\\
    (\A\circ\B)\circ\C&\to\A\circ(\B\circ\C)\label{no.5}&&*\\
    \A\circ(\B\circ\C)&\to(\A\circ\B)\circ\C\label{no.6}&&*
    \intertext{3-provable from density:}
    (\A\to\rmin\A)&\to\rmin\A\label{reductio}&&*\\
    \A\land\B&\to\A\circ\B\label{contract5}&&*\\
    (\A\to\B)&\to\rmin\A\lor\B\label{contract4}&&*\\
    \intertext{4-provable from density:}
    (\A\to(\A\to\B))&\to(\A\to\B)\label{contract2}&&*\\
    (\A\to(\B\to\C))&\to(\A\land\B\to\C)\label{contract3}&&*\\
    \intertext{3-provable from commutativity:}
    \A&\to((\A\to\B)\to\B)\label{mp}&&*\\ 
    (\A\to\rmin\B)&\to(\B\to\rmin\A)\label{contra}&&*
    \intertext{4-provable from commutativity:}
    (\A\to(\B\to\C))&\to(\B\to(\A\to\C))\label{perm}&&*\\
    (\A\to\B)&\to((\B\to\C)\to(\A\to\C))\label{suff}&&*
    \intertext{4-provable from density and commutativity:}
    (\A\to(\B\to\C))&\to((\A\to\B)\to(\A\to\C))\label{self-dist}&&*
    \intertext{2-provable from symmetry:}
    \A\land\rmin\A&\to\B\label{symm-ax}&&*
    \intertext{3-provable from symmetry:}
    \A\circ\B&\to\B\circ\A\label{comm-ax}&&*
  \end{align}
\end{theorem}
\proof[\ref{ttt}]
1-proof of $\ttt$.
\begin{align*}
  1.&&&\sep\fm\id\0\0
  &&\text{$|\id$}\\
  2.&&&\sep\fm\ttt\0\0 &&\text{1, \eqref{t=id}}
\end{align*}
\proof[\ref{lem}] 1-proof of $\A\lor\neg\A$.
\begin{align*}
  1.&&\fm\A\0\0&\sep\fm\A\0\0
  &&\text{axiom}\\
  2.&&&\sep\fm\A\0\0,\fm{\min\A}\0\0
  &&\text{1, $|\min\blank$}\\
  3.&&&\sep\fm{\A+\min\A}\0\0
  &&\text{2, $|+$}\\
  4.&&&\sep\fm{\A\lor\neg\A}\0\0 
  &&\text{3, \eqref{or}, \eqref{neg}}
\end{align*}
\proof[\ref{lem1}] 1-proof of $\A\lor\rmin\A$.
\begin{align*}
  1.&&\fm\A\0\0&\sep\fm\A\0\0
  &&\text{axiom}\\
  2.&&\fm{\conv\A}\0\0&\sep\fm\A\0\0
  &&\text{1, $\conv{}|$}\\
  3.&&&\sep\fm\A\0\0,\fm{\min{\conv\A}}\0\0
  &&\text{2, $|\min\blank$}\\
  4.&&&\sep\fm{\A+\min{\conv\A}}\0\0
  &&\text{3, $|+$}\\
  5.&&&\sep\fm{\A\lor\rmin\A}\0\0 
  &&\text{4, \eqref{or}, \eqref{rmin}}
\end{align*}
\proof[\ref{self}, \ref{E-ax}] 2-proof of $\A\to\A$ and
$((\A\to\A)\to\B)\to\B$.
\begin{align*}
  1.&&\fm\A\0\1&\sep\fm\A\0\1
  &&\text{axiom}\\
  2.&&&\sep\fm{\A\to\A}\0\0
  &&\text{1, \eqref{equiv-1}, \eqref{equiv-4}}\\
  3.&&&\sep\fm{\conv{\A\to\A}}\0\0
  &&\text{2, $|\conv{}$}\\
  4.&&\fm\B\0\1&\sep\fm\B\0\1
  &&\text{axiom}\\
  5.&&&\sep\fm{\min\B}\0\1,\fm\B\0\1
  &&\text{4, $|\min\blank$}\\
  6.&&&\sep\fm{\conv{\A\to\A}\rp\min\B}\0\1,\fm\B\0\1
  &&\text{3, 5, $|\rp$}\\
  7.&&\fm{\min{\conv{\A\to\A}\rp\min\B}}\0\1&\sep\fm\B\0\1
  &&\text{6, $\min\blank|$}\\
  8.&&\fm{(\A\to\A)\to\B}\0\1&\sep\fm\B\0\1 
  &&\text{7, \eqref{to}}
\end{align*}
\proof[\ref{top}] 2-proof of $\A\to\B\lor\neg\B$.
\begin{align*}
  1.&&\fm\A\0\1,\fm\B\0\1&\sep\fm\B\0\1
  &&\text{axiom}\\
  2.&&\fm\A\0\1&\sep\fm\B\0\1,\fm{\min\B}\0\1
  &&\text{1, $|\min\blank$}\\
  3.&&\fm\A\0\1&\sep\fm{\B+\min\B}\0\1
  &&\text{2, $|+$}\\
  4.&&\fm\A\0\1&\sep\fm{\B\lor\neg\B}\0\1
  &&\text{3, \eqref{or}, \eqref{neg}}
\end{align*}
\proof[\ref{explosion}] 2-proof of $\A\land\neg\A\to\B$.
\begin{align*}
  1.&&\fm\A\0\1&\sep\fm\A\0\1,\fm\B\0\1
  &&\text{axiom}\\
  2.&&\fm\A\0\1,\fm{\min\A}\0\1&\sep\fm\B\0\1
  &&\text{1, $\min\blank|$}\\
  3.&&\fm{\A\cdot\min\A}\0\1&\sep\fm\B\0\1
  &&\text{2, $\cdot|$}\\
  4.&&\fm{\A\land\neg\A}\0\1&\sep\fm\B\0\1
  &&\text{3, \eqref{and}, \eqref{neg}}
\end{align*}
\proof[\ref{to4}] 3-proof of
$(\A\to\B)\land(\C\to\D)\to(\A\lor\C\to\B\lor\D)$.
\begin{align*}
  1.&&\fm\A\2\0&\sep\fm\A\2\0
  &&\text{axiom}\\
  2.&&\fm\A\2\0&\sep\fm{\conv\A}\0\2
  &&\text{1, $|\conv{}$}\\
  3.&&\fm\B\2\1&\sep\fm\B\2\1
  &&\text{axiom}\\
  4.&&&\sep\fm\B\2\1,\fm{\min\B}\2\1
  &&\text{3, $|\min\blank$}\\
  5.&&\fm\A\2\0&\sep\fm\B\2\1,\fm{\conv\A\rp\min\B}\0\1
  &&\text{2, 4, $|\rp$}\\
  6.&&\fm\C\2\0&\sep\fm\C\2\0
  &&\text{axiom}\\
  7.&&\fm\C\2\0&\sep\fm{\conv\C}\0\2
  &&\text{6, $|\conv{}$}\\
  8.&&\fm\D\2\1&\sep\fm\D\2\1
  &&\text{axiom}\\
  9.&&&\sep\fm\D\2\1,\fm{\min\D}\2\1
  &&\text{8, $|\min\blank$}\\
  10.&&\fm\C\2\0&\sep\fm\D\2\1,\fm{\conv\C\rp\min\D}\0\1
  &&\text{7, 9, $|\rp$}\\
  11.&&\fm{\A+\C}\2\0\sep\fm\B\2\1,\fm\D\2\1,
  &\fm{\conv\A\rp\min\B}\0\1,\fm{\conv\C\rp\min\D}\0\1
  &&\text{5, 10, $+|$}\\
  12.&&\fm{\A+\C}\2\0\sep\fm{\B+\D}\2\1,
  &\fm{\conv\A\rp\min\B}\0\1,\fm{\conv\C\rp\min\D}\0\1
  &&\text{11, $|+$}\\
  13.&&\fm{\conv{\A+\C}}\0\2,\fm{\min{\B+\D}}\2\1
  &\sep\fm{\conv\A\rp\min\B}\0\1,\fm{\conv\C\rp\min\D}\0\1
  &&\text{12, $\conv{}|$, $\min\blank|$}\\
  14.&&\fm{\conv{\A+\C}\rp\min{\B+\D}}\0\1
  &\sep\fm{\conv\A\rp\min\B}\0\1,\fm{\conv\C\rp\min\D}\0\1
  &&\text{13, no $\2$}\\
  15.&&\fm{\min{\conv\A\rp\min\B}}\0\1,\fm{\min{\conv\C\rp\min\D}}\0\1,
  &\fm{\conv{\A+\C}\rp\min{\B+\D}}\0\1\sep
  &&\text{14, $\min\blank|$}\\
  16.&&\fm{\min{\conv\A\rp\min\B}\cdot\min{\conv\C\rp\min\D}}\0\1
  &\sep\fm{\min{\conv{\A+\C}\rp\min{\B+\D}}}\0\1
  &&\text{15, $\cdot|$, $|\min\blank$}\\
  17.&&\fm{(\A\to\B)\land(\C\to\D)}\0\1&\sep\fm{\A\lor\C\to\B\lor\D}\0\1
  &&\text{16, \eqref{or}, \eqref{and}, \eqref{to}}
\end{align*}
\proof[\ref{fus4}] 3-proof of $\A\to(\B\to\A\circ\B)$.
\begin{align*}
  1.&&\fm\A\0\1&\sep\fm\A\0\1
  &&\text{axiom}\\
  2.&&\fm\B\2\0&\sep\fm\B\2\0
  &&\text{axiom}\\
  3.&&\fm\A\0\1,\fm\B\2\0&\sep\fm{\B\rp\A}\2\1
  &&\text{1, 2, $|\rp$}\\
  4.&&\fm\A\0\1,\fm{\conv\B}\0\2,\fm{\min{\B\rp\A}}\2\1&\sep
  &&\text{3, $\conv{}|$, $\min\blank|$}\\
  5.&&\fm\A\0\1,\fm{\conv\B\rp\min{\B\rp\A}}\0\1&\sep
  &&\text{4, $\rp|$, no $\2$}\\
  6.&&\fm\A\0\1&\sep\fm{\min{\conv\B\rp\min{\B\rp\A}}}\0\1
  &&\text{5, $|\min\blank$}\\
  7.&&\fm\A\0\1&\sep\fm{\B\to\A\circ\B}\0\1
  &&\text{6, \eqref{to}, \eqref{circ}}
\end{align*}
\proof[\ref{mp-what}] 3-proof of $\A\to((\B\to\rmin\A)\to\rmin\B)$.
\begin{align*}
  1.&&\fm\A\0\1&\sep\fm\A\0\1
  &&\text{axiom}\\
  2.&&\fm\B\1\2&\sep\fm\B\1\2
  &&\text{axiom}\\
  3.&&\fm\A\0\1&\sep\fm{\min{\min{\conv\A}}}\1\0
  &&\text{1, $|\conv{}$, $\min\blank|$, $|\min\blank$}\\
  4.&&\fm\B\1\2&\sep\fm{\conv\B}\2\1
  &&\text{2, $|\conv{}$}\\
  5.&&\fm\A\0\1,
  \fm\B\1\2&\sep\fm{\conv\B\rp\min{\min{\conv\A}}}\2\0
  &&\text{3, 4, $|\rp$}\\
  6.&&\fm\A\0\1,\fm{\min{\conv\B\rp\min{\min{\conv\A}}}}\2\0,
  \fm\B\1\2&\sep
  &&\text{5, $\min\blank|$}\\
  7.&&\fm\A\0\1,\fm{\B\to\rmin\A}\2\0,\fm\B\1\2&\sep
  &&\text{6, \eqref{rmin}, \eqref{to}}\\
  8.&&\fm\A\0\1,\fm{\conv{\B\to\rmin\A}}\0\2,
  &\fm{\min{\min{\conv\B}}}\2\1\sep
  &&\text{7, $\conv{}|$, $|\min\blank$, $|\min\blank$}\\
  9.&&\fm\A\0\1,
  \fm{\conv{\B\to\rmin\A}\rp\min{\min{\conv\B}}}\0\1&\sep
  &&\text{8, $\rp|$, no $\2$}\\
  10.&&\fm\A\0\1
  &\sep\fm{\min{\conv{\B\to\rmin\A}\rp\min{\min{\conv\B}}}}\0\1
  &&\text{9, $|\min\blank$}\\
  11.&&\fm\A\0\1&\sep\fm{(\B\to\rmin\A)\to\rmin\B}\0\1
  &&\text{10, \eqref{rmin}, \eqref{to}}
\end{align*}
\proof[\ref{reflection1}] 3-proof of
$\A\circ\B\land\C\to\A\circ(\B\land\rmin\D)
\lor(\A\land\C\circ\D)\circ\B$.
\begin{align*}
  1.&&\fm\C\0\1&\sep\fm\C\0\1
  &&\text{axiom}\\
  2.&&\fm\D\2\0&\sep\fm\D\2\0
  &&\text{axiom}\\
  3.&&\fm\C\0\1,\fm\D\2\0&\sep\fm{\D\rp\C}\2\1
  &&\text{1, 2, $|\rp$}\\
  4.&&\fm\A\2\1&\sep\fm\A\2\1
  &&\text{axiom}\\
  5.&&\fm\A\2\1,\fm\C\0\1,\fm\D\2\0&\sep
  \fm{(\A\cdot\C\rp\D)}\2\1
  &&\text{3, 4, $|\cdot$}\\
  6.&&\fm\B\0\2&\sep\fm\B\0\2
  &&\text{axiom}\\
  7.&&\fm\B\0\2,\fm\A\2\1,\fm{\C&}\0\1,
  \fm\D\2\0\sep\fm{\B\rp(\A\cdot\C\rp\D)}\0\1
  &&\text{5, 6, $|\rp$}\\
  8.&&\fm\B\0\2,\fm\A\2\1,\fm\C\0\1&\sep
  \fm{\min{\conv\D}}\0\2,\fm{\B\rp(\A\cdot\C\rp\D)}\0\1
  &&\text{7, $\conv{}|$, $|\min\blank$}\\
  9.&&\fm\B\0\2,\fm\A\2\1,\fm\C\0\1&\sep
  \fm{\B\cdot\min{\conv\D}}\0\2,\fm{\B\rp(\A\cdot\C\rp\D)}\0\1
  &&\text{6, 8, $|\cdot$}\\
  10.&&\fm\B\0\2,\fm\A\2\1,\fm\C\0\1&\sep
  \fm{(\B\cdot\min{\conv\D})\rp\A}\0\1,\fm{\B\rp
    (\A\cdot\D\rp\C)}\0\1
  &&\text{4, 9, $|\rp$}\\
  11.&&\fm\B\0\2,\fm\A\2\1,\fm\C\0\1&\sep
  \fm{(\B\cdot\min{\conv\D})\rp\A+\B\rp(\A\cdot\D\rp\C)}\0\1
  &&\text{10, $|+$}\\
  12.&&\fm{\B\rp\A}\0\1,\fm\C\0\1&\sep
  \fm{(\B\cdot\min{\conv\D})\rp\A+\B\rp(\A\cdot\D\rp\C)}\0\1
  &&\text{11, $\rp|$, no $\2$}\\
  13.&&\fm{\B\rp\A\cdot\C}\0\1&\sep
  \fm{(\B\cdot\min{\conv\D})\rp\A+\B\rp(\A\cdot\D\rp\C)}\0\1 
  &&\text{12, $\cdot|$}\\
  14.&&\fm{\A\circ\B\land\C}\0\1\sep
  \fm{&\A\circ(\B\land\rmin\D)\lor(\A\cdot\C\circ\D)\circ\B}\0\1 
  &&\text{13, \eqref{and}, \eqref{rmin}, \eqref{circ}}
\end{align*}
\proof[\ref{right-id}] 3-proof of $\ttt\circ\A\to\A$.
\begin{align*}
  1.&& \fm\A\0\1&\sep\fm\A\0\1
  &&\text{axiom}\\
  2.&& \fm\A\0\2,\fm\id\2\1&\sep\fm\A\0\1
  &&\text{1, $\id|$}\\
  3.&& \fm{\A\rp\id}\0\1&\sep\fm\A\0\1
  &&\text{2, $\rp|$, no $\2$}\\
  4.&& \fm{\ttt\circ\A}\0\1&\sep\fm\A\0\1
  &&\text{3, \eqref{circ}, \eqref{t=id}}\\
\end{align*}
\proof[\ref{left-id}] 3-proof of $\A\circ\ttt\to\A$.
\begin{align*}
  1.&& \fm{\conv\A}\1\2&\sep\fm{\conv\A}\1\2
  &&\text{axiom}\\
  2.&& \fm{\conv\A}\1\0,\fm\id\0\2&\sep\fm{\conv\A}\0\1
  &&\text{1, $\id|$}\\
  3.&& \fm{\conv\id}\0\2&\sep\fm\id\0\2
  &&\text{\eqref{t*->t}}\\
  4.&& \fm{\conv\A}\1\2,\fm{\conv\id}\0\2&\sep \fm{\conv\A}\0\1
  &&\text{2, 3, Cut}\\
  5.&& \fm\id\0\2,\fm\A\2\1&\sep\fm\A\0\1
  &&\text{4, $\conv{}|$, $|\conv{}$}\\
  6.&& \fm{\id\rp\A}\0\1&\sep\fm\A\0\1
  &&\text{5, $\rp|$, no $\2$}\\
  7.&& \fm{\A\circ\ttt}\0\1&\sep\fm\A\0\1
  &&\text{6, \eqref{circ}, \eqref{t=id}}\\
\end{align*}
\proof[\ref{A->Aot*}] 3-proof of $\A\to\A\circ\star\ttt$.
\begin{align*}
  1.&&\fm\A\0\1&\sep\fm\A\0\1
  &&\text{axiom}\\
  2.&&\fm\A\0\1&\sep\fm{\conv\A}\1\0
  &&\text{1, $|\conv{}$}\\
  3.&&&\sep\fm\id\0\0
  &&\text{$|\id$}\\
  4.&&\fm\A\0\1&\sep\fm{\conv\A\rp\id}\1\0
  &&\text{2, 3, $|\rp$}\\
  5.&&\fm\A\0\1&\sep\fm{\conv{\conv\A\rp\id}}\0\1
  &&\text{4, $|\conv{}$}\\
  6.&&\fm{\conv\A}\1\2&\sep\fm{\conv\A}\1\2
  &&\text{axiom}\\
  7.&&\fm{\conv\A}\1\2&\sep\fm{\conv{\conv\A}}\2\1
  &&\text{6, $|\conv{}$}\\
  8.&&\fm\id\2\0&\sep\fm\id\2\0
  &&\text{axiom}\\
  9.&&\fm\id\2\0&\sep\fm{\conv\id}\0\2
  &&\text{8, $|\conv{}$}\\
  10.&&\fm{\conv\A}\1\2,\fm\id\2\0
  &\sep\fm{\conv\id\rp\conv{\conv\A}}\0\1
  &&\text{7, 9, $|\rp$}\\
  11.&&\fm{{\conv\A}\rp\id}\1\0
  &\sep\fm{\conv\id\rp\conv{\conv\A}}\0\1
  &&\text{10, $\rp|$, no $\2$}\\
  12.&&\fm{\conv{{\conv\A}\rp\id}}\0\1
  &\sep\fm{\conv\id\rp\conv{\conv\A}}\0\1
  &&\text{11, $\conv{}|$}\\
  13.&&\fm\A\0\1&\sep\fm{\conv\id\rp\conv{\conv\A}}\0\1
  &&\text{5, 12, Cut}\\
  14.&&\fm{\conv\id}\0\2&\sep\fm{\conv\id}\0\2
  &&\text{axiom}\\
  15.&&\fm\A\2\1&\sep\fm\A\2\1
  &&\text{axiom}\\
  16.&&\fm{\conv{\conv\A}}\2\1&\sep\fm\A\2\1
  &&\text{15, $\conv{}|$}\\
  17.&&\fm{\conv\id}\0\2,\fm{\conv{\conv\A}}\2\1
  &\sep\fm{\conv\id\rp\A}\0\1
  &&\text{14, 16, $|\rp$}\\
  18.&&\fm{\conv\id\rp\conv{\conv\A}}\0\1 &\sep\fm{\conv\id\rp\A}\0\1
  &&\text{17, $\rp|$, no $\2$}\\
  19.&&\fm\A\0\1&\sep\fm{\conv\id\rp\A}\0\1
  &&\text{13, 18, Cut}\\
  20.&&\fm\A\0\1&\sep\fm{\A\circ\star\ttt}\0\1 &&\text{19,
    \eqref{circ}, \eqref{star}, \eqref{t=id}}
\end{align*}
\proof[\ref{t->t*}] 3-proof of $\ttt\to\star\ttt$.
\begin{align*}
  1.&&\fm\id\0\1&\sep\fm{\conv\id\rp\id}\0\1
  &&\text{\eqref{A->Aot*}, $\A=\id$}\\
  2.&&\fm{\conv\id}\0\1&\sep\fm{\conv\id}\0\1
  &&\text{axiom}\\
  3.&&\fm{\conv\id}\0\2,\fm\id\2\1&\sep\fm{\conv\id}\0\1
  &&\text{2, $\id|$}\\
  4.&&\fm{{\conv\id}\rp\id}\0\1&\sep\fm{\conv\id}\0\1
  &&\text{3, $\rp|$, no $\2$}\\
  5.&&\fm\id\0\1&\sep\fm{\conv\id}\0\1
  &&\text{1, 4, Cut}\\
  6.&&\fm\ttt\0\1&\sep\fm{\star\ttt}\0\1
  &&\text{5, \eqref{star}, \eqref{t=id}}
\end{align*}
\proof[\ref{t*->t}] 3-proof of $\star\ttt\to\ttt$.
\begin{align*}
  1.&&\fm\id\1\0&\sep\fm{\conv\id}\1\0
  &&\text{\eqref{t->t*}, \eqref{equiv-5}}\\
  2.&&\fm{\conv\id}\0\1&\sep\fm{\conv{\conv\id}}\0\1
  &&\text{1, $\conv{}|$, $|\conv{}$}\\
  3.&&\fm\id\0\1&\sep\fm\id\0\1
  &&\text{axiom}\\
  4.&&\fm{\conv{\conv\id}}\0\1&\sep\fm\id\0\1
  &&\text{3, $\conv{}|$}\\
  5.&&\fm{\conv\id}\0\1&\sep\fm\id\0\1
  &&\text{2, 4, Cut}\\
  6.&&\fm{\star\ttt}\0\1&\sep\fm\ttt\0\1 
  &&\text{5, \eqref{star}, \eqref{t=id}}
\end{align*}
\proof[\ref{min-th}] 3-proof of $\ttt\land\rmin\ttt\to\A$.
\begin{align*}
  1.&&\fm\id\0\1&\sep\fm{\conv\id\rp\id}\0\1
  &&\text{\eqref{A->Aot*}}\\
  2.&&\fm{\conv\id}\0\1&\sep\fm{\conv\id}\0\1,\fm\A\0\1
  &&\text{axiom}\\
  3.&&\fm{\conv\id}\0\2,\fm\id\2\1&\sep\fm{\conv\id}\0\1,\fm\A\0\1
  &&\text{2, $\id|$}\\
  4.&&\fm{{\conv\id}\rp\id}\0\1&\sep\fm{\conv\id}\0\1,\fm\A\0\1
  &&\text{3, $\rp|$, no $\2$}\\
  5.&&\fm\id\0\1&\sep\fm{\conv\id}\0\1,\fm\A\0\1
  &&\text{1, 4, Cut}\\
  6.&&\fm\id\0\1,\fm{\min{\conv\id}}\0\1&\sep\fm\A\0\1
  &&\text{5, $\min\blank|$}\\
  7.&&\fm{\id\cdot\min{\conv\id}}\0\1&\sep\fm\A\0\1
  &&\text{6, $\cdot|$}\\
  8.&&\fm{\ttt\land\rmin\ttt}\0\1&\sep\fm\A\0\1 
  &&\text{7, \eqref{and}, \eqref{rmin}, \eqref{t=id}}
\end{align*}
\proof[\ref{adj}] 1-proof of $\A,\B\proves\A\land\B$. This is the
adjunction rule.
\begin{align*}
  1.&&&\sep\fm\A\0\0
  &&\text{sequent in an $\n$-proof}\\
  2.&&&\sep\fm\B\0\0
  &&\text{sequent in an $\n$-proof}\\
  3.&&&\sep\fm{\A\cdot\B}\0\0
  &&\text{1, 2, $|\cdot$}\\
  4.&&&\sep\fm{\A\land\B}\0\0 
  &&\text{3, \eqref{and}}
\end{align*}
\proof[\ref{modus-p}] 1-proof of $\A\to\B,\A\proves\B$. This rule is
modus ponens.  Note that \eqref{equiv-2} is applicable when $\i=\j$
and is therefore a 1-provable rule.
\begin{align*}
  1.&&&\sep\fm{\A\to\B}\0\0
  &&\text{sequent in an $\n$-proof}\\
  2.&&&\sep\fm\A\0\0
  &&\text{sequent in an $\n$-proof}\\
  3.&&\fm\A\0\0&\sep\fm\B\0\0
  &&\text{1, \eqref{equiv-2}}\\
  4.&&&\sep\fm\B\0\0 
  &&\text{2, 3, Cut}
\end{align*}
\proof[\ref{disj-syll}] 1-proof of $\A\lor\B,\rmin\A\proves\B$. This
rule is disjunctive syllogism.
\begin{align*}
  1.&&&\sep\fm{\A\lor\B}\0\0
  &&\text{sequent in an $\n$-proof}\\
  2.&&&\sep\fm{\rmin\A}\0\0
  &&\text{sequent in an $\n$-proof}\\
  3.&&&\sep\fm{\A+\B}\0\0
  &&\text{1, \eqref{or}}\\
  4.&&&\sep\fm{\min{\conv\A}}\0\0
  &&\text{2, \eqref{rmin}}\\
  5.&&\fm\A\0\0&\sep
  &&\text{4, $\min\blank|$, $\conv{}|$}\\
  6.&&\fm\B\0\0&\sep\fm\B\0\0
  &&\text{axiom}\\
  7.&&\fm{\A+\B}\0\0&\sep\fm\B\0\0
  &&\text{5, 6, $+|$}\\
  8.&&&\sep\fm\B\0\0 &&\text{3, 7, Cut}
\end{align*}
\proof[\ref{star-rule}] 1-proof of $\A\proves\A^*$.
\begin{align*}
  1.&&&\sep\fm\A\0\0
  &&\text{sequent in an $\n$-proof}\\
  2.&&&\sep\fm{\conv\A}\0\0
  &&\text{1, $|\conv{}$}\\
  3.&&&\sep\fm{\star\A}\0\0
  &&\text{2, \eqref{star}}
\end{align*}
\proof[\ref{tran-rule}] 2-proof of $\A\to\B,\B\to\C\proves
\A\to\C$. This is the transitivity rule. Note that \eqref{equiv-1}
requires $\i\neq\j$, and is therefore a 2-provable rule.
\begin{align*}
  1.&&&\sep\fm{\A\to\B}\0\0&&\text{sequent in an $\n$-proof}\\
  2.&&&\sep\fm{\B\to\C}\0\0&&\text{sequent in an $\n$-proof}\\
  3.&& \fm\A\1\0&\sep\fm\B\1\0
  &&\text{1, \eqref{equiv-2}}\\
  4.&&\fm\B\1\0&\sep\fm\C\1\0
  &&\text{2, \eqref{equiv-2}}\\
  5.&&\fm\A\1\0&\sep\fm\C\1\0  &&\text{3, 4, Cut}\\
  6.&&&\sep\fm{\A\to\C}\0\0 
  &&\text{5, \eqref{equiv-1}}
\end{align*}
\proof[\ref{contra-rule1}] 2-proof of
$\A\to\B\proves\rmin\B\to\rmin\A$. This is one form of the rule of
contraposition.
\begin{align*}
  1.&&&\sep\fm{\A\to\B}\0\0
  &&\text{sequent in an $\n$-proof}\\
  2.&&&\sep\fm{\A\to\B}\1\1
  &&\text{1, \eqref{equiv-4}}\\
  3.&&\fm\A\0\1&\sep\fm\B\0\1
  &&\text{2, \eqref{equiv-2}}\\
  4.&&\fm{\min{\conv\B}}\1\0&\sep\fm{\min{\conv\A}}\1\0
  &&\text{3, $\conv{}|$, $|\conv{}$, $\min\blank|$, $|\min\blank$}\\
  5.&&\fm{\rmin\B}\1\0&\sep\fm{\rmin\A}\1\0
  &&\text{4, \eqref{rmin}}\\
  6.&&&\sep\fm{\rmin\B\to\rmin\A}\0\0 
  &&\text{5, \eqref{equiv-1}}
\end{align*}
\proof[\ref{contra-rule2}] 2-proof of
$\A\to\rmin\B\proves\B\to\rmin\A$. This is another form of the rule of
contraposition.
\begin{align*}
  1.&&&\sep\fm{\A\to\rmin\B}\0\0
  &&\text{sequent in an $\n$-proof}\\
  2.&&\fm\A\1\0&\sep\fm{\rmin\B}\1\0
  &&\text{1, \eqref{equiv-2}}\\
  3.&&\fm\A\1\0&\sep\fm{\min{\conv\B}}\1\0
  &&\text{2, \eqref{rmin}}\\
  4.&& \fm\B\0\1&\sep\fm\B\0\1
  &&\text{axiom}\\
  5.&& \fm\B\0\1,\fm{\min{\conv\B}}\1\0&\sep
  &&\text{4, $|\conv{}$, $\min\blank|$}\\
  6.&& \fm\A\1\0,\fm\B\0\1&\sep
  &&\text{3, 5, Cut}\\
  7.&&\fm\B\0\1&\sep\fm{\min{\conv\A}}\0\1,
  &&\text{6, $\conv{}|$, $|\min\blank$}\\
  8.&&\fm\B\0\1&\sep\fm{\rmin\A}\0\1
  &&\text{7, \eqref{rmin}}\\
  9.&&&\sep\fm{\B\to\rmin\A}\1\1
  &&\text{8, \eqref{equiv-1}}\\
 10.&&&\sep\fm{\B\to\rmin\A}\0\0
  &&\text{9, \eqref{equiv-4}}
\end{align*}
\proof[\ref{cut-rule}] 2-proof of
$\A\land\B\to\C,\B\to\C\lor\A\proves\B\to\C$.  This rule is a derived
rule in Basic Logic, where it is called DR2
\cite[p.\,291]{Routleyetal1982} (derived rule number 2).
\begin{align*}
  1.&&&\sep\fm{\B\to\C\lor\A}\0\0
  &&\text{sequent in an $\n$-proof}\\
  2.&&&\sep\fm{\A\land\B\to\C}\0\0
  &&\text{sequent in an $\n$-proof}\\
  3.&&&\sep\fm{\B\to\C+\A}\0\0&&\text{1, \eqref{or}}\\
  4.&&&\sep\fm{\A\cdot\B\to\C}\0\0&&\text{2, \eqref{and}}\\
  5.&&\fm\B\1\0&\sep\fm{\C+\A}\1\0
  &&\text{3, \eqref{equiv-2}}\\
  6.&&\fm{\A\cdot\B}\1\0&\sep\fm\C\1\0
  &&\text{4, \eqref{equiv-2}}\\
  7.&&\fm\A\1\0&\sep \fm\A\1\0&&\text{axiom}\\
  8.&&\fm\B\1\0&\sep \fm\B\1\0&&\text{axiom}\\
  9.&&\fm\C\1\0&\sep\fm\C\1\0&&\text{axiom}\\
  10.&&\fm{\C+\A}\1\0&\sep\fm\A\1\0,\fm\C\1\0
  &&\text{7, 9, $+|$}\\
  11.&&\fm\B\1\0&\sep\fm\A\1\0,\fm\C\1\0&&\text{5, 10, Cut}\\
  12.&&\fm\A\1\0,\fm\B\1\0&\sep \fm{\A\cdot\B}\1\0
  &&\text{7, 8, $|\cdot$}\\
  13.&&\fm\A\1\0,\fm\B\1\0&\sep \fm\C\1\0&&\text{6, 12, Cut}\\
  14.&&\fm\B\1\0&\sep\fm\C\1\0&&\text{11, 13, Cut}\\
  15.&&&\sep\fm{\B\to\C}\0\0 
  &&\text{14, \eqref{equiv-1}}
\end{align*}
\proof[\ref{E-rule}] 2-proof of $\A\proves(\A\to\B)\to\B$.  This is
the E-rule \cite[p.\,8]{MR3728341}, also called BR1
\cite[p.\,289]{Routleyetal1982} (basic rule number 1) and R5
\cite[p.\,193]{MR3728341} (rule number 5).
\begin{align*}
  1.&&&\sep\fm\A\0\0
  &&\text{sequent in an $\n$-proof}\\
  2.&&&\sep\fm\A\1\1
  &&\text{1, \eqref{equiv-4}}\\
  3.&&&\sep\fm{\conv\A}\1\1
  &&\text{2, $|\conv{}$}\\
  4.&&\fm\B\1\0&\sep\fm\B\1\0
  &&\text{axiom}\\
  5.&&&\sep\fm{\min\B}\1\0,\fm\B\1\0
  &&\text{4, $|\min\blank$}\\
  6.&& &\sep\fm{\conv\A\rp\min\B}\1\0,\fm\B\1\0
  &&\text{3, 5, $|\rp$}\\
  7.&& \fm{\min{\conv\A\rp\min\B}}\1\0&\sep\fm\B\1\0
  &&\text{6, $\min\blank|$}\\
  8.&& \fm{\A\to\B}\1\0&\sep\fm\B\1\0
  &&\text{7, \eqref{to}}\\
  9.&&&\sep\fm{(\A\to\B)\to\B}\0\0
  &&\text{8, \eqref{equiv-1}}
\end{align*}
\proof[\ref{antilogism}] 2-proof of
$\A\land\B\to\neg\C\proves\A\land\C\to\neg\B$. This rule has been
called antilogism, a term coined by Christine
\citet{Ladd-Franklin1928}.
\begin{align*}
  1.&&&\sep\fm{\A\land\B\to\neg\C}\0\0
  &&\text{sequent in an $\n$-proof}\\
  2.&&&\sep\fm{\A\cdot\B\to\min\C}\0\0
  &&\text{1, \eqref{and}, \eqref{neg}}\\
  3.&&\fm{\A\cdot\B}\1\0&\sep\fm{\min\C}\1\0
  &&\text{2, \eqref{equiv-2}}\\
  4.&&\fm\A\1\0&\sep\fm\A\1\0
  &&\text{axiom}\\
  5.&&\fm\B\1\0&\sep\fm\B\1\0
  &&\text{axiom}\\
  6.&&\fm\A\1\0,\fm\B\1\0&\sep\fm{\A\cdot\B}\1\0
  &&\text{4, 5, $|\cdot$}\\
  7.&&\fm\A\1\0,\fm\B\1\0&\sep\fm{\min\C}\1\0
  &&\text{3, 6, Cut}\\
  8.&&\fm\C\1\0&\sep\fm\C\1\0
  &&\text{axiom}\\
  9.&&\fm\C\1\0,\fm{\min\C}\1\0&\sep
  &&\text{8, $\min\blank|$}\\
  10.&&\fm\A\1\0,\fm\B\1\0,\fm\C\1\0&\sep
  &&\text{7, 9, Cut}\\
  11.&&\fm{\A\cdot\C}\1\0,\fm\B\1\0&\sep
  &&\text{10, $\cdot|$}\\
  12.&&\fm{\A\cdot\C}\1\0&\sep\fm{\min\B}\1\0
  &&\text{11, $|\min\blank$}\\
  13.&&&\sep\fm{\A\cdot\C\to\min\B}\0\0
  &&\text{12, \eqref{equiv-1}}\\
  14.&&&\sep\fm{\A\land\C\to\neg\B}\0\0
  &&\text{13, \eqref{and}, \eqref{neg}}
\end{align*}
\proof[\ref{suff-rule}] 3-proof of
$\A\to\B\proves(\B\to\C)\to(\A\to\C)$.  This is the suffixing rule.
\begin{align*}
  1.&&&\sep\fm{\A\to\B}\0\0
  &&\text{sequent in an $\n$-proof}\\
  2.&&\fm\A\2\0&\sep\fm\B\2\0
  &&\text{1, \eqref{equiv-2}}\\
  3.&&\fm{\conv\A}\0\2&\sep\fm{\conv\B}\0\2
  &&\text{2, $\conv{}|$, $|\conv{}$}\\
  4.&&\fm{\min\C}\2\1&\sep\fm{\min\C}\2\1
  &&\text{axiom}\\
  5.&&\fm{\conv\A}\0\2,\fm{\min\C}\2\1
  &\sep\fm{\conv\B\rp\min\C}\0\1
  &&\text{3, 4, $|\rp$}\\
  6.&&\fm{\conv\A\rp\min\C}\0\1&\sep\fm{\conv\B\rp\min\C}\0\1
  &&\text{5, ${\rp}|$,~no~$\2$}\\
  7.&&\fm{\min{\conv\B\rp\min\C}}\0\1
  &\sep\fm{\min{\conv\A\rp\min\C}}\0\1
  &&\text{6, $\min\blank|$, $|\min\blank$}\\
  8.&&\fm{\B\to\C}\0\1&\sep\fm{\A\to\C}\0\1
  &&\text{7, \eqref{to}}\\
  9.&&&\sep\fm{(\B\to\C)\to(\A\to\C)}\0\0 
  &&\text{8, \eqref{equiv-1}, \eqref{equiv-4}}
\end{align*} 

\proof[\ref{pref-rule}] 3-proof of
$\A\to\B\proves(\C\to\A)\to(\C\to\B)$. This is the prefixing rule. The
prefixing axiom $(\A\to\B)\to((\C\to\A) \to(\C\to\B))$ \eqref{no.1} is
4-provable, so if $\A\to\B$ is 4-provable then
$((\C\to\A)\to(\C\to\B))$ is also 4-provable by modus ponens
\eqref{modus-p}. However, the prefixing rule only needs three
variables.
\begin{align*}
  1.&&&\sep\fm{\A\to\B}\0\0
  &&\text{sequent in an $\n$-proof}\\
  2.&&\fm\A\2\0&\sep\fm\B\2\0
  &&\text{1, \eqref{equiv-2}}\\
  3.&&\fm{\min\B}\2\0&\sep\fm{\min\A}\2\0
  &&\text{2, $\min\blank|$, $|\min\blank$}\\
  4.&&\fm{\conv\C}\1\2&\sep\fm{\conv\C}\1\2
  &&\text{axiom}\\
  5.&&\fm{\conv\C}\1\2,\fm{\min\B}\2\0
  &\sep\fm{\conv\C\rp\min\A}\1\0
  &&\text{3, 4, $|\rp$}\\
  6.&&\fm{\conv\C\rp\min\B}\1\0&\sep\fm{\conv\C\rp\min\A}\1\0
  &&\text{5, ${\rp}|$,~no~$\2$}\\
  7.&&\fm{\min{\conv\C\rp\min\A}}\1\0
  &\sep\fm{\min{\conv\C\rp\min\B}}\1\0
  &&\text{6, $\min\blank|$, $|\min\blank$}\\
  8.&&\fm{\C\to\A}\1\0&\sep\fm{\C\to\B}\1\0
  &&\text{7, \eqref{to}}\\
  9.&&&\sep\fm{(\C\to\A)\to(\C\to\B)}\0\0 
  &&\text{8, \ref{equiv-1}}
\end{align*}
\proof[\ref{affixing}] 3-provability of
$\A\to\B,\C\to\D\proves(\B\to\C)\to(\A\to\D)$.  This is the affixing
rule. Assume $\A\to\B$ and $\C\to\D$ are 3-provable.  By the
3-provable prefixing rule \eqref{pref-rule}, $(\A\to\C) \to (\A\to\D)$
is 3-provable. By the 3-provable suffixing rule \eqref{suff-rule},
$(\B\to\C)\to(\A\to\C)$ is 3-provable.  Hence, by the 2-provable
transitivity rule \eqref{tran-rule}, $(\B\to\C)\to(\A\to\D)$ is
3-provable.

\proof[\ref{fus-rule}] 3-proof of
$\A\to\B,\C\to\D\proves\A\circ\C\to\B\circ\D$.  This is the rule that
fusion is monotonic.  It preserves 4-provability because the monotonic
fusion axiom \eqref{no.4} is 4-provable. However, the monotonic fusion
rule is actually 3-provable.
\begin{align*}
  1.&&&\sep\fm{\C\to\D}\0\0
  &&\text{sequent in an $\n$-proof}\\
  2.&&&\sep\fm{\A\to\B}\0\0
  &&\text{sequent in an $\n$-proof}\\
  3.&&\fm\C\1\2&\sep\fm\D\1\2
  &&\text{1, \eqref{equiv-4}, \eqref{equiv-2}}\\
  4.&&\fm\A\2\0&\sep\fm\B\2\0
  &&\text{2, \eqref{equiv-2}}\\
  5.&&\fm\C\1\2,\fm\A\2\0&\sep\fm{\D\rp\B}\1\0
  &&\text{3, 4, $|\rp$}\\
  6.&&\fm{\C\rp\A}\1\0&\sep\fm{\D\rp\B}\1\0
  &&\text{5, $\rp|$, no $\2$}\\
  7.&&&\sep\fm{\C\rp\A\to\D\rp\B}\0\0
  &&\text{6, \eqref{equiv-1}}\\
  8.&&&\sep\fm{\A\circ\C\to\B\circ\D}\0\0 
  &&\text{7, \eqref{circ}}
\end{align*} 
\proof[\ref{cycling}] 3-proof of
$\A\to(\B\to\C)\proves\B\to(\rmin\C\to\rmin\A)$.  This is the cycling
rule.
\begin{align*}
  1.&&&\sep\fm{\A\to(\B\to\C)}\0\0
  &&\text{sequent in an $\n$-proof}\\
  2.&&\fm\A\0\2&\sep\fm{\B\to\C}\0\2
  &&\text{1, \eqref{equiv-4}, \eqref{equiv-2}}\\
  3.&&\fm\A\0\2&\sep\fm{\min{\conv\B\rp\min\C}}\0\2
  &&\text{2, \eqref{to}}\\
  4.&&\fm\B\1\0&\sep\fm\B\1\0
  &&\text{axiom}\\
  5.&&\fm\B\1\0&\sep\fm{\conv\B}\0\1
  &&\text{4, $|\conv{}$}\\
  6.&&\fm\C\1\2&\sep\fm\C\1\2
  &&\text{axiom}\\
  7.&&&\sep\fm\C\1\2,\fm{\min\C}\1\2
  &&\text{6, $|\min\blank$}\\
  8.&&\fm\B\1\0&\sep\fm\C\1\2,\fm{\conv\B\rp\min\C}\0\2
  &&\text{5, 7, $|\rp$ }\\
  9.&&\fm\B\1\0,\fm{\min{\conv\B\rp\min\C}}\0\2&\sep\fm\C\1\2
  &&\text{8, $\min\blank|$}\\
  10.&&\fm\A\0\2,\fm\B\1\0&\sep\fm\C\1\2
  &&\text{3, 9, Cut}\\
  11.&&\fm\B\1\0,\fm{\conv{\min{\conv\C}}}\1\2&
  ,\fm{\min{\min{\conv\A}}}\2\0\sep
  &&\text{10, $\conv{}|$, $|\conv{}$, $\min\blank|$, $|\min\blank$}\\
  12.&&\fm\B\1\0,\fm{\conv{\min{\conv\C}}
    &\rp\min{\min{\conv\A}}}\1\0\sep
  &&\text{11, $\rp|$,~no~$\2$}\\
  13.&&\fm\B\1\0&\sep\fm{\min{\conv{\min{\conv\C}}\rp
      \min{\min{\conv\A}}}}\1\0
  &&\text{12, $|\min\blank$}\\
  14.&&\fm\B\1\0&\sep\fm{\rmin\C\to\rmin\A}\1\0
  &&\text{13, \eqref{rmin}, \eqref{to}}\\
  15.&&&\sep\fm{\B\to(\rmin\C\to\rmin\A)}\0\0 
  &&\text{14, \eqref{equiv-1}}
\end{align*}
\proof[\ref{no.1}] 4-proof of
$(\A\to\B)\to((\C\to\A)\to(\C\to\B))$. This is the prefixing axiom.
\begin{align*}
  1.&&\fm\A\3\0&\sep\fm\A\3\0
  &&\text{axiom}\\
  2.&&&\sep\fm{\min\A}\3\0,\fm{\conv\A}\0\3
  &&\text{1, $|\min\blank$, $|\conv{}$}\\
  3.&&\fm{\min\B}\3\1&\sep\fm{\min\B}\3\1
  &&\text{axiom}\\
  4.&&\fm{\min\B}\3\1&\sep\fm{\min\A}\3\0,\fm{\conv\A\rp\min\B}\0\1
  &&\text{2, 3, $|\rp$}\\
  5.&&\fm{\conv\C}\2\3&\sep\fm{\conv\C}\2\3
  &&\text{axiom}\\
  6.&&\fm{\conv\C}\2\3,\fm{\min\B}\3\1
  &\sep\fm{\conv\C\rp\min\A}\2\0,\fm{\conv\A\rp\min\B}\0\1
  &&\text{4, 5, $|\rp$}\\
  7.&&\fm{\conv\C\rp\min\B}\2\1
  &\sep\fm{\conv\C\rp\min\A}\2\0,\fm{\conv\A\rp\min\B}\0\1
  &&\text{6, $\rp|$, no $\3$}\\
  8.&&\fm{\min{\conv\C\rp\min\A}}\2\0&,\fm{\min{\conv\A\rp\min\B}}\0\1
  \sep\fm{\min{\conv\C\rp\min\B}}\2\1
  &&\text{7, $|\min\blank$, $\min\blank|$}\\
  9.&&\fm{\C\to\A}\2\0&,\fm{\A\to\B}\0\1\sep\fm{\C\to\B}\2\1
  &&\text{8, \eqref{to}}\\
  10.&&\fm{\conv{\C\to\A}}\0\2&,\fm{\min{\C\to\B}}\2\1,
  \fm{\A\to\B}\0\1\sep
  &&\text{9, $\min\blank|$, $\conv{}|$}\\
  11.&&\fm{\conv{\C\to\A}\rp\min{\C\to\B}}\0\1&,\fm{\A\to\B}\0\1\sep
  &&\text{10, $\rp|$, no $\2$}\\
  12.&&\fm{\A\to\B}\0\1
  &\sep\fm{\min{\conv{\C\to\A}\rp\min{\C\to\B}}}\0\1
  &&\text{11, $|\min\blank$}\\
  13.&&\fm{\A\to\B}\0\1&\sep\fm{(\C\to\A)\to(\C\to\B)}\0\1 
  &&\text{12, \eqref{to}}
\end{align*}
\proof[\ref{no.2}] 4-proof of
$(\A\to(\B\to\C))\to(\A\circ\B\to\C)$. This is the axiom of bunching
hypotheses.
\begin{align*}
  1.&&\fm\A\3\0&\sep\fm\A\3\0
  &&\text{axiom}\\
  2.&&\fm\B\2\3&\sep\fm\B\2\3
  &&\text{axiom}\\
  3.&&\fm\C\2\1&\sep\fm\C\2\1
  &&\text{axiom}\\
  4.&&\fm\B\2\3&\sep\fm{\conv\B}\3\2
  &&\text{2, $|\conv{}$}\\
  5.&&&\sep\fm\C\2\1,\fm{\min\C}\2\1
  &&\text{3, $|\min\blank$}\\
  6.&&\fm\B\2\3&\sep\fm\C\2\1,\fm{\conv\B\rp\min\C}\3\1
  &&\text{4, 5, $\rp$}\\
  7.&&\fm\B\2\3
  &\sep\fm\C\2\1,\fm{\min{\min{\conv\B\rp\min\C}}}\3\1
  &&\text{6, $\min\blank|$, $|\min\blank$}\\
  8.&&\fm\A\3\0&\sep\fm{\conv\A}\0\3
  &&\text{1, $|\conv{}$}\\
  9.&&\fm\B\2\3,\fm\A\3\0
  &\sep\fm\C\2\1,\fm{\conv\A\rp\min{\min{\conv\B\rp\min\C}}}\0\1
  &&\text{7, 8, $|\rp$}\\
  10.&&\fm{\B\rp\A}\2\0
  &\sep\fm\C\2\1,\fm{\conv\A\rp\min{\min{\conv\B\rp\min\C}}}\0\1
  &&\text{9, $\rp|$, no $\3$}\\
  11.&&\fm{\conv{\B\rp\A}}\0\2,\fm{\min\C}\2\1
  &\sep\fm{\conv\A\rp\min{\min{\conv\B\rp\min\C}}}\0\1
  &&\text{10, $\conv{}|$, $\min\blank|$}\\
  12.&&\fm{\conv{\B\rp\A}\rp\min\C}\0\1
  &\sep\fm{\conv\A\rp\min{\min{\conv\B\rp\min\C}}}\0\1
  &&\text{11, $\rp|$, no $\2$}\\
  13.&&\fm{\min{\conv\A\rp\min{\min{\conv\B\rp\min\C}}}}\0\1
  &\sep\fm{\min{\conv{\B\rp\A}\rp\min\C}}\0\1
  &&\text{12, $|\min\blank$, $\min\blank|$}\\
  14.&&\fm{\A\to(\B\to\C)}\0\1&\sep\fm{\A\circ\B\to\C}\0\1
  &&\text{13, \eqref{to}, \eqref{circ}}
\end{align*}
\proof[\ref{no.3}] 4-proof of
$(\A\circ\B\to\C)\to(\A\to(\B\to\C))$. This is the converse of
bunching.
\begin{align*}
  1.&&\fm\A\2\0&\sep\fm\A\2\0
  &&\text{axiom}\\
  2.&&\fm\B\3\2&\sep\fm\B\3\2
  &&\text{axiom}\\
  3.&&\fm{\min\C}\3\1&\sep\fm{\min\C}\3\1
  &&\text{axiom}\\
  4.&&\fm\A\2\0,\fm\B\3\2&\sep\fm{\B\rp\A}\3\0
  &&\text{1, 2, $|\rp$}\\
  5.&&\fm{\conv\A}\0\2,\fm{\conv\B}\2\3
  &\sep\fm{\conv{\B\rp\A}}\0\3
  &&\text{4, $|\conv{}$, $\conv{}|$}\\
  6.&&\fm{\conv\A}\0\2,\fm{\conv\B}\2\3,\fm{\min\C}\3\1
  &\sep\fm{\conv{\B\rp\A}\rp\min\C}\0\1
  &&\text{3, 5, $|\rp$}\\
  7.&&\fm{\conv\A}\0\2,\fm{\conv\B\rp\min\C}\2\1
  &\sep\fm{\conv{\B\rp\A}\rp\min\C}\0\1
  &&\text{6, $\rp|$, no $\3$}\\
  8.&&\fm{\conv\A}\0\2,\fm{\min{\min{\conv\B\rp\min\C}}}\2\1
  &\sep\fm{\conv{\B\rp\A}\rp\min\C}\0\1
  &&\text{7, $|\min\blank$, $\min\blank|$}\\
  9.&&\fm{\conv\A\rp\min{\min{\conv\B\rp\min\C}}}\0\1
  &\sep\fm{\conv{\B\rp\A}\rp\min\C}\0\1
  &&\text{8, $\rp|$, no $\2$}\\
  10.&&\fm{\min{\conv{\B\rp\A}\rp\min\C}}\0\1
  &\sep\fm{\min{\conv\A\rp\min{\min{\conv\B\rp\min\C}}}}\0\1
  &&\text{9, $|\min\blank$, $\min\blank|$}\\
  11.&&\fm{\A\circ\B\to\C}\0\1&\sep\fm{\A\to(\B\to\C)}\0\1
  &&\text{10, \eqref{to}, \eqref{circ}}
\end{align*}
\proof[\ref{no.4}] 4-proof of
$(\A\to\B)\to(\A\circ\C\to\B\circ\C)$. This is monotonicity of fusion
in the left argument.
\begin{align*}
  1.&&\fm\A\3\0&\sep\fm\A\3\0
  &&\text{axiom}\\
  2.&&\fm\B\3\1&\sep\fm\B\3\1
  &&\text{axiom}\\
  3.&&\fm\C\2\3&\sep\fm\C\2\3
  &&\text{axiom}\\
  4.&&&\sep\fm{\min\B}\3\1,\fm\B\3\1
  &&\text{2, $|\min\blank$}\\
  5.&&\fm\C\2\3&\sep\fm{\min\B}\3\1,\fm{\C\rp\B}\2\1
  &&\text{3, 4, $|\rp$}\\
  6.&&\fm\A\3\0&\sep\fm{\conv\A}\0\3
  &&\text{1, $|\conv{}$}\\
  7.&&\fm\C\2\3,\fm\A\3\0
  &\sep\fm{\C\rp\B}\2\1,\fm{\conv\A\rp\min\B}\0\1
  &&\text{5, 6, $|\rp$}\\
  8.&&\fm{\C\rp\A}\2\0
  &\sep\fm{\C\rp\B}\2\1,\fm{\conv\A\rp\min\B}\0\1
  &&\text{7, $\rp|$, no $\3$}\\
  9.&&\fm{\conv{\C\rp\A}}\0\2,\fm{\min{\C\rp\B}}\2\1
  &\sep\fm{\conv\A\rp\min\B}\0\1
  &&\text{8, $\min\blank|$, $\conv{}|$}\\
  10.&&\fm{\conv{\C\rp\A}\rp\min{\C\rp\B}}\0\1
  &\sep\fm{\conv\A\rp\min\B}\0\1
  &&\text{9, $\rp|$, no $\2$}\\
  11.&&\fm{\min{\conv\A\rp\min\B}}\0\1
  &\sep\fm{\min{\conv{\C\rp\A}\rp\min{\C\rp\B}}}\0\1
  &&\text{10, $|\min\blank$, $\min\blank|$}\\
  12.&&\fm{\A\to\B}\0\1 &\sep\fm{\A\circ\C\to\B\circ\C}\0\1
  &&\text{11, \eqref{to}, \eqref{circ}}
\end{align*}
\proof[\ref{no.5}] 4-proof of
$(\A\circ\B)\circ\C\to\A\circ(\B\circ\C)$.  Fusion is associative in
one direction.
\begin{align*}
  1.&&\fm\B\2\3&\sep\fm\B\2\3
  &&\text{axiom}\\
  2.&&\fm\C\0\2&\sep\fm\C\0\2
  &&\text{axiom}\\
  3.&&\fm\C\0\2,\fm\B\2\3&\sep\fm{\C\rp\B}\0\3
  &&\text{1, 2, $|\rp$}\\
  4.&&\fm\A\3\1&\sep\fm\A\3\1
  &&\text{axiom}\\
  5.&&\fm\C\0\2,\fm\B\2\3,\fm\A\3\1&\sep
  \fm{(\C\rp\B)\rp\A}\0\1
  &&\text{3, 4, $|\rp$}\\
  6.&&\fm\C\0\2,\fm{\B\rp\A}\2\1&\sep\fm{(\C\rp\B)\rp\A}\0\1
  &&\text{5, $\rp|$, no $\3$}\\
  7.&&\fm{\C\rp(\B\rp\A)}\0\1&\sep\fm{(\C\rp\B)\rp\A}\0\1
  &&\text{6, $\rp|$, no $\2$}\\
  8.&&\fm{(\A\circ\B)\circ\C}\0\1&\sep\fm{\A\circ(\B\circ\C)}\0\1
  &&\text{7, \eqref{circ}}
\end{align*}
\proof[\ref{no.6}] 4-proof of
$\A\circ(\B\circ\C)\to(\A\circ\B)\circ\C$.  Fusion is associative in
the other direction.
\begin{align*}
  1.&&\fm\B\2\3&\sep\fm\B\2\3
  &&\text{axiom}\\
  2.&&\fm\A\3\1&\sep\fm\A\3\1
  &&\text{axiom}\\
  3.&&\fm\B\2\3,\fm\A\3\1&\sep\fm{\B\rp\A}\2\1
  &&\text{1, 2, $|\rp$}\\
  4.&&\fm\C\0\2&\sep\fm\C\0\2
  &&\text{axiom}\\
  5.&&\fm\C\0\2,\fm\B\2\3,\fm\A\3\1&\sep \fm{\C\rp(\B\rp\A)}\0\1
  &&\text{3, 4, $|\rp$}\\
  6.&&\fm{\C\rp\B}\0\3,\fm\A\3\2&\sep\fm{\C\rp(\B\rp\A)}\0\1
  &&\text{5, $\rp|$, no $\2$}\\
  7.&&\fm{(\C\rp\B)\rp\A}\0\1&\sep\fm{\C\rp(\B\rp\A)}\0\1
  &&\text{6, $\rp|$, no $\3$}\\
  8.&&\fm{\A\circ(\B\circ\C)}\0\1&\sep\fm{(\A\circ\B)\circ\C}\0\1
  &&\text{7, \eqref{circ}}
\end{align*}
\proof[\ref{reductio}] 3-proof of $(\A\to\rmin\A)\to\rmin\A$
from density.  This is the {\it reductio ad absurdum} axiom.
\begin{align*}
  1.&&\fm{\conv\A}\0\1&\sep\fm{\conv\A\rp\conv\A}\0\1
  &&\text{density}\\
  2.&&\fm{\conv\A}\0\2&\sep\fm{\conv\A}\0\2
  &&\text{axiom}\\
  3.&&\fm{\conv\A}\2\1&\sep\fm{\conv\A}\2\1
  &&\text{axiom}\\
  4.&&\fm{\conv\A}\2\1&\sep\fm{\min{\min{\conv\A}}}\2\1
  &&\text{3, $\min\blank|$, $|\min\blank$}\\
  5.&&\fm{\conv\A}\0\2,\fm{\conv\A}\2\1
  &\sep\fm{\conv\A\rp\min{\min{\conv\A}}}\0\1
  &&\text{2, 4, $|\rp$}\\
  6.&&\fm{\conv\A\rp\conv\A}\0\1
  &\sep\fm{\conv\A\rp\min{\min{\conv\A}}}\0\1
  &&\text{5, $\rp|$, no $\2$}\\
  7.&&\fm{\conv\A}\0\1&\sep\fm{\conv\A\rp\min{\min{\conv\A}}}\0\1
  &&\text{1, 6, Cut}\\
  8.&&\fm{\min{\conv\A\rp\min{\min{\conv\A}}}}\0\1
  &\sep\fm{\min{\conv\A}}\0\1
  &&\text{7, $\min\blank|$, $|\min\blank$}\\
  9.&&\fm{\A\to\rmin\A}\0\1&\sep\fm{\rmin\A}\0\1 
  &&\text{8, \eqref{rmin}, \eqref{to}}
\end{align*}
\proof[\ref{contract5}] 3-proof of $\A\land\B\to\A\circ\B$ from
density.  The 3-proof of \eqref{contract5} is the first one in which
we have a real need for a derived rule of the sequent calculus called
{\bf weakening}, indicated by a ``W''. Its form as a rule that could
have been included in Table \ref{rules} is
\begin{equation*}
  \boxed{W}\quad\frac{\Gamma\sep\Delta}
        {\Gamma,\Gamma'\sep\Delta,\Delta'}
\end{equation*}
It is easily proved by induction on the lengths of $\n$-proofs that
this rule can be admitted without effect on the notion of
$\n$-provability.  The idea is that any predicate that one wishes to
add later via weakening can simply be added to the previous
sequents. The form of the rules allows this; such additions take
instances of rules to instances of rules. The base step of the
induction is that weakenings of axioms are still axioms. This fact was
used in the 2-proofs of \eqref{top} and \eqref{explosion} (sequent
number 1 has a superflous formula in it). It was used in the 3-proof
of \eqref{min-th}.  The sequents corresponding to the equations of
density, commutativity, and symmetry are also closed under weakenings;
see the proof of \eqref{symm-ax}.  The proofs of \eqref{contract5},
\eqref{contract4}, and \eqref{contract3} would be unduly cluttered by
actually using this device, so we use the weakening rule instead.
\begin{align*}
  1.&&\fm{\A\cdot\B}\0\1&\sep\fm{(\A\cdot\B)\rp(\A\cdot\B)}\0\1
  &&\text{density}\\
  2.&&\fm\A\2\1&\sep\fm\A\2\1
  &&\text{axiom}\\
  3.&&\fm\B\0\2&\sep\fm\B\0\2
  &&\text{axiom}\\
  4.&&\fm\B\0\2,\fm\A\2\1&\sep\fm{\B\rp\A}\0\1
  &&\text{2, 3, $|\rp$}\\
  5.&&\fm\A\0\2,\fm\B\0\2,\fm\A\2\1,\fm\B\2\1&\sep\fm{\B\rp\A}\0\1
  &&\text{4, W}\\
  6.&&\fm{\A\cdot\B}\0\2,\fm{\A\cdot\B}\2\1&\sep\fm{\B\rp\A}\0\1
  &&\text{5, $\cdot|$}\\
  7.&&\fm{(\A\cdot\B)\rp(\A\cdot\B)}\0\1&\sep\fm{\B\rp\A}\0\1
  &&\text{6, $\rp|$, no $\2$}\\
  8.&&\fm{\A\cdot\B}\0\1&\sep\fm{\B\rp\A}\0\1
  &&\text{1, 7, Cut}\\
  9.&&\fm{\A\land\B}\0\1&\sep\fm{\A\circ\B}\0\1 
  &&\text{8, \eqref{and}, \eqref{circ}}
\end{align*}
\proof[\ref{contract4}] 3-proof of $(\A\to\B)\to\rmin\A\lor\B$ from
density.
\begin{align*}
  1.&&\fm{\A\cdot\min{\conv\B}}\1\0
  &\sep\fm{(\A\cdot\min{\conv\B})\rp(\A\cdot\min{\conv\B})}\1\0
  &&\text{density}\\
  2.&&\fm\B\2\1&\sep\fm\B\2\1
  &&\text{axiom}\\
  3.&&\fm{\min{\conv\B}}\1\2&\sep\fm{\min\B}\2\1
  &&\text{2, $|\conv{}$, $\min\blank|$, $|\min\blank$}\\
  4.&&\fm\A\2\0&\sep\fm\A\2\0
  &&\text{axiom}\\
  5.&&\fm\A\2\0&\sep\fm{\conv\A}\0\2
  &&\text{4, $|\conv{}$}\\
  6.&&\fm{\min{\conv\B}}\1\2,\fm\A\2\0&\sep\fm{\conv\A\rp\min\B}\0\1
  &&\text{3, 5, $|\rp$}\\
  7.&&\fm\A\1\2,\fm{\min{\conv\B}}\1\2&,\fm\A\2\0,
  \fm{\min{\conv\B}}\2\0 \sep\fm{\conv\A\rp\min\B}\0\1
  &&\text{6, W}\\
  8.&&\fm{(\A\cdot\min{\conv\B})}\1\2&,
  \fm{(\A\cdot\min{\conv\B})}\2\0 \sep\fm{\conv\A\rp\min\B}\0\1
  &&\text{7, $\cdot|$}\\
  9.&&\fm{(\A\cdot\min{\conv\B})&\rp(\A\cdot\min{\conv\B})}\1\0
  \sep\fm{\conv\A\rp\min\B}\0\1
  &&\text{8, $\rp|$, no $\2$}\\
  10.&&\fm{\A\cdot\min{\conv\B}}\1\0&\sep\fm{\conv\A\rp\min\B}\0\1
  &&\text{1, 9, Cut}\\
  11.&&\fm{\min{\conv\A\rp\min\B}}\0\1
  &\sep\fm{\min{\A\cdot\min{\conv\B}}}\1\0
  &&\text{10, $|\min\blank$, $\min\blank|$}\\
  12.&&\fm\A\1\0&\sep\fm\A\1\0
  &&\text{axiom}\\
  13.&&\fm\B\0\1&\sep\fm\B\0\1
  &&\text{axiom}\\
  14.&&&\sep\fm{\min{\conv\B}}\1\0,\fm\B\0\1
  &&\text{13, $\conv|$, $|\min\blank$}\\
  15.&&\fm\A\1\0&\sep\fm{\A\cdot\min{\conv\B}}\1\0,\fm\B\0\1
  &&\text{12, 14, $|\cdot$}\\
  16.&&\fm{\min{\A\cdot\min{\conv\B}}}\1\0
  &\sep\fm{\min{\conv\A}}\0\1,\fm\B\0\1
  &&\text{15, $\min\blank|$, $\conv{}|$, $|\min\blank$}\\
  17.&&\fm{\min{\A\cdot\min{\conv\B}}}\1\0&\sep\fm{\min{\conv\A}+\B}\0\1
  &&\text{16, $|+$}\\
  18.&&\fm{\min{\conv\A\rp\min\B}}\0\1&\sep\fm{\min{\conv\A}+\B}\0\1
  &&\text{11, 17, Cut}\\
  19.&&\fm{\A\to\B}\0\1&\sep\fm{\rmin\A\lor\B}\0\1 
  &&\text{18, \eqref{or}, \eqref{rmin}, \eqref{to}}
\end{align*}
\proof[\ref{contract2}] 4-proof of $(\A\to(\A\to\B))\to(\A\to\B)$ from
density.  This is the contraction axiom.
\begin{align*}
  1.&&\fm{\conv\A}\0\2&\sep\fm{\conv\A\rp\conv\A}\0\2
  &&\text{density}\\
  2.&&\fm{\conv\A}\0\3&\sep\fm{\conv\A}\0\3
  &&\text{axiom}\\
  3.&& \fm{\conv\A}\3\2&\sep\fm{\conv\A}\3\2
  &&\text{axiom}\\
  4.&& \fm{\min\B}\2\1&\sep\fm{\min\B}\2\1
  &&\text{axiom}\\
  5.&&\fm{\conv\A}\3\2,\fm{\min\B}\2\1
  &\sep\fm{\conv\A\rp\min\B}\3\1
  &&\text{3, 4, $|\rp$}\\
  6.&& \fm{\conv\A}\3\2,\fm{\min\B}\2\1
  &\sep\fm{\min{\min{\conv\A\rp\min\B}}}\3\1
  &&\text{5, $\min\blank|$, $|\min\blank$}\\
  7.&&\fm{\conv\A}\0\3,\fm{\conv\A}\3\2,\fm{\min\B}\2\1
  &\sep\fm{\conv\A\rp\min{\min{\conv\A\rp\min\B}}}\0\1
  &&\text{2, 6, $|\rp$}\\
  8.&&\fm{\conv\A\rp\conv\A}\0\2,\fm{\min\B}\2\1
  &\sep\fm{\conv\A\rp\min{\min{\conv\A\rp\min\B}}}\0\1
  &&\text{7, $\rp|$, no $\3$}\\
  9.&&\fm{\conv\A}\0\2,\fm{\min\B}\2\1
  &\sep\fm{\conv\A\rp\min{\min{\conv\A\rp\min\B}}}\0\1
  &&\text{1, 8, Cut}\\
  10.&&\fm{\conv\A\rp\min\B}\0\1
  &\sep\fm{\conv\A\rp\min{\min{\conv\A\rp\min\B}}}\0\1
  &&\text{9, $\rp|$, no $\2$}\\
  11.&&\fm{\min{\conv\A\rp\min{\min{\conv\A\rp\min\B}}}}\0\1
  &\sep\fm{\min{\conv\A\rp\min\B}}\0\1
  &&\text{10, $\min\blank|$, $|\min\blank$}\\
  12.&&\fm{\A\to(\A\to\B)}\0\1&\sep\fm{\A\to\B}\0\1 
  &&\text{11, \eqref{to}}
\end{align*}
\proof[\ref{contract3}] 4-proof of
$(\A\to(\B\to\C))\to(\A\land\B\to\C)$ from density.
\begin{align*}
  1.&&\fm{\A\cdot\B}\2\0&\sep\fm{(\A\cdot\B)\rp(\A\cdot\B)}\2\0
  &&\text{density}\\
  2.&&\fm\A\3\0&\sep\fm\A\3\0
  &&\text{axiom}\\
  3.&&\fm\A\3\0&\sep\fm{\conv\A}\0\3
  &&\text{2, $|\conv{}$}\\
  4.&&\fm\B\2\3&\sep\fm\B\2\3
  &&\text{axiom}\\
  5.&&\fm\B\2\3&\sep\fm{\conv\B}\3\2
  &&\text{4, $|\conv{}$}\\
  6.&&\fm{\min\C}\2\1&\sep\fm{\min\C}\2\1
  &&\text{axiom}\\
  7.&&\fm\B\2\3,\fm{\min\C}\2\1&\sep\fm{\conv\B\rp\min\C}\3\1
  &&\text{5, 6, $|\rp$}\\
  8.&&\fm\B\2\3,\fm{\min\C}\2\1
  &\sep\fm{\min{\min{\conv\B\rp\min\C}}}\3\1
  &&\text{7, $\min\blank|$, $|\min\blank$}\\
  9.&&\fm\B\2\3,\fm\A\3\0,\fm{\min\C}\2\1
  &\sep\fm{\conv\A\rp\min{\min{\conv\B\rp\min\C}}}\0\1
  &&\text{3, 8, $|\rp$}\\
  10.&&\fm\A\2\3,\fm\B\2\3,\fm\A\3\0,\fm\B\3\0&,
  \fm{\min\C}\2\1\sep\fm{\conv\A\rp\min{\min{\conv\B\rp\min\C}}}\0\1
  &&\text{9, W}\\
  11.&&\fm{(\A\cdot\B)}\2\3,\fm{(\A\cdot\B)}\3\0&,\fm{\min\C}\2\1
  \sep\fm{\conv\A\rp\min{\min{\conv\B\rp\min\C}}}\0\1
  &&\text{10, $\cdot|$}\\
  12.&&\fm{(\A\cdot\B)\rp(\A\cdot\B)}\2\0,\fm{\min\C}\2\1
  &\sep\fm{\conv\A\rp\min{\min{\conv\B\rp\min\C}}}\0\1
  &&\text{12, $\rp|$, no $\3$}\\
  13.&&\fm{\A\cdot\B}\2\0,\fm{\min\C}\2\1
  &\sep\fm{\conv\A\rp\min{\min{\conv\B\rp\min\C}}}\0\1
  &&\text{1, 12, Cut}\\
  14.&&\fm{\conv{\A\cdot\B}}\0\2,\fm{\min\C}\2\1
  &\sep\fm{\conv\A\rp\min{\min{\conv\B\rp\min\C}}}\0\1
  &&\text{13, $\conv{}|$}\\
  15.&&\fm{\conv{\A\cdot\B}\rp\min\C}\0\1
  &\sep\fm{\conv\A\rp\min{\min{\conv\B\rp\min\C}}}\0\1
  &&\text{14, $\rp|$, no $\2$}\\
  16.&&\fm{\min{\conv\A\rp\min{\min{\conv\B\rp\min\C}}}}\0\1
  &\sep\fm{\min{\conv{\A\cdot\B}\rp\min\C}}\0\1
  &&\text{15, $\min\blank|$, $|\min\blank$}\\
  17.&&\fm{\A\to(\B\to\C)}\0\1&\sep\fm{\A\land\B\to\C}\0\1 
  &&\text{16, \eqref{and}, \eqref{to}}
\end{align*}
\proof[\ref{mp}] 3-provability of $\A\to((\A\to\B)\to\B)$
from commutativity.  This is an axiomatic form of modus ponens.
\begin{align*}
  1.&&\fm\A\0\1&\sep\fm\A\0\1&&\text{axiom}\\
  2.&&\fm{\A\to\B}\2\0&\sep\fm{\A\to\B}\2\0&&\text{axiom}\\
  3.&&\fm\A\0\1,\fm{\A\to\B}\2\0&\sep\fm{(\A\to\B)\rp\A}\2\1
  &&\text{1, 2, $|\rp$}\\
  4.&&\fm{(\A\to\B)\rp\A}\2\1&\sep\fm{\A\rp(\A\to\B)}\2\1
  &&\text{commutativity}\\
  5.&&\fm\A\0\1,\fm{\A\to\B}\2\0&\sep\fm{\A\rp(\A\to\B)}\2\1
  &&\text{3, 4, Cut}\\
  6.&&\fm\A\2\0&\sep\fm\A\2\0&&\text{axiom}\\
  7.&&\fm\B\2\1&\sep\fm\B\2\1&&\text{axiom}\\
  8.&&\fm\A\2\0&\sep\fm{\conv\A}\0\2&&\text{6, $|\conv{}$}\\
  9.&&&\sep\fm{\min\B}\2\1,\fm\B\2\1&&\text{7, $|\min\blank$}\\
  10.&&\fm\A\2\0&\sep\fm{\conv\A\rp\min\B}\0\1,\fm\B\2\1
  &&\text{8, 9, $|\rp$}\\
  11.&&\fm\A\2\0,\fm{\min{\conv\A\rp\min\B}}\0\1&\sep\fm\B\2\1
  &&\text{10, $\min\blank|$}\\
  12.&&\fm\A\2\0,\fm{\A\to\B}\0\1&\sep\fm\B\2\1&&\text{11, \eqref{to}}\\
  13.&&\fm{\A\rp(\A\to\B)}\2\1&\sep\fm\B\2\1&&\text{12, $\rp|$, no $\0$}\\
  14.&&\fm\A\0\1,\fm{\A\to\B}\2\0&\sep\fm\B\2\1&&\text{5, 13, Cut}\\
  15.&&\fm\A\0\1,\fm{\A\to\B}\2\0,\fm{\min\B}\2\1&\sep
  &&\text{14, $\min\blank|$}\\
  16.&&\fm\A\0\1,\fm{\conv{\A\to\B}}\0\2,\fm{\min\B}\2\1&\sep
  &&\text{15, $\conv{}|$}\\
  17.&&\fm\A\0\1,\fm{\conv{\A\to\B}\rp\min\B}\0\1&\sep
  &&\text{16, $\rp|$, no $\2$}\\
  18.&&\fm\A\0\1&\sep\fm{\min{\conv{\A\to\B}\rp\min\B}}\0\1
  &&\text{17, $|\min\blank$}\\
  19.&&\fm\A\0\1&\sep\fm{(\A\to\B)\to\B}\0\1&&\text{18, \eqref{to}}
\end{align*}
\proof[\ref{contra}] 3-proof of $(\A\to\rmin\B)\to(\B\to\rmin\A)$ from
commutativity.  This is a contraposition axiom.
\begin{align*}
  1.&&\fm{\B\rp\A}\1\0&\sep\fm{\A\rp\B}\1\0
  &&\text{commutativity}\\
  2.&&\fm\A\1\2&\sep\fm\A\1\2
  &&\text{axiom}\\
  3.&&\fm\A\1\2&\sep\fm{\min{\min{\conv\A}}}\2\1
  &&\text{2, $|\conv{}$, $\min\blank|$, $|\min\blank$}\\
  4.&&\fm\B\2\0&\sep\fm\B\2\0
  &&\text{axiom}\\
  5.&&\fm\B\2\0&\sep\fm{\conv\B}\0\2
  &&\text{4, $|\conv{}$}\\
  6.&&\fm\A\1\2,\fm\B\2\0 &\sep\fm{\conv\B\rp\min{\min{\conv\A}}}\0\1
  &&\text{3, 5, $|\rp$}\\
  7.&&\fm\A\1\2,\fm\B\2\0,
  \fm{\min{\conv\B\rp\min{\min{\conv\A}}}}\0\1
  &\sep&&\text{6, $\min\blank|$}\\
  8.&&\fm\A\1\2,\fm\B\2\0,\fm{\B\to\rmin\A}\0\1&\sep
  &&\text{7, \eqref{rmin}, \eqref{to}}\\
  9.&&\fm{\A\rp\B}\1\0,\fm{\B\to\rmin\A}\0\1&\sep
  &&\text{8, $\rp|$, no $\2$}\\
  10.&&\fm{\B\rp\A}\1\0,\fm{\B\to\rmin\A}\0\1&\sep
  &&\text{1,  9, Cut}\\
  11.&&\fm\B\1\2&\sep\fm\B\1\2
  &&\text{axiom}\\
  12.&&\fm{\min{\min{\conv\B}}}\2\1&\sep\fm\B\1\2
  &&\text{11, $\conv{}|$, $|\min\blank$, $\min\blank|$}\\
  13.&&\fm\A\2\0&\sep\fm\A\2\0
  &&\text{axiom}\\
  14.&&\fm{\conv\A}\0\2&\sep\fm\A\2\0
  &&\text{13, $\conv{}|$}\\
  15.&&\fm{\conv\A}\0\2,\fm{\min{\min{\conv\B}}}\2\1
  &\sep\fm{\B\rp\A}\1\2
  &&\text{12, 14, $|\rp$}\\
  16.&&\fm{\conv\A\rp\min{\min{\conv\B}}}\0\1 &\sep\fm{\B\rp\A}\1\0
  &&\text{15, $\rp|$, no $\2$}\\
  17.&&&\sep\fm{\B\rp\A}\1\0,
  \fm{\min{\conv\A\rp\min{\min{\conv\B}}}}\0\1
  &&\text{16, $|\min\blank$}\\
  18.&&&\sep\fm{\B\rp\A}\1\0,\fm{\A\to\rmin\B}\0\1
  &&\text{17, \eqref{rmin}, \eqref{to}}\\
  19.&&\fm{\B\to\rmin\A}\0\1&\sep\fm{\A\to\rmin\B}\0\1 
  &&\text{10, 18, Cut}
\end{align*}
\proof[\ref{perm}] 4-proof of $(\A\to(\B\to\C))\to(\B\to(\A\to\C))$
from commutativity.  This is the permutation axiom.  Let
$\D=\A\to(\B\to\C)$. Then $\D=\min{\conv\A\rp
  \min{\min{\conv\B\rp\min\C}}}$ by \eqref{rmin} and \eqref{to}.
\begin{align*}
  1.&&\fm\A\3\0&\sep\fm\A\3\0
  &&\text{axiom}\\
  2.&&\fm\A\3\0&\sep\fm{\conv\A}\0\3
  &&\text{1, $|\conv{}$}\\
  3.&&\fm\B\2\3&\sep\fm\B\2\3
  &&\text{axiom}\\
  4.&&\fm\B\2\3&\sep\fm{\conv\B}\3\2
  &&\text{3, $|\conv{}$}\\
  5.&&\fm{\min\C}\2\1&\sep\fm{\min\C}\2\1
  &&\text{axiom}\\
  6.&&\fm\B\2\3,\fm{\min\C}\2\1&\sep\fm{\conv\B\rp\min\C}\3\1
  &&\text{4, 5, $|\rp$}\\
  7.&&\fm\B\2\3,\fm{\min\C}\2\1&\sep\fm{\min{\min{\conv\B\rp\min\C}}}\3\1
  &&\text{6, $\min\blank|$, $|\min\blank$}\\
  8.&&\fm\B\2\3,\fm\A\3\0,\fm{\min\C}\2\1
  &\sep\fm{\conv\A\rp\min{\min{\conv\B\rp\min\C}}}\0\1
  &&\text{2, 7, $|\rp$}\\
  9.&&\fm\B\2\3,\fm\A\3\0,\fm{\min\C}\2\1&,\fm\D\0\1\sep
  &&\text{8, $\min\blank|$, def.\ $\D$}\\
  10.&&\fm{\B\rp\A}\2\0,\fm{\min\C}\2\1,\fm\D\0\1&\sep
  &&\text{9, $\rp|$, no $\3$}\\
  11.&&\fm{\A\rp\B}\2\0&\sep\fm{\B\rp\A}\2\0
  &&\text{commutativity}\\
  12.&&\fm{\A\rp\B}\2\0,\fm{\min\C}\2\1,\fm\D\0\1&\sep
  &&\text{10, 11, Cut}\\
  13.&&\fm\A\2\3&\sep\fm\A\2\3
  &&\text{axiom}\\
  14.&&\fm{\conv\A}\3\2&\sep\fm\A\2\3
  &&\text{13, $\conv{}|$}\\
  15.&&\fm\B\3\0&\sep\fm\B\3\0
  &&\text{axiom }\\
  16.&&\fm{\conv\B}\3\0&\sep\fm\B\3\0
  &&\text{15, $\conv{}|$}\\
  17.&&\fm{\conv\B}\0\3,\fm{\conv\A}\3\2&\sep\fm{\A\rp\B}\2\0
  &&\text{14, 16, $|\rp$}\\
  18.&&\fm{\conv\B}\0\3,\fm{\conv\A}\3\2,\fm{\min\C}\2\1&,\fm\D\0\1\sep
  &&\text{12, 17, Cut}\\
  19.&&\fm{\conv\B}\0\3,\fm{\conv\A\rp\min\C}\3\1&,\fm\D\0\1\sep
  &&\text{18, $\rp|$, no $\2$}\\
  20.&&\fm{\conv\B}\0\3,\fm{\min{\min{\conv\A\rp\min\C}}}\3\1&,
  \fm\D\0\1\sep&&\text{19, $|\min\blank$, $\min\blank|$}\\
  21.&&\fm{\conv\B\rp\min{\min{\conv\A\rp\min\C}}}\0\1&,
  \fm\D\0\1\sep&&\text{20, $\rp|$, no $\3$}\\
  22.&&\fm\D\0\1&\sep\fm{\B\to(\A\to\C)}\0\1
  &&\text{21, $|\min\blank$, \eqref{to}}\\
  23.&&\fm{\A\to(\B\to\C)}\0\1&\sep\fm{\B\to(\A\to\C)}\0\1 
  &&\text{22, def.\ $\D$}
\end{align*}
\proof[\ref{suff}] 4-provability of
$(\A\to\B)\to((\B\to\C)\to(\A\to\C))$ from commutativity.  The
suffixing axiom \eqref{suff} (4-provable from commutativity) can be
derived from the permutation axiom \eqref{perm} (4-provable from
commutativity) using the 1-provable derived rule of modus ponens
\eqref{modus-p}.

\proof[\ref{self-dist}] 4-proof of
$(\A\to(\B\to\C))\to((\A\to\B)\to(\A\to\C))$ from density and
commutativity.  Let
\begin{align*}
\D&=(\A\to\B)\to((\B\to(\A\to\C))\to(\A\to(\A\to\C))),\\
\E&=(\A\to\B)\to((\B\to(\A\to\C))\to(\A\to\C)),\\
\F&=(\B\to(\A\to\C))\to((\A\to\B)\to(\A\to\C)),\\
\G&=(\A\to(\B\to\C))\to((\A\to\B)\to(\A\to\C)).
\end{align*}
Contraction $(\A\to(\A\to\C))\to(\A\to\C)$ is 4-provable from density
by \eqref{contract2}.  Apply modus ponens \eqref{modus-p} (a
1-provable rule) twice to instances of 4-provable prefixing
\eqref{no.1} to conclude that $\D\to\E$ is 4-provable from density.
$\D$ is an instance of suffixing \eqref{suff}, which is 4-provable
from commutativity.  By modus ponens \eqref{modus-p}, $\E$ is
therefore 4-provable from density and commutativity.  $\E\to\F$ is an
instance of permutation \eqref{perm}, so by modus ponens
\eqref{modus-p}, $\F$ is 4-provable from density and commutativity.
Apply the 2-provable transitivity rule \eqref{tran-rule} to $\F$ and
permutation $(\A\to(\B\to\C))\to(\B\to(\A\to\C))$ (4-provable from
commutativity) to conclude that $\G$ is 4-provable from density and
commutativity, as desired.  This proof uses more instances of density
and commutativity than are required.  The following 4-proof from
density and commutativity shows that $\A\leq\A\rp\A$ and
$\A\rp\B\equals\B\rp\A$ are sufficient. By the way, $\A\leq\A\rp\A$
and $\A\rp(\A\to\B)\equals(\A\to\B)\rp\A$ are also sufficient by
\cite[Thm 5.1(62)]{MR2641636}.
\begin{align*}
  1.&&\fm\B\2\3&\sep\fm\B\2\3
  &&\text{axiom}\\
  2.&&\fm\B\2\3&\sep\fm{\conv\B}\3\2
  &&\text{1, $|\conv{}$}\\
  3.&&\fm{\min\C}\2\1&\sep\fm{\min\C}\2\1
  &&\text{axiom}\\
  4.&&\fm\B\2\3,\fm{\min\C}\2\1&\sep\fm{\conv\B\rp\min\C}\3\1
  &&\text{2, 3, $|\rp$}\\
  5.&&\fm\B\2\3,\fm{\min\C}\2\1
  &\sep\fm{\min{\min{\conv\B\rp\min\C}}}\3\1
  &&\text{4, $\min\blank|$, $|\min\blank$}\\
  6.&&\fm\A\3\0&\sep\fm\A\3\0
  &&\text{axiom}\\
  7.&&\fm\A\3\0&\sep\fm{\conv\A}\0\3
  &&\text{6, $|\conv{}$}\\
  8.&&\fm\B\2\3,\fm\A\3\0,\fm{\min\C}\2\1
  &\sep\fm{\conv\A\rp\min{\min{\conv\B\rp\min\C}}}\0\1
  &&\text{5, 7, $|\rp$}\\
  9.&&\fm{\B\rp\A}\2\0,\fm{\min\C}\2\1
  &\sep\fm{\conv\A\rp\min{\min{\conv\B\rp\min\C}}}\0\1
  &&\text{8, $\rp|$, no $\3$}\\
  10.&&\fm\A\2\3&\sep\fm\A\2\3
  &&\text{axiom}\\
  11.&&\fm\B\3\0&\sep\fm\B\3\0
  &&\text{axiom}\\
  12.&&\fm\A\2\3,\fm\B\3\0&\sep\fm{\A\rp\B}\2\0
  &&\text{10, 11, $|\rp$}\\
  13.&&\fm{\A\rp\B}\2\0&\sep\fm{\B\rp\A}\2\0
  &&\text{commutativity}\\
  14.&&\fm\A\2\3,\fm\B\3\0&\sep\fm{\B\rp\A}\2\0
  &&\text{12, 13, Cut}\\
  15.&&\fm\A\2\3,\fm\B\3\0,\fm{\min\C}\2\1
  &\sep\fm{\conv\A\rp\min{\min{\conv\B\rp\min\C}}}\0\1
  &&\text{9, 14, Cut}\\
  16.&&\fm\B\3\0,\fm{\conv\A}\3\2,\fm{\min\C}\2\1
  &\sep\fm{\conv\A\rp\min{\min{\conv\B\rp\min\C}}}\0\1
  &&\text{15, $\conv{}|$}\\
  17.&&\fm\B\3\0,\fm{\conv\A\rp\min\C}\3\1
  &\sep\fm{\conv\A\rp\min{\min{\conv\B\rp\min\C}}}\0\1
  &&\text{16, $\rp|$, no $\2$}\\
  18.&&\fm\A\3\2&\sep\fm\A\3\2
  &&\text{axiom}\\
  19.&&\fm\B\3\0&\sep\fm\B\3\0
  &&\text{axiom}\\
  20.&&\fm\A\3\2&\sep\fm{\conv\A}\2\3
  &&\text{18, $|\conv{}$}\\
  21.&&&\sep\fm\B\3\0,\fm{\min\B}\3\0
  &&\text{19, $|\min\blank$}\\
  22.&&\fm\A\3\2&\sep\fm\B\3\0,\fm{\conv\A\rp\min\B}\2\0
  &&\text{20, 21, $|\rp$}\\
  23.&&\fm\A\3\2,\fm{\min{\conv\A\rp\min\B}}\2\0&\sep\fm\B\3\0
  &&\text{22, $\min\blank|$}\\
  24.&&\fm{\conv{\min{\conv\A\rp\min\B}}}\0\2,
  \fm{\conv\A}\2\3&\sep\fm\B\3\0
  &&\text{23, $\conv{}|$}\\
  25.&&\fm{\conv{\min{\conv\A\rp\min\B}}}\0\2,\fm{\conv\A}\2\3,
  &\fm{\conv\A\rp\min\C}\3\1
  \sep\fm{\conv\A\rp\min{\min{\conv\B\rp\min\C}}}\0\1
  &&\text{17, 24, Cut}\\
  26.&&\fm{\conv{\min{\conv\A\rp\min\B}}}\0\2,
  \fm{\conv\A\rp&(\conv\A\rp\min\C)}\2\1
  \sep\fm{\conv\A\rp\min{\min{\conv\B\rp\min\C}}}\0\1
  &&\text{25, $\rp|$, no $\3$}\\
  27.&&\fm\A\0\3&\sep\fm\A\0\3
  &&\text{axiom}\\
  28.&&\fm\A\0\3&\sep\fm{\conv\A}\3\0
  &&\text{27, $|\conv{}$}\\
  29.&&\fm{\min\C}\0\1&\sep\fm{\min\C}\0\1
  &&\text{axiom}\\
  30.&&\fm\A\0\3,\fm{\min\C}\0\1&\sep\fm{\conv\A\rp\min\C}\3\1
  &&\text{28, 29, $|\rp$}\\
  31.&&\fm\A\3\2&\sep\fm\A\3\2
  &&\text{axiom}\\
  32.&&\fm\A\3\2&\sep\fm{\conv\A}\2\3
  &&\text{31, $|\conv{}$}\\
  33.&&\fm\A\0\3,\fm\A\3\2,\fm{\min\C}\0\1
  &\sep\fm{\conv\A\rp(\conv\A\rp\min\C)}\2\1
  &&\text{30, 32, $|\rp$}\\
  34.&&\fm{\A\rp\A}\0\2,\fm{\min\C}\0\1
  &\sep\fm{\conv\A\rp(\conv\A\rp\min\C)}\2\1
  &&\text{33, $\rp|$, no $\3$}\\
  35.&&\fm\A\0\2&\sep\fm{\A\rp\A}\0\2
  &&\text{density}\\
  36.&&\fm\A\0\2,\fm{\min\C}\0\1
  &\sep\fm{\conv\A\rp(\conv\A\rp\min\C)}\2\1
  &&\text{34, 35, Cut}\\
  37.&&\fm{\conv\A}\2\0,\fm{\min\C}\0\1
  &\sep\fm{\conv\A\rp(\conv\A\rp\min\C)}\2\1
  &&\text{36, $\conv{}|$}\\
  38.&&\fm{\conv\A\rp\min\C}\2\1
  &\sep\fm{\conv\A\rp(\conv\A\rp\min\C)}\2\1
  &&\text{37, $\rp|$, no $\0$}\\
  39.&&\fm{\conv{\min{\conv\A\rp\min\B}}}\0\2,
  \fm{\conv\A\rp\min\C}\2\1
  &\sep\fm{\conv\A\rp\min{\min{\conv\B\rp\min\C}}}\0\1
  &&\text{26, 38, Cut}\\
  40.&&\fm{\conv{\min{\conv\A\rp\min\B}}}\0\2,
  \fm{\min{\min{\conv\A\rp\min\C}}}\2\1
  &\sep\fm{\conv\A\rp\min{\min{\conv\B\rp\min\C}}}\0\1
  &&\text{39, $|\min\blank$, $\min\blank|$}\\
  41.&&\fm{\conv{\min{\conv\A\rp\min\B}}\rp
    \min{\min{\conv\A\rp\min\C}}}\0\1
  &\sep\fm{\conv\A\rp\min{\min{\conv\B\rp\min\C}}}\0\1
  &&\text{40, $\rp|$, no $\2$}\\
  42.&&\fm{\min{\conv\A\rp\min{\min{\conv\B\rp\min\C}}}}\0\1
  &\sep\fm{\min{\conv{\min{\conv\A\rp\min\B}}\rp
      \min{\min{\conv\A\rp\min\C}}}}\0\1
  &&\text{41, $\min\blank|$, $|\min\blank$}\\
  43.&&\fm{\A\to(\B\to\C)}\0\1 &\sep\fm{(\A\to\B)\to(\A\to\C)}\0\1
  &&\text{42, \eqref{to}}
\end{align*}
\proof[\ref{symm-ax}] 2-proof of $\A\land\rmin\A\to\B$ from symmetry.
 \begin{align*}
   1.&&\fm\A\0\1&\sep\fm{\conv\A}\0\1,\fm\B\0\1
   &&\text{symmetry}\\
   2.&&\fm{\A\cdot\min{\conv\A}}\0\1&\sep\fm\B\0\1
   &&\text{1, $\min\blank|$}\\
   3.&&\fm{\A\land\rmin\A}\0\1&\sep\fm\B\0\1
   &&\text{2, \eqref{and}, \eqref{rmin}}
\end{align*}
\proof[\ref{comm-ax}] 3-proof of $\A\circ\B\to\B\circ\A$ from
symmetry.
 \begin{align*}
   1.&&\fm\A\2\1&\sep\fm\A\2\1
   &&\text{axiom}\\
   2.&&\fm\A\2\1&\sep\fm{\conv\A}\1\2
   &&\text{1, $|\conv{}$}\\
   3.&&\fm{\conv\A}\1\2&\sep\fm\A\1\2
   &&\text{symmetry}\\
   4.&&\fm\A\2\1&\sep\fm\A\1\2
   &&\text{2, 3, Cut}\\
   5.&&\fm\B\0\2&\sep\fm\B\0\2
   &&\text{axiom}\\
   6.&&\fm\B\0\2&\sep\fm{\conv\B}\2\0
   &&\text{5, $|\conv{}$}\\
   7.&&\fm{\conv\B}\2\0&\sep\fm\B\2\0
   &&\text{symmetry}\\
   8.&&\fm\B\0\2&\sep\fm\B\2\0
   &&\text{6, 7, Cut}\\
   9.&&\fm\B\0\2,\fm\A\2\1&\sep\fm{\A\rp\B}\1\0
   &&\text{4, 8, $|\rp$}\\
   10.&&\fm\B\0\2,\fm\A\2\1&\sep\fm{\conv{\A\rp\B}}\0\1
   &&\text{9, $|\conv{}$}\\
   11.&&\fm{\B\rp\A}\0\1&\sep\fm{\conv{\A\rp\B}}\0\1
   &&\text{10, $\rp|$, no $\2$}\\
   12.&&\fm{\conv{\A\rp\B}}\0\1&\sep\fm{\A\rp\B}\0\1
   &&\text{symmetry}\\
   13.&&\fm{\B\rp\A}\0\1&\sep\fm{\A\rp\B}\0\1
   &&\text{11, 12, Cut}\\
   14.&&\fm{\A\circ\B}\0\1&\sep\fm{\B\circ\A}\0\1 
   &&\text{13, \eqref{circ}}
\end{align*}
\endproof
\section{4-provable predicates in $\RR$ and 
  $\TR$ that are not 3-provable}
\label{sect15}
This section presents predicates that are 4-provable but not
3-provable, even in the presence of the non-logical assumptions of
density, commutativity, and symmetry.
\begin{theorem}\label{thm12}
  Predicates \eqref{no.1}--\eqref{no.6}, \eqref{perm}--\eqref{suff},
  \eqref{contract2}--\eqref{self-dist} are in $\RR$ and $\TR$ but not
  $\CT_3$. They are not 3-provable from density and symmetry.
\end{theorem}
\proof To show the predicates \eqref{no.1}--\eqref{no.6},
\eqref{perm}--\eqref{suff}, \eqref{contract2}--\eqref{self-dist}
belong to $\RR$ it suffices to check that they are valid in every
$\CR$-frame.  This is well known and will not be done here.  To show
they cannot be proved with four variables, we use a dense symmetric
(hence commutative) semi-associative relation algebra that is not
associative.  Let $\gc\K_1 = \< \K, \R, \star{}, \{0\}\>$, where $\K =
\{0,\a,\b,\c\}$, $\R\subseteq\K^3$, $\<\x,\y,\z\>\in\R$ iff
$\z\in\{\x\}\rp\{\y\}$, and $\x^*=\x$ for all $\x,\y,\z\in\K$, and
$\rp$ is defined in Table \ref{K1}.
\begin{table}
\begin{equation*}
  \gc\K_1=
\begin{array}{|c|cccc|}\hline
  \rp&\{0\}&\{\a\}&\{\b\}&\{\c\}\\\hline
  \{0\}&\{0\}&\{\a\}&\{\b\}&\{\c\}\\
  \{\a\}& \{\a\}&\{0,\a\}&\{\c\}&\{\b\}\\
  \{\b\}& \{\b\}&\{\c\}&\{0,\b\}&\{\a\}\\
  \{\c\}& \{\c\}&\{\b\}&\{\a\}&\{0,\c\}\\\hline
\end{array}
\end{equation*}
\caption{$\gc{K}_1$ is a $\CR$-frame whose complex algebra is a
  dense symmetric semi-associative relation algebra that is not
  associative.}
\label{K1}
\end{table}
Then $\Cm{\gc\K_1}\cong\gc\E_4(\{1,3\})\in\SA$ by \cite[Thm
  2.5(4)(a)]{MR662049}. Obviously $\gc\K_1$ satisfies \eqref{symm} so
it also satisfies \eqref{comm} (see Lemma \ref{SA:symm->comm}).  It is
also easy to check directly from the table that $\gc\K_1$ satisfies
\eqref{dense} and \eqref{comm}. Therefore, $\Cm{\gc\K_1}$ is dense,
commutative, and symmetric by Theorem \ref{thm7-}. For each
$\i\in\{1,2,3,4,5,6\}$ suppose $\h_\i\colon \gc\P \to \Cm{\gc\K_1}$ is
a homomorphism such that
\begin{align*} 
  \h_1(\A)&=\{\a\},&\h_1(\B)&=\{\b\},&   \h_1(\C)&=\{0,\b\},\\
  \h_2(\A)&=\{\a\},&\h_2(\B)&=\{\a,\b\},&\h_2(\C)&=\{0,\a,\b\},\\
  \h_3(\A)&=\{\a\},&\h_3(\B)&=\{\b\},&   \h_3(\C)&=\{\a,\b\},\\
  \h_4(\A)&=\{\a\},&\h_4(\B)&=\{\b\},&   \h_4(\C)&=\{\a,\c\},\\
  \h_5(\A)&=\{\a\},&\h_5(\B)&=\{\b\},&   \h_5(\C)&=\{\b\},\\
  \h_6(\A)&=\{\a\},&\h_6(\B)&=\{\b,\c\},&\h_6(\C)&=\{\b\}.
\end{align*}
The predicates are invalidated by the homomophisms according to the
following table. The symbol $\times$ means that the predicate named in
the top row is invalidated by the homomorphism listed in the leftmost
column. The symbol $\holds$ indicates that the predicate is validated
by the homomorphism. These homomorphisms were calculated with
\cite{GAP4}.
\begin{align*} 
  \begin{array}{cccccccccccc}
    &\eqref{no.1}&\eqref{no.2}&\eqref{no.3}&\eqref{no.4} 
    &\eqref{no.5}&\eqref{no.6}&\eqref{perm}&\eqref{suff} 
    &\eqref{contract2}&\eqref{contract3}&\eqref{self-dist}\\
    \h_1&\holds&\times&\times&\times&\holds&\holds
    &\times&\times&\holds&\times&\times\\
    \h_2&\holds&\times&\times&\holds&\holds&\holds
    &\times&\holds&\times&\times&\times\\
    \h_3&\holds&\holds&\times&\holds&\holds&\times
    &\times&\holds&\holds&\times&\times\\
    \h_4&\holds&\holds&\holds&\times&\times&\times
    &\holds&\times&\holds&\times&\times\\
    \h_5&\times&\holds&\times&\times&\holds&\times
    &\holds&\holds&\holds&\times&\holds\\
    \h_6&\times&\holds&\holds&\times&\times&\times
    &\holds&\holds&\holds&\holds&\holds\\
\end{array}
\end{align*}
\endproof
\section{Predicates in $\RR$ and $\TR$ that are not
  $\omega$-provable from density}
\label{sect16} 
The predicates that rely on commutativity cannot be proved without
that assumption, even in the presence of density and infinitely many
variables.
\begin{theorem}\label{thm13}
  Predicates \eqref{mp}--\eqref{suff} and \eqref{self-dist} are in
  $\RR$ and $\TR$ but not in $\CT_\omega$. They are not
  $\omega$-provable from density.  They are invalid in a $\CR$-frame
  whose complex algebra is a non-commutative dense representable
  relation algebra.
\end{theorem}
\proof Let $\gc\K_2=\<\K,\R,\star{},\{0\}\>$ be the frame determined
by $\K=\{0,\a,\b,\b^*\}$, $0^*=0$, $\a^*=\a$, $(\b^*)^*=\b$,
$\R\subseteq\K^3$, and $\<\x,\y,\z\>\in\R$ iff $\z\in\{\x\}\rp\{\y\}$,
where $\rp$ is specified in Table \ref{K2}.  The complex algebra of
$\gc\K_2$ is relation algebra number \alg{13}{37}; see
\cite[p.\,437]{MR2628352}.
\begin{table}
\begin{align*}
  \gc\K_2&=
  \begin{array}{|c|cccc|}\hline
    \rp&\{0\}&\{\a\}&\{\b\}&\{\b^*\}\\\hline
    \{0\}&\{0\}&\{\a\}&\{\b\}&\{\b^*\}\\
    \{\a\}&\{\a\}&\{0,\a,\b,\b^*\}&\{\a,\b\}&\{\a\}\\
    \{\b\}&\{\b\}&\{\a\}&\{\b\}&\{0,\a,\b,\b^*\}\\
    \{\b^*\}&\{\b^*\}&\{\a,\b^*\}&\{0,\b,\b^*\}&\{\b^*\}\\
    \hline
  \end{array}
\end{align*}
\caption{$\Cm{\gc{K}_2}$ is a non-commutative dense representable 
  relation algebra.}
\label{K2}
\end{table}
$\gc\K_2$ satisfies the conditions \eqref{left reflection},
\eqref{right reflection}, \eqref{identity}, and \eqref{Pasch}, hence
$\Cm{\gc\K_2}\in\RA$ by Theorem \ref{thm7}\eqref{thm7ii}.  $\gc\K_2$
is also satisfies \eqref{dense} and $\Cm{\gc\K_2}$ is a dense relation
algebra (see \SS\ref{sect6}).  On the other hand, $\Cm{\gc\K_2}$ is
not commutative and \eqref{comm} fails in $\gc\K_2$.  The predicates
\eqref{mp}--\eqref{suff} are invalidated in many ways, but in rather
few ways if the invalidating homomorphisms are required to map the
predicates to the empty set and the propositional variables to
singleton subsets of $\K$.  All such homomorphisms have been found by
using \cite{GAP4} and are listed in Table \ref{list} for predicates
\eqref{mp}--\eqref{suff}. If a homomophism sends $\A$, $\B$, and $\C$
to the corrsponding subsets listed in some line in the column of a
given predicate, then it invalidates that predicate by mapping it to a
set not containing $0$.
\begin{table}
\begin{align*}
\begin{array}{|cc|cc|ccc|ccc|}\hline
  \eqref{mp}&&\eqref{contra}&&&\eqref{perm}&&&\eqref{suff}&\\\hline
  \A&\B&\A&\B&\A&\B&\C&\A&\B&\C\\\hline
  \{\a\}&\{\a\}&\{\a\}&\{\b\}&\{\a\}&\{\b\}&\{\a\}&\{0\}&\{\a\}&\{\a\}\\
  \{\b\}&\{\a\}&\{\b\}&\{\b^*\}&\{\b^*\}&\{\a\}&\{\a\}&\{0\}&\{\b\}&\{\a\}\\
  &&\{\b^*\}&\{\a\}&&&&\{\b\}&\{\a\}&\{\a\}\\
  &&&&&&&\{\b\}&\{\b\}&\{\a\}\\\hline
\end{array}
\end{align*}
\caption{Assignments of propositional variables to $\Cm{\gc{K}_2}$ 
  that invalidate \eqref{mp}--\eqref{suff}.}
\label{list}
\end{table}
\par
To show that \eqref{mp}--\eqref{suff} and \eqref{symm-ax} are not in
$\CT_\omega$ it suffices, by taking $\Psi$ to be the equations of
density in Theorem \ref{thm10}\eqref{thm10-a}, to show that
$\Cm{\gc\K_2}$ is isomorphic with a dense proper relation algebra. A
{\bf finite sequence} is a function $\f$ with domain $\text{dom}(\f) =
\{1,\cdots,\n\}$ for some finite non-zero $\n\in\omega$.  Let
$\rationals$ be the set of rational numbers.  Let $\U$ be the set of
finite sequences of rational numbers.  Define a binary relation
$\B\subseteq\U\times\U$ for $\f,\g\in\U$ by $\f\B\g$ (we say $\f$ is
below $\g$ or $\f$ comes before $\g$) iff for some finite $\n>0$,
$\text{dom}(\f) = \{1,\cdots,\n\} \subseteq \text{dom}(\g)$,
$\f_\i=\g_\i$ for all $\i<\n$, and $\f_\n<\g_\n$.  Let
\begin{align*}
  \sigma(0)&=\{\<\x,\x\>\colon\x\in\U\},\\
  \sigma(\b^*)&=\B,\\
  \sigma(\b)&=\conv\B,\\
  \sigma(\a)&=(\U\times\U)\setminus(\sigma(0)\cup\B\cup\conv\B),\\
  \rho(\X)&=\bigcup_{\x\in\X}\sigma(\x)\text{ for all
    $\X\subseteq\K$.}
\end{align*}
It can be checked that $\rho$ is an isomorphism of $\Cm{\gc\K_2}$ with
a dense proper relation algebra. By Theorem
\ref{thm10}\eqref{thm10-a}, every predicate in $\CT_\omega$ is valid
in $\gc\K_2$.  Since \eqref{mp}--\eqref{suff} and \eqref{symm-ax} are
not valid in $\gc\K_2$, they are not in $\CT_\omega$.  \endproof It
has been confirmed with \cite{GAP4} that self-distribution
\eqref{self-dist} is valid in $\gc\K_2$ but is not valid in the
5-element frame whose complex algebra is the dense non-commutative
relation algebra \alg{29}{83}; see \cite[p.\,448]{MR2628352} for its
multiplication table. If \alg{29}{83} is representable (which seems
extremely likely) then \eqref{self-dist} can be added to the list of
predicates that rely on commutativity and are not $\omega$-provable
from density alone.
\section{Predicates in $\RR$ and $\TR$ that 
  are not $\omega$-provable from 
  commutativity}\label{sect16a}
The predicates that are 3-provable or 4-provable from density require
that assumption and are not even $\omega$-provable from commutativity.
\begin{theorem}\label{thm-new}
  Predicates \eqref{reductio}--\eqref{self-dist} are in $\RR$ and
  $\TR$ but not $\CT_\omega$. They are not $\omega$-provable from
  commutativity. They are invalid in the $\CR$-frame of the 2-element
  group, whose complex algebra is a commutative representable relation
  algebra.
\end{theorem}
\proof The presence of \eqref{reductio}--\eqref{self-dist} in $\RR$ is
mentioned in many sources.
\begin{table}
\begin{equation*}
  \gc\K_3=
\begin{array}{|c|cc|}\hline
  \rp&\{0\}&\{\a\}\\\hline
  \{0\}&\{\a\}&\{0\}\\
  \{\a\}& \{0\}&\{\a\}\\\hline
\end{array}
\end{equation*}
\caption{$\gc{K}_3$ is the frame of the 2-element group.}
\label{K9}
\end{table}
Let $\gc\K_3$ be the $\CR$-frame of the 2-element group, shown in
Table \ref{K9}. Note that $\gc\K_3$ does not satisfy \thetag{dense}.
For each $\i\in\{1,2,3,4,5,6\}$ suppose $\h_\i\colon \gc\P \to
\Cm{\gc\K_3}$ is a homomorphism such that
\begin{align*} 
  \h_1(\A)&=\h_2(\A)=\h_3(\A)=\h_4(\A)=\h_5(\A)=\h_6(\A)=\{\a\},\\
  \h_2(\B)&=\h_5(\B)=\{\a\},\\
  \h_3(\B)&=\h_4(\B)=\h_6(\B)=\{0\},\\
  \h_5(\C)&=\h_6(\C)=\{0\}.
\end{align*}
Then $\g(\D)=\{\a\}$ whenever $\g=\h_\i$ for some $\i\in
\{1,2,3,4,5,6\}$ and $\D$ is any one of \eqref{reductio},
\eqref{contract5}, \eqref{contract4}, \eqref{contract2},
\eqref{contract3}, or \eqref{self-dist}. Since $0\notin\g(\D)$ this
shows the six predicates are invalid in the group frame. The complex
algebra of any group is in $\RRA$, an observation first made by
J.\ C.\ C.\ McKinsey; see \cite[Thm 5.10]{MR0045086}, \cite[Thm
  233]{MR2269199}. In this case the representation is quite
simple. Map $0$ to the identity relation on a 2-element set, and map
$\a$ to the transposition that interchanges the two elements.
\endproof
\section{3-provable predicates not in $\RR$ 
  or $\CR$}\label{sect17} The frames charactistic for classical
relevant logic $\CR$ need not satisfy any of the frame properties
\eqref{center reflection}--\eqref{right reflection}. This leads to the
problem, solved in this section, of determining a predicate in the
vocabulary of $\RR$ that corresponds to these properties.
\begin{theorem}\label{thm15}
  \begin{enumerate}
  \item\label{thm15i} Predicates
    \eqref{reflection1}--\eqref{reflection1a} are in $\TT_3$ but not
    $\RR$.  Predicate \eqref{dedekind} is in $\CT_3$ but not $\CR$.
  \item\label{thm15ii} \eqref{reflection1} is valid in a frame
    $\gc\K$ satisfying \eqref{identity} iff $\gc\K$ satisfies
    \eqref{left reflection}.
      \end{enumerate}
\end{theorem}
\proof Predicates \eqref{reflection1}--\eqref{dedekind} are 3-provable
by Theorem \ref{thm11}. The use of $\star{}$ in \eqref{dedekind} puts
it in $\CT_3$.
\begin{table}
  \begin{align*}
    \gc\K_4=\begin{array}{|c|ccc|}\hline
      \rp&\{0\}&\{\a\}&\{\star\a\}\\\hline
      \{0\}&\{0\}&\{\a\}&\{\star\a\}\\
      \{\a\}&\{\a\}&\{\a,\star\a\}&\{0,\a,\star\a\}\\
      \{\star\a\}&\{\star\a\}&\{0,\a,\star\a\}&\{\star\a\}\\
      \hline
    \end{array}
    \quad
    \begin{array}{|l|l|}\hline
      \x&\star\x\\\hline
      0&0\\
      \a&\star\a\\
      \star\a&\a\\\hline
    \end{array}
  \end{align*}
  \caption{A $\CR$-frame whose complex algebra is not a 
    semi-associative relation algebra.}  \label{K3}
\end{table}
It is straightforward to check that $\gc\K_4=\<\K,\R,{}^*,\{0\}\>$ in
Table \ref{K3} is a $\CR$-frame.  Therefore, everything in $\CR$ is
valid in $\gc\K_4$.  The frame conditions \eqref{center reflection},
\eqref{left reflection}, and \eqref{right reflection} all fail in
$\gc\K_4$ because $\<\a,\a,\star\a\>\in\R$ but $\<\star\a,\star\a,\a\>
\notin\R$. Axiom \eqref{BIX} fails because $\conv{\{\a\}\rp\{\a\}} =
\{\a,\star\a\}\conv{} = \{\a,\star\a\} \neq \{\star\a\} =
\{\star\a\}\rp\{\star\a\} = \{\a\}\conv{}\rp\{\a\}\conv{}$.  Thus,
$\Cm{\gc\K_4}\notin\SA$.  If $\h\colon\gc\P\to \Cm{\gc\K_4}$ is
homomorphism such that $\h(\A)=\h(\B)=\{\a\}$ and
$\h(\C)=\h(\D)=\{\star\a\}$, then $\h$ sends \eqref{reflection1} to
the empty set, hence \eqref{reflection1} is not valid in $\gc\K_4$ and
is not in $\CR$. The other predicates can be handled similarly.

For part \eqref{thm15ii}, assume $\gc\K$ is a frame satisfying
\eqref{identity}. We first prove that if $\gc\K$ also satisfies
\eqref{left reflection} then \eqref{reflection1} is valid in $\gc\K$.
Note that \eqref{reflection1} is a implication $\F\to\G$ where
$\F=\B\circ\A\land\C$ and $\G=\B\circ(\A\land\rmin\D) \lor
(\B\land\C\circ\D)\circ\A$.  By \eqref{or}, \eqref{and}, \eqref{rmin},
and \eqref{circ},
\begin{align*}
  \F&=\A\rp\B\cdot\C,\\
  \G&=(\A\cdot\min{\conv\D})\rp\B+\A\rp(\B\cdot\D\rp\C).
\end{align*}
By definition, $\F\to\G$ is valid in $\gc\K_4$ if the equation
$\id\leq\F\to\G$ is true in $\Cm{\gc\K_4}$, \ie,
$\h(\id)\subseteq\h(\F\to\G)$ for any homomorphism
$\h:\gc\P\to\Cm{\gc\K_4}$.  Any such homomorphism must send $\id$ to
$\{0\}$, so $\F\to\G$ is valid in $\gc\K_4$ iff $0\in\h(\F\to\G)$ for
every homomorphism $\h:\gc\P\to\Cm{\gc\K}$. By \eqref{to},
\eqref{min}, and the homomorphism properties of $\h$,
$0\in\h(\F\to\G)$ iff
\begin{align*}
  0\in\h\Big(\min{\conv\F\rp\min\G}\Big)
  =\K\setminus\Big(\conv{\h(\F)}\rp(\K\setminus\h(\G))\Big).
\end{align*}
Writing this out according to \eqref{rp}, we get
\begin{align*}
  \text{for all $\x,\y\in\K$, if }\x\in\conv{\h(\F)} \text{ and }
  \y\notin\h(\G) \text{ then }\<\x,\y,0\>\notin\R.
\end{align*}
By \eqref{conv} this is equivalent to
\begin{align*}
  \text{for all $\z,\y\in\K$, if }\z\in\h(\F) \text{ and }
  \<\star\z,\y,0\>\in\R\text{ then }\y\in\h(\G).
\end{align*}
By Lemma \ref{lem8}\eqref{lem8i}, $\<\star\z,\y,0\>\in\R$ iff
$\<\z,0,y\>\in\R$, which is equivalent by \eqref{identity} to $\z=\y$.
Therefore, $0\in\h(\F\to\G)$ iff
\begin{align}\label{valid}
  \text{for all $\z\in\K$, if }\z\in\h(\F)\text{ then }\z\in\h(\G),
  \text{\ie, }\h(\F)\subseteq\h(\G).
\end{align}
We will prove \eqref{valid}. Assume $\z\in\h(\F)$. Compute
\begin{equation*}
  \h(\F)=\h(\A\rp\B\cdot\C)=\h(\A)\rp\h(\B)\cap\h(\C),
\end{equation*}
so $\z\in\h(\C)$ and $\<\x,\y,\z\>\in\R$ for some $\x\in\h(\A)$ and
$\y\in\h(\B)$ by \eqref{rp}. There are two cases.  First assume
$\x\in\h(\min{\conv\D})$. Then
\begin{equation*}
  \x\in\h(\A)\cap\h(\min{\conv\D})=\h(\A\cdot\min{\conv\D}).
\end{equation*}
From this, $\<\x,\y,\z\>\in\R$, and $\y\in\h(\B)$ we get
$\z\in\h((\A\cdot\min{\conv\D})\rp\B)$ by \eqref{rp}.  By \eqref{cup},
$\h(\G)$ is the union of this last set with
$\G(\A\rp(\B\cdot\D\rp\C))$, so $\z\in\h(\G)$, as desired.

For the second case, assume $\x\notin\h(\min{\conv\D})$.  Then
$\x\in\K\setminus\conv{\h(\D)}$ by \eqref{min}, hence
$\x\in\conv{\h(\D)}$. By \eqref{conv}, $\x=\w^*$ for some
$\w\in\h(\D)$. Then $\<\w^*,\y,\z\>\in\R$ since $\<\x,\y,\z\>\in\R$,
hence $\<\w,\z,\y\>\in\R$ by \eqref{identity}, \eqref{left
  reflection}, and Lemma \ref{lem8}\eqref{lem8i}. Therefore,
$\y\in\h(\D\rp\C)$ by $\w\in\h(\D)$, $\z\in\h(\C)$, and \eqref{rp}.
From this and $\y\in\h(\B)$ we get $\y\in\h(\B\cdot\D\rp\C)$, hence by
$\x\in\h(\A)$ and $\<\x,\y,\z\>\in\R$, we have $\z\in\h(\A\rp(\B \cdot
\D\rp\C))$.  But $\h(\G)$ is the union of this set with
$\h((\A\cdot\min{\conv\D})\rp\B)$, so $\z\in\h(\G)$, as desired.

Assume we have a frame $\gc\K=\<\K,\R,\star{},\II\>$ satisfying
\eqref{identity} in which \eqref{left reflection} fails because there
are $\x,\y,\z\in\K$ such that $\<\x,\y,\z\>\in\R$ and
$\<\star\x,\z,\y\>\notin\R.$ We will show $\id\leq\F\to\G$ is not true
in $\Cm{\gc\K}$.  Suppose $\h\colon\gc\P\to \Cm{\gc\K}$ a homomorphism
such that
\begin{align*}
  \h(\A)&=\{\x\},&\h(\B)&=\{\y\},&\h(\C)&=\{\z\},&\h(\D)&=\{\star\x\}.
\end{align*}
Since $\<\x,\y,\z\>\in\R$ we have $\h(\C)=\{\z\} \subseteq
\{\x\}\rp\{\y\} = \h(\A)\rp\h(\B) = \h(\A\rp\B),$ so $\h(\F) =
\h(\A\rp\B\cdot\C) = \h(\A\rp\B)\cap\h(\C) = \h(\C) = \{\z\}.$ We also
have $\h(\A\cdot\min{\conv\D}) = \h(\A)\cap\h(\min{\conv\D}) =
\{\x\}\cap(\K\setminus\{\x\}) = \emptyset,$ hence, independently of
the value for $\B$,
\begin{align}\label{alpha} 
  \h((\A\cdot\min{\conv\D})\rp\B) &=\h(\A\cdot\min{\conv\D})\rp\h(\B)=
  \emptyset\rp\h(\B)=\emptyset.
\end{align}
Furthermore, $\y\notin\{\star\x\}\rp\{\z\}$ since
$\<\star\x,\z,\y\>\notin\R$, hence $\h(\D\rp\C) = \h(\D)\rp\h(\C) =
\{\star\x\}\rp\{\z\} \subseteq \K \setminus \{\y\}.$ Consequently,
$\h(\B\cdot\D\rp\C) =\h(\B) \cap \h(\D\rp\C) \subseteq \{\y\} \cap
\(\K\setminus\{\y\}\) =\emptyset.$ We therefore have
\begin{align}\label{beta} 
  \h(\A\rp(\B\cdot\D\rp\C))&= \h(\A)\rp\h(\B\cdot\D\rp\C)
  =\{\x\}\rp\emptyset=\emptyset.
\end{align}
From \eqref{alpha} and \eqref{beta} we have $\h(\G) =
\h((\A\cdot\min{\conv\D}) \rp\B) \cup \h(\A\rp(\B\cdot\D\rp\C)) =
\emptyset$. The inclusion $\{\z\}=\h(\F)\subseteq\h(\G)=\emptyset$ is
false, contradicting \eqref{valid}, so $\id\leq\F\to\G$ is not true in
$\Cm{\gc\K}$ and $\F\to\G$ is not valid in $\gc\K$.  The
contrapositive of what we have just proved is that if $\F\to\G$ is
valid in a frame $\gc\K$ satisfying \eqref{identity} then $\gc\K$
satisfies \eqref{left reflection}.
\endproof
\section{Counterexample to a theorem of Kowalski}\label{sect18} 
According to \cite[Thm 8.1]{MR3289545}, $\RR$ is ``complete with
respect to square-increasing [dense], commutative, integral relation
algebras.'' On the contrary, $\RR$ does \emph{not} contain all the
formulas true in this class of algebras. By Theorem \ref{thm15},
\eqref{reflection1} is not a theorem of $\RR$ but it is in $\TT_3$ and
is therefore true in all semi-associative relation algebras, including
all dense commutative integral relation algebras.  Thus,
\eqref{reflection1} is a counterexample to \cite[Thm 8.1]{MR3289545},
which was obtained as an immediate consequence of \cite[Thm
  7.1]{MR3289545}, that every normal De Morgan monoid is embeddable in
a dense commutative integral relation algebra. However, the complex
algebra of $\gc\K_4$ is a counterexample because $\gc\K_4$ fails to
satisfy \eqref{left reflection} and it would have to satisfy
\eqref{left reflection} if $\Cm{\gc\K_4}$ were embeddable in a
relation algebra.  The difficulty seems to arise in the proof of
\cite[Lemma~5.4(1)]{MR3289545}.
\section{5-provable predicates not in $\RR$ or $\CT_4$}\label{sect19}
This section presents two examples \eqref{L''} and \eqref{M''} of
predicates that are 5-provable but not 4-provable.  By creating
infinitely many such predicates, \citet{MR2496334} proved that
$\TT_\omega$ is not finitely axiomatizable.

The original form of \eqref{L''} was a logically valid sentence C2,
due to \citet[p.\,712]{MR37278}, that could not be derived from
Tarski's axioms for the calculus of relations.  Lyndon's sentence C2
was recast as the equation (L) \cite[p.\,30]{MR2269199} by
\citet[p.\,354]{MR0043763}.  The equation (L) expresses Desargues
Theorem when the propositional variables denote points in a projective
geometry and $\A\rp\B$ is the set of points on the line passing
through points $\A$ and $\B$.  Predicate \eqref{L''} was obtained from
the equation (L) by reformulating it as a predicate that uses only the
relevance logic operators $\land$, $\rmin$, $\to$, and $\circ$ in the
form $\rp$.  Predicate \eqref{L''} is a consequence of (L) and is not
necessarily equivalent to (L).  Predicate \eqref{M''} arises in the
same way from another sentence C3 that is also due to
\citet[p.\,712]{MR37278}.
\begin{theorem}\label{5var}{\rm\cite[Thms 8.1, 8.2]{MR2641636}}
  Predicates \eqref{L''} and \eqref{M''} are in $\TT_5$ but not in
  $\CT_4$.
\begin{gather}\label{L''}
  \A\rp\B\land\C\rp\D\land\E\rp\F\to{}\\\notag
  \big((\A\land\rmin\A)\rp\B\land\C\rp\D\land\E\rp\F\big)
  \lor\big(\A\rp\B\land\C\rp(\D\land\rmin\D)\land\E\rp\F\big)\\\notag
  {}\lor\big(\A\rp\B\land\C\rp\D\land(\E\land\rmin\E)\rp\F\big)
  \lor\big(\A\rp\B\land\C\rp\D\land\E\rp(\F\land\rmin\F)\big)\\\notag
  {}\lor\A\rp\big(\A\rp\C\land\B\rp\D\land(\A\rp\E\land\B\rp\F)
  \rp(\E\rp\C\land\F\rp\D)\big)\rp\D
\end{gather}
\begin{gather}\label{M''}
  \A\land(\B\land\C\rp\D)\rp(\E\land\F\rp\G)\to{}\\\notag
  \big(\A\land(\B\land(\C\land\rmin\C)\rp\D)\rp(\E\land\F\rp\G)\big)\\
  \notag {}\lor\big(\A\land(\B\land\C\rp\D)\rp(\E\land
  \F\rp(\G\land\rmin\G))\big)\\
  \notag {}\lor\C\rp\big((\C\rp\A\land\D\rp\E)\rp\G
  \land\D\rp\F\land\C\rp(\A\rp\G \land\B\rp\F)\big)\rp\G
\end{gather}
\end{theorem}
\proof \eqref{L''} and \eqref{M''} both have the form $\H\to\J$.  By
\cite[Thm 8.1]{MR2641636}, if $\A,\B,\C,\D,\E,\F,\G$ are binary
relations on $\U$ and the operations in $\H$ and $\J$ are interpreted
acccording to Table \ref{defs-ops} then $\H\to\J$ is a binary relation
that contain the identity relation on $\U$, or, equivalently,
$\H\subseteq\J$.  In both cases a straightforward set-theoretical
proof of this fact refers to five elements of $\U$.  Two elements are
assumed to be in the relation $\H$ and there are three more elements
corresponding to the occurrences of $\rp$ in the predicate $\H$.  The
proof consists of deducing facts expressed by $\J$ about the five
elements from the assumption that the five elements are related to
each other in ways described by $\H$.  The set-theoretical proofs can
be written up as 5-proofs in the sequent calculus that are more
elaborate but very similar to the 3-proof of \eqref{reflection1} or
the 4-proofs of \eqref{no.5} and \eqref{no.6}.  It follows that
\eqref{L''} and \eqref{M''} are in $\TT_5$ by Theorem \ref{two-obs-2}.

To show \eqref{L''} and \eqref{M''} are not in $\CT_4$ let $\gc\K_5 =
\< \K, \R, {}^*, \{\id\} \>$ where $\K=\{0,\a,\b,\c\}$, $\x^*=\x$ for
every $\x\in\K$, and $\R$ is determined in Table~\ref{K28}.  This
table appears as \cite[Table 6]{MR2641636}. Recall that $\<\x,\y,\z\>
\in\R$ iff $\z\in\{\x\}\rp\{\y\}$ where $\rp$ is defined by
\eqref{rp}.
\begin{table}
  \begin{equation*}
    \gc\K_5=
    \begin{array}{|c|cccc|}\hline
      \rp&\{\id\}&\{\a\}&\{\b\}&\{\c\}\\\hline
      \{\id\}&\{\id\}&\{\a\}&\{\b\}&\{\c\}\\
      \{\a\}&\{\a\}&\{\id,\a,\c\}&\{\b,\c\}&\{\a,\b\}\\
      \{\b\}&\{\b\}&\{\b,\c\}&\{\id,\a,\b\}&\{\a,\c\}\\
      \{\c\}&\{\c\}&\{\a,\b\}&\{\a,\c\}&\{\id,\b,\c\}\\\hline
    \end{array}
  \end{equation*}
  \caption{The smallest $\KR$-frame that invalidates 
    \eqref{L''} and \eqref{M''}.}
  \label{K28}
\end{table} 
Then $\gc\K_5$ is a $\KR$-frame in which both \eqref{L''} and
\eqref{M''} fail if $\h:\gc\P\to \Cm{\gc\K_5}$ is a homomorphism such
that $\{\a\}=\h(\A)=\h(\B)= \h(\E) = \h(\G)$, $\{\c\} = \h(\C) =
\h(\F)$, and $\{\b\}=\h(\D)$.  Such homomorphisms exist if we pick
$\id\neq\A,\B,\C,\D,\E,\F,\G\in\Pi$ such that $\C,\D,\F\notin
\{\A,\B,\E,\G\}$, and $\D\notin \{\C,\F\}$.
\endproof
The proof shows there are instances of \eqref{L''} and \eqref{M''}
with only three distinct propositional variables that fail to be in
$\CT_4$. To get them, let $\A = \B = \E = \G \neq \C = \F \neq \D$.

There are 14 $\KR$-frames with four elements.  $\gc\K_5$ is the unique
4-element $\KR$-frame that invalidates both \eqref{L''} and
\eqref{M''}. It has 28 triples in its ternary relation $\R$.  Its
complex algebra is the relation algebra \alg{42}{65} in
\cite{MR2269199}.  There are exactly two other 4-element $\KR$-frames
that invalidate \eqref{L''}. However, \eqref{M''} is valid in both of
them. Their complex algebras are the relation algebras \alg{36}{65}
and \alg{50}{65}.  Among the 390 $\KR$-frames with five elements, the
number of them that invalidate both \eqref{L''} and \eqref{M''} is
58. The smallest two, both with 41 triples, are \alg{118}{3013} and
\alg{200}{3013}.
\section{Counting characteristic $\TR$-frames and $\KR$-frames}
\label{sect20}
Alasdair Urquhart wrote \cite[p.\,349]{MR1223997},
\begin{quote}
  The list of small models was enormously extended by a computer
  search using some remarkable programs written by Slaney, Meyer,
  Pritchard, Abraham, and Thistlewaite \dots\ These programs churned
  out huge quantities of $\RR$ matrices and model structures of all
  shapes and sizes. Clearly, there are lots and lots of $\RR$ model
  structures out there! But what are they like? Can we classify them
  in some intelligent fashion? Are there general constructions that
  produce interesting examples? The answer to the first two questions
  is still obscure, though clearer than it was. The answer to the last
  question is an emphatic ``yes!''.
\end{quote}
He continued on \cite[p.\,350]{MR1223997},
\begin{quote}
  The first indication that $\KR$ is indeed nontrivial came from the
  computer, which churned out reams of interesting $\KR$ matrices. In
  retrospect, this is hardly surprising, because we now know that
  $\KR$ models can be manufactured {\it ad lib} from projective
  geometries.
\end{quote}
These thoughts were on his mind when I first met Alasdair Urquhart, at
a mathematical meeting in the late 1970s. At that time, I had
constructed many finite symmetric integral relation algebras by hand
and described to him how easy it was to create them. From a chance
encounter with \citet{MR0409114} while looking up the proof by
\citet{MR0392564} of Tarski's theorem that the proof of the
associative law \eqref{BIV} requires four variables, I knew that
finite relation algebras share several properties with relevant model
structures. I mentioned this and the ease with which they could be
constructed to Alasdair. He said that explained why so many relevant
model structures were being generated by his students, colleagues, and
computers, as documented by his remarks quoted above. My own
subsequent computer investigations showed that there are 14
$\KR$-frames with four elements and 390 $\KR$-frames with five
elements. (The number of ``$\RR$ matrices and model structures'' is
much larger.) These large numbers illustrate a pattern that continues
as the number of elements increases.

Alasdair pioneered and exploited the connection between projective
geometries and $\KR$-frames. This fruitful and crucial connection does
provide a general construction that produces interesting examples, but
does not account for the large number of $\KR$-frames. One can create
$\KR$-frames by a random process, as will be shown in this section.
Alasdair's astute judgement on his first two questions, ``But what are
they like? Can we classify them in some intelligent fashion?'', that
the answer is ``still obscure'' is confirmed here by Theorem
\ref{27}\eqref{27viii}. The number of $\KR$-frames on $\n$
elements grows like 2 raised to the power of a cubic polynomial.
$\KR$-frames are roughly as numerous as ternary relations and more
numerous than graphs, whose number rises only like 2 raised to the
power of a quadratic polynomial. In fact, a $\KR$-frame can be
constructed from an arbitrary graph. By contrast, projective planes
are known only for those orders that are powers of primes.

By Theorem \ref{thm7}\eqref{thm7iii} and Lemmas \ref{five} and
\ref{lem8}\eqref{lem8i}, the complex algebra of a frame is in $\NA$
iff it satisfies \eqref{left rotation}--\eqref{involution}. Therefore,
any frame satisfying those seven conditions is called an {\bf
  $\NA$-frame}.  An $\NA$-frame is {\bf associative} if it satisfies
\eqref{Pasch}, {\bf dense} if it satisfies \eqref{dense}, {\bf
  commutative} if it satisfies \eqref{comm}, and {\bf symmetric} if it
satisfies \eqref{symm}.  A {\bf $\TR$-frame} is an associative
dense commutative $\NA$-frame. By Theorem \ref{thm10-}\eqref{charTR}
and Lemmas \ref{five} and \ref{lem8}\eqref{lem8i}, the $\TR$-frames
are characteristic for $\TR$.  By Theorem \ref{thm8a}\eqref{thm8aii},
Theorem \ref{thm7}\eqref{thm7ii}, Theorem
\ref{thm7-}\eqref{thm7-ii}\eqref{thm7-iii}, Lemma
\ref{lem8}\eqref{lem8i}, and Lemma \ref{five}\eqref{five2}, a frame is
a $\KR$-frame iff it is an associative dense symmetric $\NA$-frame.
By the definition of $\KR$, the $\KR$-frames are characteristic for
$\KR$.

An asymptotic formula for the number of isomorphism types of
$\TR$-frames and $\KR$-frames can be computed by adapting the proof of
\cite[Thm 12]{MR823016}.  The first step is to count the number of
commutative $\NA$-frames on a given finite set $\K$ with fixed
involution $\star{}\colon\K\to\K$ and $\II=\{0\}$. After that one
observes that a randomly chosen dense commutative $\NA$-frame has a
probability approaching 1 (as the number of elements increases) of
being associative and having very few automorphisms. In fact, for any
fixed dimension $3\leq\d\in\omega$, the probability that every
$\d$-provable predicate is valid in a random dense commutative
$\NA$-frame also approaches 1.  In the following theorem, a function
$\N(\n,\s)$ is said to be {\bf asymptotic to} another function
$\M(\n,\s)$ if for every real number $\r>0$ there is some
$\n_0\in\omega$ such that if $\n_0<\n\in\omega$ and $1\leq\s\leq\n$
then $|1-\N(\n,\s)/\M(\n,\s)|<\r$.
\begin{theorem}\label{27}
Assume $\n=|\K|\in\omega$, $0\in\K$, $\II=\{0\}$,
$\star{}\colon\K\to\K$, $\star0=0$, $\star\x\star{}=\x$ for all
$\x\in\K$, and $\s=|\{\x:\star\x=\x\in\K\}|$. For every ternary
relation $\R\subseteq\K^3$ let $\gc\K(\R)=\<\K,\R,\star{},\II\>$.  Let
\begin{align*}
  \F(\n,\s)&=\textstyle\frac16{(\s-1)\s(\s+1)}
  +\frac1{12}{(\n-\s)(\n-\s+1)(\n-\s+2)}
  \\ &\,{}+\textstyle\frac14{(\s-1)(\n-\s)(\n+2)},
  \\ \G(\n,\s)&=\textstyle\frac16{(\s-1)(\s-2)(\s+3)}
  +\frac1{12}{(\n-\s)(\n-\s-1)(\n-\s+4)}
  \\ &\,{}+\textstyle\frac14{(\s-1)(\n-\s)(\n+2)},
  \\\P(\n,\s)&=(\s-1)!\(\textstyle\frac12(\n-\s)\)!2^{\frac12(\n-\s)}.
\end{align*}
\begin{enumerate}
\item\label{27i} The numbers of relations $\R\subseteq\K^3$ for which
  $\gc\K(\R)$ is a commutative, symmetric, dense commutative, or dense
  symmetric $\NA$-frame are
\begin{align*}
    2^{F(n,s)}&=|\{\R:\text{$\gc\K(\R)$ is a commutative
      $\NA$-frame}\}|,
    \\ 2^{\frac16{(n-1)n(n+1)}}&=|\{\R:\text{$\gc\K(\R)$ is a
      symmetric $\NA$-frame}\}|,
    \\ 2^{G(n,s)}&=|\{\R:\text{$\gc\K(\R)$ is a dense commutative
      $\NA$-frame}\}|,
    \\ 2^{\frac16{(n-1)(n-2)(n+3)}}&=|\{\R:\text{$\gc\K(\R)$ is a
      dense symmetric $\NA$-frame}\}|.
\end{align*}
\item\label{27ii} $\P(\n,\s)$ is the number of automorphisms of
  $\<\K,\star{},\{0\}\>$.
\item\label{27iii} The number of isomorphism types of commutative
  $\NA$-frames with $\n$ elements and $\s\leq\n-2$ symmetric elements is
  asymptotic  to
  \begin{equation*}
   \textstyle\frac1{\P(\n,\s)-1}\cdot2^{F(n,s)}.
   \end{equation*}
\item\label{27iv} The number of isomorphism types of symmetric
  $\NA$-frames with $\n$ elements is asymptotic to
  \begin{equation*}
   \textstyle\frac1{(\n-1)!}\cdot2^{\frac16{(n-1)n(n+1)}}.
   \end{equation*}
\item\label{27v} The number of isomorphism types of dense commutative
  $\NA$-frames with $\n$ elements and $\s\leq\n-2$ symmetric elements
  is asymptotic to
  \begin{equation*}
   \textstyle\frac1{\P(\n,\s)-1}\cdot2^{G(n,s)}.
   \end{equation*}
\item\label{27vi} The number of isomorphism types of dense symmetric
  $\NA$-frames with $\n$ elements is asymptotic to
  \begin{equation*}
   \textstyle\frac1{(\n-1)!}\cdot2^{\frac16{(n-1)(n-2)(n+3)}}.
   \end{equation*}
\item\label{27vi.1} For every dimension $3\leq\d<\omega$, the
  probability approaches $1$ as $\n\to\infty$ that the complex algebra
  of a randomly chosen commutative, dense commutative, symmetric, or
  dense symmetric $\NA$-frame with $\n$ elements is in $\RA_\d$.
\item\label{27vii} The number of isomorphism types of $\TR$-frames
  with $\n$ elements and $\s\leq\n-2$ symmetric elements is asymptotic
  to
  \begin{equation*}
   \textstyle\frac1{\P(\n,\s)-1}\cdot2^{G(n,s)}.
   \end{equation*}
\item\label{27viii} The number of isomorphism types of $\KR$-frames
  with $\n$ elements is asymptotic to
  \begin{equation*}
   \textstyle\frac1{(\n-1)!}\cdot2^{\frac16{(n-1)(n-2)(n+3)}}.
   \end{equation*}
\item\label{27ix} For every dimension $3\leq\d<\omega$, the number of
  isomorphism types of $\TR$-frames with $\n$ elements and
  $\s\leq\n-2$ symmetric elements whose complex algebras are in
  $\RA_\d$ and in which every $\d$-provable predicate is valid is
  asymptotic to
  \begin{equation*}
   \textstyle\frac1{\P(\n,\s)-1}\cdot2^{G(n,s)}.
   \end{equation*}
\item\label{27x} For every dimension $3\leq\d<\omega$, the number of
  isomorphism types of $\KR$-frames with $\n$ elements whose complex
  algebras are in $\RA_\d$ and in which every $\d$-provable predicate
  is valid is asymptotic to
  \begin{equation*}
   \textstyle\frac1{(\n-1)!}\cdot2^{\frac16{(n-1)(n-2)(n+3)}}.
   \end{equation*}
\end{enumerate}
\end{theorem}
\proof[\ref{27i}] We will count the number of commutative $\NA$-frames
$\gc\K(\R)=\<\K,\R,\star{},\II\>$ on a given finite set $\K$ with
fixed involution $\star{}\colon\K\to\K$ and $\II=\{0\}\subseteq\K$.
Consider an arbitrary $\R\subseteq\K^3$.  If $\gc\K(\R)$ satisfies
\eqref{identity} then $\R_0\subseteq\R$ where
\begin{align*}
\R_0 =\bigcup_{\x\in\K}\{&\<0,\x,\x\>, \<\x,\star\x,0\>,
\<\star\x,0,\star\x\>, \\&\<\x,0,\x\>,
\<\star\x,\x,0\>,\<0,\star\x,\star\x\>\}\cup\{ \<0,0,0\>\}.
\end{align*}
$\R$ cannot contain any other triples with $0$ in them, lest
\eqref{identity} be falsified. To get a commutative $\NA$-frame we
must consider only those relations $\R\subseteq\K^3$ that include
$\R_0$ and have no other triples in them that contain $0$.

For all $\x,\y,\z\in\K$ let $\C(\x,\y,\z)$ be the smallest set of
triples in $\K^3$ containing $\<\x,\y,\z\>$ such that
$\gc\K(\C(\x,\y,\z)\cup\R_0)$ is a commutative $\NA$-frame.  Such a
set is called a {\bf cycle}.  The isomorphism types of cycles are
listed Table \ref{isms}.
\begin{table}
  \begin{align*}
\R_0 =\bigcup_{\x\in\K}\{&\<0,\p,\p\>, \<\p,\star\p,0\>,
\<\star\p,0,\star\p\>, \<\p,0,\p\>,
\<\star\p,\p,0\>,\<0,\star\p,\star\p\>,\<0,0,0\>\}
\\ \C(\a,\a,\a)=\{&\<\a,\a,\a\>\}
\\ \C(\a,\b,\b)=\{&\<\a,\b,\b\>, \<\b,\a,\b\>, \<\b,\b,\a\>\}
\\ \C(\a,\b,\c)=\{&\<\a,\b,\c\>, \<\a,\c,\b\>, \<\b,\a,\c\>,
\<\b,\c,\a\>, \<\c,\a,\b\>, \<\c,\b,\a\>\}
\\ \C(\p,\a,\p)=\{&\<\a,\p,\p\>, \<\a,\star\p,\star\p\>, \<\p,\a,\p\>,
\<\p,\star\p,\a\>, \<\star\p,\a,\star\p\>, \<\star\p,\p,\a\>\}
\\ \C(\p,\p,\a)=\{&\<\a,\p,\star\p\>, \<\a,\star\p,\p\>,
\<\p,\a,\star\p\>, \<\p,\p,\a\>, \<\star\p,\a,\p\>,
\<\star\p,\star\p,\a\>\}
\\ \C(\p,\a,\a)=\{&\<\a,\a,\p\>, \<\a,\a,\star\p\>, \<\a,\p,\a\>,
\<\a,\star\p,\a\>, \<\p,\a,\a\>, \<\star\p,\a,\a\>\}
\\ \C(\p,\a,\b)=\{&\<\a,\b,\p\>, \<\a,\b,\star\p\>, \<\a,\p,\b\>,
\<\a,\star\p,\b\>, \<\b,\a,\p\>, \<\b,\a,\star\p\>, \\&\<\b,\p,\a\>,
\<\b,\star\p,\a\>, \<\p,\a,\b\>, \<\p,\b,\a\>, \<\star\p,\a,\b\>,
\<\star\p,\b,\a\>\}
\\ \C(\p,\q,\a)=\{&\<\a,\p,\star\q\>, \<\a,\star\p,\q\>,
\<\a,\q,\star\p\>, \<\a,\star\q,\p\>, \<\p,\a,\star\q\>, \<\p,\q,\a\>,
\\&\<\star\p,\a,\q\>, \<\star\p,\star\q,\a\>, \<\q,\a,\star\p\>,
\<\q,\p,\a\>, \<\star\q,\a,\p\>, \<\star\q,\star\p,\a\>\}
\\ \C(\p,\p,\p)=\{&\<\p,\p,\p\>, \<\p,\star\p,\p\>,
\<\p,\star\p,\star\p\>, \<\star\p,\p,\p\>, \<\star\p,\p,\star\p\>,
\<\star\p,\star\p,\star\p\>\}
\\ \C(\p,\p,\star\p)=\{&\<\p,\p,\star\p\>, \<\star\p,\star\p,\p\>\}
\\ \C(\p,\p,\q)=\{&\<\p,\p,\q\>, \<\p,\star\q,\star\p\>,
\<\star\p,\star\p,\star\q\>, \<\star\p,\q,\p\>, \<\q,\star\p,\p\>,
\<\star\q,\p,\star\p\>\}
\\ \C(\p,\q,\p)=\{&\<\p,\star\p,\q\>, \<\p,\star\p,\star\q\>,
\<\p,\q,\p\>, \<\p,\star\q,\p\>, \<\star\p,\p,\q\>,
\<\star\p,\p,\star\q\>,\\&\<\star\p,\q,\star\p\>,
\<\star\p,\star\q,\star\p\>, \<\q,\p,\p\>, \<\q,\star\p,\star\p\>,
\<\star\q,\p,\p\>, \<\star\q,\star\p,\star\p\>\}
\\ \C(\p,\q,\r)=\{&\<\p,\q,\r\>, \<\p,\star\r,\star\q\>,
\<\star\p,\star\q,\star\r\>, \<\star\p,\r,\q\>, \<\q,\p,\r\>,
\<\q,\star\r,\star\p\>, \\&\<\star\q,\star\p,\star\r\>,
\<\star\q,\r,\p\>, \<\r,\star\p,\q\>, \<\r,\star\q,\p\>,
\<\star\r,\p,\star\q\>, \<\star\r,\q,\star\p\>\}
\end{align*}
  \caption{Cycles $C(\text{-},\text{-},\text{-})\subseteq K^3$ with
    $p,q,r,a,b,c\in K$, $\star r\neq r$, $\star q\neq q$, $\star r\neq
    r$, $\star a=a$, $\star b=b$, $\star c=c$, such that $\gc
    K(C(\text{-},\text{-},\text{-})\cup R_0)$ is a commutative
    $\NA$-frame.}
  \label{isms}
\end{table}
For each isomorphism type of cycle, the number of triples in it and
the number of cycles of that type are listed in Table \ref{choices}.
\begin{table}
  \begin{align*}
  \begin{array}{|r|c|c|c|c|l|}\hline
    \text{No.}&\text{ Relation type }&\text{ Size } &\text{ Comm?
    }&\text{and dense?  } &\text{ The number of relations to choose
      from } \\\hline 1&\C(\a,\a,\a)&\ 1&\text{yes}&\text{no}&\,\,\s-1
    \\\hline 2&\C(\a,\b,\b)&\ 3&\text{yes}&\text{yes}&(\s-1)(\s-2)
    \\3&\C(\a,\b,\c)&\ 6&\text{yes}&\text{yes}&(\s-1)(\s-2)(\s-3)/6
    \\\hline 4&\C(\p,\a,\p)&\ 6&\text{yes}&\text{yes}&(\s-1)(\n-\s)/2
    \\5&\C(\p,\p,\a)&\ 6&\text{yes}&\text{yes}&(\s-1)(\n-\s)/2
    \\6&\C(\p,\a,\a)&\ 6&\text{yes}&\text{yes}&(\s-1)(\n-\s)/2
    \\7&\C(\p,\a,\b)&12&\text{yes}&\text{yes}&(\s-1)(\s-2)(\n-\s)/4
    \\8&\C(\p,\q,\a)&12&\text{yes}&\text{yes}&(\s-1)(\n-\s)(\n-\s-2)/4
    \\\hline 9&\C(\p,\p,\p)&\ 6&\text{yes}&\text{no}&(\n-\s)/2
    \\\hline10&\C(\p,\p,\star\p)&\ 2&\text{yes}&\text{yes}&(\n-\s)/2
    \\11&\C(\p,\p,\q)&\ 6&\text{yes}&\text{yes}&(\n-\s)(\n-\s-2)/2
    \\12&\C(\p,\q,\p)&12&\text{yes}&\text{yes}&(\n-\s)(\n-\s-2)/4
    \\13&\C(\p,\q,\r)&12&\text{yes}&\text{yes}&(\n-\s)(\n-\s-2)(\n-\s-4)/12
    \\\hline
  \end{array}
  \end{align*}
  \caption{The sizes and numbers of cycles
    $C(\text{-},\text{-},\text{-})$ on an $n$-element set with $s$
    symmetric elements, where $\star p\neq p$, $\star q\neq q$, $\star
    r\neq r$, $\star a=a$, $\star b=b$, $\star c=c$.\newline Types
    4--13 disappear in the symmetric case $s=n$ and type 3
    predominates.\newline Types 1--8 disappear in the non-symmetric
    case $s=1$ and type 13 predominates.}
  \label{choices}
\end{table}\par
If $\gc\K(\R\cup\R_0)$ is a commutative $\NA$-frame then $\R$ must be
the union of cycles.  Therefore, to create such a relation one must
choose to include or exclude each cycle. The number of choices
available for each isomorphism type of cycle occurs in the rightmost
column of Table \ref{choices}. For example, there are $\s-1$ cycles of
the form $\C(\a,\a,\a)$ with $\star\a=\a$ from which to choose. Since
$0$ is a symmetric element there are only $\s-1$ other symmetric
elements that create cycles of the form $\{\<\a,\a,\a\>\}$ as
available choices.

To obtain a commutative $\NA$-frame, the total number of available
choices is obtained by adding all the numbers in the rightmost column
that are labelled ``yes'' in the column headed with ``Comm?'', \ie, do
you want a commutative $\NA$-frame? If so, add the number of choices
in the rightmost column. If you want the frame to also be dense, then
add or do not add the number in the rightmost column according to the
entry in the column headed ``and dense?". Two of those entries are
``no'' since all triples in the cycles $\C(\a,\a,\a)$ and
$\C(\x,\x,\x)$ will necessarily be included in the desired $\R$ and
are therefore not available as choices for inclusion or exclusion from
$\R$.

$\F(\n,\s)$ is the sum of all the numbers in the last column of Table
\ref{choices}. This accounts for the first equation in part
\eqref{27i}.  $\G(\n,\s)$ is the sum of all the numbers in the last
column of Table \ref{choices} that occur in rows whose entry under
``and dense?'' is ``no''. This accounts for the third equation in part
\eqref{27i}.  Note that $\s=\n$ holds iff the resulting frame is
symmetric. The second and fourth equations are therefore obtained by
setting $\s=\n$ in the first and third equations.  This completes the
proof of part \eqref{27i}.

\proof[\ref{27ii}] Let $\perm$ be the set of permutations of $\K$.
Let $\iota$ be the identity permutation on $\K$, \ie, $\iota(\x)=\x$
for every $\x\in\K$.  Let $\aut$ be the set of automorphisms of
$\<\K,\star{},\II\>$.
\begin{align*}
  \aut=\{\sigma:\sigma\in\perm,\sigma(0)=0,
  \star{\sigma(\x)}=\sigma(\star\x)\text{ for every $\x\in\K$}\}.
\end{align*}
To see that $|\aut|=\P(\n,\s)$, note first that the symmetric elements
distinct from $0$ can be arbitrarily permuted, which accounts for the
term $(\s-1)!$.  The pairs of the form $\<\x,\star\x\>$ can also be
arbitrarily permuted, accounting for the term $\big(\frac12(\n -
\s)\big)!$. When sending the pair $\<\x,\star\x\>$ to the pair
$\<\y,\star\y\>$, an automorphism can either send $\x$ to $\y$ and
$\star\x$ to $\star\y$, or the other way around. For each pair there
are two ways to send it to its target pair, so the total number of
ways of doing this is $2^{\frac12(\n-\s)}$. $\P(\n,\s)$ is the product
of these three terms.  \proof[\ref{27iii}--\ref{27vi}] Partition
$\aut$ into three sets, $\autb$, $\autf$, and $\{\iota,\star{}\}$ by
setting
\begin{align*}
  \autb&=\{\sigma:\sigma\in\aut, \x\neq\sigma(\x)\neq\star\x \text{
    for some $\x\in\K$}\},
  \\ \autf&=\aut\setminus\big(\autb\cup\{\iota,\star{}\}\big).
\end{align*}
The automorphisms in $\{\iota,\star{}\}$ are called {\bf trivial} and
the ones in $\autb\cup\autf$ are called {\bf non-trivial}.  Recall
that $\n=|\K|$ and $\s=|\{\x:\star\x=\x\in\K\}|$.  Then $1\leq\s$
because $\star0=0\in\K$, and $\gc\K(\R\cup\R_0)$ is symmetric iff
$\s=\n$. If $\s=\n$ then $\star{}=\iota$, $\aut=\perm$, $\autb=\aut
\setminus \{\star{}\}$, $\autf=\emptyset$, $|\autb| = \P(\n,\s)-1$,
and $|\autf|=0$.

Suppose $\s<\n$ and $\sigma\in\autf$.  Since $\sigma\notin\autb$,
there is no $\x\in\K$ such that $\x\neq\sigma(\x)\neq\star\x$.
Therefore, if $\sigma(\x)\neq\x$ then $\sigma(\x)=\star\x$ and
$\sigma(\star\x)=\star{(\sigma(\x))}=\star\x\star{}=\x$.  This shows
the pair $\{\x,\star\x\}$ is switched by $\sigma$ whenever $\sigma$
moves $\x$, hence every $\sigma\in\autf$ is obtained by composing at
least one (since $\sigma\neq\iota$) but not all (since
$\sigma\neq\star{}$) transpositions of the form $(\x,\star\x)$ in
cycle notation.  The number of pairs $\{\x,\star\x\}$ is
$\frac12(\n-\s)$, so $|\autb|=\P(\n,\s)-2^{\frac12(\n-\s)}$ and
$|\autf|=2^{\frac12(\n-\s)}-2$. Let
\begin{equation*}
  \W=\{\R:\text{$\gc\K(\R\cup\R_0)$ is a commutative $\NA$-frame}\}.
\end{equation*}
For any property $\varphi$, let $\prop\varphi$ be the probability that
$\R\in\W$ has property $\varphi$,
\begin{align*}
  \prop\varphi=2^{-\F(\n,\s)}|\{\R:\text{ $\R\in\W$, $\R$ has property
    $\varphi$}\}|.
\end{align*}
We will compute the probability that $\gc\K(\R\cup\R_0) =
\<\K,\R\cup\R_0,\star{},\II\>$ has a non-trivial automorphism.  For
every ternary relation $\R\subseteq\K^3$ let $\autr$ be the set of
non-trivial automorphisms of $\gc\K(\R\cup\R_0)$,
\begin{align*}
  \autr=\{\sigma:\sigma\in\autb\cup\autf,
  \R=\{\<\sigma(\x),\sigma(\y), \sigma(\z)\>:\<\x,\y,\z\>\in\R \}\}.
\end{align*}
For every non-trivial $\sigma\in\autb\cup\autf$ we will compute an
upper bound on the number of ternary relations $\R\in\W$ such that
$\sigma\in\autr$ and multiply by $|\autb|$ or $|\autf|$ to obtain an
upper bound on the number of relations in $\W$ that have a non-trivial
automorphism.  With this upper bound we can show the probability of
having a non-trivial automorphism approaches $0$ as $\n\to\infty$.
Suppose $\sigma\in\autb$. Then for some $\x\in\K$,
$0\neq\x\neq\sigma(\x)\neq \star\x$. Let
\begin{equation*}
  \X=\{0,\x, \star\x, \sigma(\x), \sigma(\star\x) \}.
\end{equation*}
Since $|\K\setminus\X|=3$ if $\x=\star\x$ and $|\K\setminus\X|=5$ if
$\x\neq\star\x$, the number ways to choose one or two elements from
$\K\setminus\X$ is at least $(\n-5)+\frac12(\n-5)(\n-6) =
\frac12(\n-4)(\n-5)$.  For each such choice $\y,\z\in\K\setminus\X$ we
also have
\begin{equation}\label{notclosed}
  \sigma(\x)\notin\Y = \{\x, \star\x, \y, \star\y, \z, \star\z\}.
\end{equation}
To see this, note that we have assumed $\sigma(\x) \neq \x,
\star\x$. Both $\sigma(\x)=\y$ and $\sigma(\x)=\z$ violate the choice
that $\y,\z\notin\X$.  Since $\sigma$ is an automorphism, $\sigma(\x)
= \star\y$ implies $\sigma{(\star\x)} = \star{\sigma(\x)} =
\star\y\star{}=\y$, contradicting $\y\notin\X$, and similarly
$\sigma(\x)=\star\z$ contradicts $\z\notin\X$.  Suppose $\C(\x,\y,\z)$
is fixed by $\sigma$, that is, $\C(\x,\y,\z) = \C(\sigma(\x),
\sigma(\y), \sigma(\z))$.  Then $\sigma$ must map $\Y$ onto itself
because $\C(\x,\y,\z) \subseteq \Y^3$ and every element of $\Y$ is in
some triple in $\C(\x,\y,\z)$, as can be easily seen in Table
\ref{isms}. But this contradicts \eqref{notclosed}.  Therefore,
$\C(\x,\y,\z)$ is moved by $\sigma$. This proves that
\begin{equation}\label{zx2215}
\text{if $\sigma\in\autb$ then at least
  $\textstyle\frac12(\n-4)(\n-5)$ cycles are moved by $\sigma$.}
\end{equation}
Suppose $\sigma\in\autf$. This is possible only if $\s\leq\n-2$.
The number of cycles moved by $\sigma$ is
\begin{align*}
  \f(\n,\s,\m)&=\textstyle\frac14\m(\n-\s-\m)(\n+\s-2),
  \text{ where $\m=|\{\x:\sigma(\x)\neq\x\}|$.}
\end{align*}
This computation was checked with GAP \cite{GAP4}. Note that $\m$ is
always an even number, since moved elements come in pairs of the form
$\{\x,\star\x\}$. Also, $2\leq\m<\n-\s$ because $\sigma$ moves
something but differs from $\star{}$.  Under these constraints the
smallest non-zero value occurs when $\m=2$ or $\m=\n-\s-2$ and is
$\f(\n,\s,2) = \f(\n,\s,\n-\s-2) = \frac12((\n-2)^2-\s^2)$.
This proves that
\begin{equation}\label{zx2215a}
\text{if $\sigma\in\autf$ then at least $\textstyle\frac12
  ((\n-2)^2-\s^2)$ cycles are moved by $\sigma$.}
\end{equation}
For every $\sigma\in\autb\cup\autf$, let $\M_\sigma$ be the number of
cycles moved by $\sigma$.  To make a relation $\R\in\W$ with
$\sigma\in\autr$, one can freely choose to include or exclude each
unmoved cycle in $\R$.  The cycles in each orbit under $\sigma$ must
be either all included in $\R$ or all excluded from $\R$.  If the
moved cycles form a single orbit then one can only include or exclude
the entire orbit, which has size $\M_\sigma$. No cycles in that orbit
are available as choices to include or exclude, so the number of
available choices is $\F(\n,\s)-\M_\sigma+1$. The number of orbits can
vary from one orbit of size $\M_\sigma$ up to $\M_\sigma/2$ orbits of
size 2.  The number of unmoved cycles is $\F(\n,\s) - \M_\sigma$. They
are all free to be included in $\R$ or excluded. The moved elements
offer somewhere between one orbit and $\M_\sigma/2$ orbits as choices
for inclusion.  The number of choices ranges from at least
$\F(\n,\s)-\M_\sigma+1$ up to at most $\F(\n,\s) - \M_\sigma +
\M_\sigma/2 = \F(\n,\s)-\M_\sigma/2$.  From \eqref{zx2215} and
\eqref{zx2215a} we know $( (\n-2)^2-\s^2 )/4 \leq \M_\sigma/2$ if
$\sigma\in\autf$ and $(\n-4)(\n-5)/4 \leq \M_\sigma/2$ if
$\sigma\in\autb$, so the number of choices is at most $\F(\n,\s)- (
(\n-2)^2-\s^2 )/4$ if $\sigma\in\autf$ and $\F(\n,\s) -
(\n-4)(\n-5)/4$ if $\sigma\in\autb$.  Thus,
\begin{align*}
  |\{\R:\sigma\in\autr\}|&<
  \begin{cases}2^{\F(\n,\s)-\frac14(\n-4)(\n-5)} &\text{ if $\sigma\in\autb$,}
    \\ 2^{\F(\n,\s)-\frac12 ((\n-2)^2-\s^2 )} &\text{ if
      $\sigma\in\autf$.}
  \end{cases}
\end{align*}
Since $|\autb|<\P(\n,\s)$ and $|\autf|<2^{\frac12(\n-\s)}$, there are
fewer than
\begin{equation*}
  \P(\n,\s)2^{\F(\n,\s)-\frac14(\n-4)(\n-5)}
\end{equation*}
relations $\R\in\W$ such that $\autr$ is not empty, and if
$\s\leq\n-2$ there are fewer than
\begin{equation*}
  2^{\frac12(\n-\s)}2^{\F(\n,\s)-\frac12((\n-2)^2-\s^2)}
\end{equation*}
relations $\R\in\W$ such that $\autr$ is not empty. When $\s\leq\n-2$,
our over-estimate of the fraction of $\W$ that has a non-trivial
automorphism in $\autr$ is obtained by adding these two numbers and
dividing by $\F(\n,\s)$.
\begin{equation*}
  \E(\n,\s)=\frac{(\s-1)!\(\frac12(\n-\s)\)!2^{\frac12(\n-\s)}}
    {2^{\frac14(\n-4)(\n-5)}} +\frac{
      2^{\frac12(\n-\s)}}{2^{\frac12((\n-2)^2-\s^2)}}.
\end{equation*}
If $\s=\n$ then $\E(\n,\s)$ is just the first term of this sum.  In
either case, a straightforward analysis of the growth rates for the
numerators and denominators shows that $\lim_{\n\to\infty}\E(\n,\s)
=0$.  Since $\E(\n,\s)$ is the probability that $\gc\K(\R\cup\R_0)$
has a non-trivial automorphism, the probability that
$\gc\K(\R\cup\R_0)$ has no non-trivial automorphisms approaches $1$ as
$\n\to\infty$. Thus, a randomly selected $\R\in\W$ will almost
certainly show up in $\P(\n,\s)-1$ ways in $\W$ if $\s\leq\n-2$ and
$\P(\n,\s)$ ways in $\W$ if $\s=\n$. To estimate the number of
isomorphism types of commutative $\NA$-frames $\gc\K(\R\cup\R_0)$ we
must therefore divide the total number $2^{\F(\n,\s)}$ of such frames
by $\P(\n,\s)-1$ if $\s\leq\n-2$ and $\P(\n,\n)$ if $\s=\n$. This
gives us the approximations in parts \eqref{27iii} and \eqref{27iv}.
Parts \eqref{27v} and \eqref{27vi} are proved in the same way, using
$\G(\n,\s)$ instead of $\F(\n,\s)$ and redefining $\W$ as the set of
$\R\subseteq\K^3$ such that $\gc\K(\R\cup\R_0)$ is a dense commutative
$\NA$-frame. The reasoning applies to both definitions of $\W$.

\proof[\ref{27vi.1}--\ref{27x}] For part \eqref{27vi.1}, we will prove
for any fixed dimension $3\leq\d<\omega$ that the probability
approaches $1$ as $\n\to\infty$ that a randomly chosen $\R\in\W$ has a
commutative $\NA$-frame $\gc\K(\R\cup\R_0)$ whose complex algebra is
in $\RA_\d$ (because it almost certainly has a much stronger property)
and hence, by Theorem \ref{key}\eqref{key4}, every $\d$-provable
predicate is valid in $\gc\K(\R\cup\R_0)$.  Part \eqref{27vi.1}
together with parts \eqref{27v} and \eqref{27vi} imply parts
\eqref{27ix} and \eqref{27x}.  Parts \eqref{27vii} and \eqref{27viii}
follow from parts \eqref{27ix} and \eqref{27x} by Theorem
\ref{key}\eqref{key3}. Indeed, associativity is equivalent to the
validity of the 4-provable predicates \eqref{no.5} and \eqref{no.6}
and is obtained with near certainty by taking $\d=4$.

For every integer $\t\geq1$, define the {\bf diamond property}
\dee\t\ of relations $\R\in\W$ by
\begin{align}\label{deet}
  \dee\t:\ \text{if }&0\neq\x_1,\dots, x_\t,\y_1, \dots, \y_\t\in\K
  \text{ then for some }0\neq\z\in\K, \\\notag
  &\bigcup_{1\leq\i\leq\t}\C(\x_\i,\y_\i,\z)\subseteq\R.
\end{align}
We will show for a fixed $1\leq\t\in\omega$ that $\dee\t$ almost
certainly holds as $\n\to\infty$.  Consider one instance of $\dee\t$,
say $0\neq\x_1,\dots,\x_t, \y_1,\dots,\y_\t \in\K$.  We want to show
there is likely to be some $\z$ that works, where
\begin{equation*}
  \text{{\bf $\z$ works} iff }
  \bigcup_{1\leq\i\leq\t}\C(\x_\i,\y_\i,\z)\subseteq\R.
\end{equation*}
Note that $\z$ works only if all the cycles $\C(\x_\i,\y_\i,\z)$ are
contained in $\R$. Each cycle is contained with probability $1/2$. If
the cycles are disjoint their inclusions in $\R$ are independent
events and the probability that $\z$ works is exactly $2^{-\t}$, but
otherwise it is more, so $\z$ works with probability at least
$2^{-\t}$ and the probability that $\z$ does not work is at most
$1-2^{-\t}$. This is almost certain if $\t$ is large, but there are
more and more $\z$'s as $\n$ increases.  One of them is bound to work
and only one is needed.  We will calculate the probability that none
of them work and see that it goes to zero.  Let
\begin{equation*}
\K^-= \K \setminus \big(\{0\} \cup \bigcup_{1\leq\i\leq\t} \{\x_\i,
\star\x_\i, \y_\i, \star\y_\i\}\big).
\end{equation*}
Then $\K^-$ is closed under $\star{}$ so we can partition $\K^-$ into
pairs $\{\z,\star\z\}$ with $\z\neq\star\z$ and singletons $\{\z\}$
for those $\z$ such that $\z=\star\z$. Let $\z_\i$, $1\leq\i\leq\q$,
be a selection of one element from each pair or singleton.  There is
at least one and at most $4\t$ elements in $\bigcup_{1\leq\i\leq\t}
\{\x_\i, \star\x_\i, \y_\i, \star\y_\i\}$, so $\n-4\t-1 \leq |\K^-|
\leq\n-1$.  The largest partition of $\K^-$ occurs when all its
elements are symmetric and the smallest when all its elements are
non-symmetric.  Together these constraints put bounds on $\q$, namely
$(\n-4\t-1)/2\leq\q\leq\n-1$.

To see that distinct $\z_\i$'s create independent events, suppose
$1\leq\j<\k\leq\q$. The three sets $\{\z_\j,\star\z_\j\}$,
$\{\z_\k,\star\z_\k\}$, and $\bigcup_{1\leq\i\leq\t} \{ \x_\i,
\star\x_\i, \y_\i, \star\y_\i\}$ are disjoint, so
$\bigcup_{1\leq\i\leq\t} \C(\x_\i,\y_i,\z_\j)$ and
$\bigcup_{1\leq\i\leq\t} \C(\x_\i,\y_i,\z_\k)$ are also disjoint
because every triple in the former set contains $\z_\j$ or
$\star\z_\j$ but no triple in the latter set does so.  Thus, the
events that $\z_\j$ works and $\z_\k$ works are independent, as are
their complements. The probability that $\z_\j$ doesn't work for every
$\j\leq\q$ is the product of the probabilities that each $\z_\j$
doesn't work.  The probability that $\z_\j$ doesn't work is at most
$1-2^{-\t}$, so the probability that none of them works is at most the
product of $\q$ copies of $1-2^{-\t}$, one for each $\z$. The number
$\q$ of $\z$'s goes to infinity with $\n$ since $(\n-4\t-1)/2\leq\q$,
so $(1-2^{-\t})^\q$ approaches $0$ because $1-2^{-\t}<1$.

This means that any instance of $\dee\t$ will eventually hold, but we
want to know they all hold. There are $(\n-1)^{2t}$ instances of
$\dee\t$. The probability that some instance fails is no more than the
sum over all instances of the probability that each instance fails.
Each instance has a probability of failing that is at most
$(1-2^{-\t})^\q$ with $\q$ depending on the instance. Since
$1-2^{-\t}<1$, bigger exponents make a smaller product and these
probabilities are largest when $\q$ is smallest. We have
$(\n-4\t-1)/2\leq\q$ for every instance, so
$(1-2^{-\t})^{(\n-4\t-1)/2}$ is an upper bound on the probability that
any particular instance of \dee\t\ fails because no $\z$ works.  The
sum of these probabilities is therefore bounded above by the product
of $(\n-1)^{2t}$, the number of instances, times the upper bound
$(1-2^{-\t})^{(\n-4\t-1)/2}$.  Thus, the probability that $\dee\t$
fails is at most $(\n-1)^{2\t}(1-2^{-\t})^{(\n-4\t-1)/2}$. This bound
goes to zero as $\n\to\infty$ because it is a polynomial in $\n$
multiplied by a constant smaller than $1$ raised to a power that is
another polynomial in $\n$. The probability of the complementary
event, that $\dee\t$ holds, therefore approaches $1$ as $\n\to\infty$.

Suppose $3\leq\d<\omega$.  We have seen that the commutative
$\NA$-frame $\gc\K(\R\cup\R_0)$ of a randomly chosen $\R\in\W$ almost
certainly satisfies $\dee{\d-2}$.  Even if it does not, its complex
algebra $\Cm{\gc\K(\R\cup\R_0)}$ is atomic and is in $\NA$.  We will
show, assuming $\dee{\d-2}$ holds, that $\Cm{\gc\K(\R\cup\R_0)}$ is in
$\RA_\d$ because it has a $\d$-dimensional relational basis.  The
atoms of the complex algebra are singletons of elements in $\K$, so
our $\d$-dimensional relational basis $\B$ will consist of
$\d\times\d$ matrices of singletons of elements of $\K$.  Let
$\R'=\R\cup\R_0$.  Let $\B$ consist of those $\x\in\{\{\a\}:\a\in\K
\}^{\d\times\d}$ such that for all $\i,\j,\k<\d$ and all
$\a,\b,\c\in\K$, $\x_{\i\i}=\{0\}$ and if $\x_{\i\j}=\{\a\}$,
$\x_{\j\k}=\{\b\}$, and $\x_{\i\k}=\{\c\}$ then $\<\a,\b,\c\>\in \R'$.

Let $\x\in\B$. To show $\B$ is a $\d$-dimensional basis, we must
verify conditions \eqref{basis1}, \eqref{basis2}, and \eqref{basis3}
in Definition \ref{basis}.  By the definition of complex algebra,
$\id=\{0\}$, but we have $\x_{\i\i}=\{0\}$ by the definition of $\B$,
so $\x_{\i\i}\subseteq\id$.  Given $\i,\j<\d$, suppose
$\x_{\i\j}=\{\a\}$, $\x_{\j\i}=\{\b\}$, and $\x_{\i\i}=\{\c\}$.  By
the definition of $\B$, we have $\c=0$ and $\<\a,\b,\c\>\in\R'$, so
$\<\a,\b,0\>\in\R'$. Then $\<\star\a,0,\b\>\in\R'$ by \eqref{left
  reflection}, so $\star\a=\b$ by \eqref{identity}. But
$\conv{\{\a\}}=\{\star\a\}$ by \eqref{conv}, so $(\x_{\i\j})\conv{} =
\conv{\{\a\}} = \{\star\a\} = \{\b\}=\x_{\j\i}$, as desired.  Given
$\i,\j,\k<\d$, suppose $\x_{\i\j}=\{\a\}$, $\x_{\j\k}=\{\b\}$, and
$\x_{\i\k}=\{\c\}$. By the definition of $\B$, $\<\a,\b,\c\>\in\R'$,
so by \eqref{rp}, $\{\a\}\rp\{\b\} \supseteq\{c\}$, \ie,
$\x_{\i\k}\subseteq\x_{\i\j}\rp\x_{\j\k}$, as desired.  This complete
the proof of part \eqref{basis1} in Definition \ref{basis}.

For part \eqref{basis2}, consider $\a\in\K$. We want $\x\in\B$ with
$\x_{01}=\{\a\}$. It is enough to define $\x$ by $\{\a\} = \x_{01} =
\x_{02}$, $\{\star\a\} = \x_{10} = \x_{20}$, and $\{0\} = \x_{00} =
\x_{11} = \x_{22} = \x_{12} = \x_{21}$, \ie,
\begin{equation*} \x=\begin{bmatrix} \{0\}
  &\{\a\}&\{\a\}\\ \{\star\a\}&\{0\} &\{0\}\\ \{\star\a\}&\{0\}
  &\{0\}\end{bmatrix}.
\end{equation*}
For part \eqref{basis3}, assume $\i,\j<\n$, $\x\in\B$, $\a,\b\in\K$,
$\x_{\i\j}\subseteq\{\a\}\rp\{\b\}$, and $\i,\j\neq\k<\n$.  We need
$\y\in\B$ such that $\y_{\i\k}=\{\a\}$, $\y_{\k\j}=\{\b\}$, and
$\x_{\ell\m}=\y_{\ell\m}$ whenever $\k\neq\ell,\m<\n$.  If $\a=0$ then
the required $\y$ is obtained directly from $\x$ by setting
\begin{equation*}
  \y_{\ell\m}=\begin{cases} \x_{\ell\m}\text{ if } \k\neq\ell,\m,
  \\ \x_{\i\m}\text{ if } \k=\ell\neq\m, \\ \x_{\ell\i} \text{ if }
  \k=\m\neq\ell, \\ \x_{\i\i}\text{ if } \k=\m=\ell.
  \end{cases}
\end{equation*}
The key frame property that shows $\y\in\B$ is \eqref{identity}. A
similar definition can be used for $\y$ when $\b=0$.  We may therefore
assume $0\neq\a,\b$.

Let $\c\in\K^{\d\times\d}$ be the matrix of elements of $\K$ whose
singletons are the entries in the matrix $\x$, so that
$\{\c_{\i\j}\}=\x_{\i\j}$ for all $\i,\j<\d$.  Without loss of
generality we may assume that $\k=\d-1$ and $0=\i\leq\j\leq1$.  We
want $\y\in\B$ such that $\y_{0\k}=\{\a\}$, $\y_{\k\j}=\{\b\}$, and
$\x_{\ell\m}=\y_{\ell\m}$ whenever $\ell,\m<\k=\d-1$.  Define $\y$ on
all arguments differing from $\k$ so that $\y$ agrees with $\x$ by
setting $\y_{\ell\m}=\x_{\ell\m}$ for all $\ell,\m<\k=\d-1$. We must
also set $\y_{\k\k}=\{0\}$, $\y_{0\k}=\{\a\}$, $\y_{\k0}=\{\star\a\}$,
$\y_{\k\j}=\{\b\}$, and $\y_{\j\k}=\{\star\b\}$. What remains is to
choose $\y_{\ell\k}$ and $\y_{\k\ell}$ whenever $\j<\ell<\k=\d-1$ in
such a way that $\y\in\B$.  Note that $\y_{\ell\k}=\{\c_{\ell\k}\}$
and $\y_{\k\ell} = \{\c_{\k\ell}\}$ whenever $\j<\ell<\k=\d-1$. We
extend $\c$ by setting $\c_{\k\k}=0$, $\c_{0\k}=\a$,
$\c_{\k0}=\star\a$, $\c_{\k\j}=\b$, and $\c_{\j\k}=\star\b$. We will
choose $\c_{\ell\k}, \c_{\k\ell}\in\K$ whenever $\j<\ell<\k$, and set
$\y_{\ell\k} = \{\c_{\ell\k}\}$ and $\y_{\k\ell} = \{\c_{\k\ell}\}$.
Note that $\dee{\d-2}$ implies $\dee\t$ whenever $1\leq\t\leq\d-2$.

Suppose $\j=0$. Then $\x_{0\j}=\x_{00}=\{0\}\subseteq\{\a\}\rp\{\b\}$,
hence $\<\a,\b,0\>\in\R'$, so $\<\star\a,0,\b\>\in\R'$ by \eqref{left
  reflection}, $\star\a=\b$ by \eqref{identity}, and finally
$\a=\star\b$ by \eqref{involution}.  Apply $\dee1$ to $\u_1=\c_{10}$
and $\v_1=\a=\star\b$ to get $0\neq\w\in\K$ such that
$\C(\u_1,\v_1,\w)=\C(\c_{10},\a,\w)\subseteq\R'$.  Set $\c_{1\k}=\w$,
$\c_{\k1}=\star\w$, $\y_{1\k}=\{\w\}$, and $\y_{\k1}=\{\star\w\}$.
The proof for the rest of the case $\j=0$ proceeds in the same way as
the case in which $\j=1$, except that we do not know (or need)
$\a=\star\b$.

Assume $\j=1$.  We are done if $\d=3$ so assume $\d>3$ and
$\k=\d-1>2$. Apply $\dee2$ to $\u_1=\c_{20}$, $\u_2=\c_{21}$,
$\v_1=\a$, and $\v_2=\star\b$ to get $0\neq\w\in\K$ such that
$\C(\u_1,\v_1,\w)=\C(\c_{20},\a,\w)=\C(\c_{20},\c_{0\k},\w)\subseteq\R'$
and $\C(\u_2,\v_2,\w) = \C(\c_{21},\star\b,\w)
=\C(\c_{21},\c_{1\k},\w) \subseteq\R'$.  Set $\c_{2\k}=\w$, $\c_{\k2}
= \star\w$, $\y_{2\k}=\{\w\}$, and $\y_{\k2}=\{\star\w\}$.  We are
done if $\d=4$ and $\k=3$ so assume $\d>4$. Apply $\dee3$ to
$\u_1=\c_{30}$, $\u_2=\c_{31}$, $\u_3=\c_{32}$, $\v_1=\a=\c_{0\k}$,
$\v_2=\star\b=\c_{1\k}$, and $\v_3=\c_{2\k}$ to get $\c_{3\k}$ such
that $\C(\c_{30},\a,\c_{3\k}) = \C(\c_{30},\c_{0\k},\c_{3\k})
\subseteq\R'$, $\C(\c_{31},\star\b,\c_{3\k}) = \C(\c_{31}, \c_{1\k},
\c_{3\k}) \subseteq\R'$, and $\C(\c_{32}, \c_{2\k}, \c_{3\k}) =
\C(\c_{32}, \c_{2\k}, \c_{3\k}) \subseteq\R'$, and set $\y_{3\k} =
\{\c_{3\k}\}$.  Continue in this way until $\dee{\d-2}$ has been
used. The conditions compiled in this process show that $\y\in\B$.

This completes the proof that if $\gc\K(\R\cup\R_0)$ has the diamond
property $\dee{\d-2}$ then $\Cm{\gc\K(\R\cup\R_0)}$ has a
$\d$-dimensional relational basis, hence $\Cm{\gc\K(\R\cup\R_0)}
\in\RA_\d$ and every $\d$-provable predicate is valid in
$\gc\K(\R\cup\R_0)$. We have shown that the commutative $\NA$-frame
corresponding to a randomly chosen relation $\R\in\W$ almost certainly
has these properties. As was observed earlier, this is enough to
conclude from parts \eqref{27iii}--\eqref{27vi} that parts
\eqref{27vii}--\eqref{27x} are also true and completes the proof of
Theorem \ref{27}. \endproof

If an equation is true in every representable relation algebra then by
Theorem \ref{key}\eqref{key0}\eqref{key3} there is a smallest
$\n\in\omega$ such that it is true in every algebra in $\RA_\n$.  Let
$\Xi$ be a finite set of equations true in $\RRA$.  Any finite subset
of $\omega$ has a largest element, so there is some $\n$ such that
every equation in $\Xi$ is true in every algebra in $\RA_\n$. By
Theorem \ref{27}\eqref{27vi.1} a randomly chosen commutative
$\NA$-frame has a complex algebra that is almost certainly in
$\RA_\n$. Therefore, a randomly chosen commutative $\NA$-frame almost
certainly validates every $\n$-provable predicate and its complex
algebra almost certainly satisfies every equation in $\Xi$.  For
example, any randomly chosen large $\KR$-frame almost certainly
validates \eqref{L''}, \eqref{M''}, and every other 5-provable
predicate.

There are 594 $\TR$-frames with five elements. In each of them,
\eqref{L''} is valid whenever \eqref{M''} is valid. Predicate
\eqref{M''} is invalid in 286 of them and \eqref{L''} is invalid in
just 73.  There are 390 $\KR$-frames among those 594 $\TR$-frames, and
\eqref{L''} and \eqref{M''} are invalid in 58 of them.  The fractions
of $\TR$-frames and $\KR$-frames in which \eqref{L''} and \eqref{M''}
are invalid shrink to zero as $\n$ increases.  A randomly selected
large $\TR$-frame or $\KR$-frame almost certainly validates
\eqref{L''} and \eqref{M''} and the numbers of such frames both grow
like $\c^{\n^3}$ for some constant $\c>1$.

\section{Questions}\label{sect21}
The results in this paper leave open or suggest a few technical
questions and raise some others of a more general nature.  The
technical questions come first.
\begin{enumerate}
\item Can an axiomatization of $\TR$ be obtained by adding
  \eqref{reflection1}, \eqref{reflection1a}, or \eqref{dedekind} to an
  axiomatization of $\CR$?
\item Which subsets of \eqref{ttt}--\eqref{no.6} axiomatize $\TT_3$,
  $\TT_4$, $\CT_3$, $\CT_4$, and $\CR$?
\item Which of the derived rules of $\CT_3$ listed in Theorem
  \ref{thm11} are either derivable, admissible, or included by
  definition in $\RR$ or $\CR$?  For example, the first two 1-provable
  rules are included in $\RR$ by definition, but the third one turned
  out to be admissible. What about all the others?
\item Are there any deductive rules of $\CT_4$ that require four
  variables?
\item Is almost every finite relation algebra representable?
\item Is relation algebra \alg{29}{83} representable? See the end of
  \SS\ref{sect16}.
\end{enumerate}
Here are some questions about logic, philosophy, and history.
\begin{enumerate}
  \setcounter{enumi}{6}
\item ``Will the real negation please stand up?''
  \cite[p.\ 174]{MR1223997}.  ``Which is the \emph{real} negation?''
  \cite[p.\ 492]{MR1223997}.  Does Table \ref{defs-ops} reveal the
  real negation? The results here suggest that Boolean negation is
  real and De Morgan negation is relevant negation.
\item Do the predicates \eqref{reflection1}, \eqref{reflection1a}, and
  \eqref{dedekind} have logical significance? Would any of them be
  proposed as an axiom for a relevance logic? Why were they never
  previously considered?
\item What are the philosophical implications of the fact for any
  fixed $\n\in\omega$ a randomly selected $\TR$-frame or $\KR$-frame
  will almost certainly validate every $\n$-provable predicate?
\item The algebra $\Re\U$ of binary relations on a set $\U$ is the
  prototypical example of a relation algebra.  Since $\Re\U$ is the
  complex algebra of the pair-frame on $\U$, could the pair-frame on
  $\U$ serve as a prototypical example of a frame for relevance logic?
\item Why do relevance logic and relation algebra overlap despite
  arising independently through the pursuit of completely different
  goals?
\item Why were Schr\"oder's studies, Tarski's axiomatization, and
  relevance logic confined to the 4-variable fragment of the calculus
  of relations?
\end{enumerate}

\end{document}